\documentclass[12pt,a4paper,oneside]{article}
\usepackage[english]{babel}%
\usepackage[round]{natbib}
\usepackage[utf8]{inputenc}
\usepackage{setspace,graphicx,amsmath,amsfonts,amssymb,amsthm}
\graphicspath{ {./images/} }
\usepackage{latexsym}
\usepackage{mathrsfs}
\usepackage{color}
\usepackage{graphicx}
\usepackage{subfigure}
\usepackage{booktabs}
\usepackage{array}       
\usepackage{tabularx}
\usepackage{dcolumn} 
\usepackage{lscape,longtable}
\usepackage{rotating}
\usepackage{multirow}
\usepackage[autostyle]{csquotes}
\usepackage{hyperref}
\usepackage{xpatch}
\usepackage{bm}
\usepackage{bbm}
\usepackage{booktabs}
\usepackage{multirow}
\usepackage{multicol}
\usepackage{amsbsy}
\usepackage[table,dvipsnames]{xcolor}
\usepackage{upgreek}
\usepackage[ruled, vlined]{algorithm2e}
\usepackage{tkz-graph}
\usetikzlibrary{arrows}
\usepackage{comment}
\SetKwInput{kwInit}{Initialize}
\SetKwInput{kwInput}{Input}
\SetKwInput{kwOutput}{Output}
\newtheorem{theorem}{Theorem}[section]

\newtheorem{remark}{Remark}[section]

\newcommand\bbone{\ensuremath{\mathbbm{1}}}
\hypersetup{
    colorlinks = {true}, 
    linkcolor = {red}, 
    urlcolor = {black}, 
    citecolor = {blue}
}

\usepackage[lmargin=1in, rmargin=1in, tmargin=1in, bmargin=1in,paperwidth=8.5in, paperheight=11in]{geometry}

\def\qmo{``}
\def\qmc{''}
\def\qmcsp{'' }

\title{\bf Bayesian Dynamic Quantile Model Averaging}
\author{M. Bernardi \\
Department of Statistical Sciences, University of Padova 
\and 
R. Casarin\\  University Ca' Foscari of Venice, Italy
\and
B.~B.~Maillet\\  University of La Reunion, Saint-Denis, France;\\ University Ca’ Foscari of Venice, Italy; Variances, Paris, France.
\and
L.~Petrella\\  Sapienza University of Rome, Italy.
}

\date{October 16, 2024}

\begin{document}

\maketitle

\begin{abstract}
\noindent This article introduces a novel dynamic framework to Bayesian model averaging for time-varying parameter quantile regressions. By employing sequential Markov chain Monte Carlo, we combine empirical estimates derived from dynamically chosen quantile regressions, thereby facilitating a comprehensive understanding of the quantile model instabilities. The effectiveness of our methodology is initially validated through the examination of simulated datasets and, subsequently, by two applications to the US inflation rates and to the US real estate market. Our empirical findings suggest that a more intricate and nuanced analysis is needed when examining different sub-period regimes, since the determinants of inflation and real estate prices are clearly shown to be time-varying. In conclusion, we suggest that our proposed approach could offer valuable insights to aid decision making in a rapidly changing environment.
\end{abstract}

\vspace{0.5cm}

\section{Introduction}
\label{sec:intro} 
%
\noindent The study of extremes, risks and quantiles of phenomena is a broad area of interest in Operations Research, since large movements of variables are of great importance in many fields related to natural disasters, namely geophysics \citep[e.g. earth-quakes prediction][]{pisarenko1997statistical}, hydrology and environmental studies \citep[such as floods, tsunami, severe droughts; e.g.][]{lin_etal.2018}, bio-medecine \citep{huang2017quantile} and, economics and finance \citep[see, e.g.][]{amedee-manesme_barthelemy.2018,benkraiem_zopounidis.2021,bellini_etal.2021}, specifically when forecasting is the ultimate objective \citep[see, e.g.][]{chalamandaris_vlachogiannakis.2018,tian_etal.2020,candila2023mixed}. \newline
\indent Among the mobilized techniques related to extremes, quantile regressions (denoted QR, hereafter) provide a simple way to model the conditional quantiles of a response variable with respect to some covariates in order to have a more satisfactory representation of its conditional distribution than one can have with traditional linear regression. Since the seminal works by \cite{koenker_basset.1978} and \cite{koenker.2005}, QR has become popular in economics and finance \citep{behr.2010,chun_etal.2012,rockafellar_royset.2013,candila2023mixed}, especially in recent years \citep[][]{lin_etal.2018,ben_ameur_etal.2020,tsionas.2020,ben_bouheni_etal.2021,tian_etal.2020}, as a simple, robust and distribution free modelling tool. 
More specifically, a QR approach is mostly appropriate not only when the underlying model is non-linear, or the innovations are non-Gaussian, but also when modelling the distribution tails is the primary interest of the analysis. See \cite{koenker.2005}, \cite{lum_gelfand.2012} and \cite{huang2017quantile} for a review of the recent advances in QR modeling and inference. However, two greatest challenges in predicting extreme events are the robustness of the forecast accuracy in a changing environment and the correct representation of the statistical features of the variable under studies, as explicated below.\newline
\indent One of the most crucial issues in QR analyses is, indeed, related to model specification. The literature primarily focuses on shrinkage penalization approaches for simultaneous handling of both the estimation and model selection. For instance, \cite{wu_liu.2009} examine the properties of variable selection based on the least absolute shrinkage and selection operator (LASSO) and smoothly clipped absolute deviation (SCAD) penalties for static models, while \cite{wang_etal.2012} investigate SCAD penalization for high-dimensional QR models. \cite{noh_etal.2012} advocate a model selection method for time-varying QR, based on coefficient basis expansion and penalization. \cite{lee_etal.2014} proposed instead a simpler alternative, based on Bayesian information criterion. A concern in the QR estimation process emerges from the fact that the objective function is not differentiable with respect to the parameters, which makes the derivation of the asymptotics of QR estimators almost unfeasible. Thus, we followed a Bayesian model averaging (BMA) approach to overcome these difficulties \citep[see, e.g.][]{ji_etal.2012,mamatzakis2021making}. \newline
\indent From a computational point of view, Bayesian model selection requires the estimation and the post-processing of the $2^M$ models, where $M\in\mathbb{N}$ being the number of regressors. Badly, it is almost impossible to be run as soon as $M$ is even only moderately large. For this main reason, a computationally lighter version of the stochastic search variable selection (SSVS) by \cite{george_mcculloch.1993}, has been proposed for QR models in \cite{meligkotsidou_etal.2009}, \cite{reed_etal.2009} and \cite{alhamzawi_yu.2012}. In this paper, we propose an alternative approach, extending the Dynamic Model Averaging (DMA) by \cite{raftery_etal.2010} to QR models. Recently, DMA for linear models has gained significant attention in econometrics \citep{koop_korobilis.2012,koop_korobilis.2013,koop_tole.2013}. To the best of our knowledge, it remains the only computationally feasible algorithm when handling a large number of regressors. Several extensions of the DMA technique have already been proposed. \cite{mccormick_etal.2012} introduced a DMA algorithm for binary outcome regression models. \cite{belmonte_koop.2014} proposed an alternative model-based mechanism for the time-varying model switching. However, an extension to more general models, such as QR models, is still absent in the literature, primarily due to the computational cost associated with increased model complexity. In this article, we demonstrate how this cost can be significantly reduced by implementing the Sequential MCMC (SMCMC) algorithm recently proposed by \cite{guhaniyogi_etal.2018}, which is an efficient and self-tuning approach for simulating both parameters and latent variables.\newline
\indent In summary, our novel Bayesian Dynamic Quantile Model Averaging (BDQMA) combines three complementary overlays as an attempt for mitigating the specification and estimation errors inherent to our operational research problem: {\it (1)} a time-varying parameter QR setting \citep{koenker.2005}; {\it (2)} a BMA framework \citep{steel2020model}; and {\it (3)} a DMA approach \citep{raftery_etal.2010} with the use of SMCMC \citep{guhaniyogi_etal.2018}. The contribution to the literature is manifolds. First, we extend the model selection procedure for QR models by \citep{kim.2007} to time-varying parameters QR models. Second, we extend the classical DMA for linear models by \cite{raftery_etal.2010} to conditionally linear models. Our article makes a significant contribution also to the model combination literature \cite{stone.1961,hall_mitchell.2007,billio_etal.2012,fawcett_etal.2015,gneiting_ranjan.2013,bassetti_etal.2018,casarin_etal.2013,busetti.2017,billio_etal.2013,aastveit_etal_2014,aastveit_etal_2018,casarin2023flexible,mcalinn_west.2019}.\newline
\indent We investigate the effectiveness of the time-varying parameter quantile regression model, supported by a sequential estimation procedure, in extracting quantile signals from synthetic datasets designed to replicate real-world features such as abrupt and smooth changes in the relationship between the dependent variables and the covariates. This validation highlights the model's robustness and applicability to dynamic financial and macroeconomic time series, where capturing extreme quantiles and handling time-varying volatility are particularly important. We illustrate further the potential of our methodology with an application to two real-world datasets. Financial and macroeconomic time series are widely documented in the literature as exhibiting distinctive features, such as heavy tails in the conditional distribution and time-varying volatility \citep[e.g.][]{engle.1982,bollerslev.1986}. Addressing these characteristics is particularly important when the study focuses on extreme quantiles \citep[e.g.][]{taylor2020forecast,taylor2022forecasting,bonaccolto_etal.2022,Bel21,candila2023mixed}. Accordingly, our first application examines inflation \citep[following][]{koop_korobilis.2012}, while the second focuses on the U.S. real estate market \citep[following][]{risse_kern.2016,amedee-manesme_barthelemy.2018}.\newline
\indent The rest of the article is structured as follows. Section \ref{sec:DMA_quantile}  introduces our BDQMA approach. Section \ref{sec:simulated_examples} first underscores the effectiveness of our methodology using synthetic data. Section \ref{sec:empirical_application_inflation} and \ref{sec:empirical_application_estate} provide two empirical applications, whilst Section \ref{sec:conclu} encapsulates our concluding remarks. 
%
\section{Bayesian Dynamic Quantile Model Averating}
\label{sec:DMA_quantile} 
\noindent The building blocks of our BDQMA approach are: {\it (i)} time-varying parameter quantile regression (QR) modelling within a Bayesian framework; {\it (ii)} dynamic model averaging (DMA) within a Bayesian framework, and {\it (iii)} sequential Markov chain Monte Carlo (SMCMC) for sequential estimation of model risk and prediction of extreme events. Fig. \ref{fig:blocks} provides a graphical illustration of the three blocks. In the following subsections, we provide a detailed explanation of each building block and how they contribute to the overall methodology.
 \tikzstyle{frame} = [
        rectangle, draw, 
        text width=8em, text centered,
        minimum height=2em
    ]
 \tikzstyle{frameLarge} = [
        rectangle, draw, 
        text width=20em, text centered,
        minimum height=1em
    ]
 \tikzstyle{frameCirc} = [
        circle, draw, 
        text width=5em, text centered,
        minimum height=2em
    ]
    
\tikzstyle{line} = [draw, -latex']
\begin{figure}[t]
\centering
\begin{tikzpicture}[node distance = 4cm]
    \node [frame] (pop) {\small Quantile regressions};
    \node [frame, right of=pop] (B) {\small Dynamic Model Averaging};
    \node [frame, right of=B] (SMCMC) {\small Monte Carlo methods};
    \node [frameLarge, below of=B, yshift=1.5cm] (BDQMA) {\small Bayesian Dynamic Quantile Model Averaging};
    \node [frameCirc, above of=pop, yshift=-1.5cm,xshift=-0.5cm] (TVP) {\small Time-varying parameters};
    \node [frameCirc, right of=TVP, xshift=-1.5cm] (BF) {\small Bayesian framework};
    \node [frameCirc, right of=BF, xshift=0cm] (SF) {\small Sequential framework};    
    \node [frameCirc, above of=SMCMC, yshift=-1.5cm, xshift=0.5cm] (MCMC) {\small Markov chains};       
    \path [line] (TVP) -- (pop);
    \path [line] (BF) -- (pop);
    \path [line] (BF) -- (B);
    \path [line] (SF) -- (B);
    \path [line] (SF) -- (SMCMC);
    \path [line] (MCMC) -- (SMCMC);
    \path [line] (pop) -- (BDQMA);
    \path [line] (B) -- (BDQMA);
    \path [line] (SMCMC) -- (BDQMA);
\end{tikzpicture}
\caption{The three main modeling and inference blocks (boxes) of our Bayesian Dynamic Quantile Model Averaging (BDQMA) framework together with the main features of the blocks (circles).}
\label{fig:blocks}
\end{figure}
%
%
\subsection{Time-varying quantile regressions}
Regarding the first block, the dynamic QR model used in this paper extends the one in \cite{bernardi_etal.2015} to the case all the parameters evolve randomly over time. Let $\mathbf{x}_{t}=\left(1,x_{2,t},\ldots
,x_{M,t}\right)^{\top}$, $t=1,2,\ldots,T$ be the collection of $M$ exogenous variables. The dependent variable, $y_{t}$, is a linear function of $\mathbf{x}_{t}$, and, as in \cite{harvey.1989}, the parameter vector $\boldsymbol{\beta}_{t}$ follows a multivariate
random walk: 
\begin{align}
\label{eq:TV_qreg_measurement}
y_{t}& =\mathbf{x}_{t}^{\top}\boldsymbol{\beta}_{t}+\xi _{t}\\
 \label{eq:TV_qreg_transition} 
\boldsymbol{\beta}_{t}& =\boldsymbol{\beta}_{t-1}+\boldsymbol{\zeta}_{t}, \\
 \label{eq:TV_qreg_init}
\boldsymbol{\beta}_{0}& \sim \mathsf{MVN}_{M}\left( \boldsymbol{\beta}
_{0\vert 0},\mathbf{P}_{0|0}\right),  
\end{align}
where $\xi _{t}\sim\mathsf{AL}\left( \tau ,0,\sigma \right) $, $t=1,2,\ldots ,T$ are i.i.d.  variables with centred Asymmetric Laplace (AL)
distribution,  with $\tau\in\left(0,1\right)$ the quantile level, $\sigma\in\mathbb{R}^+$ the scale parameter and $\boldsymbol{\beta}_{0}\in\mathbb{R}^M$ the initial state of the evolution of $\boldsymbol{\beta}_t$ defined by the transition Eq. \eqref{eq:TV_qreg_transition}. We assume that $\boldsymbol{\beta}_{0}\in\mathbb{R}^M$ has null mean and diffuse variance-covariance matrix $\mathbf{P}_{0\vert 0}=\kappa 
\mathbf{I}_{M}$, with $\kappa \to +\infty $ and that the error terms $
\boldsymbol{\zeta }_{t}\sim\mathsf{MVN}_{M}\left( \mathbf{0},\boldsymbol{\Omega} \right) $,  are i.i.d.
and independent of the measurement equation errors $\xi _{s}$, $\forall t,s=1,2,\ldots ,T$. Hereafter, $\mathsf{MVN}_{M}(\cdot,\cdot)$ denotes the $M$-dimensional Gaussian random variable, while $\mathsf{N}(\cdot,\cdot)$ is exclusively reserved for representing univariate Gaussian random variables.\newline
\indent The linear state space model introduced in Eq. \eqref{eq:TV_qreg_measurement}-\eqref{eq:TV_qreg_init} for representing time-varying conditional quantiles deviates from Gaussianity due to the assumptions imposed on the measurement innovation terms. Consequently, traditional optimal filtering techniques, such as those relying on Kalman filter recursions, cannot be applied \citep[see,][]{durbin_koopman.2012}. Nonetheless, drawing upon the insights of \cite{kozumi_kobayashi.2011} and \cite{bernardi_etal.2015}, it is possible to leverage the representation of the centered Asymmetric Laplace probability density function, $f\left(\xi _{t}\vert\tau,\sigma\right)$, as a location-scale continuous Normal mixture. This approach yields:
\begin{equation}
f\left(\xi _{t}\vert\tau,\sigma\right)=\int_{0}^{\infty }\frac{1}{\sqrt{2\pi \delta
\sigma \omega }}\exp \left\{ -\frac{1}{2\delta \sigma \omega }(\xi
_{t}-\lambda \omega )^{2}\right\} \frac{1}{\sigma} \exp \left\{ -\frac{1}{\sigma}
\omega \right\}d\omega,
\label{eq:ald_stoch_repr}
\end{equation}
for $t=1,2,\dots ,T$. Setting $\lambda =\frac{1-2\tau }{\tau (1-\tau )}$ and $
\delta =\frac{2}{\tau (1-\tau )}$ ensures that the $\tau $-level
quantile of the measurement error $\xi _{t}$ becomes zero. Consequently, it becomes evident that our non-Gaussian state space model admits the following conditionally Gaussian and linear state space representation:
\begin{align}
\label{eq:ApproxGaussianSSM_measurement}
y_{t}& =\mathbf{x}_{t}^{\top}\boldsymbol{\beta}_{t}+\lambda \omega
_{t}+\varepsilon _{t},\quad \varepsilon _{t}\overset{\mathsf{iid}}{\sim }\mathsf{N}
\left( 0,\delta \sigma \omega _{t}\right),\\
\label{eq:ApproxGaussianSSM_transition}
\boldsymbol{\beta}_{t}& =\boldsymbol{\beta}_{t-1}+\boldsymbol{\zeta }
_{t},\quad \quad \quad \quad \boldsymbol{\zeta }_{t}\overset{\mathsf{iid}}{\sim }
\mathsf{MVN}_{M}\left( 0,\boldsymbol{\Omega} \right),\\
\boldsymbol{\beta}_{0}& \sim \mathsf{MVN}_{M}\left( \boldsymbol{\beta}
_{1|0},\mathbf{P}_{1|0}\right) ,
\end{align}
where $\omega _{t}$, $t=1,2,\ldots ,T$, are i.i.d. with an exponential
distribution of parameter $\sigma ^{-1}$, that is $\omega _{t}\sim \mathsf{Exp
}\left( \sigma ^{-1}\right)$. We further assume the
following prior distributions for the parameters $\sigma $ and $\boldsymbol{\Omega}$: 
\begin{equation}
\sigma  \sim \mathsf{IG}\left( a_{0},b_{0}\right),\qquad
\boldsymbol{\Omega}  \sim \mathsf{IW}\left( c_{0},\mathbf{C}_{0}\right),
\label{eq:prior_trans_eq}
\end{equation}
which are Inverse Gamma and Inverse Wishart distributions, with probability density functions $p\left(\sigma\vert a_0,b_0\right) \propto \sigma ^{-\left( a_{0}+1\right)}\exp
\left\{ -\frac{b_{0}}{\sigma }\right\}$ and $p \left(\boldsymbol{\Omega}\vert c_0,\boldsymbol{\Omega}_0\right)  \propto \vert\mathbf{C}_{0}\vert^{c_{0}}\vert\boldsymbol{\Omega} |^{-\left(
c_{0}+\frac{M+1}{2}\right) }\exp \left\{ -\mathrm{trace}\left( \mathbf{C}_{0}\boldsymbol{\Omega}
^{-1}\right) \right\}$, respectively, where $\left(a_0,b_0,c_0,\boldsymbol{\Omega}_0\right)$ are given hyperparameters.
%
%
\subsection{Bayesian dynamic model averaging}
For the second building block of our method, we propose a dynamic procedure for selecting and combining dynamic quantile regression (QR) models. The challenge of determining which regressors to include in the model can be addressed from two perspectives. One approach is to use a variable dimension model \citep[see, e.g.][]{marin_robert.2007}. However, in this study, we take a different route, opting for a model selection approach that requires estimating all possible sub-models. Drawing inspiration from the seminal work of \cite{raftery_etal.2010}, we implement a recursive model selection process over time to continually evaluate model adequacy. Expanding on the framework introduced by \cite{raftery_etal.2010}, we introduce a time-varying model index $L_t$, which takes values in $\left\{1,2,\ldots,K\right\}$, representing the selected model at each time point $t$. We assume it has a Markov dynamics with transition matrix 
\begin{equation}
\mathsf{P}(L_{t}=l\vert L_{t-1}=j)=\pi_{lk}, \quad\qquad l,j\in\left\{1,2,\ldots,K\right\}.
\label{eq:model_transition_mat}
\end{equation}
Let $\mathbf{x}
_{kt}=(1,x_{i_{1}t},\ldots,x_{i_{m^{(k)}}t})^{\top}$ be the set of $
(m^{(k)}+1)$ variables in the model $k$, with $i_{j}\in\{2,3,\ldots,M\}$, $
j=1,\ldots,m^{(k)}$, $i_{j}\neq i_{l}$ $\forall\,l,j$ and $
m^{(k)}\in\{0,1,\ldots,M-1\}$. 
Then, we posit that each model can be expressed in the form of Eq. \eqref{eq:TV_qreg_measurement}-\eqref{eq:TV_qreg_init}, where the dynamic regressor parameter vector $\boldsymbol{\beta}_t^{(k)}$ and the covariate vector $\mathbf{x}_t^{(k)}$ are unique to each model $k=1,2,\dots,M$. This assumption implies that at each time $t$ and for each model $k=1,2,\dots,K$, the system variables represented by the dynamic QR parameters $\boldsymbol{\beta}_t$ and the model indicator variable $L_t$, i.e. $\left\{\boldsymbol{\beta}_t,L_t\right\}$, transition to a new state according to the transition matrix specified in Eq. \eqref{eq:model_transition_mat} and the kernel defined by the dynamics of latent states in Eq. \eqref{eq:TV_qreg_transition}. The updated states comprise the vector of predicted values at time $t+1$ for each model-specific quantile regression parameter vector $\boldsymbol{\beta}_{t+1}^{(k)}$ and the updated time-varying probability $\pi_{t+1}^{(k)}$ associated with that model.\newline
\indent As anticipated in the Introduction, the dynamic quantile regression model specified in Eq. \eqref{eq:multimodel_representation_trans_eq} exhibits non-Gaussian characteristics, rendering the iterative application of the linear Kalman filter and smoother impractical for obtaining updated estimates of the latent dynamics over time. To address this limitation, we once again leverage the representation of the Asymmetric Laplace error term in the measurement equation (Eq. \eqref{eq:TV_qreg_measurement}) as a location-scale mixture of Normals. This representation hinges on the introduction of an additional latent factor $\omega_t^{(k)}$, unique to each model $k$, following an i.i.d. Exponential distribution with shape parameter $\sigma^{(k)}$, i.e. $\omega_{t}^{(k)}\sim\mathsf{Exp}(1/\sigma^{(k)})$, i.i.d. $\forall t$.
\newline
\indent For the model selection purposes, we derive from Eq. 
\eqref{eq:ApproxGaussianSSM_measurement}-\eqref{eq:ApproxGaussianSSM_transition} the following:
\begin{equation}
\begin{aligned}
y_{t}&=\sum_{k=1}^{K}\bbone_{\{k\}}(L_{t})\big(\mathbf{x}_{t}^{(k)^{\top}}
\boldsymbol{\beta}_{t}^{(k)}+\lambda\omega_{t}^{(k)}+\sqrt{
\delta\sigma^{(k)}\omega_{t}^{(k)}}\varepsilon_{t}^{(k)}\big)\\
\boldsymbol{\beta}_{t}^{(k)}&=\boldsymbol{\beta}_{t-1}^{(k)}+\boldsymbol{
\zeta}_{t}^{(k)},  \label{eq:multimodel_representation_trans_eq}
\end{aligned}
\end{equation}
where $\varepsilon_{t}^{(k)}\overset{\mathsf{i.i.d.}}{\sim}\mathsf{N}(0,1)$ and $
\boldsymbol{\zeta}_{t}^{(k)}\overset{\mathsf{i.i.d.}}{\sim}\mathsf{MVN}_{k+1}(\mathbf{0
},\boldsymbol{\Omega}^{(k)})$. The parameters $\lambda$, $\delta$, and $\tau$ are defined
as above. The description of the model is completed by the prior specification for the initial state given in Eq. \eqref{eq:TV_qreg_init}. 
%
%
\subsection{Sequential Markov chain Monte Carlo}
In the third block, we adopt a Bayesian approach to inference for QR models, and introduce auxiliary variables to represent the original model in Eq. \eqref{eq:TV_qreg_measurement}-\eqref{eq:TV_qreg_init} as a conditionally Gaussian and linear state space model. Let $\mathbf{z}_{1:t}=\left(
\mathbf{z}_{1},\mathbf{z}_{2},\ldots,\mathbf{z}_{t}\right)$ denote the sequence of vectors $\mathbf{z}_u$ up to time $t$. Then, the complete-data likelihood of the unobservable components $\boldsymbol{\beta}_{1:T}$ and $\boldsymbol{\omega}_{t:T}$ and all parameters $\boldsymbol{\gamma}=\left(\sigma,\boldsymbol{\Omega}\right)$ is given by:
\begin{align*}
\label{eq:compllete_likelihood}
&\mathcal{L}\left(\boldsymbol{\gamma},\boldsymbol{\omega}_{1:T},\boldsymbol{
\beta}_{1:T}\vert\mathbf{y}_{1:T},\mathbf{x}_{1:T}\right)\propto
f\left(\mathbf{y}_{1:T}\vert\mathbf{x}_{1:T},\boldsymbol{
\beta}_{1:T},\boldsymbol{\omega}_{1:T},\boldsymbol{\gamma}\right)f\left(\boldsymbol{
\beta}_{1:T},\boldsymbol{\omega}_{1:T}\vert \boldsymbol{\gamma}\right)\pi(\boldsymbol{\gamma})\\
&\propto\prod_{t=1}^Tf\left(y_{t}\vert\boldsymbol{\beta}
_t,\omega_{t},\sigma,\mathbf{x}_t\right)
f\left(\omega_{t}\vert\sigma\right)f\left(\boldsymbol{\beta}_1
\right)\prod_{t=2}^{T}f\left(\boldsymbol{\beta}_{t}\vert\boldsymbol{
\beta}_{t-1},\boldsymbol{\Omega}\right)  \notag \\
&\propto\prod_{t=1}^T\left(\sigma\times\omega_{t}\right)^{-\frac{1}{2}
} \exp\left\{-\frac{1}{2\sigma\delta}\sum_{t=1}^T\frac{\left(y_t-\lambda
\omega_t-\mathbf{x}_t^{\top}\boldsymbol{\beta}_t\right)^2}{\omega_t}\right\}\nonumber\\
&\qquad\times\exp\left\{-\frac{\sum_{t=1}^T\omega_{t}}{\delta}
\right\}\exp\left\{-\frac{\boldsymbol{\beta}
_1^{\top}\boldsymbol{\beta}_1}{2\kappa}\right\}  \notag \\
&\qquad\times \vert\boldsymbol{\Omega}\vert^{-\frac{T-1}{2}}\exp\left\{-\frac{1}{2}\sum_{t=2}^{T}\left(\boldsymbol{\beta}
_{t}-\boldsymbol{\beta}_{t-1}\right)^{\top}\boldsymbol{\Omega}^{-1}\left(\boldsymbol{\beta
}_{t}-\boldsymbol{\beta}_{t-1}\right)\right\}.
\end{align*}
The Gaussian scale mixture representation poses challenges for sequential inference of latent variables and parameters. In this model, the DMA approach introduced by \cite{raftery_etal.2010} may become computationally burdensome, as exact estimation of latent variables and parameters across different models requires posterior approximation. Standard Markov chain Monte Carlo (MCMC) schemes for posterior approximation are not feasible in such cases.
To address this issue, we propose combining a sequential MCMC procedure (SMCMC) \citep[][]{dunson_yang.2013} with the DMA approach to create a feasible sequential DMA procedure for latent variable models. Additionally, we introduce enhancements to reduce computational time. Firstly, we demonstrate how the forgetting factor, as introduced by \cite{raftery_etal.2010} for the state covariance, can expedite computations without compromising the validity of the SMCMC procedure. Secondly, we illustrate that computing costs, which typically increase linearly over time, can be stabilized using fixed-lag backward sampling within the SMCMC transition kernel. Subsequent subsections present the SMCMC for the single-model scenario, its extension to the multi-model scenario, and strategies for mitigating computational complexity.
%
\subsubsection{Sequential model estimation}
\label{sec:sequential_model_estim} 
%
\noindent Let $\boldsymbol{\theta}_{t}=\left(\sigma,\boldsymbol{\Omega},\boldsymbol{\omega}_{1:t},
\boldsymbol{\beta}_{1:t}\right)$, $t\in\mathbb{N}$, denote the time series of augmented parameter vectors, taking values in measurable spaces $\left(\mathbb{R}
^{d_{t}},\mathcal{B}\left(\mathbb{R}^{d_{t}}\right)\right)$, where the dimension non-decreasingly grows as $
d_{t}=d_{t-1}+d$, for $t\geq 1$. We assume that the augmented parameter
vector can be decomposed as $\boldsymbol{\theta}_{t}=\left(\boldsymbol{\theta}_
{t-1},\boldsymbol{\eta}_{t}\right)$, where $\boldsymbol{\eta}_{t}=\left(\omega_{t},
\boldsymbol{\beta}_{t}\right)$ represents the latent variable vector of dimension $d$ associated with the observation $y_{t}$.\newline
\indent The proposed quantile regression model with time-varying parameters possesses a prior distribution $\pi_{t}\left(\boldsymbol{\theta}_{t}\right)=\pi\left(\boldsymbol{\gamma
}\right)\pi\left(\boldsymbol{\omega}_{1:t},\boldsymbol{\beta}_{1:t}\right)$, which satisfies the compatibility condition:
\begin{equation}
\pi_{t+1}\left(\boldsymbol{\theta}_{t}\right)=\int_{\mathbb{R}^{d}} \pi_{t}\left(
\boldsymbol{\theta}_{t},\boldsymbol{\eta}_{t+1}\right)d\boldsymbol{\eta}_{t+1},
\end{equation}
enabling a simplification of notation. Furthermore, we denote the posterior distribution at time $t$, with respect to the Lebesgue measure of $\boldsymbol{\theta}_{t}$, as $\pi_t(\boldsymbol{\theta}_{t})=\pi\left(\boldsymbol{\theta}_{t}|\mathbf{y}_
{1:t}\right)\propto \pi\left(\boldsymbol{\gamma}\right)\mathcal{L}\left(\boldsymbol{\gamma},
\boldsymbol{\omega}_{1:t},\boldsymbol{\beta}_{1:t}\vert\mathbf{y}_{1:t},
\mathbf{x}_{1:t}\right)$. Henceforth, in the notation, we omit the conditioning on the dependent variables $\mathbf{y}_{1:t}$ and the covariates $\mathbf{x}_{1:t}$ from the posterior.\newline
\indent In the SMCMC algorithm, $L$ parallel inhomogeneous Markov chains are employed to generate samples $\boldsymbol{\theta}_
{t}^{(l,j)}$ with $j=1,\ldots,m_{t}$, $l=1,\ldots,L$ and $t=1,2,\ldots,T$
from the sequence of posterior distributions $\pi_{t}$, $t=1,2,\ldots,T$. Each
Markov chain within the population is defined by a sequence of transition
kernels $\mathcal{K}_{t}\left(\boldsymbol{\theta},A\right)$, $t\in\mathbb{N}$, which are operators
from $\left(\mathbb{R}^{d_{t-1}},\mathcal{B}\left(\mathbb{R}^{d_{t-1}}\right)\right)$ to $\left(\mathbb{
R}^{d_{t}},\mathcal{B}\left(\mathbb{R}^{d_{t}}\right)\right)$, such that $\mathcal{K}_{t}\left(\boldsymbol{
\theta},\cdot\right)$ is a probability measure for all $\boldsymbol{\theta}\in
\mathbb{R}^{d_{t-1}}$, and $\mathcal{K}_{t}\left(\cdot,A\right)$ is measurable for all $A\in
\mathcal{B}\left(\mathbb{R}^{d_{t}}\right)$.\newline
\indent The kernel $\mathcal{K}_{t}\left(\boldsymbol{\theta},A\right)$ has stationary distribution $\pi_{t}$ and results from the composition of a jumping kernel, $\mathcal{J}_{t}$, and a transition kernel, $\mathcal{T}_{t}$, expressed as:
\begin{equation*}
\mathcal{K}_{t}(\boldsymbol{\theta},A)=\mathcal{J}_{t}\circ \mathcal{T}_{t}^{m_{t}}\left(\boldsymbol{\theta}
,A\right)=\int_{\mathbb{R}^{d_{t}}}\mathcal{J}_{t}\left(\boldsymbol{\theta},d\boldsymbol{\theta}
^{\prime}\right)\mathcal{T}_{t}^{m_{t}}\left(\boldsymbol{\theta}^{\prime},A\right),
\end{equation*}
where the fixed dimension transition is defined as: 
\begin{equation}
\mathcal{T}_{t}^{m_{t}}\left(\boldsymbol{\theta},A\right)=\mathcal{T}_{t}\circ \mathcal{T}_{t}^{m_{t}-1}\left(\boldsymbol{
\theta},A\right)=\int_{\mathbb{R}^{d_{t}}}\mathcal{T}_{t}(\boldsymbol{\theta},d\boldsymbol{
\theta}^{\prime})\mathcal{T}_{t}^{m_{t}-1}\left(\boldsymbol{\theta}^{\prime},A\right),
\end{equation}
with $m_{t}\in\mathbb{N}$, and $\mathcal{T}^{0}$ is the identity kernel.
We assume that the jumping kernel satisfies 
\begin{equation}
\mathcal{J}_{t+1}\left(\boldsymbol{\theta}_{t},\boldsymbol{\theta}_{t+1}\right)=\mathcal{J}_{t+1}\left(
\boldsymbol{\theta}_{t},\boldsymbol{\eta}_{t}\right)\delta_{\boldsymbol{\theta}
_{t}}(\boldsymbol{\theta}_{t+1}),
\end{equation}
where $\mathcal{J}_{t+1}\left(\boldsymbol{\theta}_{t},\boldsymbol{\eta}_{t+1}\right)=\mathcal{J}_{t+1}\left(
\boldsymbol{\theta}_{t},\left(\boldsymbol{\theta}_{t},\boldsymbol{\eta}_{t+1}\right)\right)$. This condition ensures that the error propagation through the jumping kernel can be effectively managed across the iterations of SMCMC. This result, presented in Theorem \ref{th1}, is a straightforward extension of the findings established in the $L_{1}$-norm by \cite{guhaniyogi_etal.2018} and related works. Let us introduce the $v$-norm, denoted as $||\mu(\cdot)||_v=\sup{|f|\leq v}|\mu(f)|$, where $v:\mathbb{R}^{n}\rightarrow[
1,\infty]$ and $\mu$ represents a signed measure \citep[see][Ch. 14]{meyn_tweedie.1993}.
\begin{theorem}\label{th1}
For any probability density function $p\left(\cdot\right)$, and for $\boldsymbol{\theta}_{t-1}\in\mathbb{R}^{d_{t-1}}$ the next inequality is met:
\begin{equation}
||\pi_{t}-\mathcal{J}_{t}\circ p\Vert_{v}\leq \underset{\boldsymbol{\theta}_{t-1}\in\mathbb{R}^{d_{t}-1}}{\sup}||\pi_{t}\left(\cdot\vert\boldsymbol{\theta}_{t-1}\right)-\mathcal{J}_{t}(\boldsymbol{\theta}_{t-1},\cdot)||_{\widetilde{v}}+||\pi_{t}- p\Vert_{v},
\end{equation}
where $\widetilde{v}=\int_{\mathbb{R}^{d_{t-1}}}v\left(\left(\boldsymbol{\theta}_{t-1},\boldsymbol{\eta}_{t}\right)\right)d\boldsymbol{\theta}_{t-1}$.
\end{theorem}
\begin{proof}
See Appendix \ref{subsec:Proof}.
\end{proof}
\noindent The following findings outlined in \cite{guhaniyogi_etal.2018} demonstrate the convergence of the chain with transition kernel $\mathcal{J}_{t}\circ \mathcal{T}^{m_{t}}$ towards the target distribution, under highly general conditions.
\begin{theorem}
\label{th2}
Let us assume the following conditions hold:
\begin{itemize}
\item[{(i)}] \textit{(Geometric ergodicity)} For each $t$, there exists a function $V_t:\mathbb{R}^{d_{t}}\to\left[1,\infty\right)$, $C>0$ and $\rho_{t}\in\left(0,1\right)$ such that:
\begin{enumerate}
\item[(a)]$\int_{\mathbb{R}^{d_{t}}}V\left(\boldsymbol{\theta}_{t}\right)^{2}\pi_{t}\left(\boldsymbol{\theta}_{t}\right)d\boldsymbol{\theta}_{t}\leq C$;
\item[(b)]$\int_{\mathbb{R}^{d_t}}V\left(\boldsymbol{\theta}_{t}\right)^{2}\pi_{t}(\boldsymbol{\eta}_{t}\vert\boldsymbol{\theta}_{t})d\boldsymbol{\eta}_{t}=V_{t-1}\left(\boldsymbol{\theta}_{t-1}\right)$, where $\boldsymbol{\theta}_{t}=\left(\boldsymbol{\theta}_{t-1},\boldsymbol{\eta}_{t}\right)$;
\item[(c)] for all $\boldsymbol{\theta}_{t}\in\mathbb{R}^{d_{t}}$,  $||\mathcal{T}_{t}\left(\boldsymbol{\theta}_{t},\cdot\right)-\pi_{t}(\cdot)||_{v_{t}}\leq V_{t}\left(\boldsymbol{\theta}_{t}\right)\rho_{t}$.
\end{enumerate}
\item[{(ii)}]\textit{(Stationary convergence)} The stationary distribution $\pi_t$ of $\mathcal{T}_{t}$ satisfies
\begin{equation}
\alpha_{t}=2\sqrt{C}d_{H}\left(\pi_{t},\pi_{t-1}\right)\rightarrow 0,\nonumber
\end{equation}
where $\pi_t$ is the marginal posterior of $\boldsymbol{\theta}_{t-1}$ at time $t$.
\item[{(iii)}] \textit{(Jumping consistency)} For a sequence of $\lambda_{t}\rightarrow 0$ the following holds:
\begin{equation}
\underset{\boldsymbol{\theta}_{t-1}\in\mathbb{R}^{d_{t-1}}}{\sup}
||\mathcal{J}_{t}(\boldsymbol{\theta}_{t-1},\cdot)-\pi_{t}\left(\cdot\vert\boldsymbol{\theta}_{t-1}\right) ||_{\widetilde{v}_{t}},\nonumber
\end{equation}
where $\widetilde{v}_{t}=\int_{\mathbb{R}^{d}}v\left(\left(\boldsymbol{\theta},\boldsymbol{\eta}\right)\right)d\boldsymbol{\theta}_{t-1}$.
\end{itemize}
Let $\varepsilon_{t}=\rho_{t}^{m_{t}}$. Then for any initial distribution $\pi_{0}$,
\begin{equation}
||\mathcal{K}_{t}\circ\cdots\circ \mathcal{K}_{1}\circ \pi_{0}-\pi_{t}||_{v_{t}}
\leq \sum_{s=1}^{t}\left(\prod_{u=s}^{t}\varepsilon_{u}\right)\left(\alpha_{s}+\lambda_{s}\right).\nonumber
\end{equation}
\end{theorem}
\begin{proof}
See Appendix \ref{subsec:Proof}.
\end{proof}
\noindent 
Regarding the assumptions of geometric ergodicity and jumping consistency, several works have investigated convergence rates of Gibbs samplers for Bayesian models. For instance, \cite{roman_etal.2012} establish geometric ergodicity for a range of improper priors in Bayesian linear models. Additionally, studies by \cite{jones_hobert.2004} and \cite{papaspiliopoulos_gareth.2008} are relevant in this context. In the present paper, as we employ a Gibbs sampler for a conditionally linear and Gaussian model with proper conditionally conjugate prior distributions, these assumptions are met in accordance with \cite{roman_etal.2012}.\newline
\indent To employ SMCMC, one must specify the transition kernel $\mathcal{T}_{t+1}$ and the jumping kernel $\mathcal{J}_{t+1}$ at the iteration $t+1$. The transition kernel $\mathcal{T}_{t}$ at the iteration $t$ enables each parallel chain to explore the sample space of dimension $d_t$ and generate samples $\boldsymbol{\theta}_{t}^{(l,j)}$ from the posterior distribution $\pi_t$, given the previously generated samples. The collapsed Gibbs sampler algorithm is detailed in Algorithm \ref{alg:main_algorithm}. Within this algorithm, we opt for a transition kernel akin to a blocked (or multi-move) kernel.
\begin{algorithm}[!t]
\KwData{$\mathbf{y}_{1:T}$, $\mathbf{x}_{1:T}$} 
\For {$t=1,2,\dots,T$, and $l=1,2,\dots,L$}	{
draw $\sigma^{(l)}$ from $f(\sigma\vert\boldsymbol{\beta}^{(l)}_{1:t},
\mathbf{y}_{1:t},\mathbf{x}_{1:t})$ marginalizing the Gibbs with respect to the $\boldsymbol{
\omega}_{1:t}$;\\
draw $\omega_{s}^{(l)}$ from $f(\omega_{s}\vert\sigma^{(l)},
\boldsymbol{\Omega}^{(l)},\boldsymbol{\beta}_{1:t}^{(l)},\mathbf{y}_{1:t},\mathbf{x}_{1:t})$, $
s=1,2,\ldots,t $;\\	
draw $\boldsymbol{\beta}_{1:t}^{(l)}$ from $f(\boldsymbol{\beta}
_{1:t}\vert\sigma^{(l)},\boldsymbol{\Omega}^{(l)}, \boldsymbol{\omega}^{(l)}_{1:t},\mathbf{y}
_{1:t},\mathbf{x}_{1:t})$ with a multi-move proposal distribution;\\
draw $\boldsymbol{\Omega}^{(l)}$ from $f(\boldsymbol{\Omega}\vert\boldsymbol{\beta}
^{(l)}_{1:t})$\;}
\caption{Collapsed Gibbs sampler.}
\label{alg:main_algorithm}
\end{algorithm}
The sampling strategy for the full conditional distributions is elaborated in Appendix \ref{subsec:TV_qreg_MCMCT}. Leveraging the stochastic representation of the AL distribution in Eq. \eqref{eq:ald_stoch_repr}, we implement a partially collapsed SMCMC strategy based on Gibbs-type updating with data augmentation \citep{liu.1994,van_dyk_park.2008}. The core concept of the complete collapsed Gibbs-type simulation scheme is to avoid simulations from the full conditional distributions of blocks of model parameters whenever possible by analytically marginalizing them out. As demonstrated by \cite{park_van_dyk.2009}, this approach offers several advantages over systematic sampling, such as reduced computational time and enhanced convergence rate of the sampler. In our scenario, we can only partially integrate out the latent variables $\omega^{(l)}_{1:s}$, for $s=1,2,\dots,t$, from the full conditional of the scale parameters $\sigma^{(l)}$ in the preceding step 3, for $l=1,2,\dots,L$, thereby reducing the algorithm to a partially collapsed Gibbs sampler type move. Concerning the convergence of partially collapsed Gibbs moves, it is notable that updating the parameters in the specified order ensures that the posterior distribution at each time point $t=1,2,\dots,T$ corresponds to the stationary distribution of the generated Markov chain. In fact, combining steps 1 and 2 produces sample draws from $\pi(\sigma^{(l)},\boldsymbol{\omega}_{1:t}^{(l)}\vert\boldsymbol{\beta}^{(l)}_{1:t},\mathbf{y}_{1:t},\mathbf{x}_{1:t})$, for $t=1,2,\dots,T$ and $l=1,2,\dots,L$, i.e.  the conditional posterior distribution. The partially collapsed Gibbs sampler thus represents a blocked version of the ordinary Gibbs sampler \citep{van_dyk_park.2008,park_van_dyk.2009}.\newline 
\indent To initialize the Gibbs sampling, we begin by generating a random draw from the joint prior distribution of the parameters defined in Eq. \eqref{eq:prior_trans_eq}. Conditional on this draw, we simulate the initial values of the augmented variables $\omega^{(l)}_{1}$, for $l=1,2,\dots,L$. Subsequently, at iteration $t+1$, the jumping kernel $\mathcal{J}_{t+1}$ facilitates the transition of chains from a space of dimension $d_{t}$ to one of dimension $d_{t+1}$. We designate the jumping kernel of the $l$-th parallel chain as the transition kernel of a Gibbs sampler, characterized by the following full conditional distributions:
\begin{enumerate}
\item draw $\omega_{t+1}^{(l)}$ from $
f(\omega_{t+1}\vert\sigma^{(l)},\boldsymbol{\beta}_{t+1}^{(l)},y_{t+1},\mathbf{x}_{t+1})$;
\item draw $\boldsymbol{\beta}_{t+1}^{(l)}$ from $f(\boldsymbol{
\beta}_{t+1}\vert\sigma^{(l)},\boldsymbol{\Omega}^{(l)},\boldsymbol{\beta}_{t}^{(l)},
\omega_{t+1}^{(l)},y_{t+1},\mathbf{x}_{t+1})$.
\end{enumerate}
\noindent The specifics of the sampling strategy for the full conditional distributions are detailed in Appendix \ref{subsec:TV_qreg_MCMCJ}. Thanks to the location-scale mixture representation of the AL distribution, all the full conditional distributions involved in the transition and jumping kernels have a known closed-form representation. This feature is especially advantageous when $T$ is large or when multiple models are estimated simultaneously, as in the current scenario, because the availability of conditional sufficient statistics that can be tracked helps mitigate the increase in storage and computational burden over time.
%
\subsubsection{Monitoring convergence}
\label{sec:monitoring_convergence} 
%
\noindent The number of iterations of the Sequential MCMC sampler at each time $m_{t}$ is determined based on the cross-chain correlation \citep[see][]{guhaniyogi_etal.2018}. Specifically, we set the number of iterations at time $t$, denoted as $m_{t}$, to be the smallest integer $s$ such that $r_{t}\left( s\right) \leq 1-\epsilon $, where $r_{t}\left( s\right) $ represents the rate function associated with the transition kernel $\mathcal{T}_{t}$ and $\epsilon $ is a predefined threshold level. An upper bound for the rate function can be estimated by the lag-$s$ chain autocorrelation, which is calculated online using information provided by the output of all the parallel chains. This estimation is performed as follows:
\begin{equation}
\widehat{r}_{t}\left( s\right) =\max_{j=1,2,\dots ,p}\frac{\sum_{l=1}^{L}(
\theta _{j}^{(s+1,t,l)}-\bar{\theta}_{j}^{(s+1,t)}) (\theta
_{j}^{(1,t,l)}-\bar{\theta}_{j}^{(1,t)}) }{\left( \sum_{l=1}^{L}(
\theta _{j}^{(s+1,t,l)}-\bar{\theta}_{j}^{(s+1,t)})^{2}\right)^{
\frac{1}{2}}\left( \sum_{l=1}^{L}( \theta _{j}^{(1,t,l)}-\bar{\theta}
_{j}^{(1,t)}) ^{2}\right) ^{\frac{1}{2}}},
\end{equation}
where $\theta _{j}^{(s,t,l)}$, is the $j$-th element of the vector $
\boldsymbol{\theta }^{(s,t,l)}$ staking the fixed parameters $\boldsymbol{
\gamma }^{(l)}$ and the latent states generated up to time $t$, i.e. $( 
\boldsymbol{\omega}_{1:t}^{(l)},\boldsymbol{\beta}_{1:t}^{(l)}) $, by the $l$-th chain at the the $s$-th iteration. Moreover, $\bar{\theta}_{j}^{(s,t)}=L^{-1}
\sum_{l=1}^{L}\theta _{j}^{(s,t,l)}$ is the average of the draws of the $s$-th iteration over the $L$ parallel chains. 
%
\subsubsection{Sequential model selection}
\label{sec:sequential_model_selection} 
%
\noindent Consider scenarios where multiple models are evaluated and estimated simultaneously, especially when there is uncertainty about the correct specification of the data-generating process. We proceed under the assumption that the entire array of models can be represented as in Eq. \eqref{eq:multimodel_representation_trans_eq}, where the indicator variable $L_{t}\in \left\{ 1,2,\ldots ,K\right\} $ has a Markovian transition kernel, selecting the operative model at each time step $t=1,2,\ldots,T$. Furthermore, it is assumed that the divergence among models primarily arises from the presence or absence of specific explanatory variables within the predictor vector $\mathbf{x}_{t}$. Assuming there are $M-1$ potential explanatory variables, with the intercept being a constant across all models, the total number of models becomes $K=2^{M-1}$. Even with a moderate number of variables, such as $M=10$, the count of models can increases substantially, reaching, for instance, $1,024$, thereby making the formulation of the transition matrix a formidable task. Following the approach outlined by \cite{raftery_etal.2010}, we adopt a forgetting factor methodology, as described in the following.\newline
\indent Let $\pi _{t+1|t}^{(k)}=\mathsf{P}( L_{t+1}=k\vert \mathbf{y}
_{1:t},\mathbf{x}^{(k)}_{t+1})$ and $\pi _{t+1|t+1}^{(k)}=\mathsf{P}( L_{t+1}=k\vert\mathbf{y}
_{1:t+1},\mathbf{x}^{(k)}_{t+1})$, where again $\mathbf{y}_{1:t}=\left(y_{1},y_{2},\ldots
,y_{t}\right)$. The approximated Hamilton filter recursion involves the following two steps. At each time point $t$, the first step predicts and updates the model probabilities as follows:
\begin{align}
\pi _{t+1|t}^{(k,l)}& =\frac{(\pi _{t|t}^{(k,l)})^{\alpha }+\upsilon }{
\sum_{j=1}^{K}(\pi _{t|t}^{(j,l)})^{\alpha }+\upsilon },\qquad \qquad \forall
l=1,2,\ldots ,L,  \label{eq:prob_pred} \\
\pi _{t+1|t+1}^{(k)}& \propto \sum_{l=1}^{L}\pi _{t+1\vert
t}^{(k,l)}f(y_{t+1}\vert \mathbf{y}_{1:t},\mathbf{x}_{t+1}^{(k)},L_{t}=k,\omega _{t+1}^{(k,l)}).
\label{eq:prob_update}
\end{align}%
In Eq. \eqref{eq:prob_pred}, the exponent $\alpha$ serves as a forgetting factor, typically set close to one, while the constant $\upsilon >0$ is a small number introduced to mitigate rounding errors due to machine precision, with $\upsilon$ set to $\upsilon =0.001/K$, following the methodology outlined in \cite{raftery_etal.2010}.

To circumvent the necessity of explicitly specifying and computing the entire transition matrix of the Markov chain $L_{t}$, which comprises $K(K-1)$ entries, we introduce the forgetting factor $\alpha$. This factor allocates prior weight to information preceding time $t$, with a higher $\alpha$ parameter implying reduced weight assigned to subsequent observations at time $t+1$. Essentially, higher values of $\alpha$ indicate less prior confidence in the informational significance of new observations during the model probability update process.

As argued by \cite{raftery_etal.2010}, although this approach lacks an explicit specification of the transition matrix, it does not pose a problem for inference in model probabilities as long as the data provides sufficient information regarding which models better capture the data dynamics. Regarding the updating of model probabilities in Eq. \eqref{eq:prob_update}, it is obtained as the Monte Carlo average of the approximated predictive distributions $f(y_{t+1}\vert \mathbf{y}_{1:t},\mathbf{x}^{(k)}_{t+1},L_{t}=k,\omega _{t+1}^{(k,l)})$ over the entire population of independent Markov chains, effectively marginalizing out the simulated latent factor $\omega _{t+1}$. The predictive distribution of the observable at time $t+1$ is obtained as a by-product of the Kalman filtering updating equations for the approximated Gaussian model:
\begin{equation}
f(y_{t+1}\vert \mathbf{y}_{1:t},\mathbf{x}^{(k)}_{t+1},L_{t}=k,\omega _{t+1}^{(k,l)})\approx \mathsf{N}
( y_{t+1}\vert \widehat{y}_{t+1|t}^{(k,l)},V_{t+1}^{(k,l)}),
\end{equation}
where $\widehat{y}_{t+1|t}^{(k,l)}=\lambda \omega _{t+1}^{(k,l)}+{\mathbf{x}^{(k)}
_{t+1}}^{\top}\boldsymbol{\beta}_{t+1|t}^{(k,l)}$ is the predictive mean of $y_{t+1}$
for the $l$-th chain and model $k\in \left\{ 1,2,\dots ,K\right\} $ and the
variance-covariance matrix equal to $V_{t+1}^{(k,l)}=\big( \delta \sigma
^{(k,l)}\omega^{(k,l)}+{\mathbf{x}_{t+1}^{(k)}}^{\top}\mathbf{P}_{t+1|t}^{(k,l)}\mathbf{x}^{(k)}_{t+1}\big)
^{-1}$, from the Kalman filter recursion given in Appendix \ref{subsec:TV_qreg_MCMCT}.\newline
\indent The updated model probabilities are normalised as follows: 
\begin{equation}
\widetilde{\pi}_{t+1|t+1}^{(k,l)}=\frac{\pi _{t+1|t+1}^{(j,l)}}{
\sum_{j=1}^{K}\pi _{t+1|t+1}^{(j,l)}},\qquad \forall l=1,2,\ldots ,L,
\end{equation}%
to get a proper probability measure.\newline
\indent The second step involves updating the latent factors $\boldsymbol{\beta}_{t+1}^{(l,k)}$ and $\omega _{t+1}^{(l,k)}$, for $l=1,2,\ldots ,L$, and this can be easily accomplished by running the SMCMC algorithm defined in the previous section sequentially for each model $k=1,2,\ldots ,K$. This process is then iterated as a new sample observation becomes available. The sequence is initialized at $t=0$, with the initial model probabilities $\pi _{0|0}^{(k,l)}=\frac{1}{K}$ set uniformly for all the models $k=1,2,\ldots ,K$, and for all the parallel chains of the SMCMC sampler $l=1,2,\ldots ,L$. Static model parameters $\left( \sigma ^{(k,l)},\boldsymbol{\Omega} ^{(k,l)}\right) $ are initialized by drawing from their prior distributions defined in Eq. \eqref{eq:prior_trans_eq}. 
%
\subsubsection{Dynamic quantile averaging}
\label{sec:quantile_model_averaging} 
%
\noindent Forecasting with SMCMC methods offers a streamlined approach to integrating 
parameter uncertainty. The following steps delineate the process of generating 
one-step-ahead quantile forecasts from all competing models at each time point $t$, 
and subsequently merging them using updated model weights. These procedures are 
executed at the conclusion of each SMCMC iteration within the sampling phase, 
leveraging the complete parameter set of the model denoted as $\boldsymbol{\theta }_{t}^{(k)}$. 
At every time instance during the estimation period, an approximation 
of the predictive mean for the $\tau$-level quantile function $q_{\tau}( \mathbf{x},\boldsymbol{\beta})={\mathbf{x}^{\top}\boldsymbol{\beta}}$, can be derived for the $k$-th model, as follows: 
\begin{align}
\widehat{q}_{\tau ,t+1}^{(k)}( \mathbf{x}_{t+1}^{(k)},\boldsymbol{\beta}
_{t+1}^{(k)}) & =\mathsf{E}\Big[q_{\tau}( \mathbf{x}
_{t+1}^{(k)},\boldsymbol{\beta}_{t+1}^{(k)}) \vert \mathbf{y}_{1:t},
\mathbf{x}_{t+1}^{(k)},\boldsymbol{\beta}_{1:t}^{(k)}\Big]   \nonumber \\
& ={\mathbf{x}_{t+1}^{(k)}}^{\top}\mathsf{E}\Big[ \boldsymbol{\beta}
_{t+1}^{(k)}\vert \mathbf{y}_{1:t},\mathbf{x}_{t+1}^{(k)},\boldsymbol{\beta}
_{1:t}^{(k)}\Big].  \label{eq:quantile_forecast_single_model}
\end{align}
Here, $\widehat{\boldsymbol{\beta}}_{t+1}^{(k)}=\mathsf{E}\Big[
\boldsymbol{\beta}_{t+1}^{(k)}\vert \mathbf{y}_{1:t},\mathbf{x}_{t+1}^{(k)},
\boldsymbol{\beta}_{1:t}^{(k)}\Big]$ represents the mean of the marginal posterior distribution of the regression parameters for model $k$, i.e.  $\boldsymbol{\beta}_{t+1}^{(k)}$. This can be obtained by numerically marginalizing out the nuisance parameters $\left( \sigma ,\boldsymbol{\Omega} \right)
^{(k)}$ and the latent factors $\omega _{1:t+1}^{(k)}$ using the population
of  chains:
\begin{align}
\widehat{\boldsymbol{\beta}}_{t+1}^{(k)}& =\mathsf{E}\Big( \boldsymbol{\beta}_{t+1}^{(k)}\vert \mathbf{y}_{1:t},\mathbf{x}_{t+1}^{(k)},\boldsymbol{\beta}_{1:t}^{(k)}\Big)   \nonumber \\
& \approx \frac{1}{L}\sum_{l=1}^{L}\mathsf{E}\Big( \boldsymbol{\beta}
_{t+1}^{(k)}\vert \mathbf{y}_{1:t},\mathbf{x}_{t+1}^{(k)},\boldsymbol{\beta}
_{1:t}^{(k,l)},\boldsymbol{\omega }_{1:t}^{(k,l)},\sigma ^{(k,l)},\boldsymbol{\Omega}
^{(k,l)}\Big).
\end{align}
In the right term, the expectation is analytically provided by the Kalman filter predictive Eq. \eqref{eq:kf_updating_mean}, averaged over the entire set of parameters $(\sigma ^{(k,l)},\boldsymbol{\Omega} ^{(k,l)}, \boldsymbol{\omega }_{1:t}^{(k,l)}) _{l=1}^{L}$ generated by the parallel chains, as follows:
\begin{equation}
\mathsf{E}\left( \boldsymbol{\beta}_{t+1}^{(k)}\vert \mathbf{y}_{1:t},
\mathbf{x}_{t+1}^{(k)},\boldsymbol{\beta}_{1:t}^{(k)}\right) \approx \frac{1
}{L}\sum_{l=1}^{L}\left( \boldsymbol{\beta}
_{t-1|t-2}^{(k,l)}+\mathbf{P}_{t-1|t-2}^{(k,l)}\mathbf{x}_{t-1}V_{t-2}^{(k,l)}\nu
_{t-2}^{(k,j)}\right) ,
\end{equation}
where all elements of the expression are given in Appendix \ref{subsec:TV_qreg_MCMCT}. The resulting Rao-Blackwellized estimate of the quantile function is more efficient than a simple average over the samples from the full conditional distribution of $\boldsymbol{\beta} _{t+1}^{(k)}$, as argued by \citep[][pp. 130--134]{robert_casella.2004}.\newline
\indent The model-averaged one-step-ahead $\tau $-level quantile prediction is then formed by combining previous forecasts using the predicted model probabilities $\pi _{t+1|t}^{(k)}$, for each competing model $k=1,2,\ldots,K$:
\begin{align}
\label{eq:quantile_dma_generico}
\widehat{q}_{\tau ,t+1}^{\mathsf{DMA}}( \mathbf{x}_{t+1}^{(k)},\boldsymbol{
\beta }_{t+1}^{(k)}) & =\sum_{k=1}^{K}\pi _{t+1|t}^{(k)}\widehat{q}_{\tau
,t+1}^{(k)}( \mathbf{x}_{t+1}^{(k)},\boldsymbol{\beta}
_{t+1}^{(k)}) \\
 \label{eq:quantile_dma_specifico}
& =\sum_{k=1}^{K}\pi _{t+1|t}^{(k)}\mathbf{x}_{t+1}^{(k)}\widehat{
\boldsymbol{\beta}}_{t+1}^{(k)}, 
\end{align}
where the predicting model probabilities $\pi _{t+1|t}^{(k)}$, for $
k=1,2,\ldots ,K$ have been obtained by averaging single model predictive
probabilities defined in Eq. \eqref{eq:prob_pred} over the $L$ parallel
chains: 
\begin{equation}
\pi _{t+1|t}^{(k)}=\frac{1}{L}\sum_{l=1}^{L}\pi _{t+1|t}^{(k,l)},\qquad
\forall k=1,2,\dots ,K,  \label{eq:model_pred_prob}
\end{equation}%
and Eq. \eqref{eq:quantile_dma_specifico} can be obtained by
substituting for the definition of the predicted model specific $\tau$-level quantile in Eq. \eqref{eq:quantile_forecast_single_model} into
Eq. \eqref{eq:quantile_dma_generico}. 
The multi-model predictions of the $\tau$-level quantile function of the response variable $y_{t}$ at time $t$ are obtained through a weighted average of the model-specific $\tau$-level quantile predictions $\widehat{q}_{\tau ,t}^{(k)}(\mathbf{x}_{t}^{(k)}, \boldsymbol{\beta}_{t}^{(k)})$. These weights correspond to the posterior predictive model probabilities for sample $t$, $\pi _{t+1|t}^{(k)}$, derived as in Eq. \eqref{eq:model_pred_prob}. Furthermore, since model predictive probabilities in Eq. \eqref{eq:model_pred_prob} are obtained by averaging $\pi _{t+1|t}^{(k,l)}$ over the $L$ independent parallel chains, we can provide approximated confidence intervals for the estimated probabilities by calculating their Monte Carlo variance as follows:
\begin{equation}
\mathsf{var}( \pi _{t+1|t}^{(k)}) =\frac{1}{L-1}\sum_{l=1}^{L}(\pi
_{t|t-1}^{(k,l)}-\pi _{t|t-1}^{(k)})^{2}.
\label{eq:model_pred_prob_var}
\end{equation}
%
\subsubsection{Reducing the computational complexity}
\label{sec:reducing_complexity} 
%
\noindent Fixed-lag backward sampling is an effective approach to alleviate computational demands, particularly when dealing with a large number of observations. This involves replacing the fixed-interval backward smoothing algorithm outlined in Appendix \ref{subsec:TV_qreg_FixLagSmo} with a fixed-lag smoother. The fixed-lag smoothing procedure updates the dynamic latent states within a specified lag $h\leq t$ surrounding the current observation $y_{t}$, while keeping the remaining states from $1$ to $h-1$ unchanged. At each iteration, the resulting sampling procedure only simulates the quantile regression coefficient $\boldsymbol{\beta}_{t-h\vert t}$ at time $t$, for a given lag $h>0$. For further details, refer to \cite{anderson_moore.1979} and \cite{simon.2006}.
\begin{remark} The fixed-lag smoother can be viewed as a form of local-dynamic quantile regression, which extends the non-parametric approach introduced by \cite{yu_jones.1998}.
\end{remark}
%
\section{BDQMA efficiency on simulated series}
\label{sec:simulated_examples} 
%
\noindent We now focus on simulated examples specifically designed to assess the effectiveness of the SMCMC algorithm and the BDQMA approach in accurately recovering the true parameter values of the underlying Data Generating Process (DGP). These simulations are designed to reflect real-world challenges, such as abrupt and smooth changes in the relationships between dependent variables and covariates. These simulations provide crucial insights into the robustness and adaptability of the methods in handling complex dynamics, further validating their performance in practical scenarios.\newline
\indent Throughout all simulated scenarios, the sample length remains fixed at $T=200$ observations, aligning with the sample size of the real data examples presented in subsequent sections. The covariates are generated from a uniform distribution on $\left( -\frac{T}{2}, \frac{T}{2}\right)$, denoted as $x_{i,t}\sim \mathsf{U}\left( -\frac{T}{2},\frac{T}{2}\right)$, i.i.d. for $i=1,2,\ldots ,M$ and $t=1,2,\ldots,T$. The innovation term follows a Gaussian distribution, i.e.  $\varepsilon _{t}\sim \mathsf{N}\left( 0,\nu _{t}^{2}\right)$, independently for $t=1,2,\ldots,T$, with heteroskedastic variance. Furthermore, the true DGP is specified as a time-varying parameters regression model:
\begin{equation}
y_{t}=\mathbf{x}_{t}^{\top}\boldsymbol{\beta}_{t}^{\ast }+\varepsilon _{t},\qquad \quad
\forall t=1,2,\dots ,T,
\end{equation}%
where $M=2$, $\mathbf{x}_{t}=\left( 1,x_{1,t},x_{2,t}\right) $ and the quantile
regression parameters dynamics $\boldsymbol{\beta}_{t}^{\ast }=\left( \beta _{1,t}^{\ast
},\beta _{2,t}^{\ast },\beta _{3,t}^{\ast }\right)^{\top}$ are defined as follows. 
%
\subsection{Smooth change in quantiles}
\label{ex:sim_ex_1} 
%
\noindent We assume the true quantile parameters exhibit the following dynamics: the constant term $\beta_{1,t}^*$ remains fixed at $\beta_1^{*}=-2.5$, the coefficient $\beta_{2,t}^*$ undergoes a change in slope at a specific time point $t=100$, and the third parameter $\beta_{3,t}^*$ follows a smooth, sinusoidal transition between two distinct levels:
\begin{align*}
\beta_{2,t}^{*}&= \bigg(0.6-\frac{0.4t}{100}\bigg)\bbone_{(-\infty,100]}(t)+\bigg(\frac{0.4t}{100}-0.2\bigg)\bbone_{(100,\infty)}(t)
\notag \\
\beta_{3,t}^{*}&=a+b\left(1+\exp\left\{\frac{-c\left(2t-T-2\right)}{T}
\right\}\right)^{-1}  \notag \\
\nu_t^2&=\bbone_{(-\infty,100]}(t)+0.25\bbone_{(100,\infty)}(t),
\end{align*}
with $a=0.2$, $b=1$ and $c=5$. 
\begin{figure}[!t]
\begin{center}
\resizebox*{0.9\textwidth}{!}{\subfigure[$\beta_{1,t}$]{\label{fig:SMCMC_Example_1_Reg1}
\includegraphics[trim={0 1cm 1cm 2cm},clip,width=0.45\textwidth]{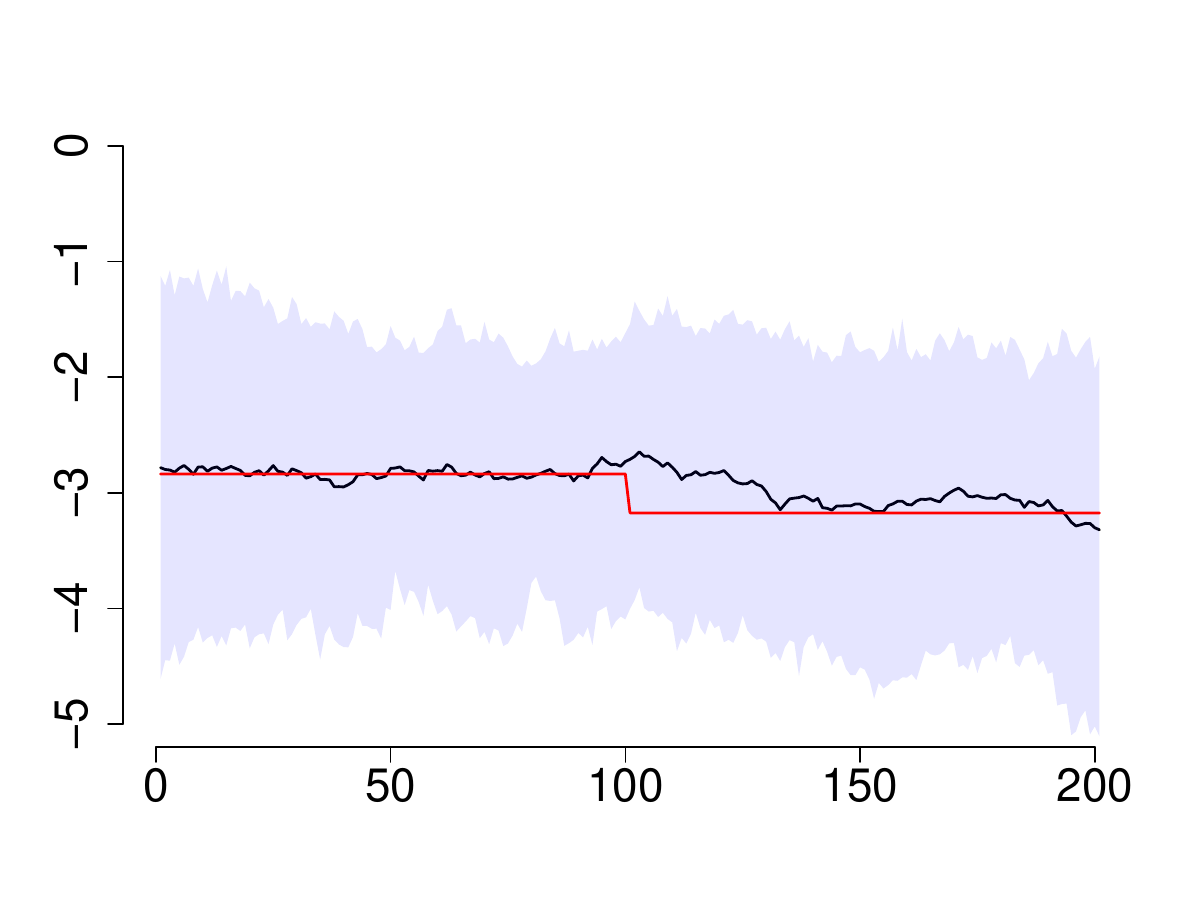}}\qquad
\subfigure[$\beta_{1,t}$]{\label{fig:SMCMC_Example_2_Reg1}
\includegraphics[trim={0 1cm 1cm 2cm},clip,width=0.45\textwidth]{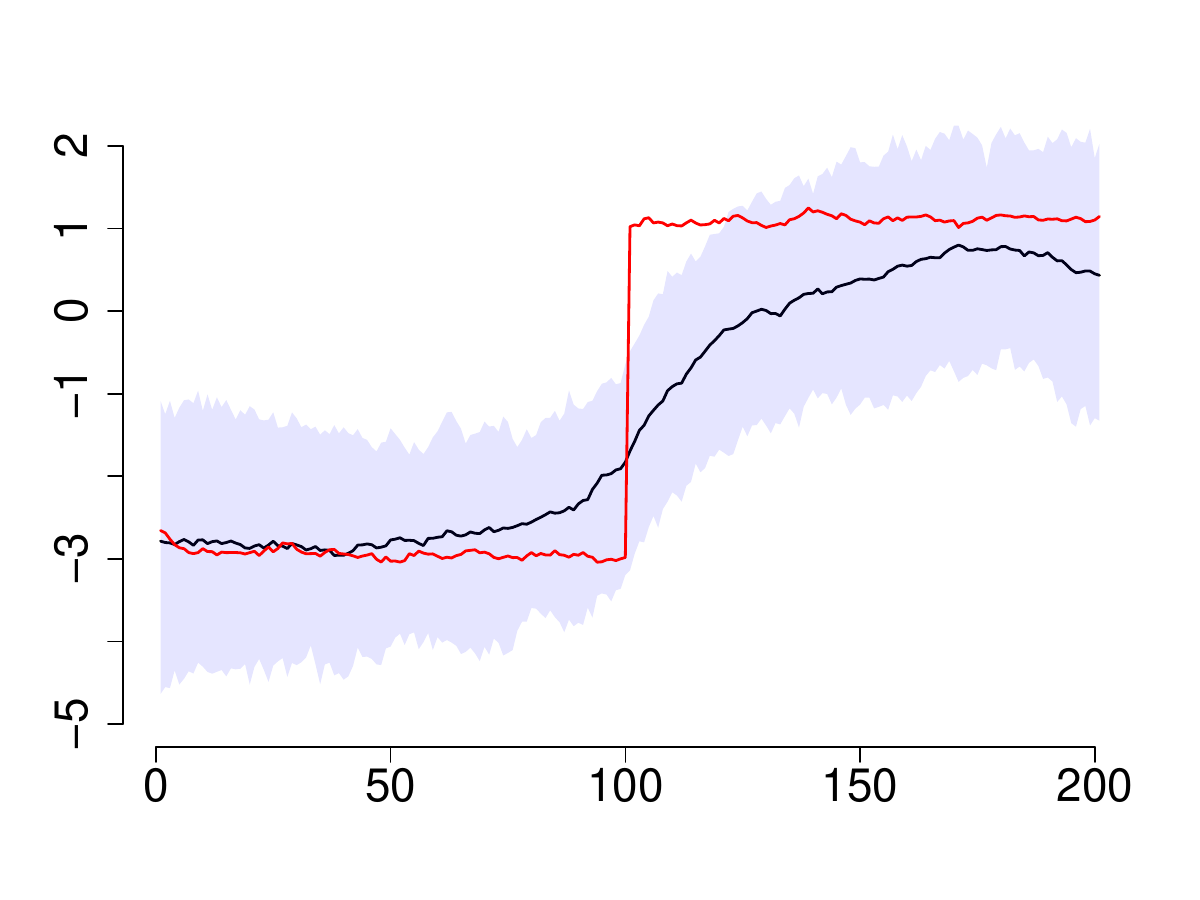}}}
\resizebox*{0.9\textwidth}{!}{\subfigure[$\beta_{2,t}$]{\label{fig:SMCMC_Example_1_Reg2}
\includegraphics[trim={0 1cm 1cm 2cm},clip,width=0.45\textwidth]{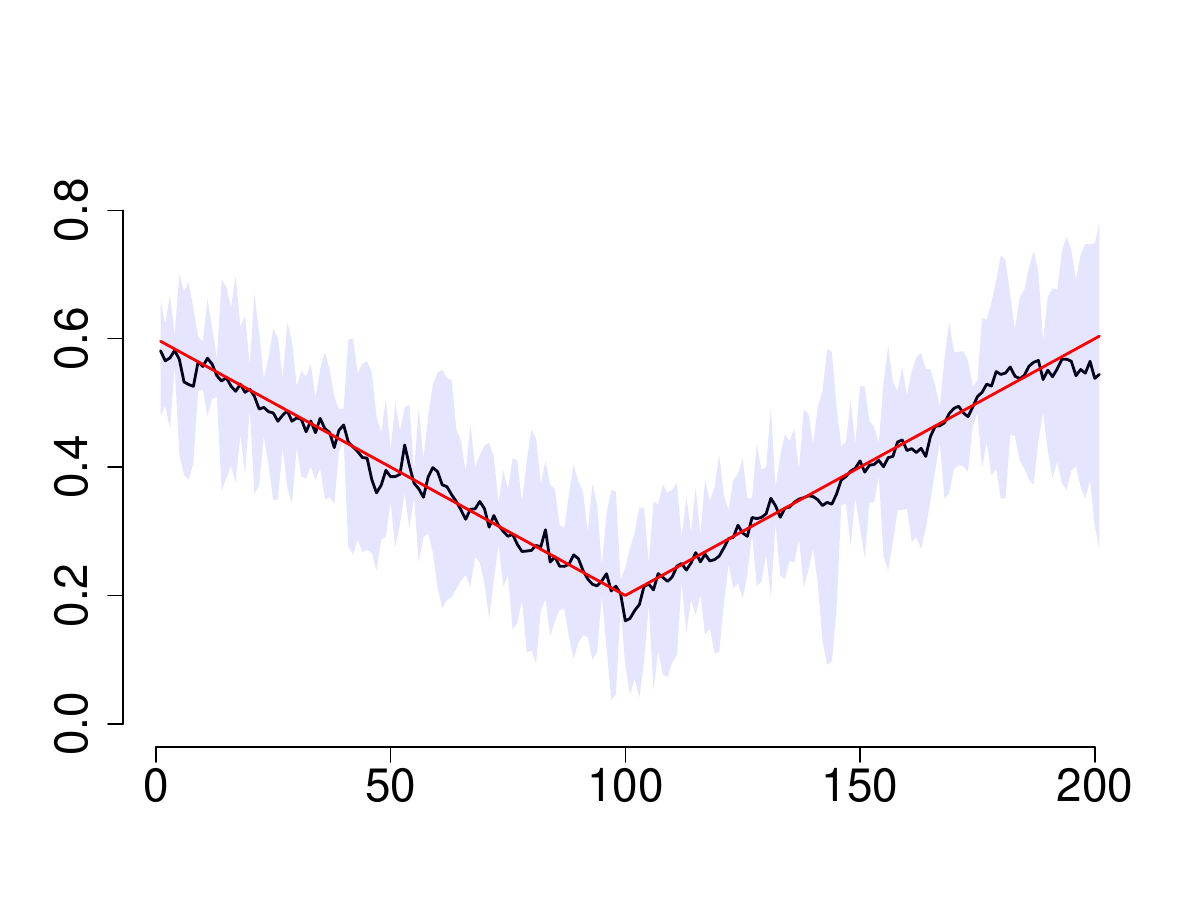}}\qquad
\subfigure[$\beta_{2,t}$]{\label{fig:SMCMC_Example_2_Reg2}
\includegraphics[trim={0 1cm 1cm 2cm},clip,width=0.45\textwidth]{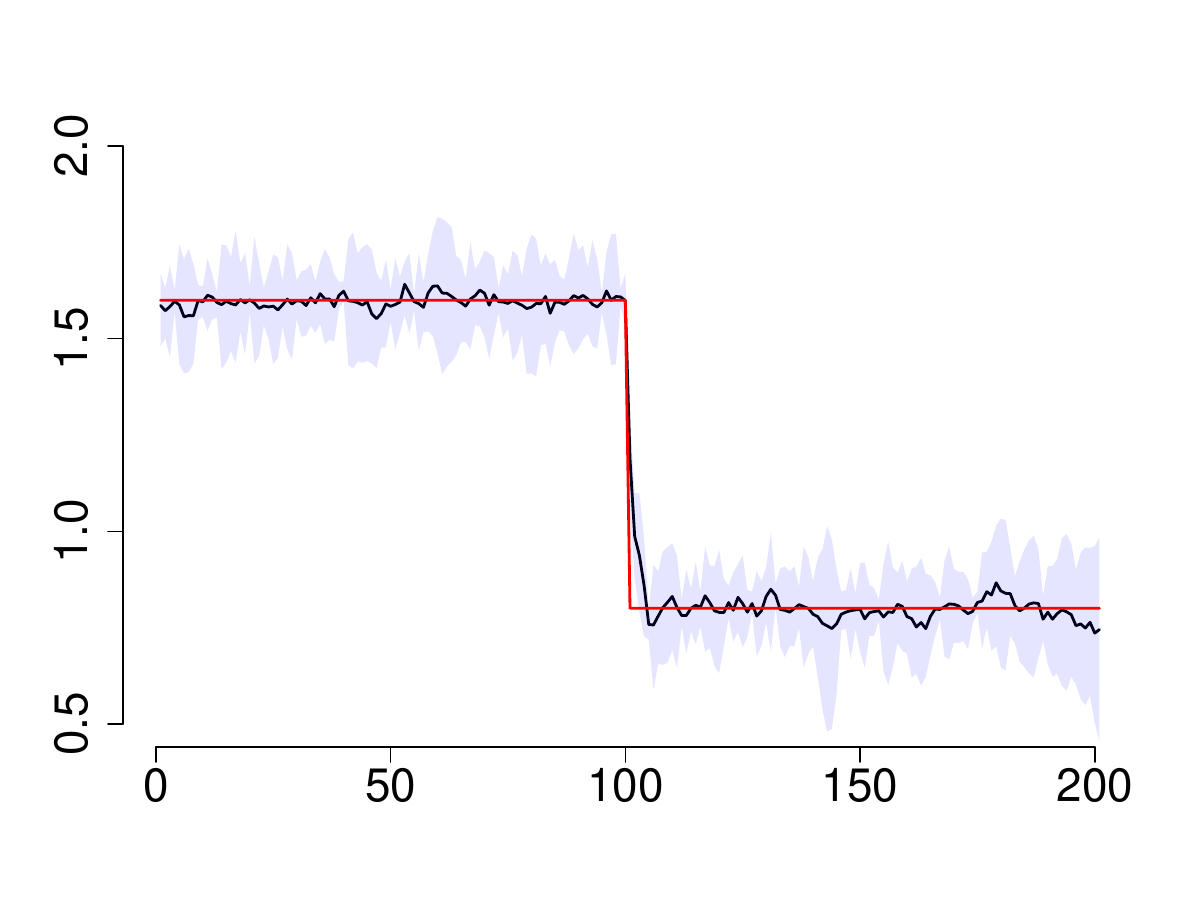}}}
\resizebox*{0.9\textwidth}{!}{\subfigure[$\beta_{3,t}$]{\label{fig:SMCMC_Example_1_Reg3}
\includegraphics[trim={0 1cm 1cm 2cm},clip,width=0.45\textwidth]{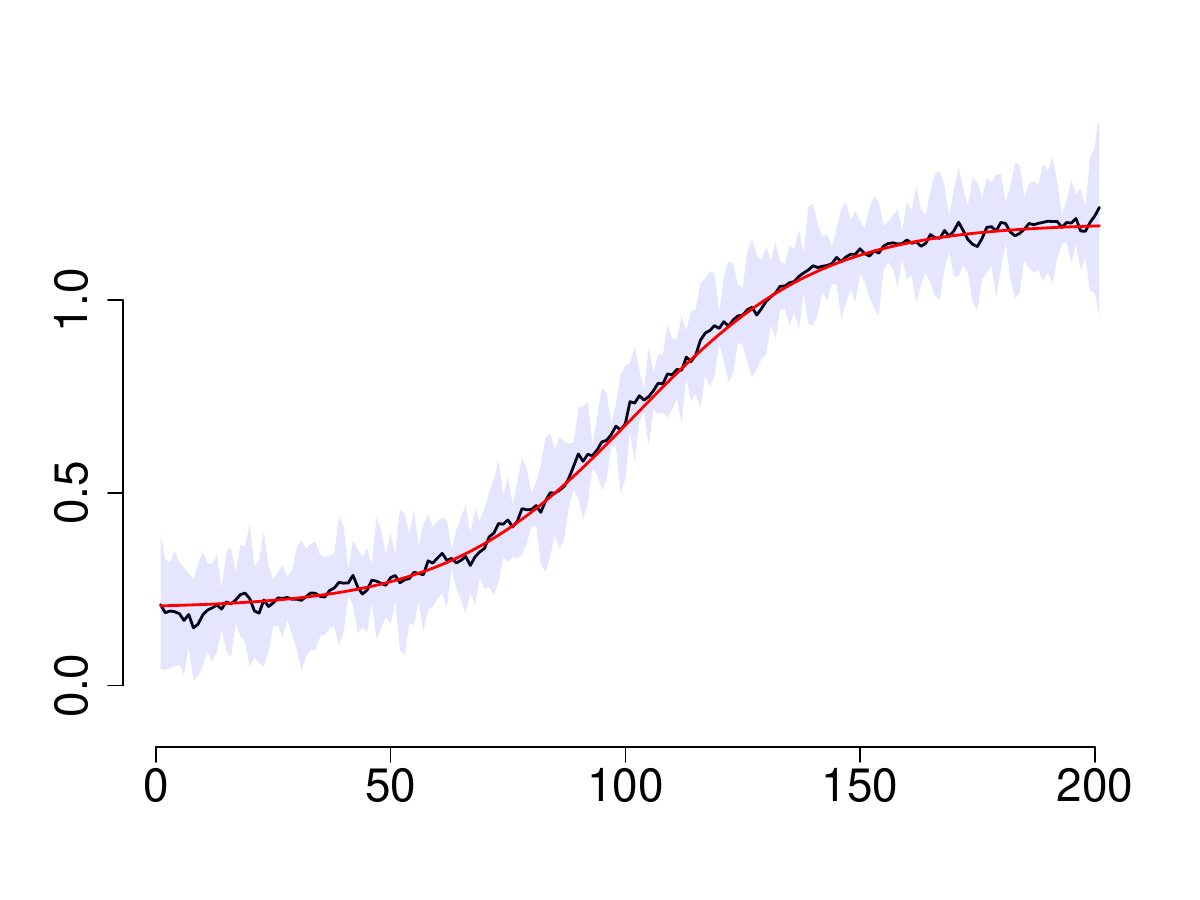}}\qquad
\subfigure[$\beta_{3,t}$]{\label{fig:SMCMC_Example_2_Reg3}
\includegraphics[trim={0 1cm 1cm 2cm},clip,width=0.45\textwidth]{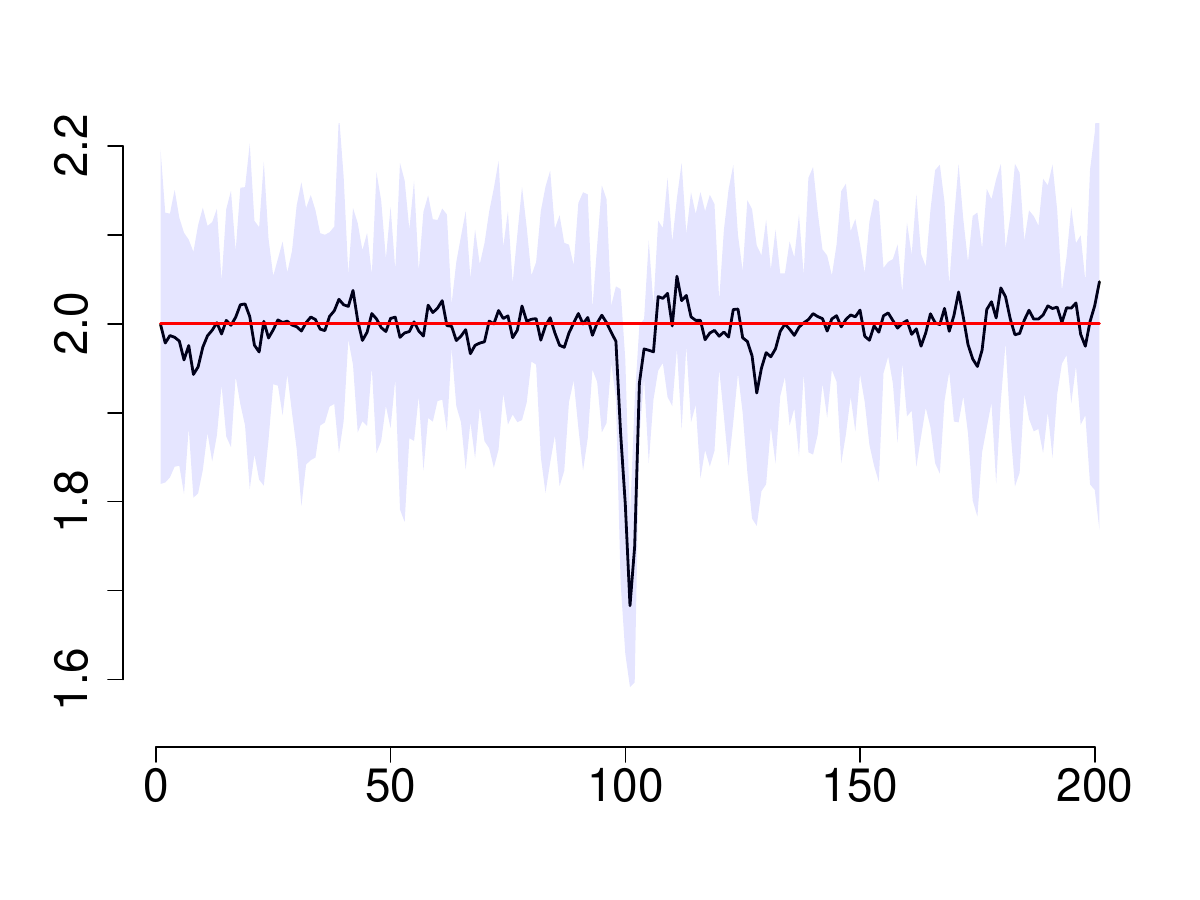}}}
\caption{Posterior mean of the regression
parameters $\boldsymbol{\beta}_{t}=\left( \protect\beta _{1,t},\protect\beta _{2,t},
\protect\beta _{3,t}\right) $, $t=1,2,\ldots ,T$, for the simulated data in Example 
\protect\ref{ex:sim_ex_1} \textit{(left panel)} and \protect\ref{ex:sim_ex_2}
\textit{(right panel)}, with quantile level $\protect\tau =0.25$ and $N=100$
parallel chains. In each plot: true parameters (red), posterior medians (dark) and $95\%$ HPD regions (gray areas).}
\label{fig:example_SMCMC_sim}
\end{center}
\end{figure}
%
\subsection{Abrupt change in quantiles}
\label{ex:sim_ex_2} 
%
\noindent We assume an abrupt change at time $t=100$ in the constant term $\beta_{1,t}^*$ and the coefficient $\beta_{2,t}^*$, and a GARCH(1,1) dynamics for the conditional volatility of the innovation term $\nu_t$:
\begin{align}
\beta _{1,t}^{\ast }& =-2\bbone_{(-\infty,100]}(t)+2\bbone_{(100,\infty)}(t)
\notag \\
\beta _{2,t}^{\ast }& =1.6\bbone_{(-\infty,100]}(t)+0.8\bbone_{(100,\infty)}(t)\notag\\
\nu _{t}^{2}& =a+b\nu _{t-1}^{2}+c\varepsilon _{t-1}^{2},  \notag
\end{align}%
for $t=1,\dots,200$, and $\beta _{3,t}^{\ast }=2,\forall t=1,2,\dots ,T$, with $a=0.05$, $b=0.9$
and $c=0.05$ and $\varepsilon_0\sim\mathsf{N}(0,\frac{a}{1-b})$. It is worth noting that the relationship between the coefficients of the DGP and the coefficients $\boldsymbol{\beta}_{t}$ of the quantile regression is as follows: $\beta_{1,t} = \beta_{1,t}^{*} + \nu_{t}\Phi^{-1}(\tau)$, $\beta_{2,t} = \beta_{2,t}^{*}$, and $\beta_{3,t} = \beta_{3,t}^{*}$, where $\tau$ represents the quantile level and $\Phi(\cdot)$ denotes the cumulative density function of a standard Normal distribution. The prior hyperparameter settings are specified as follows: $a_0 =b_0=10^{-3}$, $c_0 = M+1$, and $C_0 = 0.01\mathbf{I}_3$, where $\mathbf{I}_3$ denotes the $3$-dimensional identity matrix.\newline 
\indent In Fig. \ref{fig:example_SMCMC_sim}, we present the sequential estimates for the two simulated examples. The figure illustrates how the quantile regression model effectively captures both abrupt shifts and smooth transitions in the intercept $\beta_{1,t}$ (first row) in presence of both homoskedastic and heteroskedastic observation noise (Panel \ref{fig:SMCMC_Example_1_Reg1} and \ref{fig:SMCMC_Example_2_Reg1}, respectively). Moreover, the second and third rows demonstrate the ability of our BDQMA procedure to accurately detect both abrupt changes (Panel \ref{fig:SMCMC_Example_1_Reg2}-\ref{fig:SMCMC_Example_2_Reg2}) and smooth changes (Panel \ref{fig:SMCMC_Example_1_Reg3}) in the relationship between the dependent variable and the covariates.
%
\section{Predicting US inflation when causality changes}
\label{sec:empirical_application_inflation} 
%
%
\noindent Inflation has been a central concern for economists, as expectations around it play a critical role in shaping the decisions of economic agents and significantly impact both the economic and social development of nations \citep[][]{phillips.1958,stock_watson.1999, KumWessel2024}. Our empirical exercise serves two primary objectives: {\it 1)} to explore the relevance of covariates beyond the unemployment rate and lagged inflation in predicting current inflation at various quantile levels; and {\it 2)}  to assess whether the predictors for \qmo high\qmcsp and \qmo low\qmcsp inflation rates differ or remain consistent, as reflected by dynamically time-varying inclusion probabilities, $\pi_{t+1|t}^{(k)}$. In this study, we use a slightly modified version of the dataset presented in \cite{koop_korobilis.2012}. Detailed descriptions of the dataset and initial analyses can be found in Appendixes \ref{sec:AppData} and \ref{sec:US_infl_additional_results} of the supplementary materials.
\subsection{A Generalized Phillips curve}
\begin{figure}[p]
\begin{center}
\resizebox*{0.9\textwidth}{!}{\subfigure[$\tau=0.10$]{\label{fig:CPIAUCSL_incl_prob010}
\includegraphics[trim={0 1cm 0.5cm 1.5cm},clip,width=0.8\textwidth]{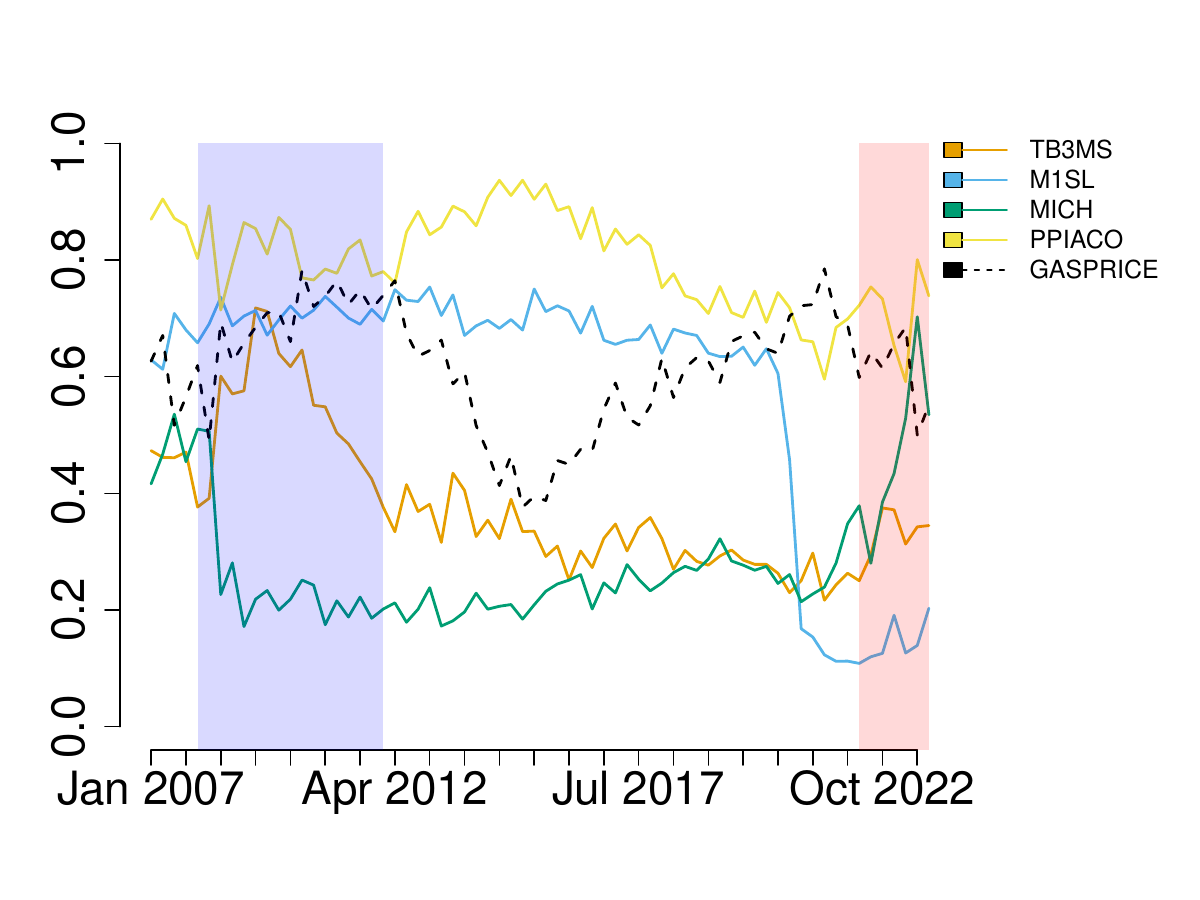}}
\subfigure[$\tau=0.25$]{\label{fig:CPIAUCSL_incl_prob025}
\includegraphics[trim={0 1cm 0.5cm 1.5cm},clip,width=0.8\textwidth]{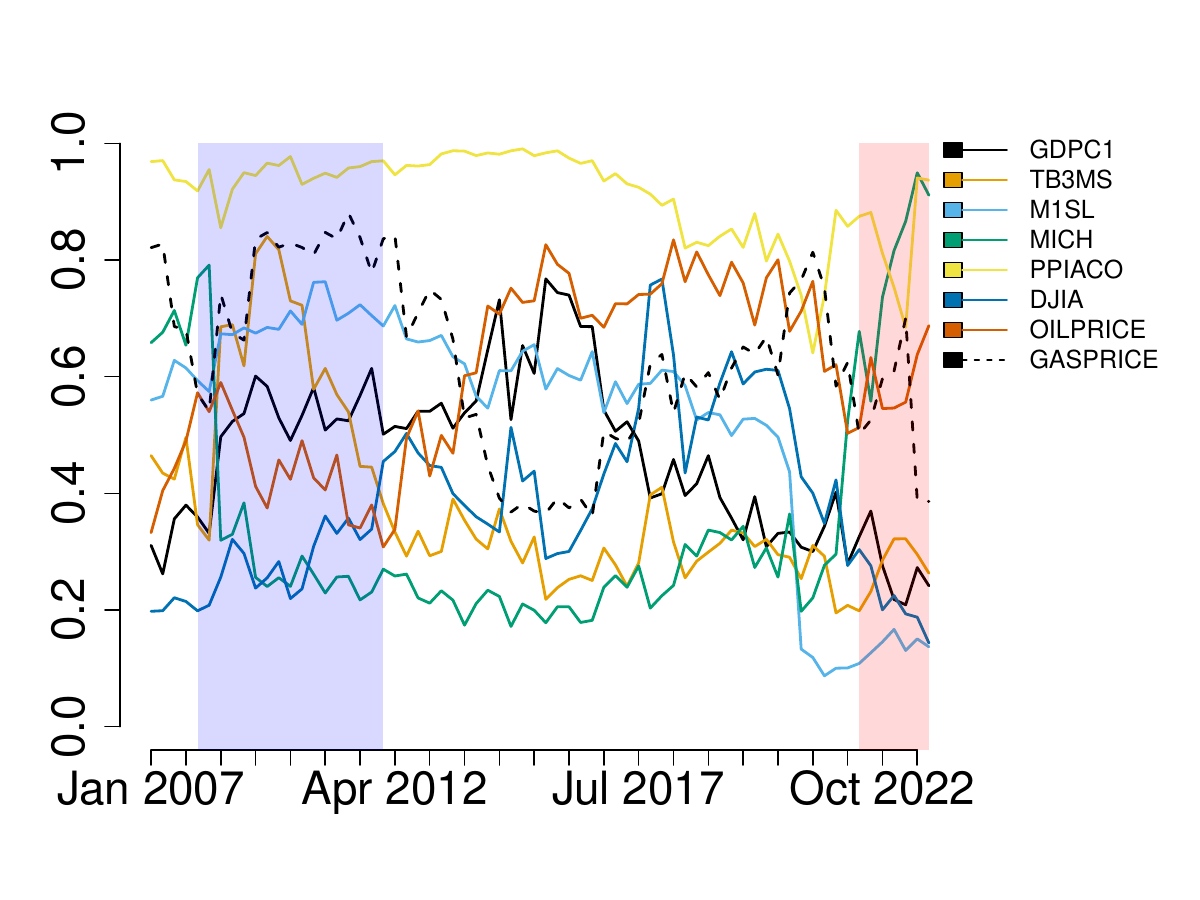}}}
\resizebox*{0.9\textwidth}{!}{\subfigure[$\tau=0.50$]{\label{fig:CPIAUCSL_incl_prob050}
\includegraphics[trim={0 1cm 0.5cm 1.5cm},clip,width=0.8\textwidth]{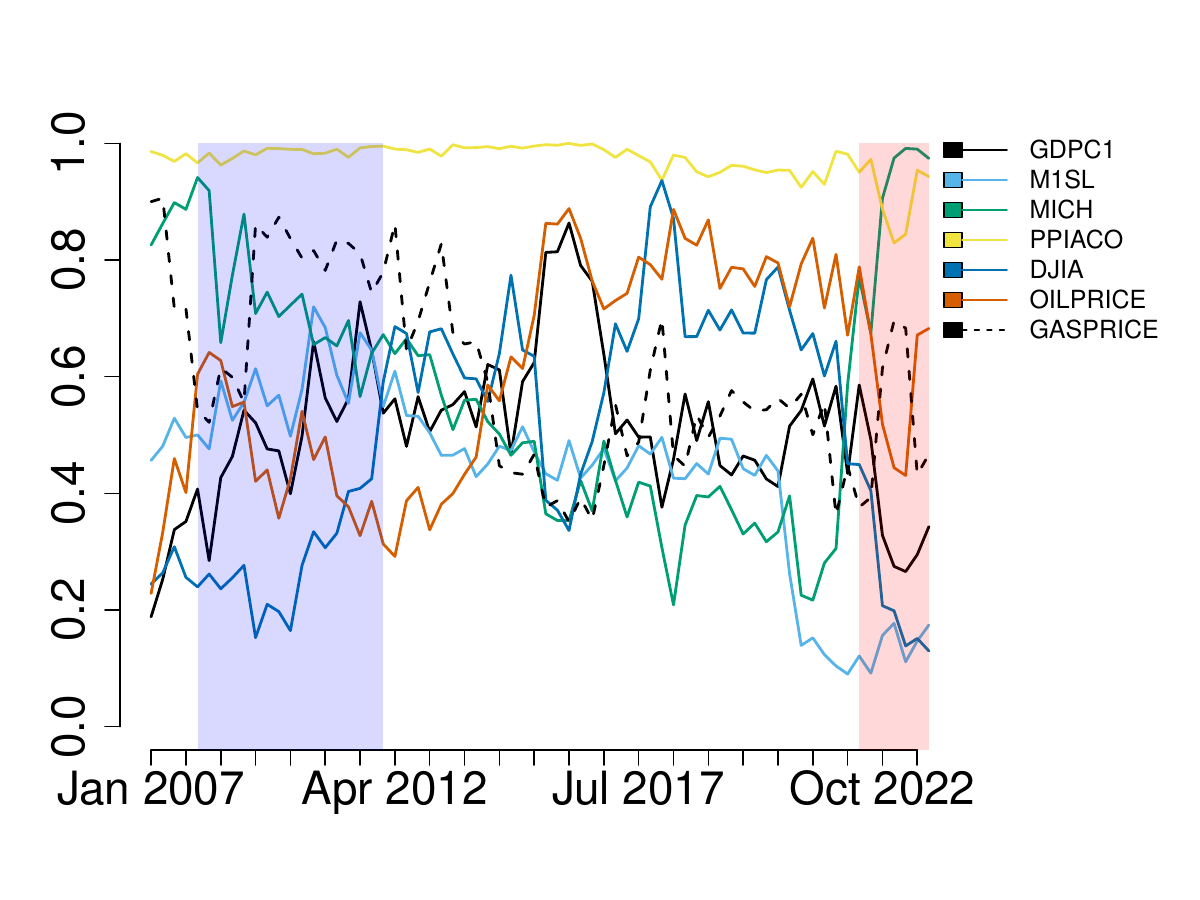}}
\subfigure[$\tau=0.75$]{\label{fig:CPIAUCSL_incl_prob075}
\includegraphics[trim={0 1cm 0.5cm 1.5cm},clip,width=0.8\textwidth]{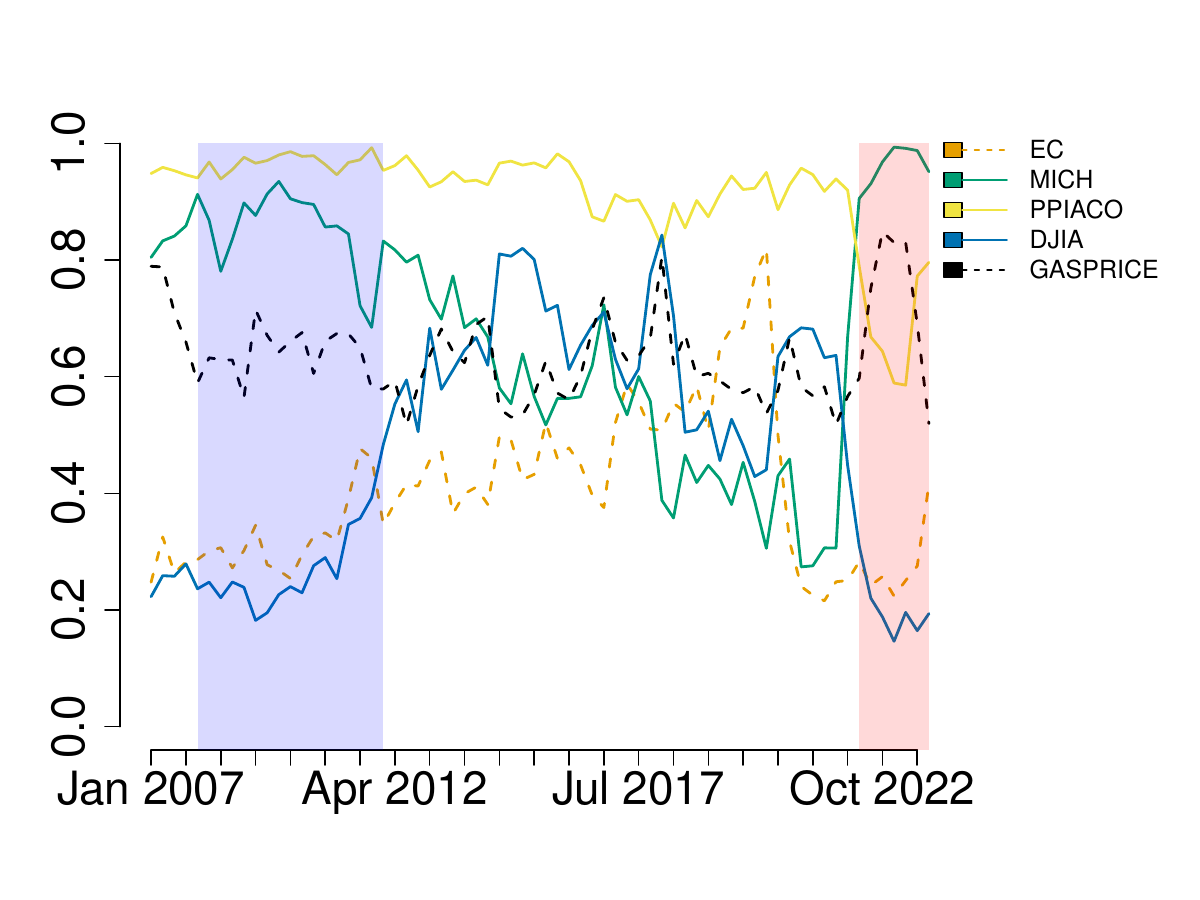}}}
\resizebox*{0.9\textwidth}{!}{\subfigure[$\tau=0.90$]{\label{fig:CPIAUCSL_incl_prob090}
\includegraphics[trim={0 1cm 0.5cm 1.5cm},clip,width=0.8\textwidth]{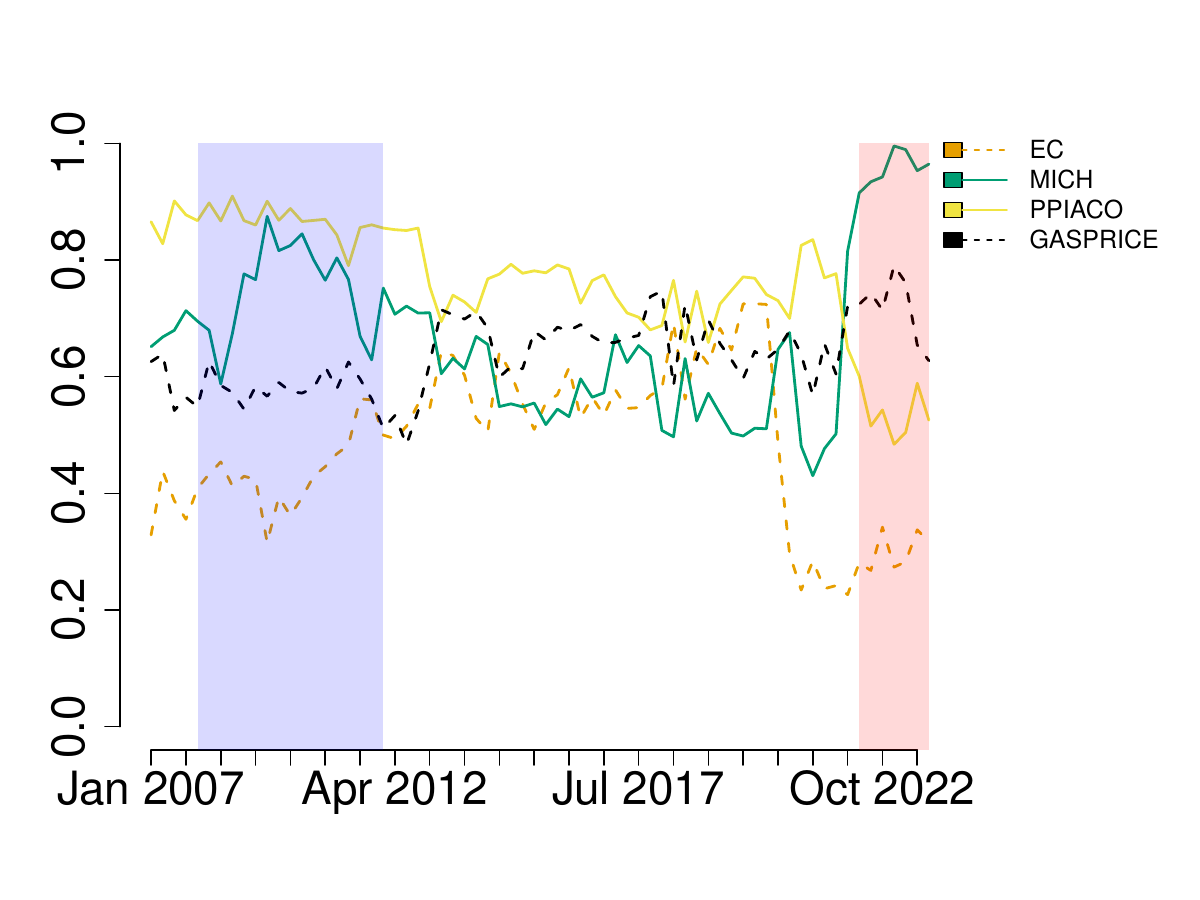}}
\subfigure[Mean regression]{\label{fig:CPIAUCSL_incl_prob}
\includegraphics[trim={0 1cm 0.5cm 1.5cm},clip,width=0.8\textwidth]{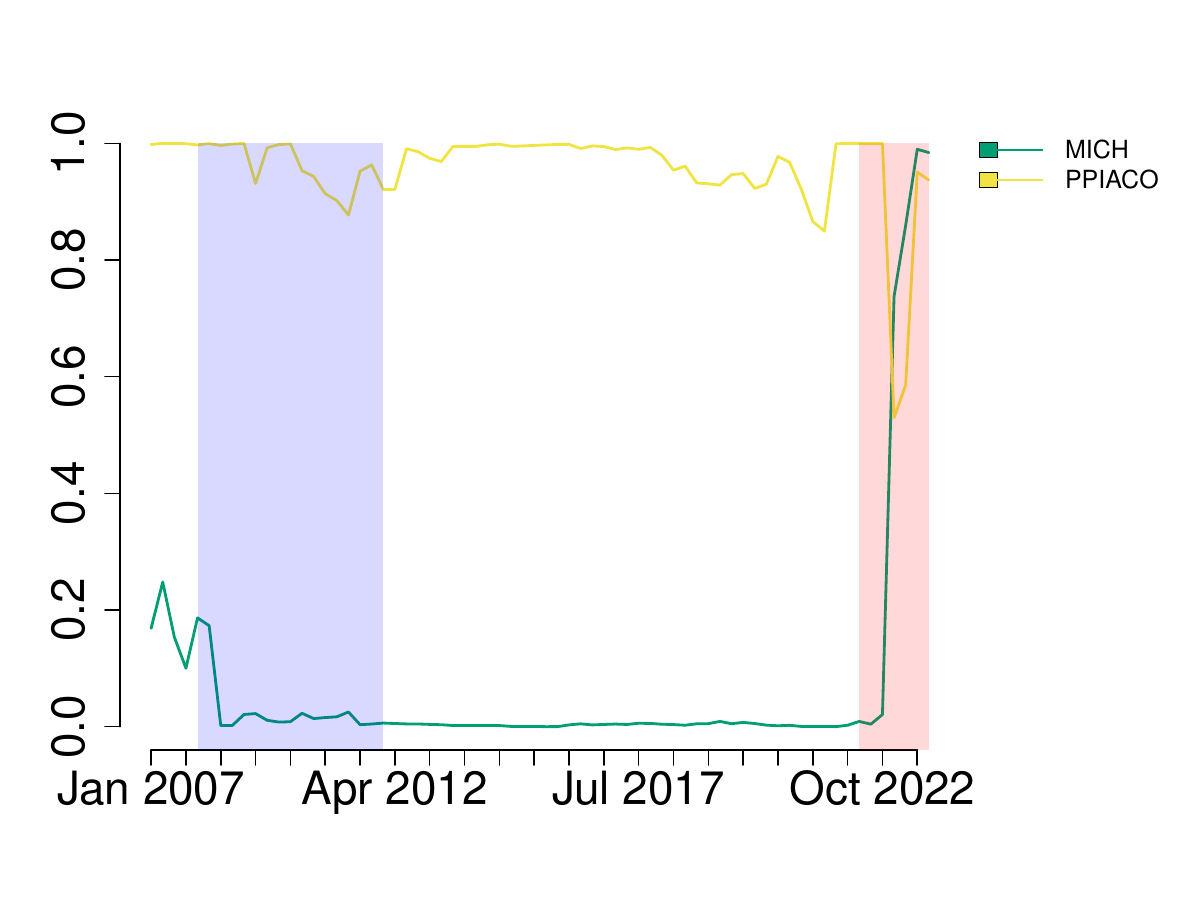}}}
\caption{US inflation dataset. Sequential update of the predicted inclusion
probabilities $\protect\pi _{t|t-1}$ by DMA (for the mean regression \ref{fig:CPIAUCSL_incl_prob}) and BDQMA (for the quantile regression \ref{fig:CPIAUCSL_incl_prob010}-\ref{fig:CPIAUCSL_incl_prob090}) for the CPIAUCSL (Consumer Price Index for All Urban Consumers: All Items in U.S. City Average). For each quantile level $\tau$ the corresponding figure only reports those parameters having inclusion probability larger or equal to $0.7$ for at least one quarter. See Table \ref{tab:table_US_inflation_data_CPIAUCSL_ss} for a summary of the relevant covariate. The shaded areas identify two significant periods: the US Great Financial Crisis from 2007-Q4 to 2011-Q4 {\it (blue)}, and the Russian-Ukraine crisis {\it (red)}.}
\label{fig:inflation_InclProb}
\end{center}
\end{figure}

\noindent  We employ our BDQMA approach to select the best subset of predictors for forecasting US inflation within a generalized Phillips curve framework. Our BDQMA approach is well-suited to address these issues due to its ability to consider dynamically evolving linear relationships between covariates and quantiles of the explained variable. Additionally, it can capture and reproduce structural breaks commonly observed in the evolution of economic variables such as inflation \citep[see also][]{primiceri.2005,koop_onorante.2012,stock_watson.2007}. Furthermore, focusing on quantiles of the predicted variable aids in identifying periods characterized by different economic conditions, particularly those featuring low or high inflation levels. 
Building on the conditionally Gaussian representation of the BDQMA model, and drawing inspiration from the seminal works on inflation by \cite{engle.1982} and \cite{bollerslev.1986}, we incorporate a conditional heteroskedastic volatility error term. This modification greatly enhances the model's flexibility in capturing asymmetric conditional distributions.\newline
\indent Regarding the model specification, we extend the autoregressive model of order $p$ with exogenous covariates (ARX$(p)$) developed by \cite{stock_watson.1999}, and
previously considered by \cite{koop_korobilis.2012}, within our BDQMA framework. The quantile model is therefore formulated as follows:
\begin{equation}
q_{\tau }\left(\mathbf{x}_{t},\boldsymbol{\beta}\right) =\mathbf{z}_{t-1}^{\top}\boldsymbol{\psi}
+\sum_{j=1}^{p}\phi _{j}y_{t-j},  \label{eq:quantile_model}
\end{equation}
where $\boldsymbol{\beta}=(\boldsymbol{\psi}^\top,\phi_1,\dots,\phi_p)^\top$ is the parameter vector, $y_{t}$ represents inflation, defined as the percentage logarithmic change of the price index $P_{t}$, that is: $y_{t}=100\log (P_{t}/P_{t-1})$. The covariate set $\mathbf{x}_{t}$ includes lagged predictors $\mathbf{z}_{t-1}$, and lagged inflation, $\left(y_{t-1},\ldots ,y_{t-p}\right)$.\newline
%
\begin{figure}[p]
\begin{center}
\resizebox*{0.9\textwidth}{!}{\subfigure[$\tau=0.10$]{\label{fig:CPIAUCSL_beta010_relevant}
\includegraphics[trim={0 1cm 0.5cm 1.5cm},clip,width=0.8\textwidth]{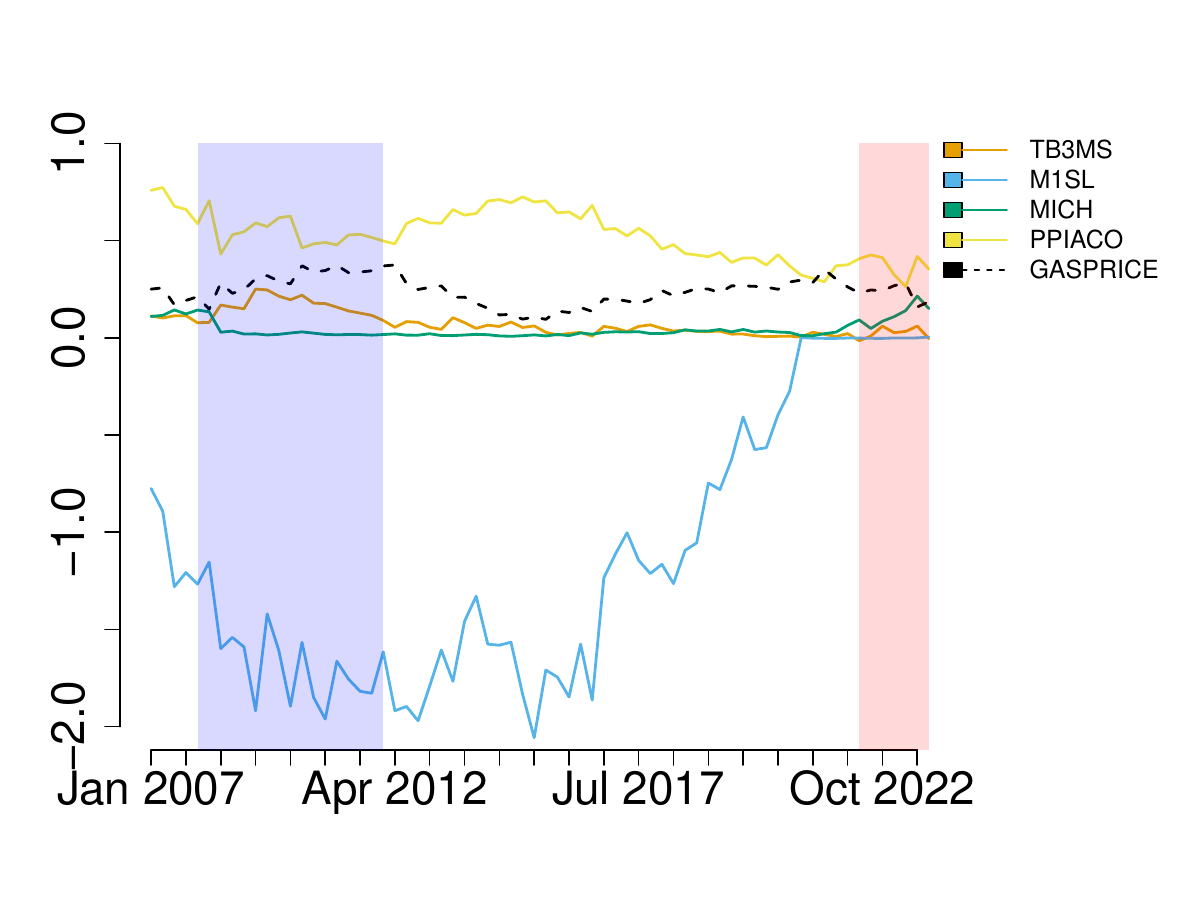}}
\subfigure[$\tau=0.25$]{\label{fig:CPIAUCSL_beta025_relevant}
\includegraphics[trim={0 1cm 0.5cm 1.5cm},clip,width=0.8\textwidth]{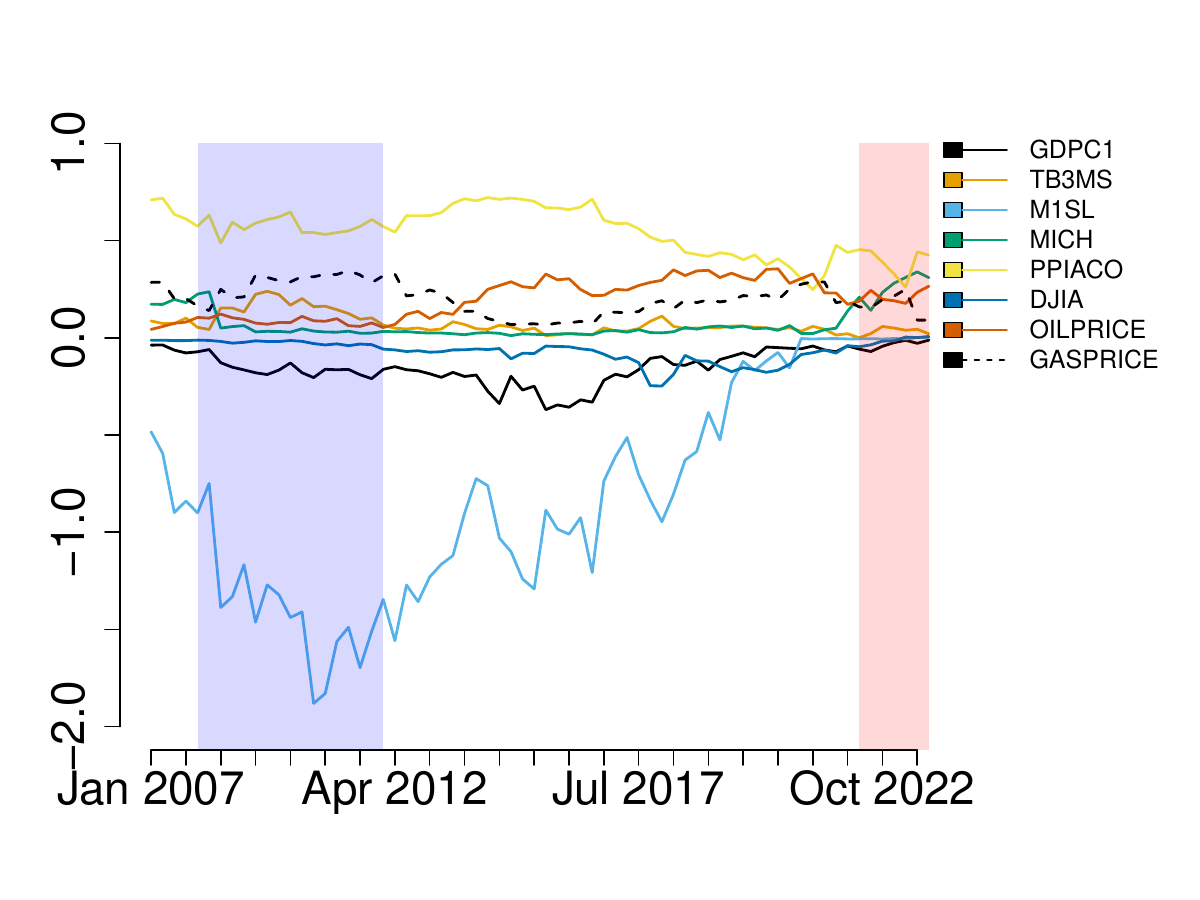}}}
\resizebox*{0.9\textwidth}{!}{\subfigure[$\tau=0.50$]{\label{fig:CPIAUCSL_beta050_relevant}
\includegraphics[trim={0 1cm 0.5cm 1.5cm},clip,width=0.8\textwidth]{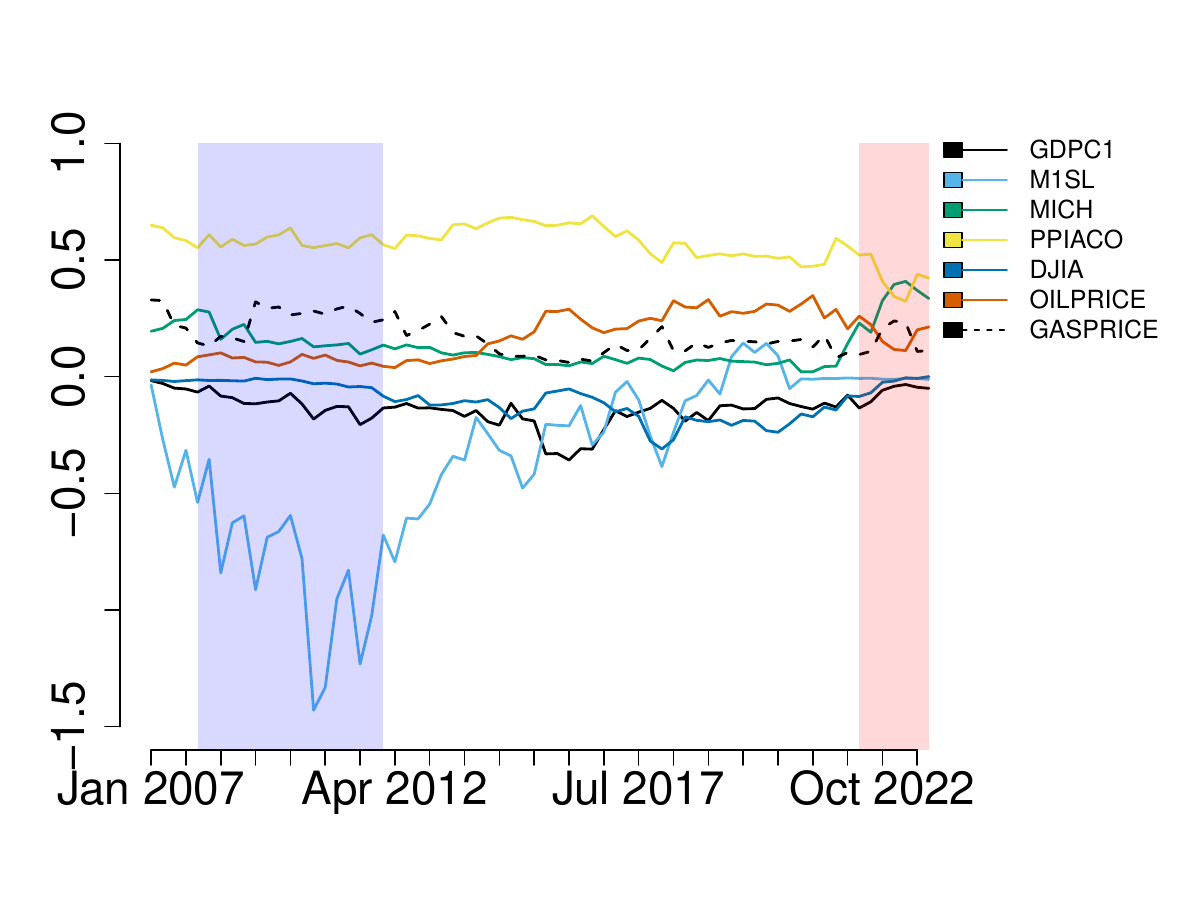}}
\subfigure[$\tau=0.75$]{\label{fig:CPIAUCSL_beta075_relevant}
\includegraphics[trim={0 1cm 0.5cm 1.5cm},clip,width=0.8\textwidth]{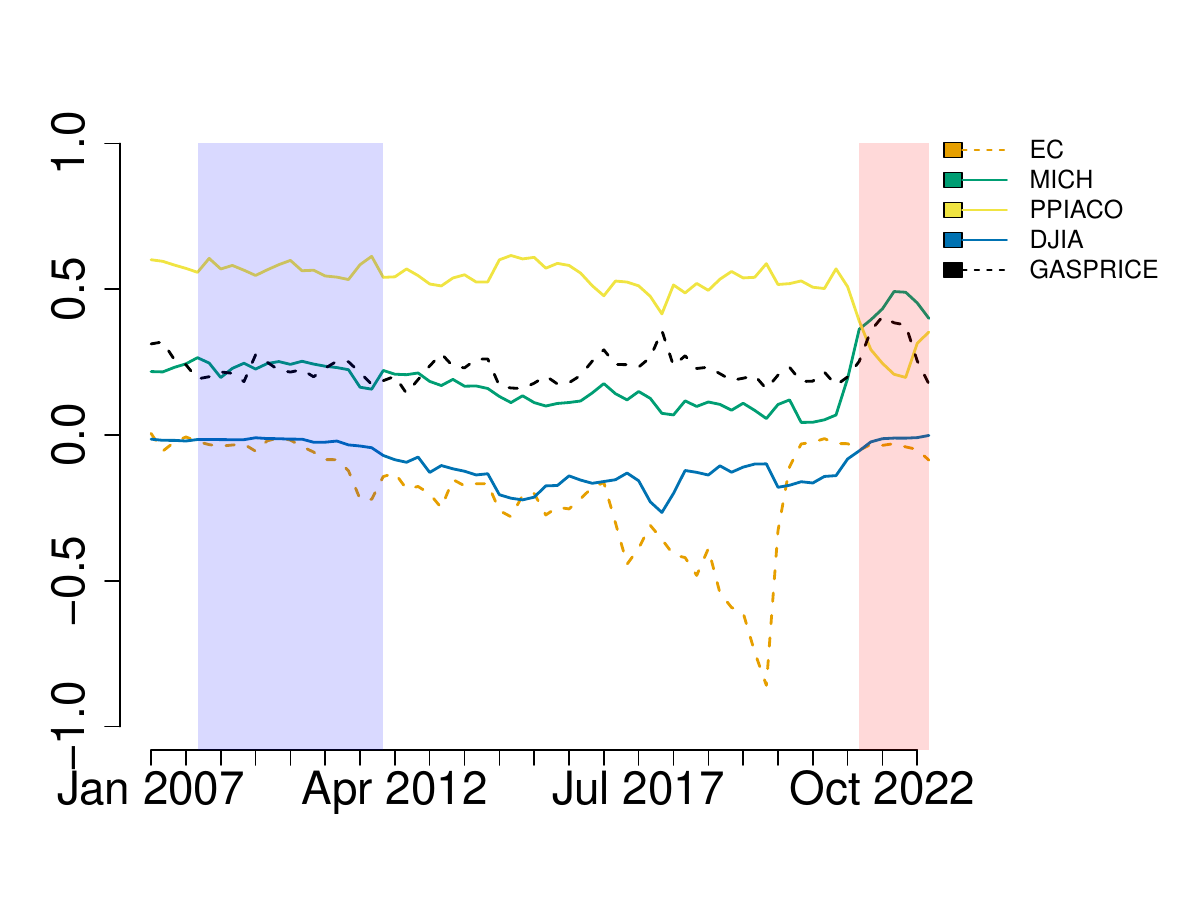}}}
\resizebox*{0.9\textwidth}{!}{\subfigure[$\tau=0.90$]{\label{fig:CPIAUCSL_beta090_relevant}
\includegraphics[trim={0 1cm 0.5cm 1.5cm},clip,width=0.8\textwidth]{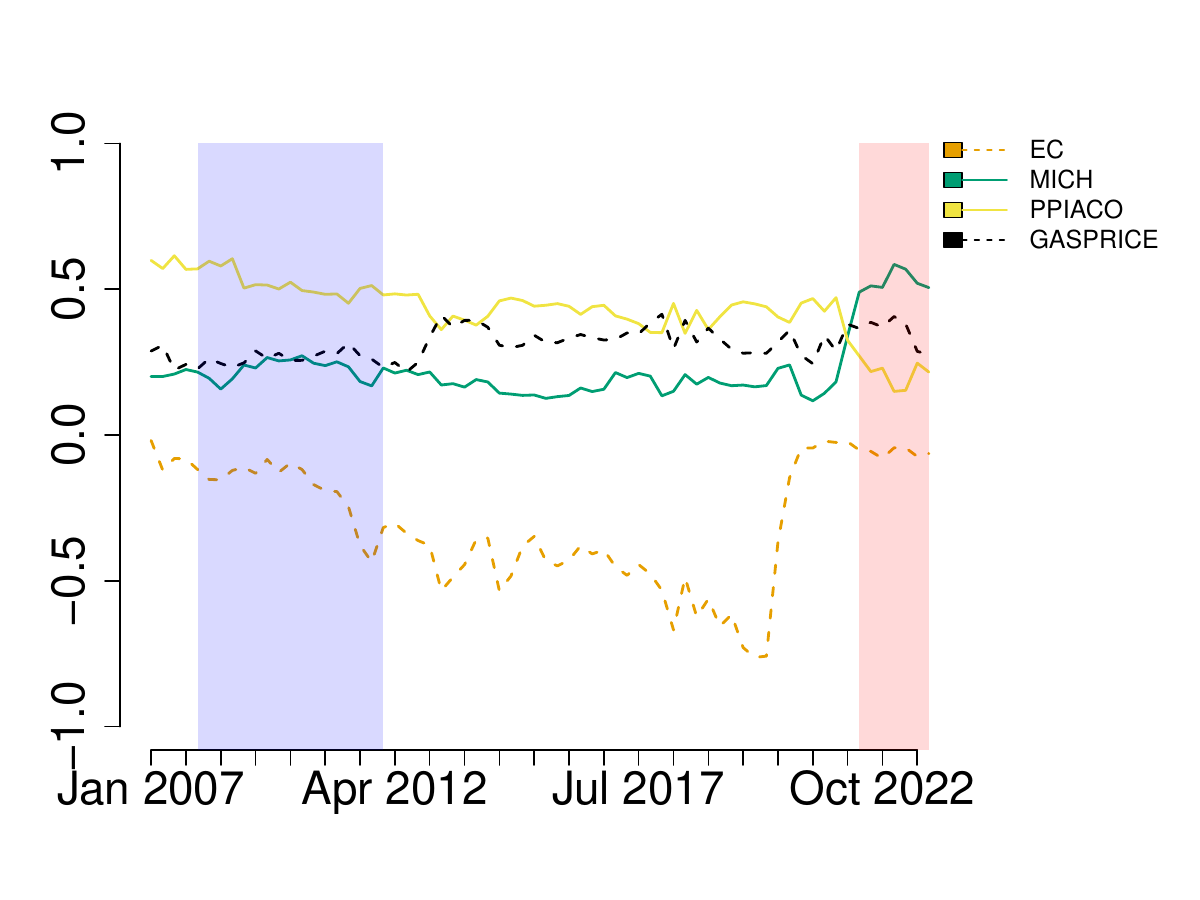}}
\subfigure[Mean regression]{\label{fig:CPIAUCSL_beta_relevant}
\includegraphics[trim={0 1cm 0.5cm 1.5cm},clip,width=0.8\textwidth]{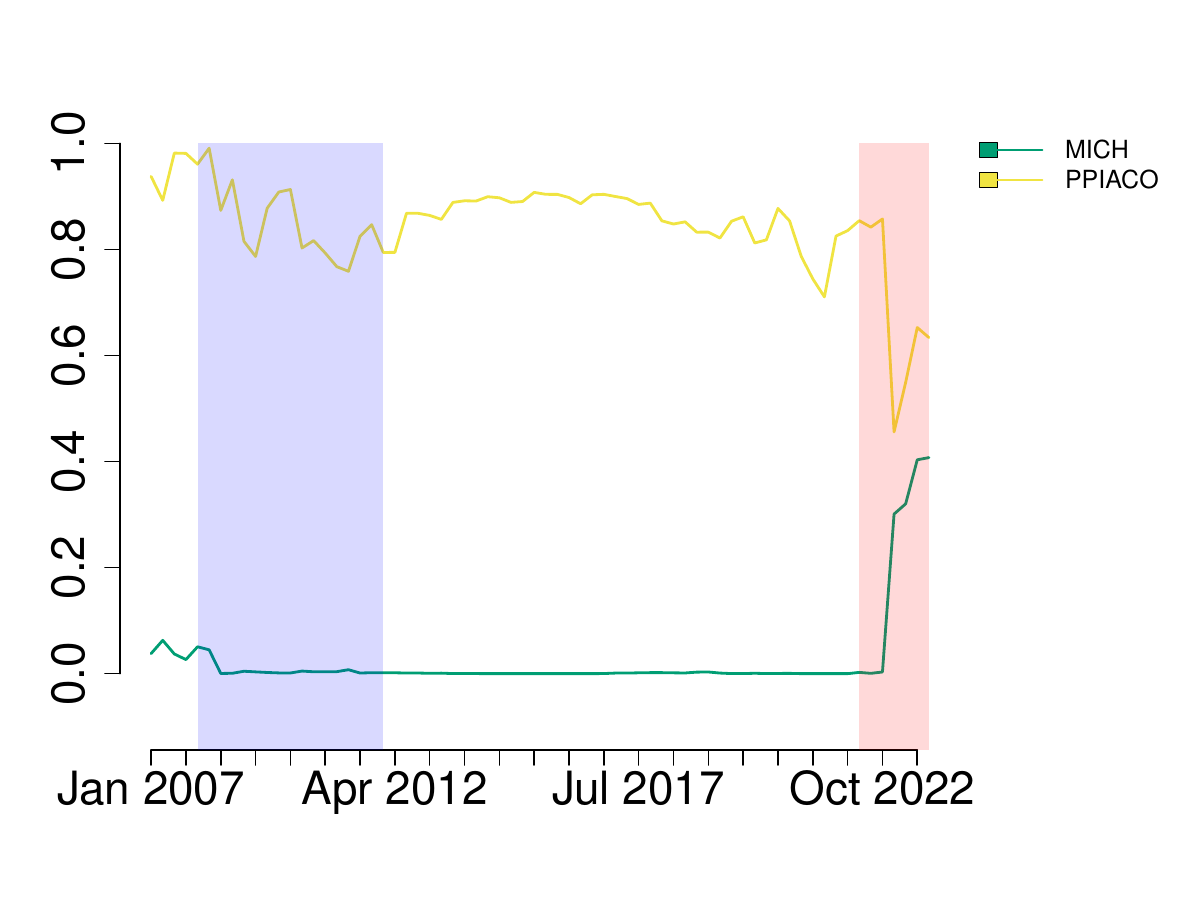}}}
\caption{US inflation data. Sequential estimates of the regression parameters by DMA (for the mean regression \ref{fig:CPIAUCSL_beta_relevant}) and BDQMA (for the quantile regression \ref{fig:CPIAUCSL_beta010_relevant}-\ref{fig:CPIAUCSL_beta090_relevant}) for the CPIAUCSL (Consumer Price Index for All Urban Consumers: All Items in U.S. City Average). For each quantile level $\tau$ the corresponding figure only reports those parameters having inclusion probability larger or equal to $0.7$ for at least one quarter. See Table \ref{tab:table_US_inflation_data_CPIAUCSL_ss} for a summary of the relevant covariate. The shaded areas identify two significant periods: the US Great Financial Crisis from 2007-Q4 to 2011-Q4 {\it (blue)}, and the Russian-Ukraine crisis {\it (red)}.}
\label{fig:inflation_RegDyn}
\end{center}
\end{figure}
\indent Drawing from standard literature on quantile regression analysis, we concentrate on five quantile levels: $\tau =0.10$ and $\tau =0.25$ for lower quantile analyses, $\tau =0.5$ for median analysis, and $\tau =0.75$ and $\tau =0.90$ for upper quantile analyses. These levels correspond to five distinct inflation thresholds at each time step. Due to space constraints, we present in Fig. \ref{fig:inflation_RegDyn}-\ref{fig:inflation_InclProb} the time-varying parameters and inclusion probabilities derived from our BDQMA model for all the quantile levels $\tau=\left(0.10,0.25,0.5,0.75,0.90\right)$, focusing solely on the response variable CPIAUCSL. 
For comparison, we also include Fig. \ref{fig:CPIAUCSL_beta_relevant} and Fig. \ref{fig:CPIAUCSL_incl_prob}, which display the results for the Gaussian version of the BDQMA model. This Gaussian counterpart follows the Dynamic Model Averaging (DMA) approach introduced by \cite{raftery_etal.2010}, while incorporating the findings of \cite{koop_korobilis.2013} and \cite{koop_onorante.2012}, adapted to the updated dataset.
Similar results have been provided for the response variable CPILFESL are provided in Appendix \ref{sec:US_infl_additional_results_CPILFESL} of the supplementary materials (Fig. \ref{fig:inflation_RegDyn_CPILFESL} and \ref{fig:inflation_InclProb_CPILFESL}). 
To maintain readability, we include in the figures only those variables relevant for predicting the different conditional quantiles, i.e.  those variables with inclusion probabilities exceeding 0.7. For comparison, we present a summary of the coefficient values obtained from a static mean and quantile regressions over the entire sample in Tab. \ref{tab:inflation_posterior_summary}--\ref{tab:US_inflation_mean_regression} and Tab. \ref{tab:US_inflation_quantile_regression_010}--\ref{tab:US_inflation_quantile_regression_090} of the Appendix \ref{sec:US_infl_additional_results} of the supplementary materials accompanying the paper.\newline 
%
\indent The median analysis provides a robust version of the time-varying regression analysis pioneered by \cite{koop_korobilis.2013} and \cite{koop_onorante.2012}. Consistent with their findings, our results (see Fig. \ref{fig:inflation_RegDyn}-\ref{fig:inflation_InclProb}, $\tau=0.5$) show that coefficients of the regressors and their probabilities of inclusion in the model change over time. Moreover, the level of inflation persistence, indicated by the sum of the autoregressive coefficients (not reported), is approximately 0.5, consistent with the findings of \cite{koop_korobilis.2013}. The primary disparity between mean regression and its robust related median version, extends beyond the quantity of pertinent regressors to elucidate inflation; it also encompasses their efficacy. Indeed, beyond the two regressors shared with the mean regression, i.e. the changes in the NAPM commodities price index (PPIACO) and expectations on future inflation (MICH), the significance of the remaining variables fluctuates over time in response to prevailing economic conditions. For instance, the prominence of oil and gas prices in explaining heightened inflation levels during the recent two years, catalyzed by the Russian-Ukraine crisis, contrasts with their lesser relevance during the Great Financial Crisis. See Appendix \ref{sec:US_infl_additional_results} for a detailed description of the results.\newline
%
%
\indent The analysis of lower and upper quantiles reveals significant variations in the coefficients of the regressors, which bear important policy implications. Specifically, in the static quantile regression, we observe that the autoregressive terms $\left(\phi _{1},\phi _{2}\right)$ (not reported) increase as we move towards higher quantile levels. For instance, at low inflation levels (e.g. $\tau =0.10$, see Tab. \ref{tab:inflation_posterior_summary}), the persistence is approximately $0.24$, whereas for high inflation levels (e.g. $\tau=0.90$, see Tab. \ref{tab:inflation_posterior_summary}), we find evidence of higher persistence, reaching around $0.61$. The dynamic regression reveals a similar pattern toward the end of the period (not reported), indicating greater inflation persistence at higher quantile levels. Notably, during the two inflation peaks in 1975 and 1981, the parameter $\phi_{2}$ follows a predominantly negative trajectory, suggesting a reduction in inflation persistence. \newline
%
%
\indent The imperative of adopting a time-varying parameters analysis becomes abundantly clear upon scrutiny of Fig. \ref{fig:inflation_InclProb}. Consider, for instance, if we were to exclusively rely on the findings of a static regression model fitted across the entire sample, meticulously detailed in Tab. \ref{tab:inflation_posterior_summary}. Such an approach might lead us to erroneously conclude that six variables out of the total 18 (EC, HOUST, USPRIV, GS10, T10Y3MM, M1SL, and PMI) hold no relevance across any inflation level at any given point in time (refer to Tab. \ref{tab:inflation_posterior_summary}). While this assertion largely holds true for most scenarios, even under the BDQMA framework, which consistently excludes HOUST, USPRIV, GS10, T10Y3MM, and PMI across all quantile levels and time spans (as gleaned from the inclusion probabilities summary statistics in Tab. \ref{tab:table_US_inflation_data_CPIAUCSL_ss}), the same cannot be said for EC and M1SL. Moreover, it becomes evident that money supply (M1SL) and real personal consumption expenditures (EC) emerge as highly informative for quantiles below the median and the median itself, respectively, and for the quantile exceeding the median. In contrast, variables such as UNRATE, PRFI, T10YFMM, and NAPMSDI, though included in the static regression for at least one quantile, are consistently deemed irrelevant by BDQMA. While the explanatory power of inflation by the NAPM vendor deliveries index remains contentious in the literature, economic theory suggests that the remaining variables---encompassing investment, the spread of long-term government bonds, the federal funds rate, and the unemployment rate---are pivotal determinants of inflation rates and are typically integrated into the Phillips curve. To shed light on this discrepancy, it is worth noting that for the static quantile regression, UNRATE and PRFI are relevant solely for the lower and two lower quantiles, respectively. However, the interest rate on US short-term federal funds (TB3MS) is included in both the static quantile regression and BDQMA for quantiles ranging from $0.1$ to $0.5$, thereby tempering the impact of other variables associated with federal funds interest rates. This nuanced analysis helps to unravel the complexity surrounding these variables. While some exhibit sporadic relevance across different regression techniques and quantiles, others, like the federal funds rate, consistently play a moderating role, highlighting the intricate interplay of economic factors influencing inflation dynamics. To conclude our exploration of the relevance of the proposed methodology in understanding the determinants of inflation dynamics, let us juxtapose the results obtained from static and dynamic quantiles for three pivotal variables: the growth rate of real GDP (GDPC1) and energy prices (OILPRICE and GASPRICE). or a detailed discussion of the static quantile regression results, please refer to Appendix \ref{sec:US_infl_additional_results_static}.
These covariates, as dictated by economic theory, are overwhelmingly relevant and are certainly included in both static and dynamic quantile estimates, albeit for varying quantile levels. However, BDQMA outperforms its static counterpart by a considerable margin in terms of efficacy and accuracy of results. BDQMA not only offers insights into the presence of a relationship between the covariate and the response variable but also illuminates the period during which this relationship is most pronounced. For instance, in the case of energy prices, BDQMA reveals that GDP is influential during periods of low inflation levels, particularly for quantiles $\tau = 0.25$ and $\tau = 0.5$, with the relationship being strongest in the middle of the sample period. Conversely, for energy prices, the association is more robust during the latter part of the sampling period (see Fig. \ref{fig:inflation_InclProb}).\newline
\indent Regarding the sign and magnitude of the impact of the most relevant regressors, as measured by the estimated coefficients in Fig. \ref{fig:inflation_RegDyn}, the dynamic regression largely corroborates the findings of the static regression for variables such as MICH, PPIACO, and energy prices. However, notable deviations are observed, particularly in the case of TB3MS, where the different impact is clearly evident when examining the results in Tab. \ref{tab:inflation_posterior_summary}. Furthermore, a crucial variable like money supply (M1SL), which was not included in the static regression, displays a consistent negative sign across all periods and quantile levels equal to or below the median. Notably, its impact exhibits an intriguing U shape, underscoring the pronounced negative effect during the Great Financial Crisis. The dynamic regression results also reveals a larger number of significant covariates for most or part of the period. More specifically, some covariates, such as M1SL and TB3SL for instance, are relevant only for the first and last quartile regressions.\newline
%
\begin{figure}[t]
\begin{center}
\includegraphics[trim={0 1cm 0 1cm},clip,width=0.8\textwidth]{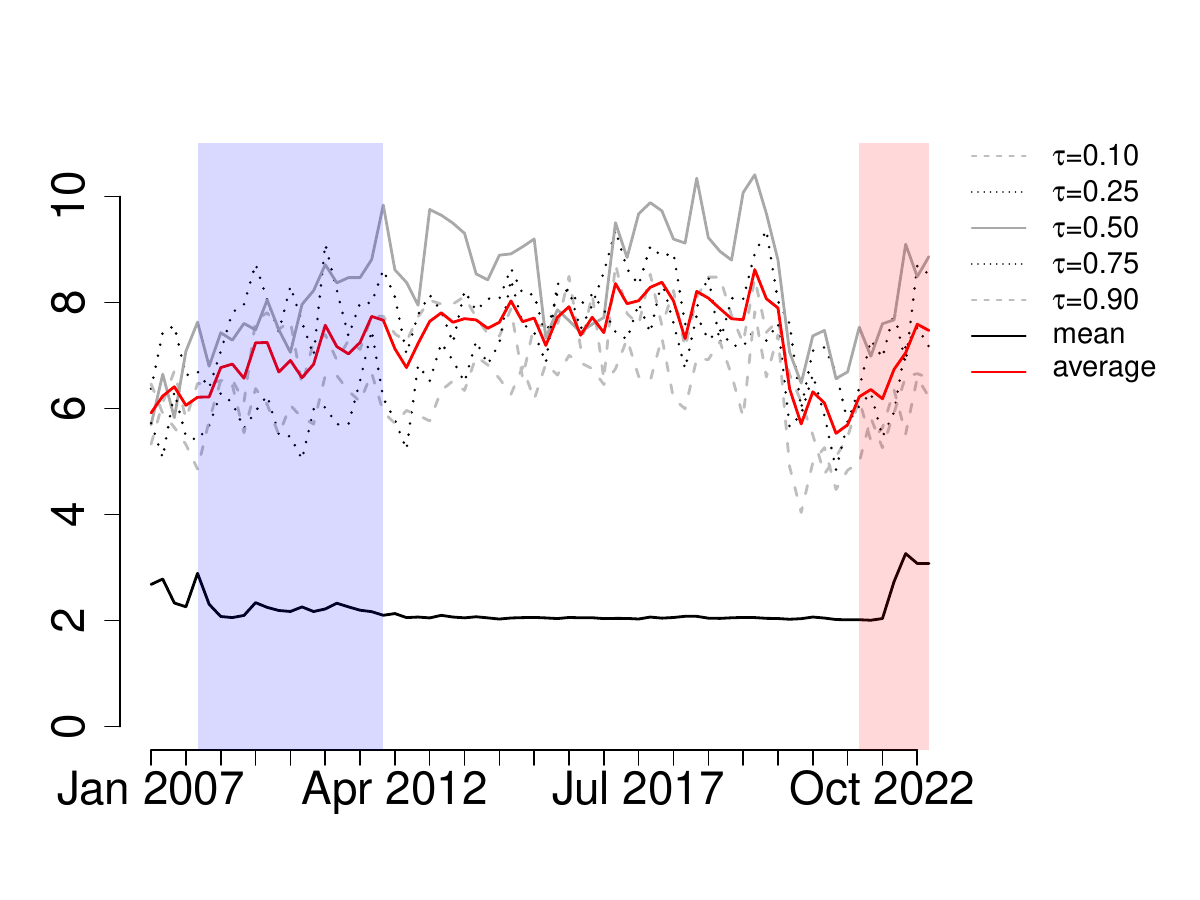} 
\end{center}
\caption{US Inflation dataset. Expected number of predictors over time for each quantile
confidence levels. The quantiles are denoted as follows: $\tau=0.10$ {\it (gray dashed line)}, $\tau=0.25$ {\it (black dotted line)}, $\tau=0.50$ {\it (gray solid line)}, $\tau=0.75$ {\it (black dotted line)}, $\tau=0.90$ {\it (gray dashed line)}, and the mean {\it (black solid line)}. The mean (Gaussian regression) is depicted by the black solid line. Additionally, the red line indicates the average number of predictors across all quantiles. The shaded areas identify two significant periods: the US Great Financial Crisis from 2007-Q4 to 2011-Q4 {\it (blue)}, and the Russian-Ukraine crisis {\it (red)}.}
\label{fig:inflation_ExpNumbReg}
\end{figure}
%
\subsection{Inclusion probabilities}
\noindent One of the significant advantages of BDQMA over traditional static model selection methods lies in its dynamic updating of posterior inclusion probabilities $\pi_{t|t-1}$. However, the superiority of BDQMA extends beyond dynamic model selection. Effective model selection procedures are also evaluated based on their ability to accurately describe the characteristics of the response variable while maintaining simplicity. When model selection is conducted for forecasting purposes, simplicity translates into avoiding overfitting and typically results in smaller forecasted confidence intervals. As suggested by \cite{koop_korobilis.2012} for dynamic model averaging (DMA), if the approach prioritizes models with fewer predictors, it can preserve out-of-sample forecasting performance without compromising goodness-of-fit properties. This emphasizes the importance of parsimony in model selection for forecasting, where simpler models often yield more reliable predictions. The expected number of predictors
included by the BDQMA model selection procedure can be analytically evaluated
using the predicted inclusion probabilities $\pi _{t\vert t-1}^{k}$, for $
t=1,2,\dots ,T$ and $k=1,2,\dots ,K,$ such as: 
\begin{equation}
\mathsf{E}\left( S_{t}^{\tau }\right) =\sum_{k=1}^{K}\pi _{t\vert
t-1}^{(k)}S_{t}^{\tau ,k}.
\end{equation}
The expected number of regressors $\mathsf{E}\left( S_{t}\right)$, where $S_{t}^{k}$ denotes the number of regressors included in model $k=1,2,\dots ,K$ at each point in time $t=1,2,\dots ,T$, provides insight into the average number of predictors included by BDQMA at time $t$. Fig. \ref{fig:inflation_ExpNumbReg} illustrates $\mathsf{E}\left( S_{t}\right)$ for all the quantile confidence levels $\tau =\left(0.10,0.25,0.50,0.75,0.90\right)$. It is important to note that the average number of predictors takes into account both significant and irrelevant regressors, while Figures \ref{fig:inflation_RegDyn}-\ref{fig:inflation_InclProb} only consider regressors that should be included at least once in the dynamic regressor set. Hence, the expected value $\mathsf{E}\left(S_{t}\right)$ tends to be slightly larger on average than the number of relevant regressors included in Figures \ref{fig:inflation_InclProb}-\ref{fig:inflation_RegDyn}. Examining Fig. \ref{fig:inflation_ExpNumbReg}, we observe that the shrinkage increases as the quantile levels deviate further from the median. This observation aligns with intuition, indicating that the number of predictors relevant for explaining higher and lower inflation levels is, on average, lower than those for the median. Furthermore, we notice a noteworthy difference in the expected number of relevant predictors ($\mathsf{E}\left( S_{t}^{0.5}\right)$) between median and mean regression for the same dataset and period. Specifically, $\mathsf{E}\left( S_{t}^{0.5}\right)$ is considerably larger for the median regression. This disparity can be attributed to the robustness properties of the median dynamic regression, which underweights observations in the extreme tails, thereby reducing the variance that affects the Kalman filter predicting equations. Consequently, the proposed probability estimates are expected to exhibit greater efficiency compared to those derived from conditional mean regression. It is noteworthy that for all considered quantile levels, $\mathsf{E}\left( S_{t}^{0.25}\right)$ varies over time and eventually converges to a stable level after a few periods. This convergence underscores the dynamic nature of predictor relevance and highlights the importance of considering temporal dynamics in regression analysis. In conclusion, it is noteworthy that the average number of relevant predictors follows a distinct pattern. It increases, reaching its peak around mid-2019, then begins to decline until the end of 2022, and subsequently rises again during the Russian-Ukraine crisis for all quantile levels. Conversely, it remains relatively stable for the mean. This observation suggests that the BDQMA method effectively captures the intricate dynamics of the relationships between explanatory variables and the response variable over time, surpassing the capabilities of traditional mean-based approaches.\newline
%
\begin{table}[t]
\centering
\setlength{\tabcolsep}{5 pt}
\caption{Backtesting results for the US inflation dataset. For each quantile level $\tau=(0.1,0.25,0.5,0.75,0.9)$, the table reports the $p$-values for three common back-testing methods: the conditional coverage test of Kupiec, (CC), the unconditional coverage test of (UC) and the CaViaR test of \cite{engle_manganelli.2004} (DQ).}
\begin{tabular}{crrr}
  \toprule
 $\tau$& UC & CC & DQ \\ 
  \midrule
\multirow{1}{*}{$0.10$} & 0.07 & 0.11 & 0.63 \\ 
\multirow{1}{*}{$0.25$} & 0.28 & 0.44 & 0.61 \\ 
\multirow{1}{*}{$0.50$} & 0.69 & 0.64 & 0.30 \\ 
\multirow{1}{*}{$0.75$} & 0.13 & 0.23 & 0.56 \\ 
\multirow{1}{*}{$0.90$} & 0.07 & 0.14 & 0.60 \\ 
   \bottomrule
\end{tabular}
\label{tab:inflation_backtesting}
\end{table}
%
\indent The flexibility of the BDQMA method results in good forecasting performances. Table \ref{tab:inflation_backtesting} presents the outcomes of three common backtesting procedures for one-step ahead quantile forecasts of US inflation data. It is notable that both the conditional coverage and unconditional coverage, as well as the DQ test, fail to reject the null hypothesis of good quantile forecast ability at the $5\%$ significance level.
%
%
\section{Real estate forecasting with time-varying market
conditions}
\label{sec:empirical_application_estate} 
\noindent We now turn to a second application where we seek to explain the monthly values of the REIT net-of-S\&P500 return from January 1991 to September 2023, building upon the dataset utilized in \cite{ling_etal.2000}. The REIT net-of-S\&P500 return represents the difference between the monthly NAREIT equity index return and the return to the S\&P500 index for that month. A detailed dataset description is provided in Appendix \ref{sec:AppData} of the supplementary materials accompanying the paper. 
\subsection{A time-varying parameter prediction model}
In this analysis, we employ the BDQMA model, e.g. a quantile-based dynamic regression approach which substantially differs in several aspects from the static mean-regression analysis with macro fundamental and financial variables utilized by \cite{ling_etal.2000}. This dynamic quantile regression method allows us to continuously update the sets of relevant regressors for predicting future REIT levels. While \cite{ling_etal.2000} employed a combination of a rolling-window approach and a stochastic search variable selection method, our approach offers a more dynamic and flexible framework for analyzing the relationship between predictors and REIT returns.\newline
\indent Our analysis extends the work of \cite{ling_etal.2000} in two key aspects. Firstly, we investigate the predictability of moderately large/low REIT returns. Secondly, we examine the median return as a function of the same set of explanatory variables as utilized by \cite{ling_etal.2000}. Intuitively, we anticipate that the determinants of high/low returns would exhibit broad similarities, yet substantial differences across various sub-periods due to evolving market conditions, underlying mechanisms, and economic forces. Indeed, our approach is inherently more robust and dynamic. In our BDQMA framework, the inclusion probabilities of relevant variables are continuously updated over time, and the parameters of the regressors follow a latent dynamic process. This dynamic updating of inclusion probabilities at different quantile levels offers insights into how evolving economic conditions influence housing bubble bursts, a phenomenon commonly observed in the market. The quantile regression model we employ is defined in Eq. \eqref{eq:quantile_model}, with lagged endogenous variables excluded from the set of covariates. This exclusion helps ensure the independence of current observations from past values, maintaining the integrity of the analysis.\newline
%
%
\begin{figure}[p]
\begin{center}
\resizebox*{0.9\textwidth}{!}{\subfigure[$\tau=0.10$]{\label{fig:REITMKT_beta010_relevant}
\includegraphics[trim={0 1cm 0.5cm 1.5cm},clip,width=0.8\textwidth]{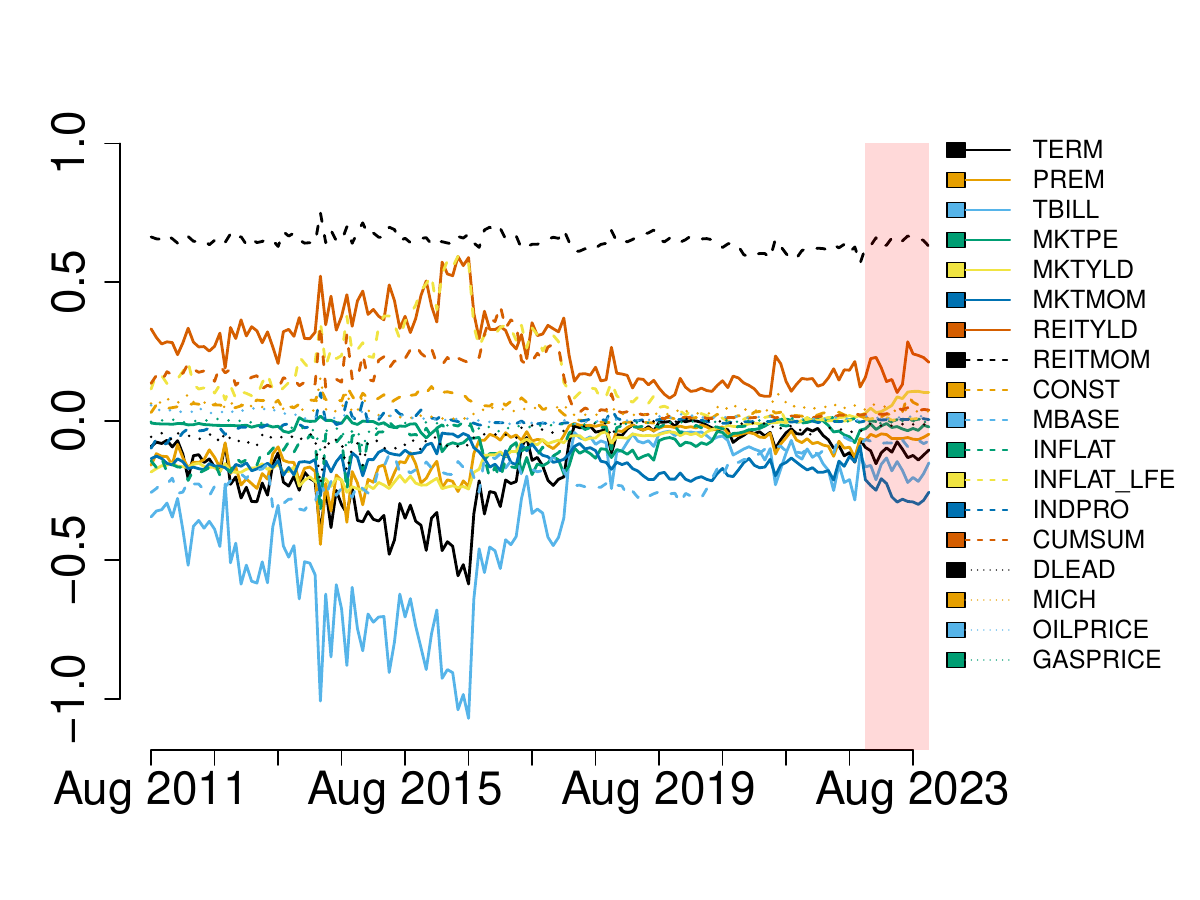}}
\subfigure[$\tau=0.25$]{\label{fig:REITMKT_beta025_relevant}
\includegraphics[trim={0 1cm 0.5cm 1.5cm},clip,width=0.8\textwidth]{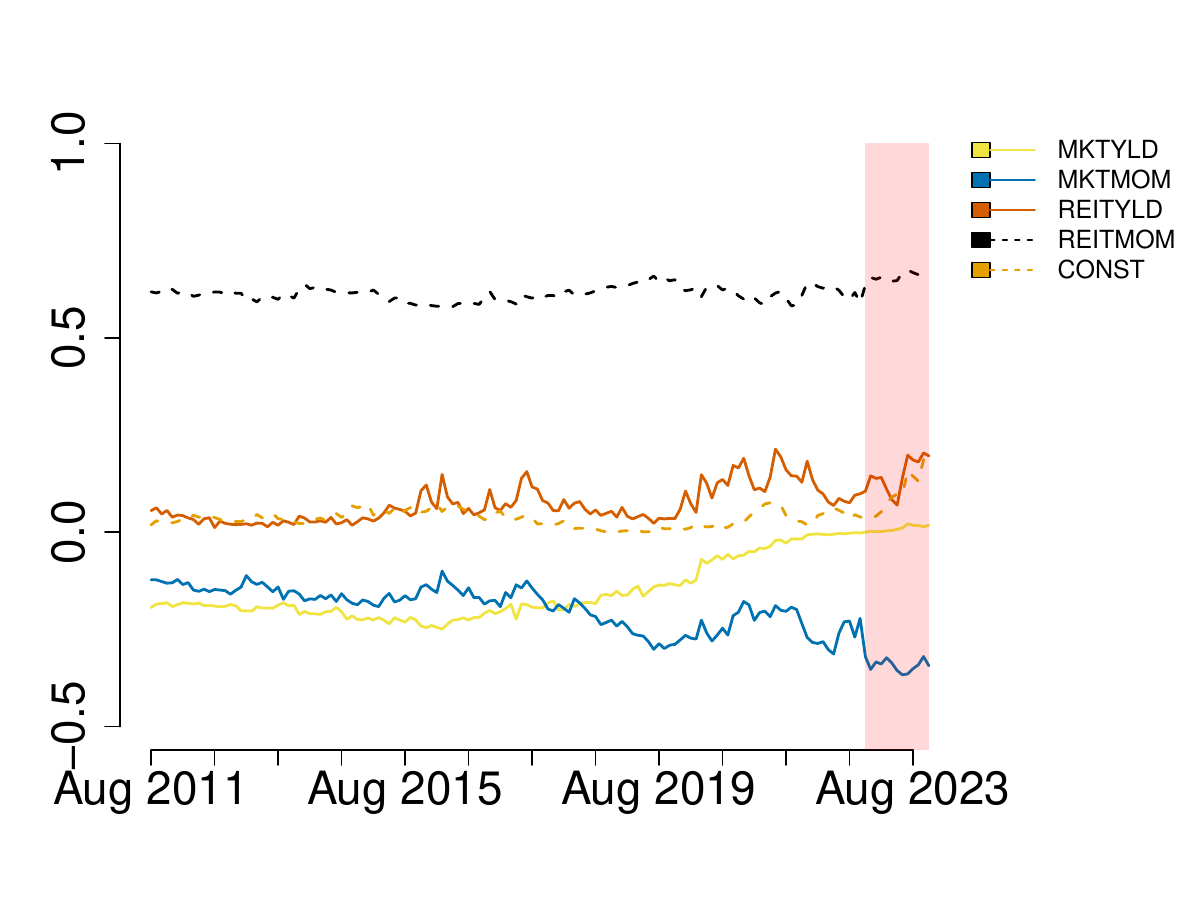}}}
\resizebox*{0.9\textwidth}{!}{\subfigure[$\tau=0.50$]{\label{fig:REITMKT_beta050_relevant}
\includegraphics[trim={0 1cm 0.5cm 1.5cm},clip,width=0.8\textwidth]{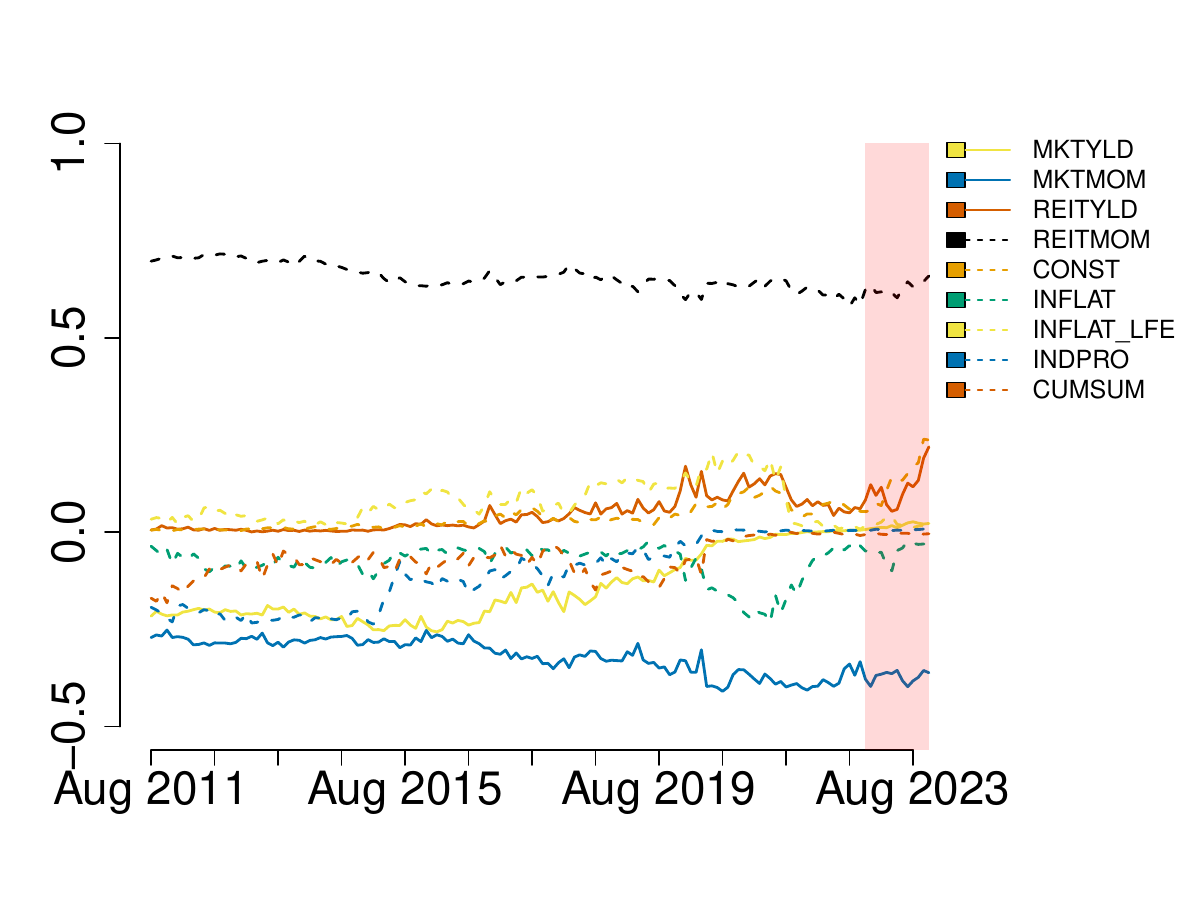}}
\subfigure[$\tau=0.75$]{\label{fig:REITMKT_beta075_relevant}
\includegraphics[trim={0 1cm 0.5cm 1.5cm},clip,width=0.8\textwidth]{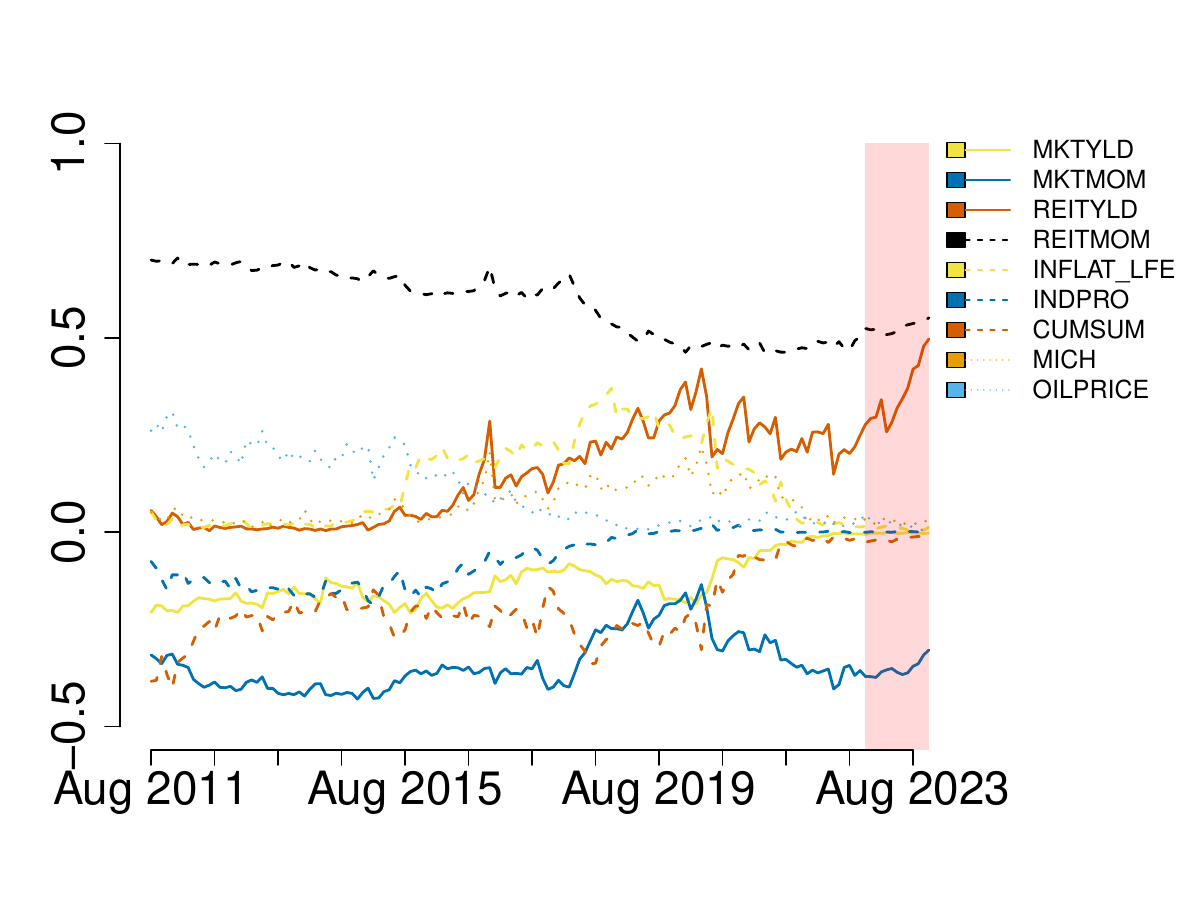}}}
\resizebox*{0.9\textwidth}{!}{\subfigure[$\tau=0.90$]{\label{fig:REITMKT_beta090_relevant}
\includegraphics[trim={0 1cm 0.5cm 1.5cm},clip,trim={0 1cm 0 1cm},clip,width=0.8\textwidth]{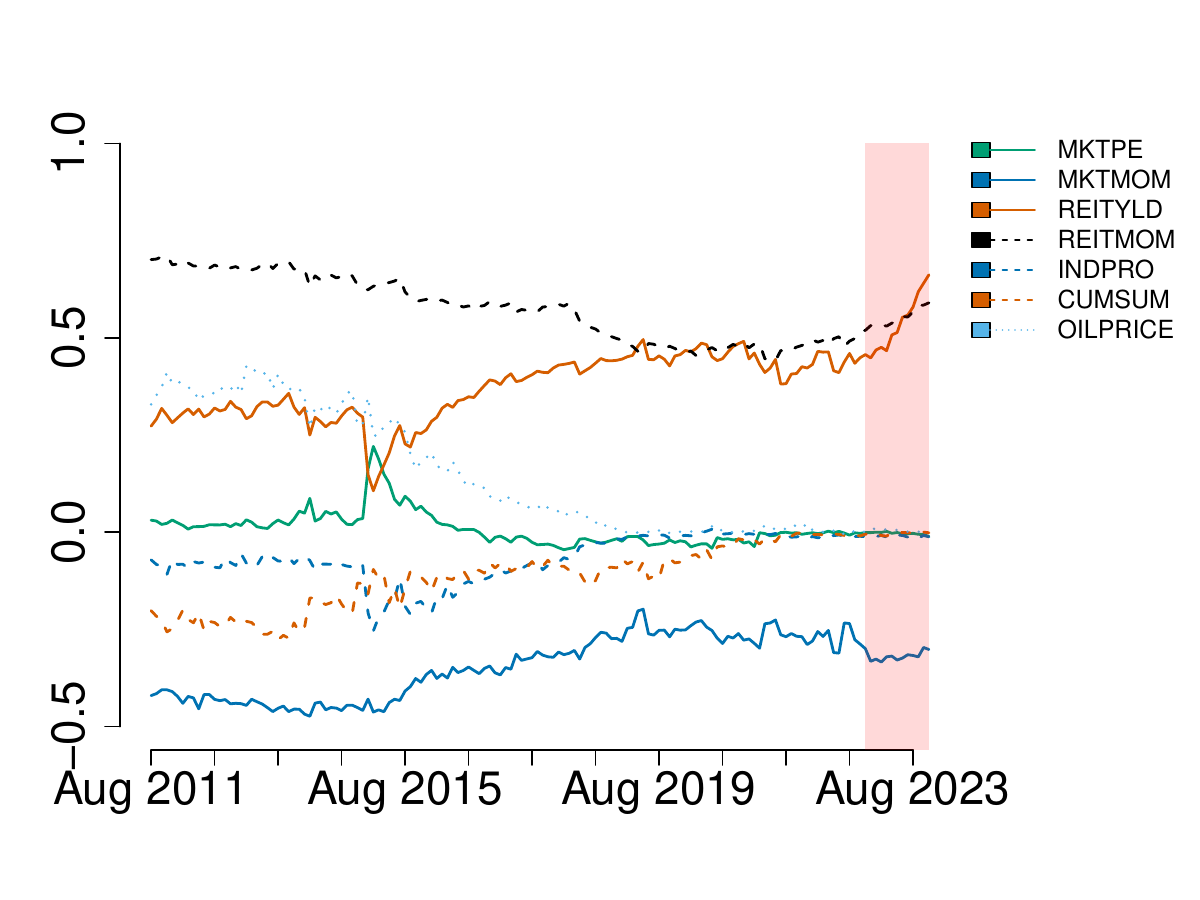}}
\subfigure[Mean regression]{\label{fig:REITMKT_beta_relevant}
\includegraphics[trim={0 1cm 0.5cm 1.5cm},clip,width=0.8\textwidth]{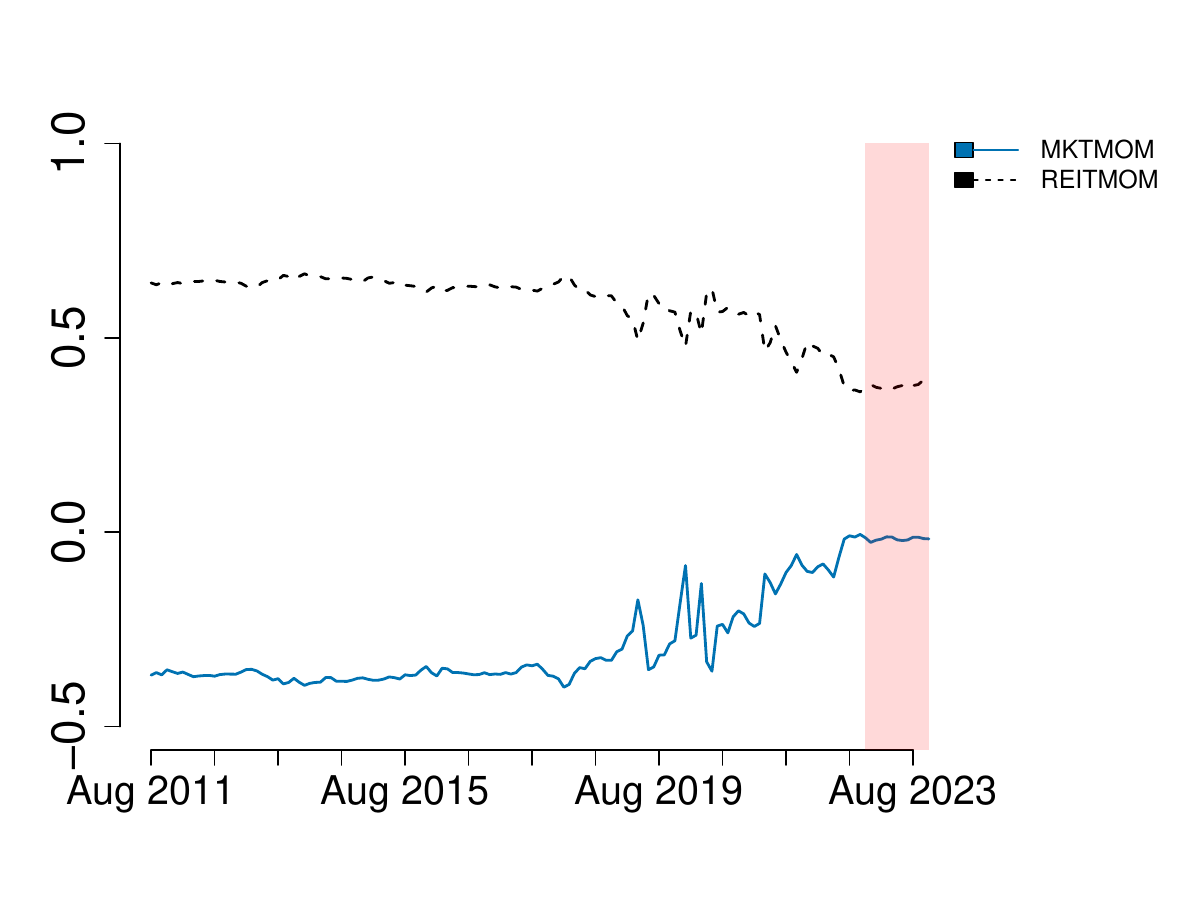}}}
\caption{US real estate data. Sequential estimates of the regression parameters by DMA (for the mean regression \ref{fig:REITMKT_beta_relevant}) and BDQMA (for the quantile regression \ref{fig:REITMKT_beta010_relevant}-\ref{fig:REITMKT_beta090_relevant}) for the REITMKT (monthly NAREIT equity return in excess of the monthly return on the S\&P 500 stock index). For each quantile level $\tau$ the corresponding figure only reports those parameters having inclusion probability larger or equal to $0.7$ for at least one quarter. See Table \ref{tab:table_US_inflation_data_CPIAUCSL_ss} for a summary of the relevant covariate. The shaded area identifies the period from 2022-02 to 2023-10 of the Russian-Ukraine crisis.}
\label{fig:realestate_RegDyn}
\end{center}
\end{figure}
%
\begin{figure}[p]
\begin{center}
\resizebox*{0.9\textwidth}{!}{\subfigure[$\tau=0.10$]{\label{fig:REITMKT_incl_prob010}
\includegraphics[trim={0 1cm 0.5cm 1.5cm},clip,width=0.8\textwidth]{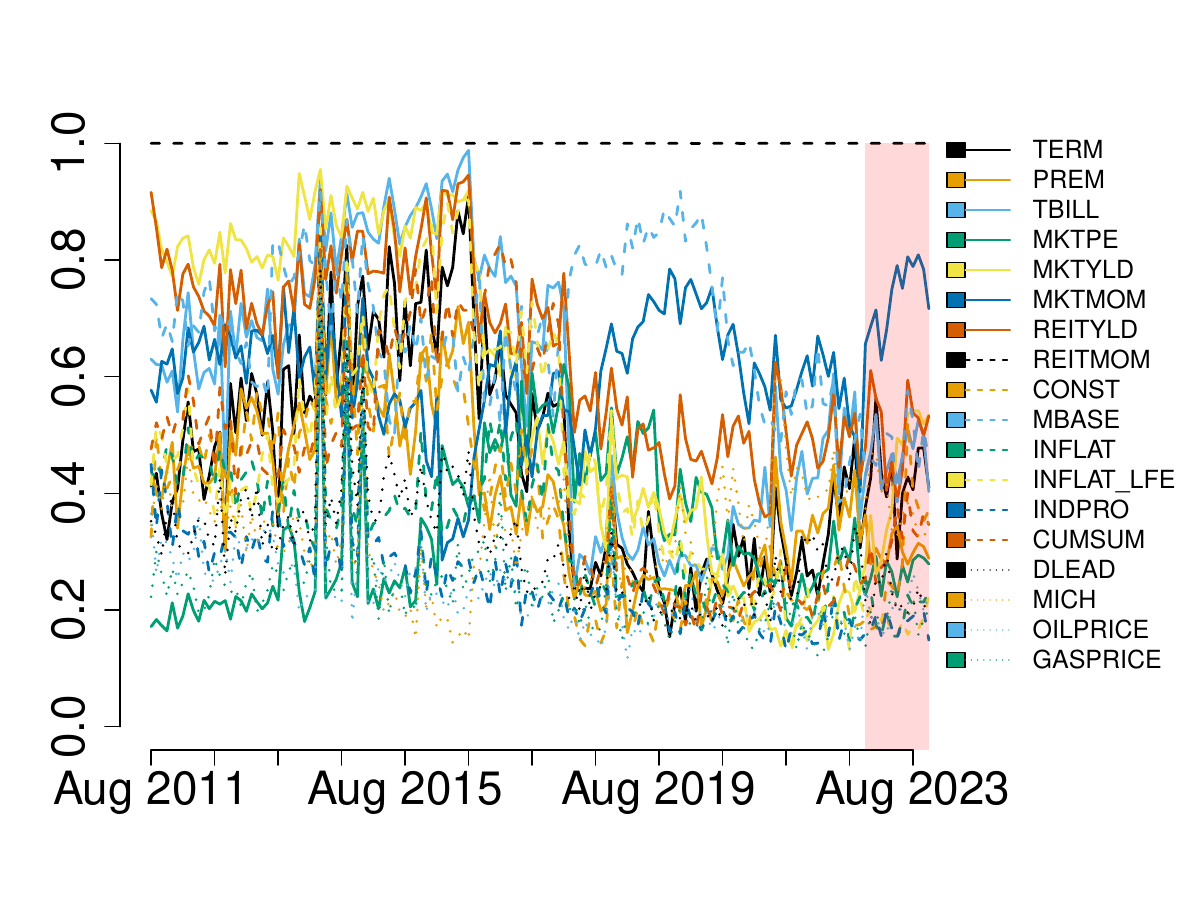}}
\subfigure[$\tau=0.25$]{\label{fig:REITMKT_incl_prob025}
\includegraphics[trim={0 1cm 0.5cm 1.5cm},clip,width=0.8\textwidth]{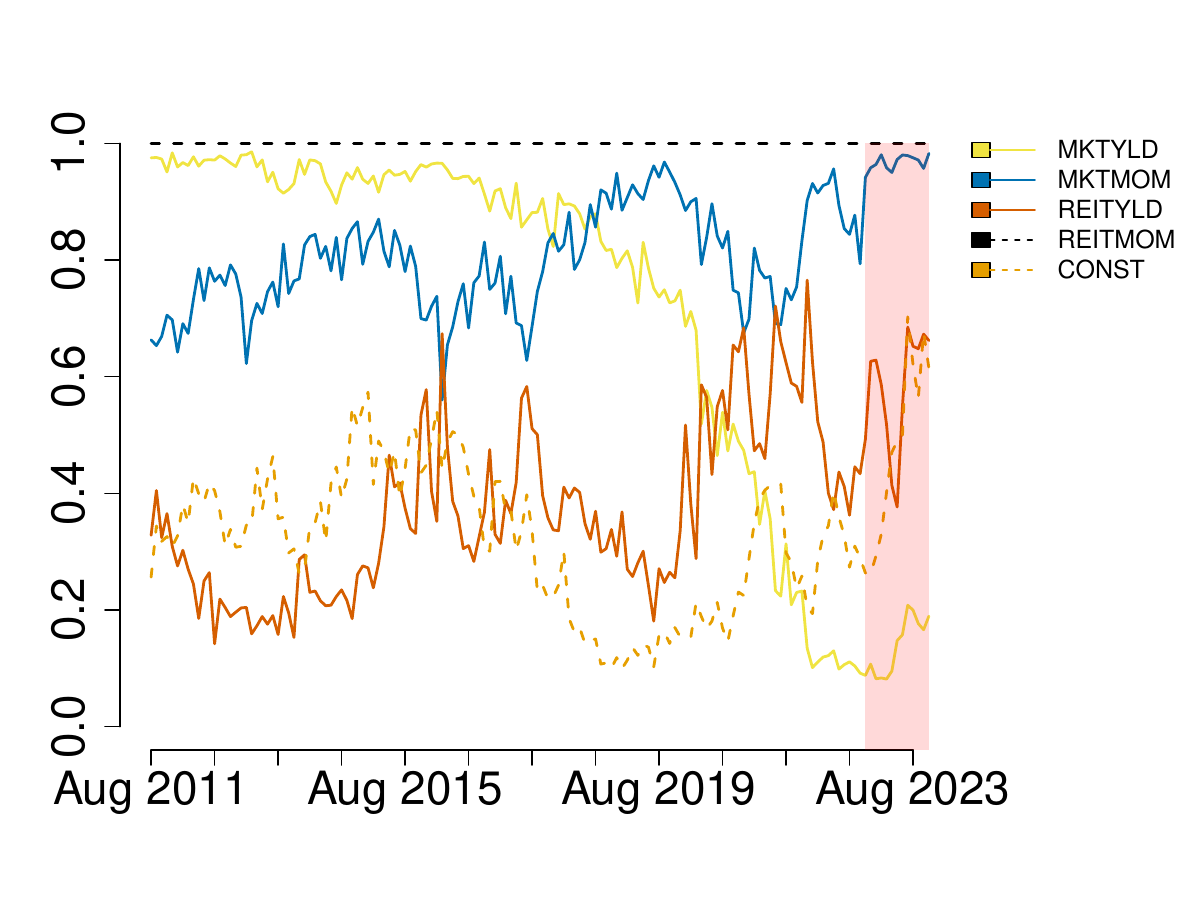}}}
\resizebox*{0.9\textwidth}{!}{\subfigure[$\tau=0.50$]{\label{fig:REITMKT_incl_prob050}
\includegraphics[trim={0 1cm 0.5cm 1.5cm},clip,width=0.8\textwidth]{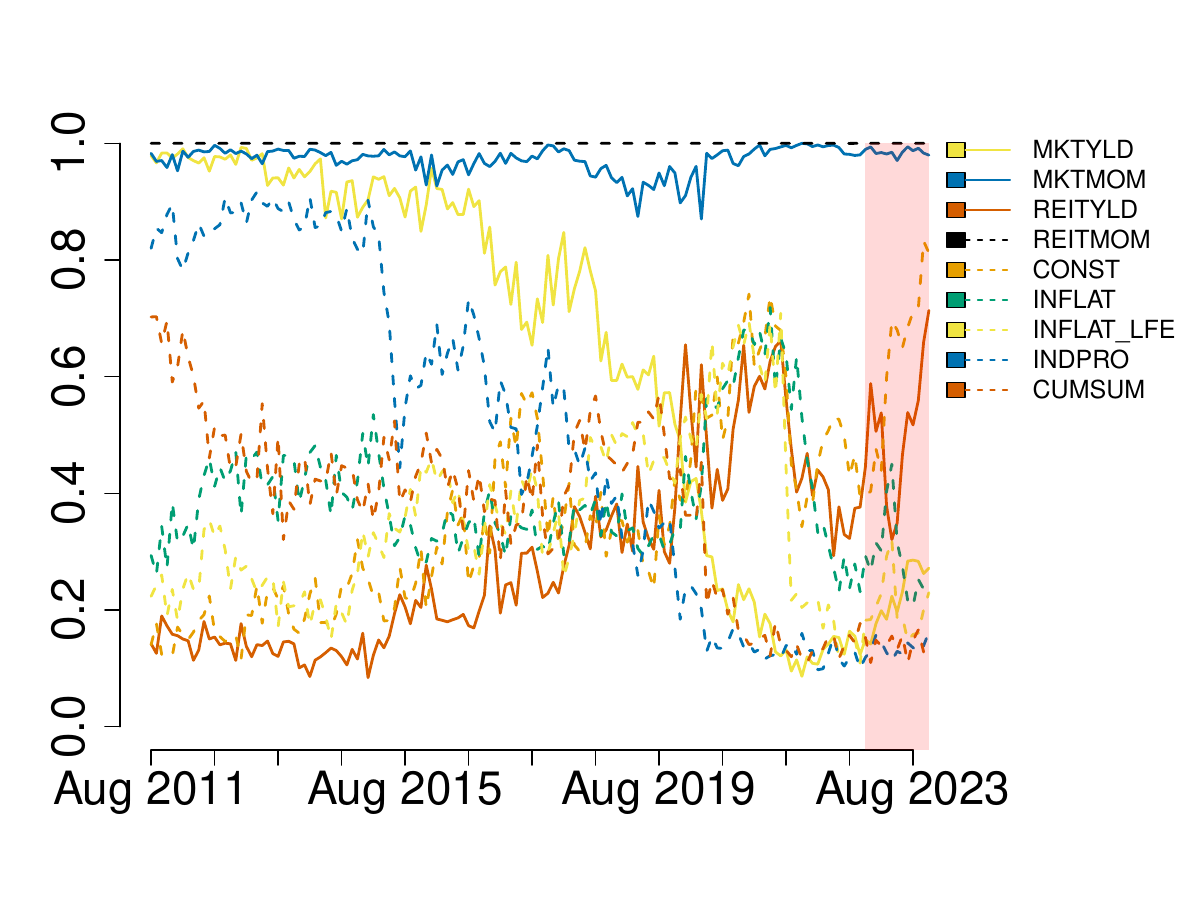}}
\subfigure[$\tau=0.75$]{\label{fig:REITMKT_incl_prob075}
\includegraphics[trim={0 1cm 0.5cm 1.5cm},clip,width=0.8\textwidth]{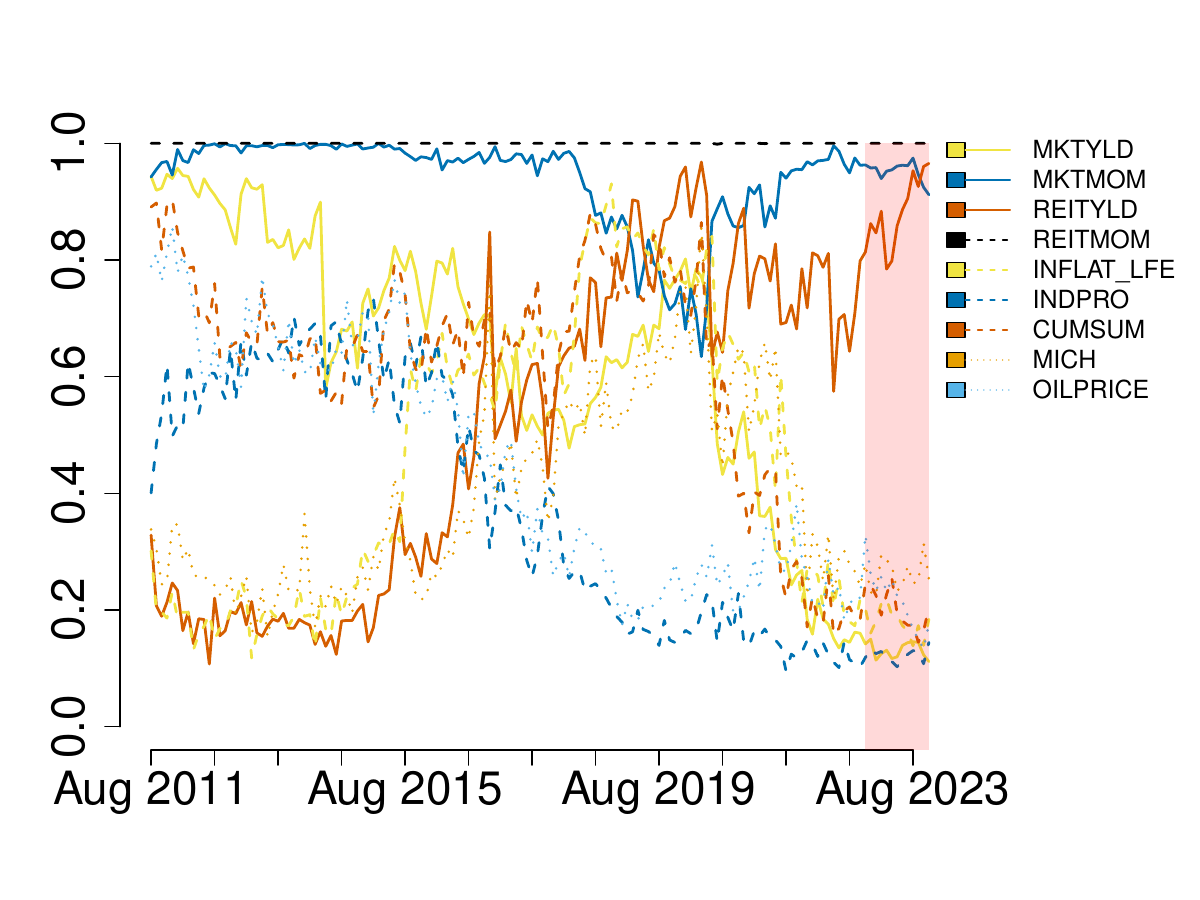}}}
\resizebox*{0.9\textwidth}{!}{\subfigure[$\tau=0.90$]{\label{fig:REITMKT_incl_prob090}
\includegraphics[trim={0 1cm 0.5cm 1.5cm},clip,width=0.8\textwidth]{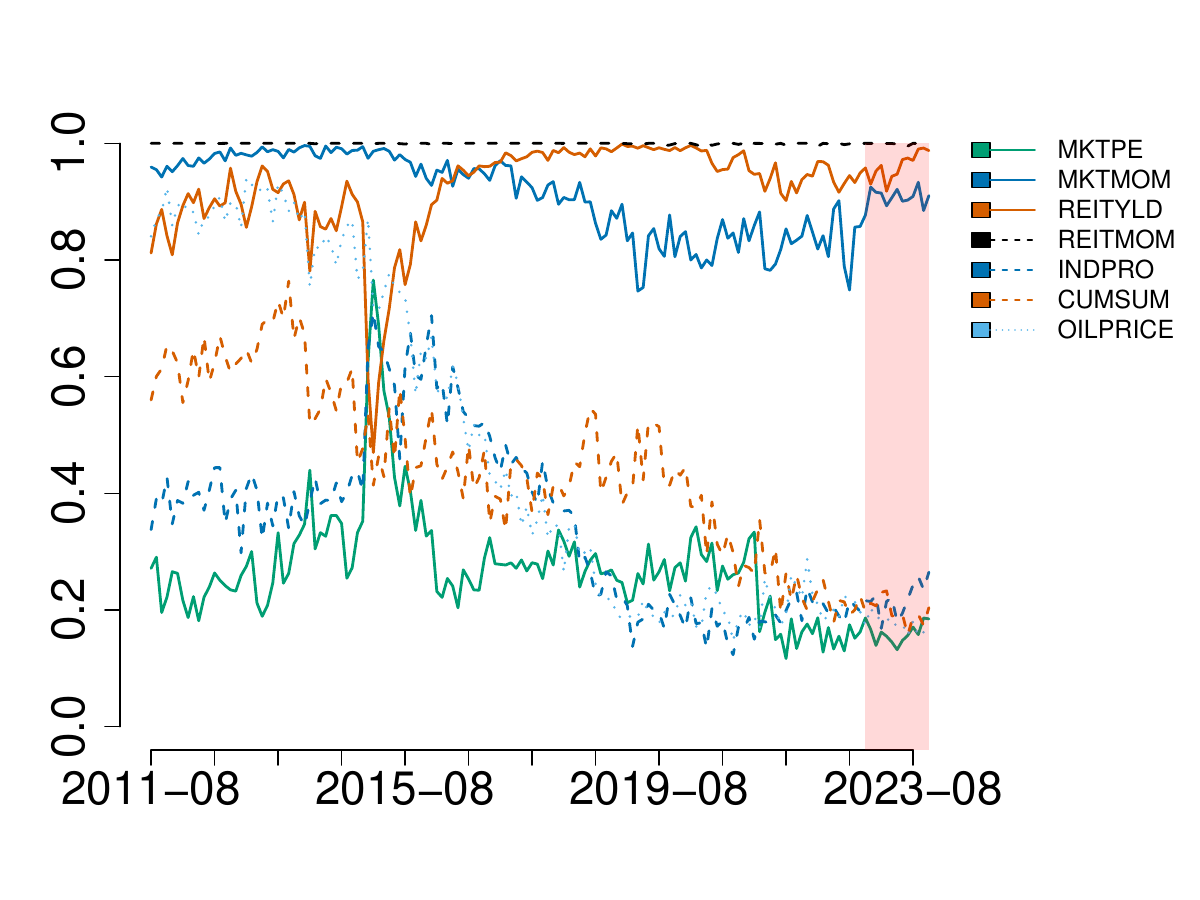}}
\subfigure[Mean regression]{\label{fig:REITMKT_incl_prob}
\includegraphics[trim={0 1cm 0.5cm 1.5cm},clip,width=0.8\textwidth]{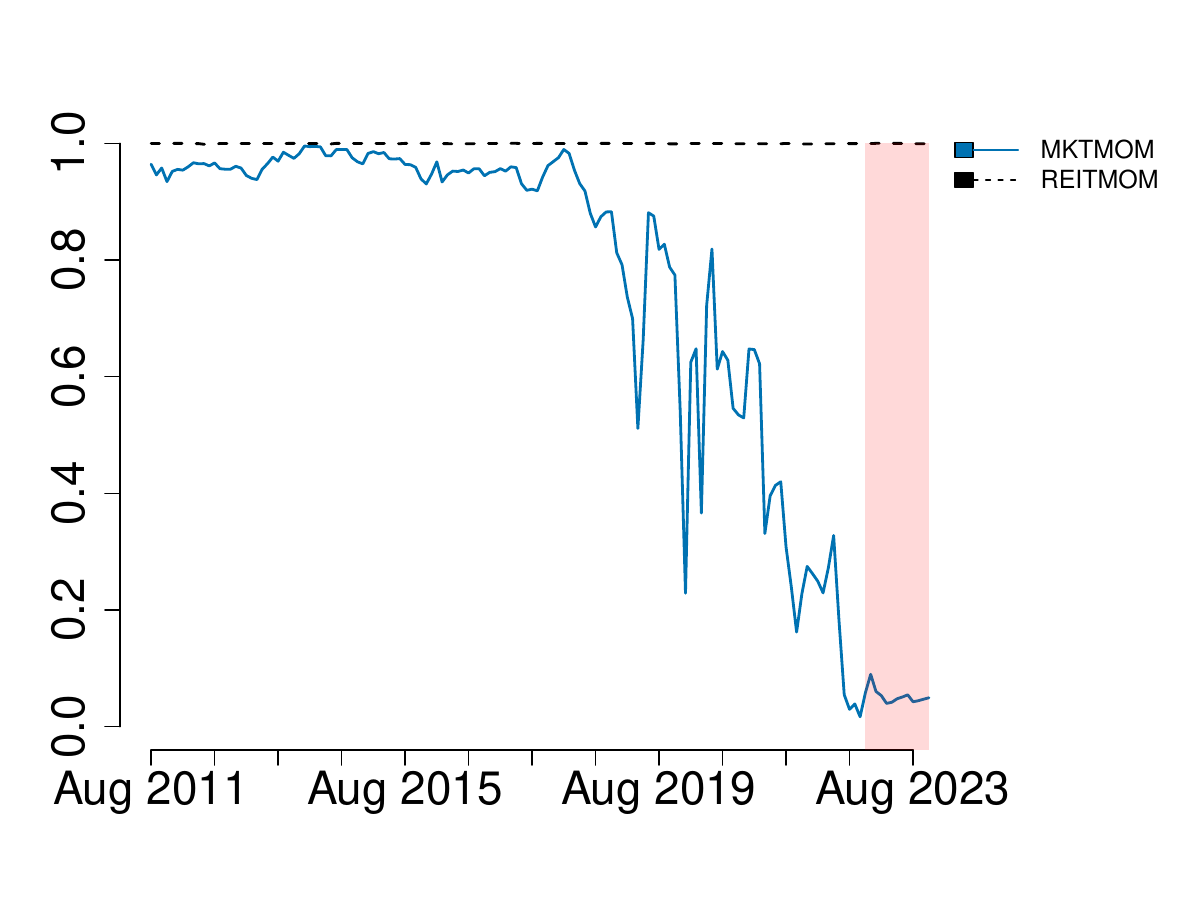}}}
\caption{US real estate data. Sequential update of the predicted inclusion
probabilities $\protect\pi _{t|t-1}$ by DMA (for the mean regression \ref{fig:REITMKT_incl_prob}) and BDQMA (for the quantile regression \ref{fig:REITMKT_incl_prob010}-\ref{fig:REITMKT_incl_prob090}) for the REITMKT (monthly NAREIT equity return in excess of the monthly return on the S\&P 500 stock index). For each quantile level $\tau$ the corresponding figure only reports those parameters having inclusion probability larger or equal to $0.7$ for at least one quarter. See Table \ref{tab:table_US_inflation_data_CPIAUCSL_ss} for a summary of the relevant covariate. The shaded area identifies the period from 2022-02 to 2023-10 of the Russian-Ukraine crisis.}
\label{fig:DQMA_RealEstate_data_results}
\end{center}
\end{figure}
%
\indent The outcomes of the initial static quantile model applied across the complete dataset are showcased in Table \ref{tab:RE_static_qreg}, Appendix \ref{sec:US_realestate_additional_results}. Upon examination, it becomes apparent that although numerous variables are pertinent for elucidating high returns, as anticipated, only a select few delineate the conditional quantile above median confidence levels. Additionally, it is noteworthy that only three variables, namely MKTYLD, MKTMOM, and REITMIN, influence all quantiles. 
A visual inspection of Fig. \ref{fig:realestate_RegDyn} reveals compelling evidence of model variations over time across all confidence levels. For instance, TBILL, which gauges the opportunity cost of real estate investment, notably contributes to the prediction of lower real estate returns surpassed only during the period spanning from late 2012 to early 2018 only for the lowest quantile level ($\tau=0.1$), after which its relevance diminishes. This suggests a temporal aspect to its influence on real estate returns. At times, however, we observe that the inclusion probability of a specific regressor exhibits a dynamic pattern that remains largely consistent as we transition from lower to upper quantiles. This pattern is notably observed with the two most influential regressors (REITMOM, MKTMOM), which consistently appear across all quantile levels. These variables are also featured in the static quantile regression in Tab. \ref{tab:RE_static_qreg} across all quantile levels, indicating their robustness in capturing real estate return dynamics. \newline
\indent From Fig. \ref{fig:DQMA_RealEstate_data_results_LOESS}, it becomes evident that certain variables, such as INFLAT (which measures the impact of inflation), are included as predictors solely for $\tau = 0.50$, indicating their relevance primarily in scenarios where real estate returns are moderate. Conversely, INFLAT\_LFE (which measures inflation excluding food and energy prices) emerges as significant in predicting high real estate returns, even extending to the third quartile ($\tau = 0.75$). This observation suggests that accounting for price dynamics, particularly those related to essential commodities, may not be crucial for forecasting future house price levels. This finding prompts further inquiry into the underlying drivers of real estate returns and the factors that shape their dynamics. It raises questions about the relative importance of macroeconomic indicators versus specific inflation components in influencing housing market performance. Moreover, it underscores the nuanced nature of predictive modeling in real estate economics, where the inclusion of certain variables may vary across different quantiles, reflecting distinct market conditions and dynamics. Overall, these insights contribute to a deeper understanding of the complex relationship between economic factors and real estate returns, informing more nuanced and effective forecasting models. By discerning the varying impacts of different predictors across quantiles, researchers and policymakers can refine their strategies for mitigating risks and maximizing opportunities in the dynamic real estate market landscape.

\subsection{Inclusion probabilities}
Fig. \ref{fig:DQMA_RealEstate_data_results} illustrates the predicted inclusion probabilities $\pi _{t|t-1}$ of the regressors across various quantile levels. The dynamic progression of these probabilities elucidates the significance of variables in forecasting future quantile levels of real estate returns at each time point $t$. Fig. \ref{fig:DQMA_RealEstate_data_results} exclusively portrays those regressors deemed relevant at least once over time, indicated by a posterior inclusion probability of at least $0.7$ for some date $t$. Unfortunately, Fig. \ref{fig:DQMA_RealEstate_data_results} contains an excessive number of regressors and considerable noise, making it challenging to discern meaningful patterns. Therefore, we have chosen to present a smoothed version in Fig. \ref{fig:DQMA_RealEstate_data_results_LOESS} in Appendix \ref{sec:US_realestate_additional_results} of the supplementary materials accompanying the paper, where each inclusion probability has undergone local polynomial regression to mitigate noise and improve interpretability.\newline
\indent Interestingly, two significant variables in Tab. \ref{tab:RE_static_qreg}, REITYLD and MKTYLD, demonstrate inclusion probability dynamics characterized by a sinusoidal shape with comparable amplitude and frequency across quantiles, transitioning to lower levels as the quantile increases. This sinusoidal behavior suggests a cyclical influence on real estate returns, warranting further investigation into its underlying drivers. Overall, the smoothed version of the inclusion probabilities provides a clearer understanding of the evolving dynamics of regressors in predicting real estate returns across different quantile levels, shedding light on both consistent influences and temporal variations in model predictors. Moreover, the percentage change in the monetary base (MBASE), which is also linked to inflation, exhibits low levels of inclusion probability across all quantile confidence levels. This suggests that the price dynamics are largely independent of the real estate market. \newline
\indent Another significant finding is that, BDQMA captured both gradual and abrupt changes in inclusion probabilities. For example, there is an abrupt change in the inclusion probability of the monetary base CUMSUM (which accounts for the growth rate of consumption expenditures for nondurable goods) from nearly $0.2$ to $0.8$ for $\tau =0.75$. However, it is worth noting that there are also many instances where the inclusion probability of a specific predictor evolves smoothly over time. Our findings underscore the importance of adaptive policymaking strategies that are responsive to changing economic dynamics. Ultimately, the dynamic selection of relevant predictors enhances the accuracy of real estate return forecasts and provides invaluable insights for policymakers seeking to optimize their interventions in alignment with evolving economic conditions.
%
\section{Conclusion}
\label{sec:conclu} 
%
In this article, we propose a new Bayesian Dynamic Quantile Model Averaging (BDQMA) approach, which for the first time to our best knowledge, combines in a tractable way Bayesian time-varying quantile regressions and dynamic model averaging. The approach dynamically accounts for model risk and parameters uncertainty in the quantiles of the response variable in presence of time-varying features such as heteroskedasticity and volatility clustering, non-linearities, breaks and jumps, leptokurticity and unconditional non-normality. These input-data features represent a major challenge when predicting the variable of interest. The tractability of BDQMA relies on an original combination of the location-scale Gaussian mixture representation of the quantile regression error terms, a conditionally Normal state-space representation of the time-varying parameter model and a sequential Monte Carlo Markov chain sampler for posterior approximation.  We demonstrate that our proposed combination of transition and jumping kernel, based on a series of Kalman filtering and smoothing steps, enables convergence to the posterior distribution. This framework is further enriched by incorporating an additional layer over the discrete space of regressors' combinations. While this resembles a conventional batch processing approach, it is typically impractical in its classical form, especially for small samples.\newline
\indent After rigorous testing with various simulation experiments, our approach has demonstrated effectiveness across different real-world data characteristics such as abrupt and smooth changes in the relationship between response and explanatory variables, in presence of heteroschedasticity in the obervation errors. We illustrate the potential of BDQMA  through an application to predicting US inflation and real estate returns. In the two applications we found evidence of disparity between mean regression and its robust median version. Indeed, beyond the disparity in the coefficients of some regressors usually considered in the literature, the impact of energy prices became significant in the recent years in explaining low and high quantiles of inflation and real estates returns. Furthermore, comparing the results obtained from static quantiles with those from BDQMA, for most of the relevant variables, the significance, the magnitude and the sign of their coefficients fluctuate over time in response to prevailing economic conditions. When compared with the mean results, BDQMA also returns a larger number of significant covariates for a large part of the sample period. This disparity can be attributed to the robustness properties of BDQMA, which dynamically underweights observations in the extreme tails. The flexibility of the BDQMA approach results in good forecasting performances as confirmed by standard backtesting procedures.\newline
\indent Finally, we may think about several extensions of our study. Relating to the method, nontrivial interesting extensions include modelling and forecasting with multiple frequencies \citep{Ghysels2007,candila2023mixed}, horizons \citep{shackleton2010multi,aastveit_etal_2014}, and quantiles \citep{meng2020estimating,taylor2022forecasting,merlo2021forecasting,ZHANG2023102498}. Furthermiore, BDQMA can be combined with an ex-ante \qmo hair-cut\qmcsp correction \citep{boucher2014risk,lazar2019model}, in order to correct on an a priori basis the output estimates in a Bayesian setting. Finally, Machine Learning techniques can be used to embed  textual media-news data \citep{huang2023improving} or to enhance the predictive power of the models \citep{medeiros2021forecasting,barkan2023forecasting,lenza2023density,lenza2023forecasting}. Regarding the applications, other geographical areas, such as Europe, or markets, such as the bond and stock markets, can be considered. Tackling the problem of predicting key economic variables, namely for instance probability of financial distress and bankruptcies \citep{huang2023improving} might be also of interest.
%

\section*{Acknowledgement}
We thank conference and seminar participants at the: SIRE Econometrics Workshop (Glasgow, 2014), the 8th CSDA International Conference on Computational and Financial Econometrics (Pisa, 2014), the 6th Italian Congress of Econometrics and Empirical Economics, and the FEM2022 Conference (Paris, 2022). An early-stage version of this article previously circulated under the title \qmo Dynamic Model Averaging for Bayesian Quantile Regression\qmc. We also thank Jean-Luc Prigent for his positive comments and encouragement on the preliminary draft of this article. The usual disclaimer applies.

\section*{Funding}
This study was funded by the MUR -- PRIN project \qmo Discrete random structures for Bayesian learning and prediction\qmcsp under g.a. n. 2022CLTYP4; the Next Generation EU -- \qmo GRINS -- Growing Resilient, INclusive and Sustainable\qmcsp project (PE0000018), National Recovery and Resilience Plan (NRRP) -- PE9 -- Mission 4, C2, Intervention 1.3; the BERN BIRD2222 01 - BIRD 2022 grant from the University of Padua, and by the European Union - Next Generation EU, Mission 4 Component 2 - CUP C53D23002580006 via the MUR-PRIN grant 2022SMNNKY. The views and opinions expressed are only those of the authors and do not necessarily reflect those of the European Union or the European Commission. Neither the European Union nor the European Commission can be held responsible for them. This research used the HPC-VERA multiprocessor cluster system at University Ca' Foscari of Venice.
\appendix
\newpage
\begin{center}
\LARGE{Supplementyary Material of:\\
\Large Bayesian Dynamic Quantile Model Averaging\\ by Bernardi, M., Casarin, R., Maillet, B.B. and Petrella, L.}
\end{center}

\noindent This supplementary materials accompanying the paper \qmo Bayesian Dynamic Quantile Model Averaging\qmc, by Bernardi, M., Casarin, R. Maillet, B.B. and Petrella, L. is organized as follows. Section \ref{subsec:Proof} provides the main Proofs (Theorem \ref{th1} and \ref{th2}). Section \ref{sec:compdet} provides further computational details for the implementation of the SMCMC algorithm for the BDQMA method. Section \ref{sec:AppData} provides a detailed description of the datasets used in the two applications: the US inflation dataset and the US real estate dataset. Section \ref{sec:US_infl_additional_results} and Section \ref{sec:US_realestate_additional_results} provide further useful results on the US inflation and real estate analysis, respectively.
%
\section{Proofs of propositions}
\label{subsec:Proof} 
\subsection{Proof of Theorem \ref{th1} and \ref{th2}}
\begin{proof}[Proof of Theorem \ref{th1}]
By the assumption on the jumping kernel it follows that
\begin{align*}
&||\pi_{t}-\mathcal{J}_{t}\circ p||_{V}=\\
&=\quad\underset{|f|\leq V}{\sup} \left|\int_{\mathbb{R}^{d_{t}}} \left(\int_{\mathbb{R}^{d_{t-1}}}\left(p(\boldsymbol{\theta}_{t-1})-
\pi_{t}(\boldsymbol{\theta}_{t-1})\mathcal{J}_{t}(\boldsymbol{\theta}_{t-1},\boldsymbol{\theta}_{t})\right)d\boldsymbol{\theta}_{t-1}\right)f(\boldsymbol{\theta}_{t})d\boldsymbol{\theta}_{t}\right.\\
&\quad+\left.\int_{\mathbb{R}^{d_{t}}} \left(\int_{\mathbb{R}^{d_{t-1}}}\left(\mathcal{J}_{t}(\boldsymbol{\theta}_{t-1},\boldsymbol{\theta}_{t})-
\pi_{t}(\boldsymbol{\theta}_{t})\pi_{t-1}(\boldsymbol{\theta}_{t-1})\right)d\boldsymbol{\theta}_{t-1}\right)f(\boldsymbol{\theta}_{t})d\boldsymbol{\theta}_{t} \right|\nonumber\\
&=\quad\underset{|f|\leq V}{\sup} \left|\int_{\mathbb{R}^{d_{t}}\times \mathbb{R}^{d}} \left(p(\boldsymbol{\theta}_{t-1})-
\pi_{t}(\boldsymbol{\theta}_{t-1})\mathcal{J}_{t}(\boldsymbol{\theta}_{t-1},\boldsymbol{\eta}_{t})\right)f(\boldsymbol{\theta}_{t})d\boldsymbol{\theta}_{t-1}d\boldsymbol{\eta}_{t}\right.\\
&\quad+\left.\int_{\mathbb{R}^{d_{t}}} \left(\int_{\mathbb{R}^{d_{t-1}}}\left(\mathcal{J}_{t}(\boldsymbol{\theta}_{t-1},\boldsymbol{\theta}_{t})-
\pi_{t}(\boldsymbol{\theta}_{t})\pi_{t-1}(\boldsymbol{\theta}_{t-1})\right)d\boldsymbol{\theta}_{t-1}\right)f(\boldsymbol{\theta}_{t})d\boldsymbol{\theta}_{t} \right|,\nonumber
\end{align*}
where $\boldsymbol{\theta}_{t}=(\boldsymbol{\theta}_{t-1},\boldsymbol{\eta}_{t})$. 
Thus one obtains
\begin{align}
&||\pi_{t}-\mathcal{J}_{t}\circ p||_{V}\leq\underset{|f|\leq V}{\sup} \left|\int_{\mathbb{R}^{d_{t}}} \left(\pi_{t}(\boldsymbol{\theta}_{t})-
\pi_{t}(\boldsymbol{\theta}_{t-1})\mathcal{J}_{t}(\boldsymbol{\theta}_{t-1},\boldsymbol{\theta}_{t})\right)f(\boldsymbol{\theta}_{t})d\boldsymbol{\theta}_{t}\right|\nonumber\\
&+\underset{|f|\leq V}{\sup} \left|\int_{\mathbb{R}^{d_{t}}} 
\left(\pi_{t}(\boldsymbol{\theta}_{t-1})\mathcal{J}_{t}(\boldsymbol{\theta}_{t-1},\boldsymbol{\theta}_{t})-p(\boldsymbol{\theta}_{t-1})\mathcal{J}_{t}(\boldsymbol{\theta}_{t-1},\boldsymbol{\theta}_{t})\right)f(\boldsymbol{\theta}_{t})d\boldsymbol{\theta}_{t}
\right|\nonumber\\
&\leq\int_{\mathbb{R}^{d_{t-1}}}\pi_{t}(\boldsymbol{\theta}_{t-1}) \underset{\boldsymbol{\theta}_{t-1}\in \mathbb{R}^{d_{t-1}}}{\sup}\underset{|f|\leq V}{\sup} \left|\int_{\mathbb{R}^{d}}\left(\pi_{t}(\boldsymbol{\eta}_{t}|\boldsymbol{\theta}_{t-1})-
\mathcal{J}_{t}(\boldsymbol{\theta}_{t-1},\boldsymbol{\theta}_{t})\right)f(\boldsymbol{\theta}_{t})d\boldsymbol{\eta}_{t}\right|d\boldsymbol{\theta}_{t-1}\nonumber\\
&+\underset{|f|\leq V}{\sup} \left|\int_{\mathbb{R}^{d_{t-1}}}\int_{\mathbb{R}^{d}} 
\left(\pi_{t}(\boldsymbol{\theta}_{t-1})-p(\boldsymbol{\theta}_{t-1})\right)\mathcal{J}_{t}(\boldsymbol{\theta}_{t-1},\boldsymbol{\theta}_{t})f(\boldsymbol{\theta}_{t})d\boldsymbol{\theta}_{t-1}d\boldsymbol{\eta}_{t}
\right|\nonumber\\
&\leq\int_{\mathbb{R}^{d_{t-1}}}\pi_{t}(\boldsymbol{\theta}_{t-1}) \underset{\boldsymbol{\theta}_{t-1}\in \mathbb{R}^{d_{t-1}}}{\sup}\underset{|f|\leq V}{\sup} \left|\int_{\mathbb{R}^{d}}\left(\pi_{t}(\boldsymbol{\eta}_{t}|\boldsymbol{\theta}_{t-1})-
\mathcal{J}_{t}(\boldsymbol{\theta}_{t-1},\boldsymbol{\theta}_{t})\right)f(\boldsymbol{\theta}_{t})d\boldsymbol{\eta}_{t}\right|d\boldsymbol{\theta}_{t-1}\nonumber\\
&+\underset{|f|\leq V}{\sup} \left|\int_{\mathbb{R}^{d_{t-1}}} 
\left(\pi_{t}(\boldsymbol{\theta}_{t-1})-p(\boldsymbol{\theta}_{t-1})\right)\left(\int_{\mathbb{R}^{d}}f((\boldsymbol{\theta}_{t-1},\boldsymbol{\eta}_{t}))d\boldsymbol{\eta}_{t}\right)d\boldsymbol{\theta}_{t-1}
\right|\nonumber\\
&=\underset{\boldsymbol{\theta}_{t-1}\in \mathbb{R}^{d_{t-1}}}{\sup}|| \pi_{t}(\boldsymbol{\eta}_{t}|\mathbb{\theta}_{t-1})-\mathcal{J}_{t}(\boldsymbol{\theta}_{t-1},\boldsymbol{\eta}_{t})||_{v}+||\pi_{t}-p||_{V}.
\end{align}
\end{proof}
\begin{proof}[Proof of Theorem \ref{th2}]
See \cite{dunson_yang.2013}.
\end{proof}
%
\section{Computational details}
\label{sec:compdet}
\subsection{Full conditional distributions of the transition kernel}
\label{subsec:TV_qreg_MCMCT} 
%
\noindent In this appendix we detail the kernel we use to simulate from the
augmented joint posterior distribution of the parameters and latent states
at time $t$, $\pi_t\left(\boldsymbol{\theta}_t\right)$. The transition
kernel consists of the following full conditional distributions. %
%
%
\begin{enumerate}
%
\item The full conditional distribution, $f\left(\sigma\vert\boldsymbol{\beta}
^{(l)}_{1:t},\mathbf{y}_{1:t}\right)$, of the measurement scale parameter $\sigma$
is obtained by collapsing the augmented latent variables up to time $t$, $
\boldsymbol{\omega}_{1:t}$. The distribution is an inverted gamma, $\mathsf{
IG}\left(a,b_l\right)$, with parameters 
\begin{equation}
a=a^0+t,\qquad b=b^0+\sum_{s=1}^{t}\rho_\tau\left(y_{s}-\mathbf{x}_{s}^\top
\boldsymbol{\beta}_s^{(l)}\right).\notag
\end{equation}
where $\rho_{\tau}\left(y\right)=y\left(\tau-\bbone\left(y<0\right)\right)$
is the $\tau$--th level quantile loss function. 
%
\item The full conditional distribution of $\omega_s^{-1}$, $
f\left(\omega_{s}^{-1}\vert\sigma^{(l)},\boldsymbol{\Omega}^{(l)},\boldsymbol{\beta}
_{1:t}^{(l)},\mathbf{y}_{1:t},\mathbf{x}_{1:t}\right)$, for $s=1,2,\dots,t$ is Inverted Normal, $
\mathsf{IN}\left(\psi_{s},\phi\right)$, with parameters 
\begin{equation}
\psi_{s}=\sqrt{\frac{\lambda^2+2\delta^2}{\left(y_{s}-\mathbf{x}_s^\top
\boldsymbol{\beta}_s\right)^2}},\qquad \phi=\frac{\lambda^2+2\delta^2}{
\delta^2\sigma}.  \notag
\end{equation}
%
%
\item The sequence of full conditional distributions of the time-varying
regression parameters $\boldsymbol{\beta}_{1:t}^{(l)}$ are obtained from the
conditional Gaussian state space representation by running the Kalman filter
prediction and filtering equations forward up to time $t$ and then
processing backward the observations using the Kalman smoothing equations.
At each time $s=1,2,\dots,t$, we have the conditional distribution 
\begin{equation}
f\left(\boldsymbol{\beta}_{1:t}\vert\sigma^{(l)},\boldsymbol{\Omega}^{(l)},\boldsymbol{\omega}
^{(l)}_{1:t},\mathbf{y}_{1:t},\mathbf{x}_{1:t}\right)=\prod_{s=1}^t f\left(\boldsymbol{\beta}_s\vert
\boldsymbol{\beta}_{s\vert t}^{(l)},\mathbf{P}_{s\vert t}^{(l)}\right),
\end{equation}
where $f\left(\boldsymbol{\beta}_s\vert\boldsymbol{\beta}_{s\vert
t}^{(l)},\mathbf{P}_{s\vert t}^{(l)}\right)$ is the density of a Normal distribution
with mean and variance 
\begin{align}
\boldsymbol{\beta}_{s\vert t}^{(l)}&=\mathsf{E}\left(\boldsymbol{\beta}
_{s}^{(l)}\vert \boldsymbol{\omega}_{1:t}^{(l)},\mathbf{y}_{1:t},\mathbf{x}_{1:t}\right) \\
\mathbf{P}_{s\vert t}^{(l)}&=\mathsf{var}\left(\boldsymbol{\beta}_{s}^{(l)}\vert \boldsymbol{
\omega}_{1:t}^{(l)},\mathbf{y}_{1:t},\mathbf{x}_{1:t}\right),
\end{align}
with $\mathbf{y}_{1:t}=\left(y_{1},\ldots,y_{t}\right)$, which are obtained by
running the following Kalman smoother recursion 
\begin{align}
\boldsymbol{\beta}_{s\vert t}^{(l)}&=\boldsymbol{\beta}_{s\vert
s}^{(l)}+\mathbf{P}_{s\vert s}^{(l)}(\mathbf{P}_{s+1\vert s}^{(l)})^{-1}\left(\boldsymbol{\beta}_{s+1\vert t}^{(l)}-\boldsymbol{\beta}_{s+1\vert s}^{(l)}\right)  \notag \\
\mathbf{P}_{s\vert t}^{(l)}&=\mathbf{P}_{s\vert s}^{(l)}+\mathbf{P}_{s\vert s}^{(l)}(\mathbf{P}_{s+1\vert
s}^{(l)})^{-1}\left(\mathbf{P}_{s+1\vert t}^{(l)}-\mathbf{P}_{s+1\vert
s}^{(l)}\right)(\mathbf{P}_{s+1\vert s}^{(l)})^{-1}\mathbf{P}_{s\vert s}^{(l)},  \notag
\end{align}
and $\left(\boldsymbol{\beta}_{s\vert s-1}^{(l)},\mathbf{P}_{s\vert s-1}^{(l)}\right)$
and $\left(\boldsymbol{\beta}_{s\vert s}^{(l)},\mathbf{P}_{s\vert s}^{(l)}\right)$
are obtained through the following Kalman predictive and filtering
equations: 
\begin{align}
\boldsymbol{\beta}_{s\vert s-1}^{(l)}&=\boldsymbol{\beta}_{s-1\vert
s-1}^{(l)}  \label{eq:kf_predicting_mean} \\
\mathbf{P}_{s\vert s-1}^{(l)}&=\mathbf{P}_{s-1\vert s-1}^{(l)}+\boldsymbol{\Omega}
\label{eq:kf_predicting_var} \\
\boldsymbol{\beta}_{s\vert s}^{(l)}&=\boldsymbol{\beta}_{s\vert
s-1}^{(l)}+\mathbf{P}_{s\vert s-1}^{(l)}\mathbf{x}_{s}V_{s}^{(l)}\nu^{(l)}_{s}
\label{eq:kf_updating_mean} \\
\mathbf{P}_{s\vert s}^{(l)}&=\mathbf{P}_{s\vert s-1}^{(l)}-\mathbf{P}_{s\vert s-1}^{(l)}\mathbf{x}
_{s}V_{s}^{(l)}\mathbf{x}_{s}^\top \mathbf{P}_{s\vert s-1}^{(l)},
\label{eq:kf_updating_var}
\end{align}
where $\nu^{(l)}_{s}=y_{s}-\widehat{y}^{(l)}_{s\vert s-1}$, with $\widehat{y}^{(l)}_{s\vert
s-1}=\lambda\omega_{s}^{(l)}+\mathbf{x}^\top_{s}\boldsymbol{\beta}
_{s\vert s-1}^{(l)}$, is the prediction error at time $s$, and $
V_{s}^{(l)}=\left(\delta\sigma^{(l)}\omega_{s}^{(l)}+\mathbf{x}
^\top_{s}\mathbf{P}_{s\vert s-1}^{(l)}\mathbf{x}_{s}\right)^{-1}$ is the variance
of the prediction error.
%
%
\item The full conditional distribution, $f(\boldsymbol{\Omega}\vert\boldsymbol{\beta}
^{(l)}_{1:t},\mathbf{y}_{1:t})$, of the transition variance-covariance
matrix $\boldsymbol{\Omega}^{(l)}$ is Inverted Wishart, $\mathsf{IW}\left(c,\mathbf{C}
\right)$ with parameters 
\begin{equation}
c=c_0+\frac{\left(t-1\right)}{2},\qquad\mathbf{C}=\mathbf{C}_0+\frac{1}{2}
\sum_{s=1}^{t-1}\left(\boldsymbol{\beta}_{s+1}-\boldsymbol{\beta}
_s\right)\left(\boldsymbol{\beta}_{s+1}-\boldsymbol{\beta}_s\right)^\top.
\notag
\end{equation}
\end{enumerate}
%
%
\subsection{Full conditional distributions of the jumping kernel}
\label{subsec:TV_qreg_MCMCJ} 
%
\noindent In this appendix we detail the jumping kernel used to sequentially
update the posterior latent states, $\left(\omega_{t+1}^{(l)},\boldsymbol{
\beta}_{t+1}^{(l)}\right)$, at each time $t$ and for the $l$-th chain of the
population, when new observations become available, $\forall l=1,2,\dots,L$.
For the easy of exposition we assume that the Jumping kernel is applied as a
new observation arrives, but the procedure can be easily extended to include
updating of block of observations.
\begin{enumerate}
\item The full conditional of $\omega_{t+1}$, $f\left(\omega_{t+1}\vert
\sigma^{(l)},\boldsymbol{\beta}_{t+1}^{(l)},y_{t+1}\right)$, is Inverted
Normal $\mathsf{IN}\left(\psi_{t+1}^{(l)},\phi^{(l)}\right)$ with
parameters 
\begin{equation}
\psi_{t+1}^{\left(l\right)}=\sqrt{\frac{\lambda^2+2\delta^2}{\left(y_{t+1}-
\mathbf{x}_{t+1}^\top\boldsymbol{\beta}_{t+1}^{\left(l\right)}\right)^2}}
,\qquad \phi^{(l)}=\frac{\lambda^2+2\delta^2}{\delta^2\sigma^{(l)}}.  \notag
\end{equation}
%
\item The full conditional distribution $f(\boldsymbol{\beta}%
_{t+1}|\sigma^{(l)},\boldsymbol{\Omega}^{(l)},\boldsymbol{\beta}_{1:t}^{(l)},\boldsymbol{\omega}
_{1:t+1}^{(l)},\mathbf{y}_{1:t+1})$ of the dynamic regression parameters $
\boldsymbol{\beta}_{t+1}$ is a Normal, $\mathsf{N}(\boldsymbol{\beta}
_{t+1\vert t+1},\mathbf{P}_{t+1\vert t+1})$, with parameters $\boldsymbol{\beta}
_{t+1\vert t+1}^{(l)}$, and $\mathbf{P}_{t+1\vert t+1}^{(l)}$ which are the filtered
mean and variance obtained by iterating the Kalman filter updating equations 
\eqref{eq:kf_predicting_mean}--\eqref{eq:kf_updating_var} for one step from
time $t$ to time $t+1$.
\end{enumerate}
%
\subsection{Fixed-lag smoother for dynamic quantile regression}
\label{subsec:TV_qreg_FixLagSmo} 
%
\noindent In this section, we describe the fixed-lag smoother for and the corresponding fixed-lag simulation smoother that are required to simulate draws from the distribution of the latent regression parameters at time $t-h$, $\boldsymbol{\beta}_{t-h}$ conditional to information up to time $t$, i.e. $\pi\left(\boldsymbol{\beta}_{t-h}\vert\mathbf{y}_{1:t},\mathbf{x}_{1:t},\boldsymbol{\omega}_{1:t},\sigma,\boldsymbol{\Omega}\right)$, where $h$ is the fixed lag and $t$ varies from $t=1,2,\dots,T$. For $t=1,2,\dots,T$, we have
\begin{equation}
f\left(\boldsymbol{\beta}_{1:t}\vert\sigma^{(l)},\boldsymbol{\Omega}^{(l)}, \boldsymbol{\omega}
^{(l)}_{1:t},\mathbf{y}_{1:t},\mathbf{x}_{1:t}\right)=\prod_{s=1}^t f\left(\boldsymbol{\beta}_s\vert
\boldsymbol{\beta}_{s-h\vert s}^{(l)},\mathbf{P}_{s-h\vert s}^{(l)}\right),
\end{equation}
where $f\left(\boldsymbol{\beta}_{s-h}\vert\boldsymbol{\beta}_{s-h\vert
s}^{(l)},\mathbf{P}_{s-h\vert s}^{(l)}\right)$ is the density of a Normal distribution
with mean and variance 
\begin{align}
\boldsymbol{\beta}_{s-h\vert s}^{(l)}&=\mathsf{E}\left(\boldsymbol{\beta}
_{s-h}^{(l)}\vert \boldsymbol{\omega}_{1:s}^{(l)},\mathbf{y}_{1:s},\mathbf{x}_{1:s}\right) \\
\mathbf{P}_{s-h\vert s}^{(l)}&=\mathsf{var}\left(\boldsymbol{\beta}_{s-h}^{(l)}\vert \boldsymbol{
\omega}_{1:s}^{(l)},\mathbf{y}_{1:s},\mathbf{x}_{1:s}\right),
\end{align}
with $\mathbf{y}_{1:s}=\left(y_{1},y_{2},\dots,y_{s}\right)$, which are obtained by
running the following Kalman fix-lag smoother recursion 
\begin{align}
\boldsymbol{\beta}_{s-h\vert s}^{(l)}&=\boldsymbol{\beta}_{s-h\vert
s-1}^{(l)}+\mathbf{P}_{s\vert s-1}^{h, (l)}\mathbf{x}_s{V^{(l)}_{s}}^{-1}\nu^{(l)}_s \\
\mathbf{P}_{s-h\vert s}^{(l)}&=\mathbf{P}_{s-h\vert s-1}^{(l)}+\mathbf{P}^{h,(l)}_{s\vert s-1}\mathbf{x}_s{V_s^{(l)}}^{-1}\mathbf{x}^\top_s\mathbf{P}_{s\vert s-1}^{h,(l)},
\end{align}
where $\nu_s^{(l)}$ and $V_s^{(l)}$ are obtained by running the Kalman filter recursion for the augmented state space model and $\mathbf{P}_{s\vert s-1}^{h+1,(l)}=\mathbf{P}_{s\vert s-1}^{h,(l)}\ell^{(l)}_s$ with $\ell^{(l)}_s=\mathbf{I}_K-\mathbf{P}_{s\vert s-1}^{h,(l)}\mathbf{x}_s{V^{(l)}_s}^{-1}\mathbf{x}_s^\top$ initialised by $\mathbf{P}_s^{0,(l)}=\mathbf{P}^{(l)}_{s\vert s-1}$, for $l=1,2,\dots,L$. Simulation from the posterior distribution of the regression parameter can be easily accomplished by exploiting joint normality of the states.

\section{Data appendix}
\label{sec:AppData} 
\begin{figure}[!t]
\begin{center}
\subfigure[Inflation]{\label{fig:re_inflation_data_infl}
\includegraphics[trim={0 1cm 1cm 1cm},clip,width=0.45\textwidth]{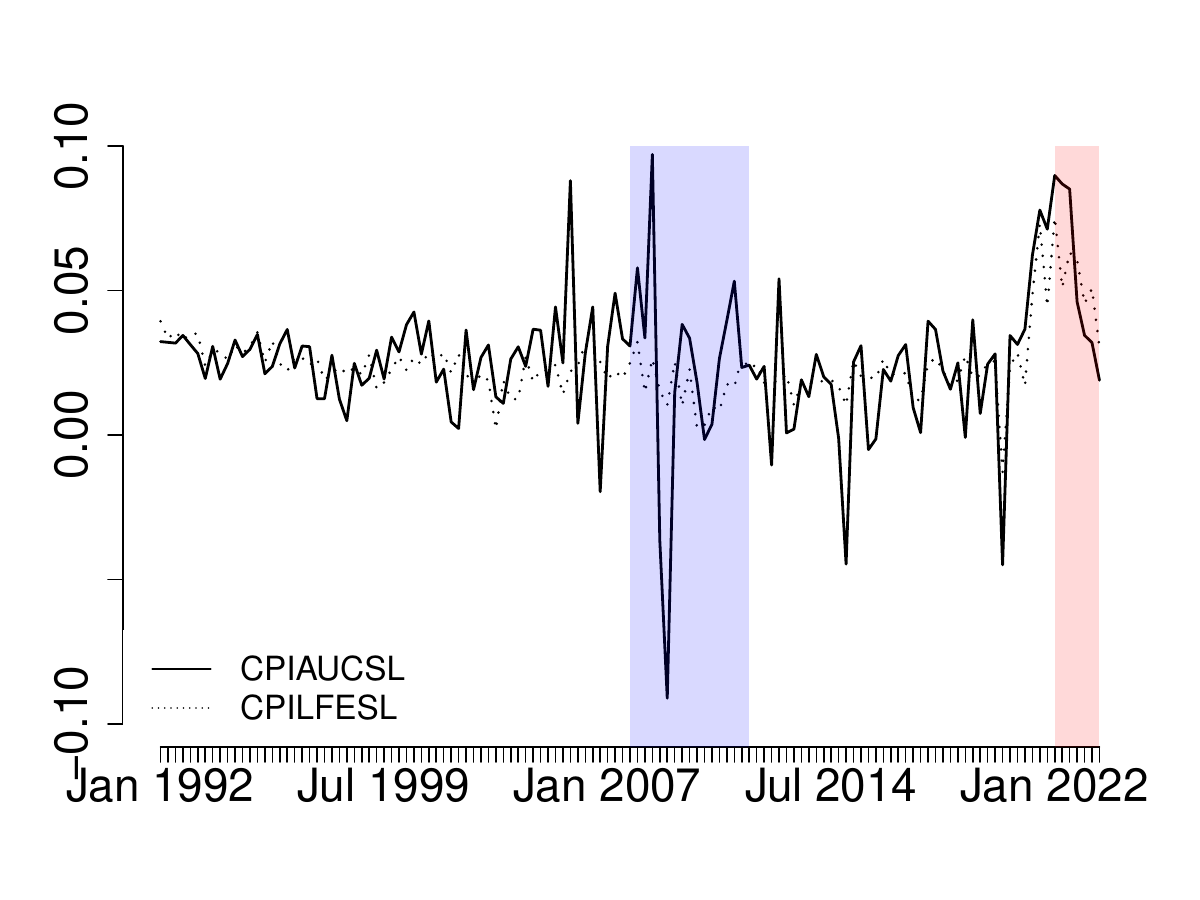}}\quad
\subfigure[Real estate]{\label{fig:re_inflation_data_re}
\includegraphics[trim={0 1cm 1cm 1cm},clip,width=0.45\textwidth]{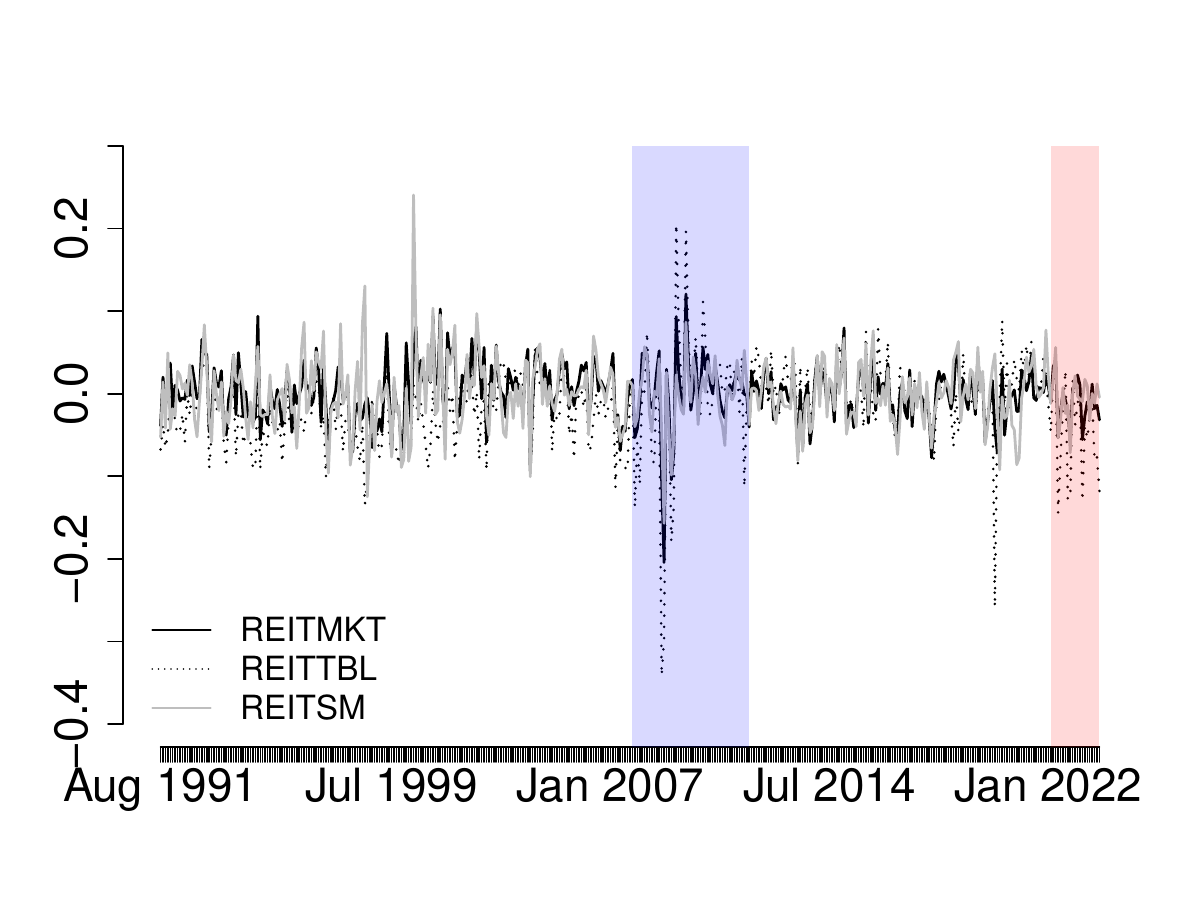}}
%
\caption{Inflation and real estate datasets. The two definitions of inflation (\qmo CPIAUCSL\qmc, \qmo CPILFESL\qmc) are provided in Table \ref{tab:table_US_inflation_data}, while the three definitions of real estate returns (\qmo REITMKT\qmc, \qmo REITTBL\qmc, \qmo REITSM\qmc) are provided in Table \ref{tab:table_US_realestate_data}. The two shaded area identifies the period from 2007-Q4 to 2011-Q4 of the US great financial crisis \textit{(blue)} and the Ukraine-Russian war (since 2022-02) \textit{(red)}.}
\label{fig:re_inflation_data}
\end{center}
\end{figure}

%
\noindent In this section, we describe the two datasets used in the empirical analyses: the US inflation dataset and the US real estate dataset. The response variables for the two analyses are plotted in Figure \ref{fig:re_inflation_data}.
%
\subsection{Inflation dataset}
\label{sec:app_inflation_dataset} 
The primary disparity between the datasets utilized by \cite{koop_korobilis.2012} and the one employed in the following analysis lies in two fundamental aspects. We will focus, as a measure of inflation, on the Consumer Price Index For All Urban Consumers: all Items in US city average (denoted \qmo CPIAUCSL\qmcsp in short) and on the Consumer Price Index For All Urban Consumers: all item excluding food and energy, US city average (denoted \qmo CPILFESL\qmcsp in short). As a complementary information, our dataset
also incorporates two additional covariates pertaining to the significance of oil and
gas prices, which were not accounted for in the \cite{koop_korobilis.2012} study. By integrating these additional variables, we aim to capture and assess the nuanced effects of oil and gas price fluctuations on the observed inflation dynamics, thereby enhancing the depth and accuracy of our analysis. Additionally, the dataset is extended in time, spanning here from 1992--Q1 to 2023--Q2, and providing a more extensive temporal scope for analysis. In particular, the analysis
includes the period 2007-2011 of the Great Financial Crisis (GFC), as well as the
period of the still on-going Russian-Ukraine confrontation (from February 2022 to
the end of the sample period). The inclusion of the two covariates accounting for
the oil and gas prices has the main purpose of providing a more comprehensive
understanding of inflation dynamics by capturing the influence of energy market
fluctuations on price levels and inflationary pressures over time. The covariates used in the empirical analysis (for both definitions of inflation) are provided in Tab. \ref{tab:inflation_posterior_summary}.\newline
\indent 
In order to make our results comparable with those obtained by \cite{koop_korobilis.2012} we consider the following predictors: the unemployment rate (UNRATE); the percentage change in real personal consumption expenditures (EC); the percentage change in private residential fixed investment (PRFI); the percentage change in real GDP (GDPC1); the logarithmic transformation of the housing starts measured as total new privately owned housing units (HOUST); the percentage change in employment measured as all employees total private industries, seasonally adjusted (USPRIV); the three month Treasury Bill rate (TB3MS); the market yield on US Treasury securities at 10-year constant maturity (GS10); the spread between the 10-year and 3-month rates on Treasury constant maturity (T10Y3MM); the spread between the 10 year and Federal Funds rates (T10YFFM); the percentage change in the money supply given by the M1 variable (M1SL); the University of Michigan measure of inflation expectations (MICH);  the change in the NAPM commodities price index (PPIACO); the percentage change in the Dow Jones Industrial Average index (DJIA); the change in the purchasing manager's composite index provided by the Institute of Supply Management (PMI); the change in the NAPM vendor deliveries index (NAPMSDI). We also accounted for the sharp upward pressure on household energy bills after the Russian invasion of the Ukraine by including two variables: oil price (OILPRICE) and gas price (GASPRICE). Further details on the variables used and their sources can be found in the Data Appendix of \cite{koop_korobilis.2013}. Table \ref{tab:table_US_inflation_data} provides a complete and detailed description of the dataset, the source of the dataset and the transformations applied to make them stationary.
%
\begin{table}[!t]
\setlength{\tabcolsep}{5 pt}
\caption{Inflation dataset. The column tcode denotes the following data transformation for a series $x$: (1) no transformation; (2) $\Delta x_t$; (3) $\Delta^2 x_t$; (4) $\log(x_t)$; (5) $\Delta\log (x_t)$; (6) $\Delta^2\log (x_t)$; (7) Percentage change $\Delta x_t/x_{t-1}$, (8) $\Delta (x_t/x_{t-1}-1)$. The FRED column gives mnemonics in FRED followed by a short description.} 
\begin{center}\resizebox{0.95\columnwidth}{!}{\begin{tabular}{lllllll}\\
\toprule
Id& Name&Type& tcode & FRED & Description\\
\cmidrule(lr){1-1}\cmidrule(lr){2-2}\cmidrule(lr){3-3}\cmidrule(lr){4-4}\cmidrule(lr){5-5}\cmidrule(lr){6-6}
%
1	&	\textsf{DATE}	&	--	&--&--&	date	\\
\multirow{2}{*}{2}	 &	\multirow{2}{*}{\textsf{CPIAUCSL}}&\multirow{2}{*}{Index}& \multirow{2}{*}{7}& \multirow{2}{*}{\textsf{CPIAUCSL}}		&	Consumer Price Index for All Urban Consumers:	\\
&&&&& All Items in U.S. City Average\\
\multirow{2}{*}{3}	&	\multirow{2}{*}{\textsf{CPILFESL}} &\multirow{2}{*}{Index}&\multirow{2}{*}{7}& \multirow{2}{*}{\textsf{CPILFESL}}		&	Consumer Price Index for All Urban Consumers: 	\\
&&&&& All Items Less Food and Energy in U.S. City Average \\
4	&	\textsf{UNRATE} &Percent	&1&\textsf{UNRATE}	&		Unemployment Rate\\
5	&	\textsf{EC}	&	Index	&7&	\textsf{DPCERA3M086SBEA}	&	Real Personal Consumption Expenditures\\
6	&	\textsf{PRFI}	&	Level	&7&	\textsf{PRFI}	&	Private Residential Fixed Investment\\
7	&	\textsf{GDPC1}	&	Level	&7&	\textsf{GDPC1}	&	Real Gross Domestic Product 	\\
8	&	\textsf{HOUST}	&	Level	&4&	\textsf{HOUST}	&		New Privately-Owned Housing Units Started: Total Units\\
9	&	\textsf{USPRIV}&	Level	&7&	\textsf{USPRIV}	&	Employees, Total Private	\\
10	&	\textsf{TB3MS}	&	Percent	&1&	\textsf{TB3MS}	&	3-Month Treasury Bill Secondary Market Rate 	\\
11	&	\textsf{GS10}	&	Percent	&1&	\textsf{GS10}	&	Market Yield on U.S. Treasury Securities at 10-Year Constant Maturity	\\
12	&	\textsf{T10Y3MM}&	Percent	&1&	\textsf{T10Y3MM}	&	10-Year Treasury Constant Maturity Minus 3-Month Treasury Constant Maturity 	\\
13	&	\textsf{T10YFFM}&	Percent	&1&	\textsf{T10YFFM}	&	10-Year Treasury Constant Maturity Minus Federal Funds Rate 	\\
14	&	\textsf{M1SL}	&	Level	&7&\textsf{M1SL}	&		Money supply - M1\\
15	&	\textsf{MICH}	&	Percent	&1&\textsf{MICH}	&		University of Michigan: Inflation Expectation\\
16	&	\textsf{PPIACO}	&	Index	&2&	\textsf{PPIACO}	&		Producer Price Index by Commodity: All Commodities \\
17	&	\textsf{DJIA}	&	Index	&7&\textsf{DJIA}&		Dow Jones Industrial Average Index\\
18	&	\textsf{PMI}	&	Index	&2&\textsf{PMI}	&	Purchasing Manager's composite index (Institute of Supply Management)\\
19	&	\textsf{NAPMSDI}&	Index	&2&	\textsf{NAPMSDI}	&		NAPM vendor deliveries index\\
20	&	\textsf{OILPRICE}&	Index	&5&	\textsf{WTISPLC}	&		Spot Crude Oil Price: West Texas Intermediate (WTI)  \\
21	&	\textsf{GASPRICE}&	Index	&5&	\textsf{GASREGCOVM}	&		US Regular Conventional Gas Price \\
%
\bottomrule 
\end{tabular}} 
\label{tab:table_US_inflation_data} 
\end{center} 
\end{table} 
%
%
%
\subsection{Real estate dataset}
\label{sec:app_real_estate_dataset} 
%
%
\begin{table}[t]
\setlength{\tabcolsep}{5 pt}
\caption{Real estate dataset. The name column is followed by a short description.} 
\begin{center}\resizebox{0.95\columnwidth}{!}{\begin{tabular}{lll}\\
\toprule
Id& Name & Description\\
\cmidrule(lr){1-1}\cmidrule(lr){2-2}\cmidrule(lr){3-3}
%
1	&	\textsf{DATE}		&	date	\\
2	&	\textsf{REITTBL} 	&	monthly NAREIT equity return in excess of the one-month T-bill\\
3	&	\textsf{REITMKT} 	&	monthly NAREIT equity return in excess of the monthly return on the S\&P 500 stock index\\
\multirow{2}{*}{4}	&	\multirow{2}{*}{\textsf{REITSM}} 		&	monthly NAREIT equity return in excess of the monthly return 	\\
&& on small-capitalization stocks from Ibbotson and Associates \\
5	&	\textsf{TERM} 		&	lagged monthly yield spread between long-term government bonds and one-month Treasury bills	\\
6	&	\textsf{PREM} 		&	lagged monthly yield spread between BAA rated bonds and government bonds\\
7	&	\textsf{TBILL} 		&	monthly T-bill yield\\
8	&	\textsf{MKTPE} 	&	lagged monthly price to earnings ratio on the S\&P 500 stock index \\
9	&	\textsf{MKTYLD} 	&	lagged monthly dividend yield on the S\&P 500 stock index \\
10	&	\textsf{MKTMOM}	&	monthly compounded return on the S\&P 500 stock index over the previous six months \\
11	&	\textsf{REITYLD} 	&	lagged dividend yield on the NAREIT equity index \\
12	&	\textsf{REITMOM} 	&	monthly compounded return on the NAREIT equity index over the previous six months \\
\multirow{2}{*}{13}	&	\multirow{2}{*}{\textsf{CONST}} 		&	growth rate (over the previous six months $t-2$ to $t-8$) of construction contracts 	\\
&& (floor space in millions of square feet) for commercial and industrial buildings \\
14	&	\textsf{MBASE} 	&	growth rate (over the previous six months $t-2$ to $t-8$) of the monetary base \\
15	&	\textsf{INFLAT} 		&	inflation (over the previous six months $t-2$ to $t-8$) rate	 \\
16	&	\textsf{INFLAT\_LFE} 		&	inflation without food and energy (over the previous six months $t-2$ to $t-8$) rate	 \\
17	&	\textsf{INDPRD} 	&	growth rate (over the previous six months $t-2$ to $t-8$) of industrial production\\
18	&	\textsf{CONSUM} 	&	growth rate (over the previous six months $t-2$ to $t-8$) of consumption expenditures for nondurable goods \\
19	&	\textsf{DLEAD} 	&	 growth rate of the index of leading economic indicators\\
20	&	\textsf{MICH} 		&	University of Michigan inflation expectations  \\
21	&	\textsf{OILPRICE}&				Spot Crude Oil Price: West Texas Intermediate (WTI)  \\
22	&	\textsf{GASPRICE}&			US Regular Conventional Gas Price \\
%
\bottomrule 
\end{tabular}} 
\label{tab:table_US_realestate_data} 
\end{center} 
\end{table}
\noindent We now provide a brief overview of the dataset used in our analysis of the U.S. real estate market. The objective is to explain the monthly values of the REIT net-of-S\&P500 return (denoted \qmo REITMKT\qmcsp in short) from January 1991 to September 2023, building upon the dataset originally used in \cite{ling_etal.2000}. The REIT net-of-S\&P500 return measures the difference between the monthly return on the NAREIT equity index and the corresponding return on the S\&P500 index. Similar to our approach in the U.S. inflation analysis, we modify the dataset by extending the sample period to September 2023 and adding two variables that capture oil and gas price dynamics. A brief description of the macroeconomic variables is provided below and summarized in Tab. \ref{tab:table_US_realestate_data}.\newline
%
%
\noindent The macroeconomic variables used for the analysis of US real estate data are: the current one-month T-bill rate
(TBILL); the spread between the yield-to-maturity (YTM) on a 30-year
government bond and the T-bill rate (TERM); the spread between the YTM on
AAA corporate bonds and the YTM on 30-year government bonds (PREM); the
percentage change in the industrial production index (INDPRD); the
percentage change in the leading economic indicators (DLEAD); the percentage
change in construction starts (CONST); the percentage change in the consumer
price index (INFL); the percentage change in nondurable consumption
(CONSUM); the percentage change in the monetary base (MBASE). 
Changes in these macroeconomic variables over the prior six-month period to avoid noise
and to decrease the impact of historical data revisions on the results. The
changes are measured from month $t-8$ to month $t-2$ for predicting month $t$
because there are reporting delays in the noninterest-rate variables. The
financial variables are: the dividend yield on the S\&P 500 (MKTYLD); the
dividend yield on the NAREIT Index (REITYLD); the S\&P500 price-earnings
(PE) ratio (MKTPE); the compounded
return to the S\&P500 during the previous six months (MKTMOM) and the
compounded return on the equity NAREIT index over the previous six months
(REITMOM). We also included the inflation without food and energy (over the previous six months $t-2$ to $t-8$) rate (INFL\_LSE)
and the University of Michigan measure of inflation expectations (MICH). As for inflation data, we included the two variables: oil price (OILPRICE) and gas price (GASPRICE) to account for the sharp upward pressure on household energy bills after the Russian invasion of the Ukraine.
See \cite{ling_etal.2000} for further details.

%
%
\section{Additional results: US inflation analysis}
\label{sec:US_infl_additional_results}
%
\subsection{Preliminary analysis}
\begin{figure}[t]
\begin{center}
\resizebox*{1\textwidth}{!}{
\subfigure[]{\label{fig:US_inflation_m}
\includegraphics[trim={0 1cm 0 1cm},clip,width=0.4\textwidth]{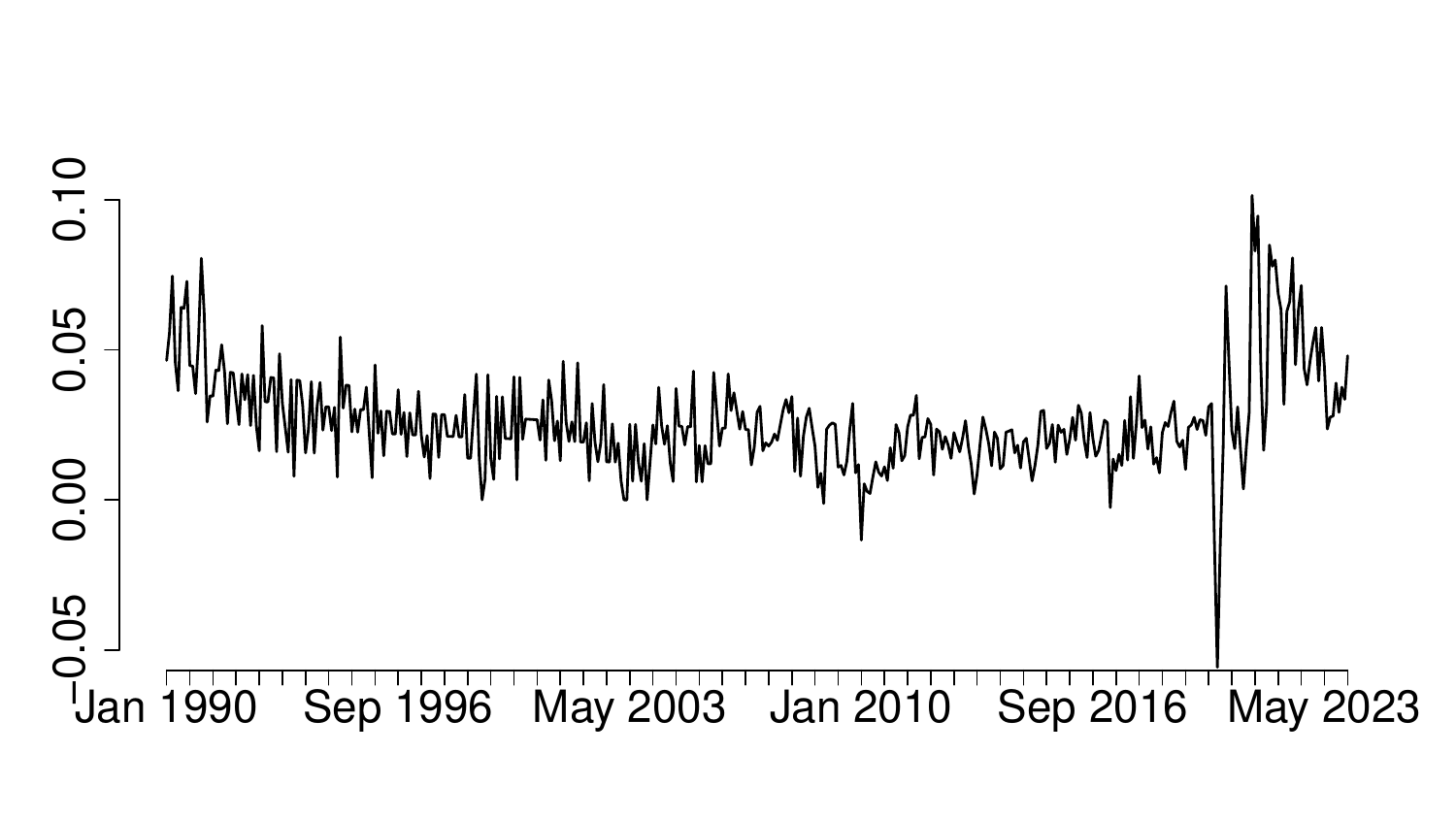}}
\subfigure[]{\label{fig:brent}
\includegraphics[trim={0 1cm 0 1cm},clip,width=0.4\textwidth]{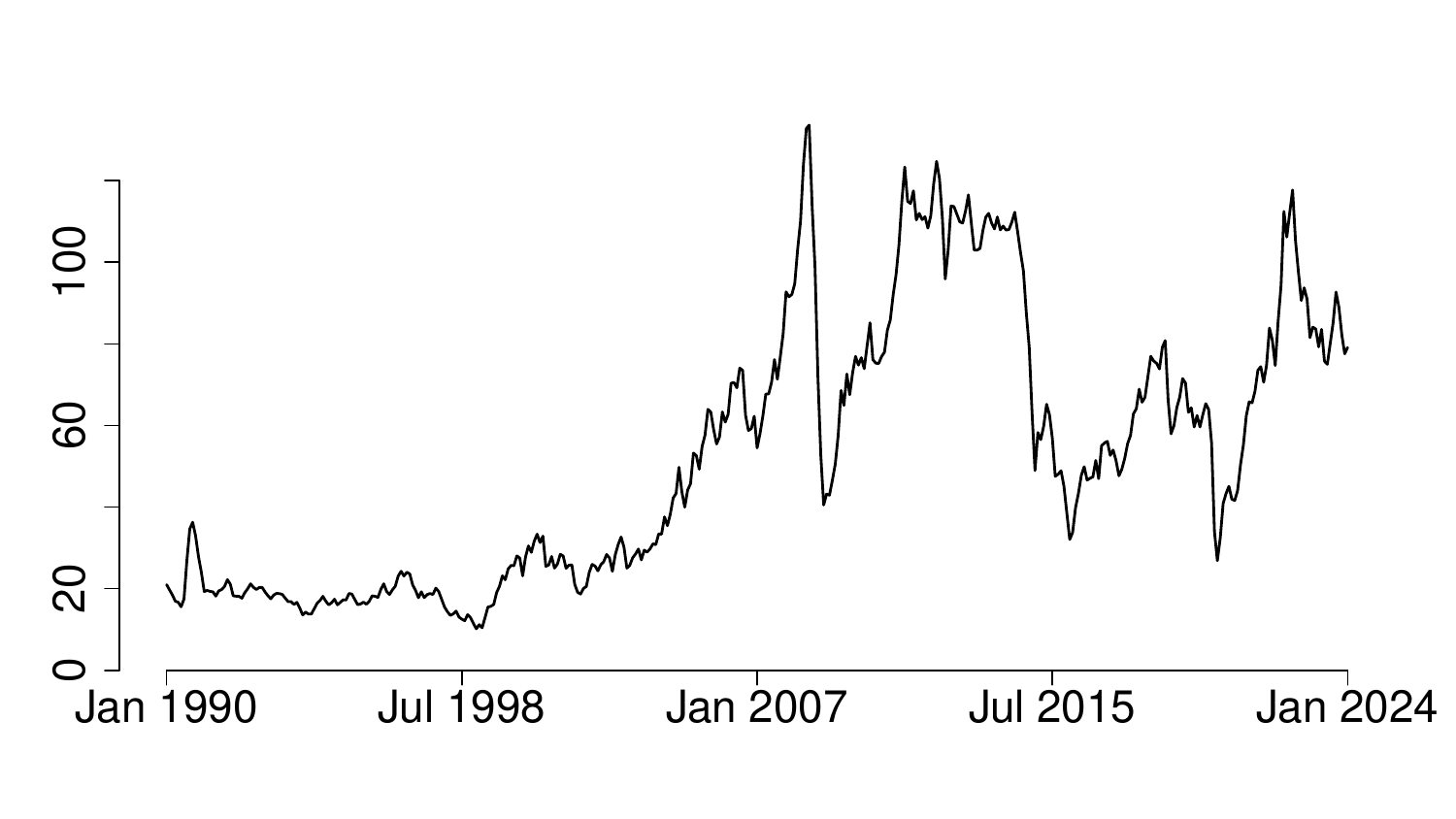}}}
\resizebox*{1\textwidth}{!}{\subfigure[]{\label{fig:google}
\includegraphics[trim={0 1cm 0 1cm},clip,width=0.4\textwidth]{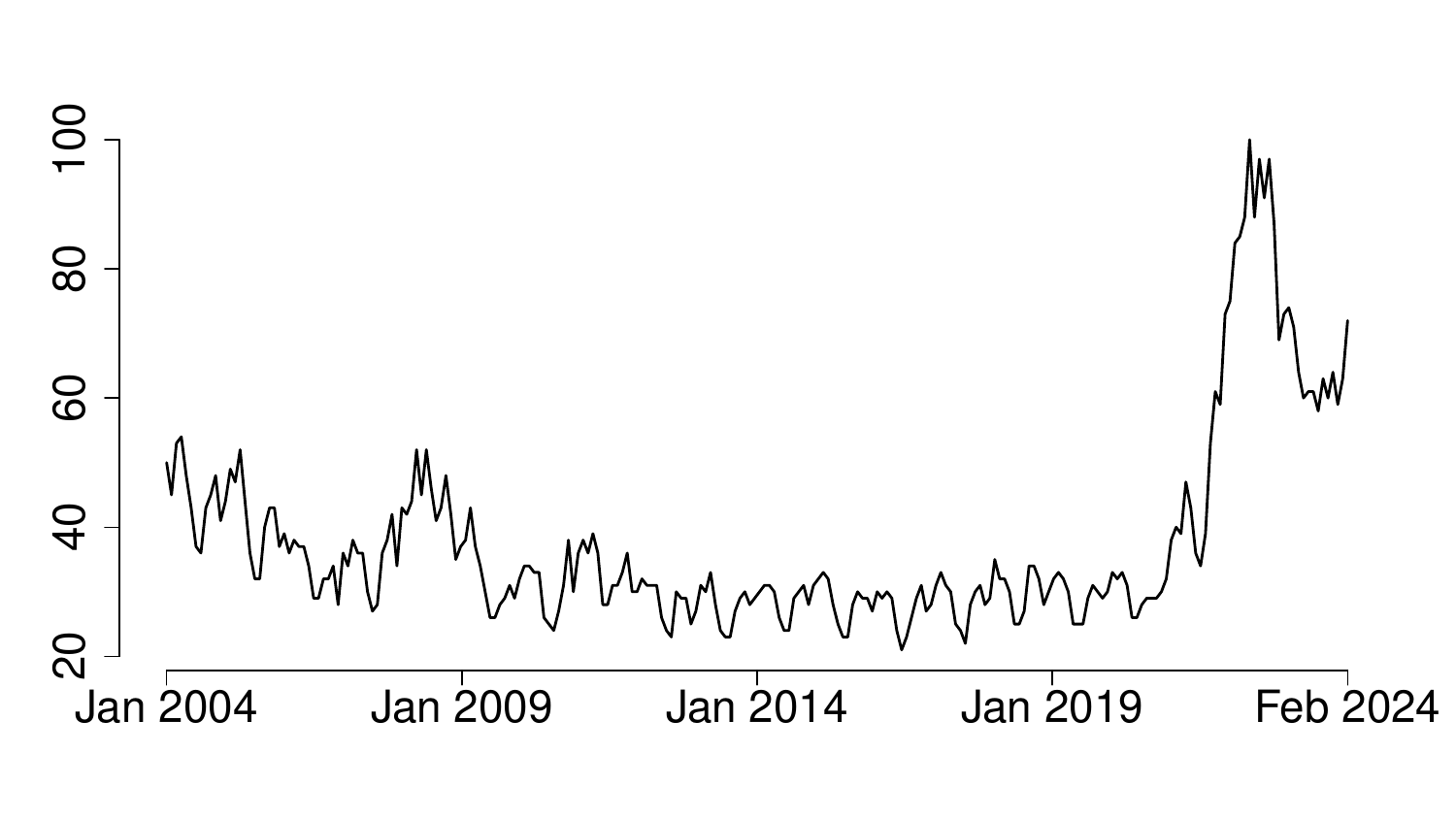}}
\subfigure[]{\label{fig:relevantIntro}
\includegraphics[trim={0 1cm 0 1cm},clip,width=0.4\textwidth]{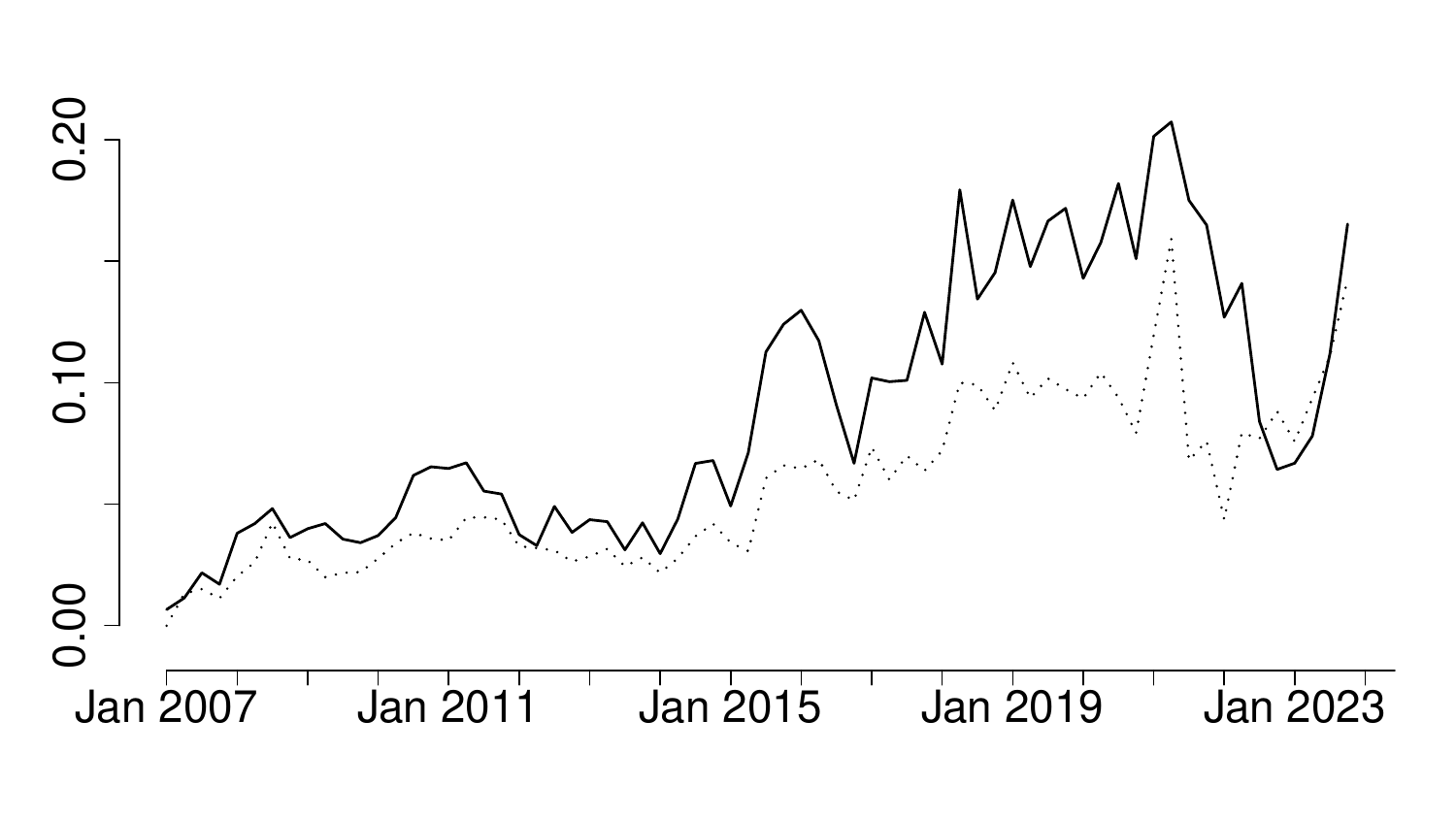}}}
\end{center}
\caption{\ref{fig:US_inflation_m} Monthly US inflation excluding food \& energy sectors. \ref{fig:brent} Quarterly Brent price in USD. \ref{fig:google} Google trend of the interest in the term ``inflation", has jumped in the US from 0 in December 2019 to 100 in December 2023. \ref{fig:relevantIntro} The impact of oil price on  ``high inflation" (solid line) and ``very high inflation" (dashed line).}
\label{fig:infl}
\end{figure}

The current inflationary surge is palpable to all consumers (see Figure \ref{fig:US_inflation_m}), as we are now witnessing its severe repercussions on a multitude of economic entities across numerous developed economies, such as ordinary citizens (and voters), consumers, investors, central, public, and private banks, international authorities, pension funds, and lastly, regulators, politicians, and decision-makers. As depicted in Figure \ref{fig:brent}, it is immediately discernible, from a simple eye-ball analysis, that Brent has recently demonstrated a heightened trend compared to recent historical levels, as well as considerably fluctuated in the past few years, especially in the aftermath of the outbreak of the Global Financial Crisis (GFC), the Pandemic COVID-19, and the last episode Energy crisis following the extreme tensions in eastern Europe. Concurrently, the Consumer Price Index (CPI) inclusive of food and energy has been on an exponential rise, notably following the COVID-19 pandemic. This ostensibly appears, on prima-facie grounds, to be the main driver of the general inflation rate.\newline
\indent As a result, reports and ``(good) old chestnuts" media-news on energy and general inflation have received considerable attention, as demonstrated by the increase in the Google searches for inflation (see Figure \ref{fig:google}), and have contributed to instigate a contraction in economic growth, accompanied by a decrease in the aggregate level of consumption and investment. This also appears to have negatively influenced the psychology and confidence of householders, as mentioned in \cite{jawadi2023revisiting}, thereby emphasizing again the need for efficient state-of-the-art operational research tools for predicting inflation dynamics, as elaborated and detailed subsequently.\newline
%
\indent On the one hand, anchoring inflation around a low level is a clear and explicit target since the seventeens,  for most (independent) main central banks around the world, since it has been largely documented that a high inflation rate is a ``threat" for the economy \citep[as named by][]{castaneda2020inflation}, destroying the signal contained in prices, and finally damaging the economy.  Since then, however, monetary aggregates' controls have gradually lost their relevance, and some endogenous or exogenous events have gained some explanatory power when explaining the inflation level \citep[see, e.g.][]{stock_watson.2007}, the increases in the money stock in several countries is still prominent in economic and political debates 	\citep[see, e.g.][]{castaneda2020inflation}. Concerning the recent COVID-19 post-pandemic period, various factors have been shown to be linked to the level of inflation, such as GDP growth, changes in energy and commodity prices, exchange and interest rates, taxes, trade balances and unemployment rates, as well as expected (anticipated) future inflation rate \citep[see, e.g.][]{goodhart2021corporate,baltensperger2023return}. \newline

\noindent Recent market tensions have once again drawn the attention of consumers, investors, and central bankers to future inflation risks, both on the downside and upside \citep[see, e.g.][]{banerjee2020inflation,apergis2021inflation,bobeica2023covid}. In the context of inflation targeting (IT) policies implemented by central banks \citep[see][]{caputo2017following,coulibaly2019inflation,huang2019inflation}, the combination of low (but rising) interest rates, potential future economic slowdowns, and adverse output-inflation trade-offs resulting from overly aggressive interventions complicates monetary policy decisions \citep{huang2019inflation}. The actions and communication strategies of central banks have become increasingly complex and sensitive, influencing inflation expectations as perceived by media outlets, consumers, investors, and market participants \citep[][]{gimeno_ibanez.2018,picault_etal.2022}. Thanks to foundational work by \cite{klein1977econometrics}, \cite{engle.1982}, \cite{bollerslev.1986}, \cite{engle1987co}, \cite{Dic79}, \cite{Dic81}, \cite{phillips1988testing}, and \cite{kwiatkowski1992testing}, we have long understood that prices and inflation exhibit specific statistical characteristics. These include the non-stationarity of prices and the (weak) stationarity of inflation, along with the autoregressive nature of price changes, which leads to heteroskedasticity. This, in turn, results in leptokurtic non-Gaussian distributions, both unconditionally and conditionally. Moreover, inflation and economic output exhibit co-integration relationships. Additional stylized facts include non-linearities, time-varying conventional moments (mean, volatility, skewness, and excess kurtosis), volatility clustering, outliers, structural breaks, jumps, threshold effects, asymmetric persistence, inertia, and sticky prices \citep[see, e.g.][]{lahiri1988interest,lahiri1987normality,culver1997there,henry2004there,canova2007structural,gambetti2008structural,hassler2014detecting,wu2018does,bataa2013structural,di2020price,harding2023understanding}. In summary, inflation is a complex and multifaceted risk with numerous causes and manifestations. The prospect of future significant inflation figures holds substantial importance for economic agents worldwide, across different regions \citep{ball2024weighted}. Figure \ref{fig:relevantIntro}  illustrates the intricate dynamics of inflation, highlighting the role of energy in shaping both \qmo high inflation\qmcsp (solid line) and \qmo very high inflation\qmcsp (dashed line). The magnitude and even the direction of the coefficients can shift across different periods and inflation regimes.
%
\subsection{Gaussian and quantile regression on the whole sample}
\label{sec:US_infl_additional_results_static}
%
\noindent In this section, we report the posterior summary statistics for the Gaussian and quantile regression models estimated on the whole datasets. Each table reports the posterior mean (mean), standard deviation (sd), $95\%$ high posterior credible sets (HPD$_{95\%}$) (reported in italics) as well as the MCMC convergence statistics: Geweke convergence test \citep{geweke.1992} (CD), Rhat statistic \citep{gelman_rubin.1992} (RHAT) and effective sample size, ESS. Summary statistics are based on a sample of $10,000$ draws post burn-in (of $25,000$ draws) for $4$ independent chains. Gray shaded areas highlight those parameters for which the corresponding HPD$_{95\%}$ does not include zero. The regression parameters associated to the lagged values of the endogenous variable $(\phi _{1},\phi _{2})$ and the scale parameter $\sigma$ are not reported to save space. They are considered as nuisance parameters.
%
%
The complete set of tables containing posterior estimates for the Gaussian regression and for quantile regression, for each quantile level $\tau=(0.1,0.25,0.5,0.75,0.9)$, are reported in Tables \ref{tab:US_inflation_mean_regression}--\ref{tab:US_inflation_quantile_regression_090}.\newline 
%
%
\indent Concerning the Phillips curve, the unemployment rate (UNRATE) shows a low probability of inclusion in the median static quantile regression analysis (see Tab. \ref{tab:inflation_posterior_summary}). However, in the dynamic regression model (see Fig. \ref{fig:inflation_InclProb} and the discussion in Section \ref{sec:empirical_application_inflation}), its influence can only be dismissed before 1973. After 1973, UNRATE exerts a negative effect on inflation (INFL), in line with the findings of \cite{koop_onorante.2012} (see Fig. \ref{fig:inflation_InclProb}). In the dynamic regression framework, several other variables are included in the median model, consistent with the static regression results. These include the percentage change in real gross domestic product (GDPC1), the percentage change in money supply (M1SL), the University of Michigan's inflation expectations (MICH), the first difference of the producer price index by commodity (PPIACO), returns on the U.S. financial index (DJIA), and variables related to energy prices, such as oil prices (OILPRICE) and gas prices (GASPRICE). Among these, only PPIACO and MICH are found significant in the dynamic mean regression (see Fig. \ref{fig:CPIAUCSL_incl_prob}). For the dynamic mean and median static regressions, additional variables like short-term interest rates (TB3MS) and gas prices (GASPRICE) are also significant. In the median static regression, short-term interest rates (TB3MS) and changes in the NAPM vendor deliveries index (NAPMSDI) are also included (see Tab. \ref{fig:CPIAUCSL_incl_prob}).\newline
\indent The producer price index of commodities (PPIACO) stands out as the sole
explanatory variable that remains consistently relevant across the entire sample.
This holds true for both the median regression and all other quantile levels, as
well as for the mean. Such a finding is in line with expectations. It is noteworthy
that the significance of the producer price index (PPIACO) diminishes towards the
latter portion of the sample period, coinciding with the onset of the Russian-Ukraine armed conflict. This
holds true across all quantiles and for the mean, albeit with varying degrees of
magnitude. While these variables contribute to the understanding inflation rates,
their relevance fluctuates over the reference period. We note, however that the
producer price index of commodity (PPIACO) is maintaining consistent
significance. Notably, expectations regarding future inflation levels (MICH) initially exhibit a high inclusion probability exceeding 0.70 at the period's outset. However, this probability declines to below 0.50 shortly after the onset of the GFC, only to sharply increase again during the Russian-Ukraine conflict in February 2022. This sudden increase in the inclusion probability of expectations regarding future inflation levels (MICH) aligns with findings from Gaussian models. However, expectations regarding future inflation levels (MICH) is completely irrelevant
before February 2023 and gradually rises to reach a probability level of one by
the end of the sample period. This highlights the dynamic nature of inflation
forecasting and the importance of understanding shifts in predictive variables over
time.\newline
\indent Another interesting finding relates to the relevance of money supply
(M1SL) and returns on the Dow Jones index (DJIA), which display similar
inclusion probability dynamics. Both variables begin to serve as significant
predictors of inflation after the onset of the GFC, but experience a sharp decrease
in relevance during the last period. This pattern suggests that certain
covariates become relevant for forecasting inflation only during periods of low
inflation volatility, indicating the complex interplay between economic indicators
and inflation dynamics. Moreover, to explain the high inflation levels experienced
after the Russian-Ukraine crisis, variables related to energy prices (OILPRICE
and GASPRICE) emerge as particularly relevant. While these variables generally
maintain significance, they exhibit distinct dynamic behaviors. Specifically, gas
price (GASPRICE) inclusion probabilities are high during the GFC, decrease
immediately after (around mid-2015), then begin to rise again until the start of the
Russian-Ukraine crisis in February 2022, when they experience a sharp decrease.
Conversely, oil price (OILPRICE) inclusion probabilities remain relatively low
during the GFC while starting to increase to very high levels immediately after to
the begin of the Russian-Ukraine crisis, after which they decrease again to very
low levels (about $0.20$). Furthermore, to elucidate the surge in inflation levels
following the Russian-Ukraine begin of the last crisis, variables associated with
energy prices (notably OILPRICE and GASPRICE), emerge as pivotal factors.
Although these variables generally maintain significance, they manifest distinct
dynamic behaviors.\newline
\indent More specifically, GASPRICE inclusion probabilities soar during the GFC,
diminish shortly after (around mid-2015), then resume an upward trajectory until the onset of the Russian-Ukraine crisis in February 2022, prompting a sharp
decline. In contrast, price of oil (OILPRICE) inclusion probabilities remain
relatively subdued during the GFC but escalate sharply immediately after the
onset of the Russian-Ukraine extreme tension, before subsequently dwindling to
minimal levels (approximately $0.20$). 
%
\begin{table}[!ht]
\setlength{\tabcolsep}{5 pt}
\caption{US inflation dataset (CPIAUCSL). Significant regression parameters (selected by using the 95\% HPD credible intervals) and their posterior mean, estimated on the whole sample for different quantile levels (columns). HPD$_{95\%}$ are in italics. The complete set of tables containing posterior estimates for the Gaussian regression and for quantile regression, for each quantile level $\tau=(0.1,0.25,0.5,0.75,0.9)$, are reported in Tables \ref{tab:US_inflation_mean_regression}--\ref{tab:US_inflation_quantile_regression_090}.}
\begin{center}\resizebox{0.9\columnwidth}{!}{\begin{tabular}{lcccccc}\\
\toprule
&\multirow{2}{*}{Mean regression}&\multicolumn{5}{c}{Quantile regression}\\
Name&& $\tau=0.10$ & $\tau=0.25$ & $\tau=0.50$ & $\tau=0.75$ & $\tau=0.90$\\
\cmidrule(lr){1-1}\cmidrule(lr){2-2}\cmidrule(lr){3-3}\cmidrule(lr){4-4}\cmidrule(lr){5-5}\cmidrule(lr){6-6}\cmidrule(lr){7-7}
	\multirow{2}{*}{\textsf{intercept}} 	&&$-0.56$& $-0.30$&&$0.29$&$0.53$\\
	&&$\it(-0.63, -0.49)$&$\it(-0.38, -0.22)$ &&$\it (0.21, 0.36)$&$\it(0.46, 0.60)$\\
	\multirow{2}{*}{\textsf{UNRATE}} 	&&-0.16 &&&&\\
	&&$\it(-0.31, -0.01)$ &&&&\\
	\multirow{2}{*}{\textsf{EC}}		& &&&&&\\
	& &&&&&\\
	\multirow{2}{*}{\textsf{PRFI}}		&&$0.29$ &$0.15$&&&\\
	&&$\it(0.12, 0.45)$ &$\it(0.01, 0.28)$&&&\\
	\multirow{2}{*}{\textsf{GDPC1}}		& &&&&$-0.20$&$-0.25$\\
	& &&&&$\it(-0.38, -0.02)$&$\it(-0.49, -0.01)$\\
	\multirow{2}{*}{\textsf{HOUST}}	& &&&&&\\
	& &&&&&\\
	\multirow{2}{*}{\textsf{USPRIV}}	& &&&&&\\
	& &&&&&\\
	\multirow{2}{*}{\textsf{TB3MS}}		&$0.15$&$0.28$ &$0.22$&$0.16$&&\\
	&$\it(0.03, 0.28)$&$\it(0.16, 0.39)$ &$\it(0.12, 0.32)$&$\it(0.05, 0.26)$&&\\
	\multirow{2}{*}{\textsf{GS10}}		& &&&&&\\
	& &&&&&\\
	\multirow{2}{*}{\textsf{T10Y3MM}}	& &&&&&\\
	& &&&&&\\
	\multirow{2}{*}{\textsf{T10YFFM}}&	& $0.44$&&&&\\
	&&$\it(0.03, 0.85)$  &&&&\\
	\multirow{2}{*}{\textsf{M1SL}}		& &&&&&\\
	&&&&&&\\
	\multirow{2}{*}{\textsf{MICH}}		&$0.33$&$0.27$ &$0.28$&$0.33$&$0.35$&$0.40$\\
	&$\it(0.18, 0.48)$&$\it(0.13, 0.41)$ &$\it(0.14, 0.42)$&$\it(0.18, 0.48)$&$\it(0.20, 0.50)$&$\it(0.25, 0.55)$\\
	\multirow{2}{*}{\textsf{PPIACO}}	&$0.42$&$0.46$ &$0.40$&$0.49$&$0.36$\\
	&$\it(0.23, 0.62)$&$\it(0.23, 0.68)$&$\it(0.22, 0.58)$&$\it(0.32, 0.66)$&$\it(0.14, 0.57)$&\\
	\multirow{2}{*}{\textsf{DJIA}}		&& $-0.10$&&&&\\
	&&$\it(-0.19, -0.00)$ &&&&\\
	\multirow{2}{*}{\textsf{PMI}}		& &&&&\\
	& &&&&&\\
	\multirow{2}{*}{\textsf{NAPMSDI}}	&& $-0.19$&&$0.13$&$0.14$&\\
	&&$\it(-0.34, -0.04)$ &&$\it(0.01, 0.25)$&$\it(0.01, 0.26)$&\\
	\multirow{2}{*}{\textsf{OILPRICE}}	&&$0.19$ &$0.17$&&&\\
	&&$\it(0.02, 0.36)$ &$\it(0.02, 0.31)$&&&\\
	\multirow{2}{*}{\textsf{GASPRICE}}&$0.21$ &&$0.22$&&$0.29$&$0.33$\\
	&$\it(0.02, 0.40)$&&$\it(0.06, 0.39)$ &&$\it(0.08, 0.50)$&$\it(0.04, 0.62)$\\
\bottomrule 
\end{tabular}} 
\label{tab:inflation_posterior_summary}
\end{center} 
\end{table} 
This nuanced understanding of inflation
dynamics underscores the importance of continuously reassessing predictive
variables and their relevance over different economic contexts. Researchers have extensively explored the dynamic nature of inflation forecasting, shedding light on the intricate relationships between economic indicators and inflationary pressures \citep[see, e.g.][]{faust2013forecasting}. For instance, \cite{feldkircher2019global} conducted a comprehensive analysis of inflation forecasting models, highlighting the evolving role of variables such as expectations on future inflation levels and financial market indices. Similarly, \cite{abdallah2023large} examined the impact of energy
prices on inflation dynamics, emphasizing the need to incorporate energy-related
variables into inflation forecasting models, especially during periods of geopolitical
uncertainty. These contributions offer valuable insights for policymakers and
researchers seeking to enhance the accuracy and robustness of inflation
forecasts in a constantly evolving economic landscape.\newline
\indent An intriguing aspect of the analysis involves comparing the inclusion probabilities of the Gaussian counterpart of the BDQMA methods, as depicted in Fig. \ref{fig:CPIAUCSL_incl_prob}, with those for the median illustrated in Fig. \ref{fig:CPIAUCSL_incl_prob050} (see Fig. \ref{fig:inflation_InclProb} and the discussion in Section \ref{sec:empirical_application_inflation}). The most notable discovery is that while only two variables are deemed relevant for the mean regression (MICH and PPIACO), BDQMA incorporates an additional five variables. Although the dynamic behavior of PPIACO remains consistent between the mean regression and BDQMA, the dynamics for MICH align only towards the end of the sampling period, particularly during the Russian-Ukraine crisis. Throughout the rest of the period, the inclusion probabilities for MICH are nearly negligible for the mean regression, whereas for BDQMA (median regression), they hover below $0.5$ only from the conclusion of the great financial crisis to the onset of 2021. This observation underscores the nuanced differences in variable selection between regression methods, emphasizing the importance of considering alternative robust modeling approaches in economic analysis.
%
\begin{table}[!t]
\caption{Inflation dataset. Summary statistics of the posterior distribution from the Gaussian regression model, based on two inflation measures: CPIAUCSL on the left, and CPILFESL on the right.}
\begin{center}
\tabcolsep=2.0mm
\resizebox{0.95\columnwidth}{!}{\begin{tabular}{rrrrrrrr|rrrrrrr}
\toprule
&\multicolumn{7}{c}{CPIAUCSL} & \multicolumn{7}{c}{CPILFESL}\\
\cmidrule(lr){1-1}\cmidrule(lr){2-8}\cmidrule(lr){9-15}
 & mean & sd & \multicolumn{2}{c}{HPD$_{95\%}$} & CD & RHAT & ESS & mean & sd &\multicolumn{2}{c}{HPD$_{95\%}$} & CD & RHAT & ESS \\ 
\midrule
intercept & -0.00 & 0.04 & -0.08 & 0.08 & 0.94 & 1.00 & 40000.00 & -0.00 & 0.06 & -0.12 & 0.11 & 0.36 & 1.00 & 40000.00 \\ 
  UNRATE & -0.08 & 0.08 & -0.23 & 0.08 & -0.54 & 1.00 & 40000.00 & -0.13 & 0.11 & -0.35 & 0.10 & 1.54 & 1.00 & 40000.00 \\ 
  EC & -0.02 & 0.19 & -0.39 & 0.35 & 0.51 & 1.00 & 40000.00 & 0.03 & 0.27 & -0.50 & 0.56 & -0.93 & 1.00 & 40000.00 \\ 
  PRFI & 0.11 & 0.06 & -0.01 & 0.23 & -0.74 & 1.00 & 40000.00 & 0.04 & 0.09 & -0.14 & 0.22 & -1.04 & 1.00 & 40000.00 \\ 
  GDPC1 & -0.14 & 0.11 & -0.36 & 0.07 & 1.23 & 1.00 & 40000.00 & -0.11 & 0.16 & -0.42 & 0.20 & -0.24 & 1.00 & 40000.00 \\ 
  HOUST & 0.07 & 0.07 & -0.08 & 0.22 & 0.74 & 1.00 & 40000.00 & 0.08 & 0.10 & -0.12 & 0.29 & 0.25 & 1.00 & 40000.00 \\ 
  USPRIV & -0.08 & 0.22 & -0.52 & 0.35 & 0.73 & 1.00 & 40000.00 & 0.19 & 0.32 & -0.43 & 0.82 & -2.34 & 1.00 & 40000.00 \\ 
TB3MS & \cellcolor{Gray}0.15 & \cellcolor{Gray}0.06 & \cellcolor{Gray}0.03 & \cellcolor{Gray}0.28 & \cellcolor{Gray}0.57 & \cellcolor{Gray}1.00 & \cellcolor{Gray}40000.00 & 0.09 & 0.09 & -0.09 & 0.26 & 0.04 & 1.00 & 40000.00 \\ 
  GS10 & -0.03 & 0.06 & -0.14 & 0.08 & -0.72 & 1.00 & 40000.00 & -0.11 & 0.08 & -0.26 & 0.05 & 1.08 & 1.00 & 40000.00 \\ 
  T10Y3MM & -0.24 & 0.24 & -0.71 & 0.24 & 0.02 & 1.00 & 40000.00 & -0.54 & 0.36 & -1.25 & 0.16 & -1.07 & 1.00 & 40000.00 \\ 
  T10YFFM & 0.29 & 0.25 & -0.20 & 0.77 & -0.89 & 1.00 & 40000.00 & 0.61 & 0.37 & -0.12 & 1.34 & 2.13 & 1.00 & 40000.00 \\ 
  M1SL & -0.03 & 0.10 & -0.22 & 0.17 & 0.20 & 1.00 & 40000.00 & -0.14 & 0.14 & -0.41 & 0.14 & 0.07 & 1.00 & 40000.00 \\ 
\rowcolor{Gray}
MICH & \cellcolor{Gray}0.33 & \cellcolor{Gray}0.07 &\cellcolor{Gray} 0.18 &\cellcolor{Gray} 0.48 &\cellcolor{Gray} -1.51 &\cellcolor{Gray} 1.00 &\cellcolor{Gray} 40000.00 & \cellcolor{Gray}\cellcolor{Gray}0.30 &\cellcolor{Gray} 0.10 &\cellcolor{Gray} 0.11 & \cellcolor{Gray}0.50 & \cellcolor{Gray}1.17 & \cellcolor{Gray}1.00 & \cellcolor{Gray}40000.00 \\ 
PPIACO & \cellcolor{Gray}0.42 & \cellcolor{Gray}0.10 & \cellcolor{Gray}0.23 & \cellcolor{Gray}0.62 & \cellcolor{Gray}-0.23 & \cellcolor{Gray}1.00 & \cellcolor{Gray}40000.00 & 0.14 & 0.14 & -0.13 & 0.41 & -0.62 & 1.00 & 40000.00 \\ 
  DJIA & -0.09 & 0.06 & -0.20 & 0.02 & -0.40 & 1.00 & 40000.00 & -0.02 & 0.08 & -0.18 & 0.13 & 0.42 & 1.00 & 40000.00 \\ 
  NAPMPMI & 0.11 & 0.07 & -0.03 & 0.25 & 1.25 & 1.00 & 40000.00 & 0.13 & 0.10 & -0.07 & 0.32 & -2.03 & 1.00 & 40000.00 \\ 
  NAPMSDI & -0.12 & 0.07 & -0.25 & 0.01 & -0.18 & 1.00 & 40000.00 & -0.17 & 0.10 & -0.36 & 0.01 & -1.28 & 1.00 & 40000.00 \\ 
  OILPRICE & 0.16 & 0.09 & -0.01 & 0.34 & 0.04 & 1.00 & 40000.00 & 0.02 & 0.13 & -0.23 & 0.27 & 0.06 & 1.00 & 40000.00 \\ 
GASPRICE & \cellcolor{Gray}0.21 & \cellcolor{Gray}0.10 &\cellcolor{Gray} 0.02 & \cellcolor{Gray}0.40 & \cellcolor{Gray}0.30 & \cellcolor{Gray}1.00 & \cellcolor{Gray}40000.00 & -0.10 & 0.13 & -0.36 & 0.15 & 2.16 & 1.00 & 40000.00 \\ 
\bottomrule
\end{tabular}}
\label{tab:US_inflation_mean_regression}
\end{center}
\end{table}
%
\begin{table}[!ht]
\caption{Inflation dataset. Summary statistics of the posterior distribution from the quantile regression model at $\tau=0.10$, based on two inflation measures: CPIAUCSL on the left, and CPILFESL on the right.}
\begin{center}
\tabcolsep=2.0mm
\resizebox{0.95\columnwidth}{!}{\begin{tabular}{rrrrrrrr|rrrrrrr}
\toprule
&\multicolumn{7}{c}{CPIAUCSL} & \multicolumn{7}{c}{CPILFESL}\\
\cmidrule(lr){1-1}\cmidrule(lr){2-8}\cmidrule(lr){9-15}
 & mean & sd & \multicolumn{2}{c}{HPD$_{95\%}$} & CD & RHAT & ESS & mean & sd &\multicolumn{2}{c}{HPD$_{95\%}$} & CD & RHAT & ESS \\ 
\midrule
\rowcolor{Gray}
intercept & \cellcolor{Gray}-0.56 & \cellcolor{Gray}0.04 & \cellcolor{Gray}-0.63 & \cellcolor{Gray}-0.49 & \cellcolor{Gray}-0.82 & \cellcolor{Gray}1.00 & \cellcolor{Gray}27366.25 & \cellcolor{Gray}-0.69 & \cellcolor{Gray}0.04 &\cellcolor{Gray}-0.78 & \cellcolor{Gray}-0.61 & \cellcolor{Gray}-0.01 & \cellcolor{Gray}1.00 & \cellcolor{Gray}28366.70 \\ 
\rowcolor{Gray}
  UNRATE & \cellcolor{Gray}-0.16 & 0.08 & -0.31 & -0.01 & -0.30 & 1.00 & 24722.99 & \cellcolor{Gray}-0.35 & 0.10 & -0.55 & -0.15 & 1.54 & 1.00 & 22817.07 \\ 
  EC & -0.27 & 0.19 & -0.64 & 0.11 & 0.35 & 1.00 & 23608.88 & -0.25 & 0.33 & -0.88 & 0.40 & -2.20 & 1.00 & 16784.56 \\ 
  \rowcolor{Gray}
  PRFI & \cellcolor{Gray}0.29 & 0.08 & 0.12 & 0.45 & -0.49 & 1.00 & 17323.31 & \cellcolor{Gray}0.29 & 0.10 & 0.10 & 0.47 & 0.54 & 1.00 & 19750.56 \\ 
  GDPC1 & -0.11 & 0.13 & -0.35 & 0.13 & 1.50 & 1.00 & 21507.94 & -0.13 & 0.13 & -0.36 & 0.13 & 0.13 & 1.00 & 23027.53 \\ 
  HOUST & 0.06 & 0.07 & -0.08 & 0.20 & 1.22 & 1.00 & 23954.61 & -0.12 & 0.08 & -0.28 & 0.03 & -0.13 & 1.00 & 24121.38 \\ 
  USPRIV & 0.00 & 0.23 & -0.45 & 0.45 & -0.21 & 1.00 & 24045.89 & 0.31 & 0.34 & -0.36 & 0.96 & -0.46 & 1.00 & 18860.69 \\ 
  \rowcolor{Gray}
  TB3MS & \cellcolor{Gray}0.28 & 0.06 & 0.16 & 0.39 & -1.09 & 1.00 & 25158.08 &\cellcolor{Gray} 0.31 & 0.09 & 0.13 & 0.48 & 0.36 & 1.00 & 18524.26 \\ 
  GS10 & -0.02 & 0.05 & -0.12 & 0.08 & -0.15 & 1.00 & 26219.00 & -0.11 & 0.06 & -0.23 & 0.01 & 0.99 & 1.00 & 24201.00 \\ 
  T10Y3MM & -0.33 & 0.20 & -0.73 & 0.07 & 0.34 & 1.00 & 29358.33 & \cellcolor{Gray}-0.62 & \cellcolor{Gray}0.28 & \cellcolor{Gray}-1.16 & \cellcolor{Gray}-0.07 &\cellcolor{Gray} 0.04 & \cellcolor{Gray}1.00 & \cellcolor{Gray}24964.81 \\ 
  \rowcolor{Gray}
  T10YFFM & \cellcolor{Gray}0.44 & 0.21 & 0.03 & 0.85 & 0.93 & 1.00 & 28939.50 & \cellcolor{Gray}0.74 & 0.30 & 0.15 & 1.32 & -0.03 & 1.00 & 23728.46 \\ 
  M1SL & -0.03 & 0.10 & -0.22 & 0.17 & 1.52 & 1.00 & 26079.26 & -0.17 & 0.15 & -0.46 & 0.13 & -0.94 & 1.00 & 20385.75 \\ 
  \rowcolor{Gray}
  MICH & \cellcolor{Gray}0.27 & 0.07 & 0.13 & 0.41 & 0.62 & 1.00 & 24685.45 & \cellcolor{Gray}0.28 & 0.07 & 0.13 & 0.42 & -0.34 & 1.00 & 26059.00 \\ 
  PPIACO & \cellcolor{Gray}0.46 & \cellcolor{Gray}0.11 & \cellcolor{Gray}0.23 & \cellcolor{Gray}0.68 & \cellcolor{Gray}-0.87 & \cellcolor{Gray}1.00 & \cellcolor{Gray}21223.60 & 0.12 & 0.12 & -0.11 & 0.34 & 0.30 & 1.00 & 24538.64 \\ 
  \rowcolor{Gray}
  DJIA & \cellcolor{Gray}-0.10 & 0.05 & -0.19 & -0.00 & -0.98 & 1.00 & 26670.11 & \cellcolor{Gray}-0.18 & 0.05 & -0.28 & -0.08 & 0.82 & 1.00 & 29279.32 \\ 
  NAPMPMI & 0.06 & 0.07 & -0.08 & 0.20 & -0.82 & 1.00 & 21662.08 & \cellcolor{Gray}0.20 & \cellcolor{Gray}0.08 & \cellcolor{Gray}0.04 & \cellcolor{Gray}0.36 & \cellcolor{Gray}0.02 & \cellcolor{Gray}1.00 & \cellcolor{Gray}22910.32 \\ 
  \rowcolor{Gray}
  NAPMSDI & \cellcolor{Gray}-0.19 & 0.08 & -0.34 & -0.04 & -1.10 & 1.00 & 20167.22 & \cellcolor{Gray}-0.17 & 0.09 & -0.34 & -0.01 & 0.51 & 1.00 & 20923.06 \\ 
  OILPRICE & \cellcolor{Gray}0.19 & \cellcolor{Gray}0.09 & \cellcolor{Gray}0.02 & \cellcolor{Gray}0.36 & \cellcolor{Gray}-0.29 & \cellcolor{Gray}1.00 & \cellcolor{Gray}22176.25 & 0.15 & 0.12 & -0.07 & 0.38 & 1.19 & 1.00 & 19902.56 \\ 
  GASPRICE & 0.17 & 0.11 & -0.03 & 0.39 & 1.66 & 1.00 & 18651.59 & -0.25 & 0.13 & -0.52 & 0.01 & 0.08 & 1.00 & 20446.24 \\ 
\bottomrule
\end{tabular}}
\label{tab:US_inflation_quantile_regression_010}
\end{center}
\end{table}
%
\begin{table}[!ht]
\caption{Inflation dataset. Summary statistics of the posterior distribution from the quantile regression model at $\tau=0.25$, based on two inflation measures: CPIAUCSL on the left, and CPILFESL on the right.}
\begin{center}
\tabcolsep=2.0mm
\resizebox{0.95\columnwidth}{!}{\begin{tabular}{rrrrrrrr|rrrrrrr}
\toprule
&\multicolumn{7}{c}{CPIAUCSL} & \multicolumn{7}{c}{CPILFESL}\\
\cmidrule(lr){1-1}\cmidrule(lr){2-8}\cmidrule(lr){9-15}
 & mean & sd & \multicolumn{2}{c}{HPD$_{95\%}$} & CD & RHAT & ESS & mean & sd &\multicolumn{2}{c}{HPD$_{95\%}$} & CD & RHAT & ESS \\ 
\midrule
  \rowcolor{Gray}
intercept & -0.30 & 0.04 & -0.38 & -0.22 & 2.86 & 1.00 & 32432.19 & -0.42 & 0.05 & -0.52 & -0.33 & -0.96 & 1.00 & 33491.22 \\ 
  UNRATE & -0.07 & 0.07 & -0.22 & 0.07 & -0.12 & 1.00 & 38796.60 & \cellcolor{Gray}-0.22 &\cellcolor{Gray}0.10 & \cellcolor{Gray}-0.41 & \cellcolor{Gray}-0.04 &\cellcolor{Gray} 0.38 & \cellcolor{Gray}1.00 & \cellcolor{Gray}36319.70 \\ 
  EC & -0.02 & 0.20 & -0.41 & 0.38 & -2.10 & 1.00 & 31075.51 & -0.15 & 0.28 & -0.71 & 0.40 & -0.21 & 1.00 & 30564.73 \\ 
   \rowcolor{Gray} 
  PRFI & 0.15 & 0.07 & 0.01 & 0.28 & -0.52 & 1.00 & 30685.85 & 0.19 & 0.09 & 0.01 & 0.36 & -1.94 & 1.00 & 30880.58 \\ 
  GDPC1 & -0.17 & 0.12 & -0.39 & 0.08 & -0.76 & 1.00 & 29412.90 & -0.06 & 0.13 & -0.29 & 0.20 & -1.03 & 1.00 & 33839.81 \\ 
  HOUST & 0.07 & 0.07 & -0.06 & 0.21 & 0.92 & 1.00 & 35537.70 & -0.08 & 0.08 & -0.24 & 0.07 & 1.02 & 1.00 & 38016.39 \\ 
  USPRIV & -0.08 & 0.22 & -0.53 & 0.35 & 0.43 & 1.00 & 33258.49 & 0.25 & 0.31 & -0.34 & 0.85 & -0.51 & 1.00 & 31771.43 \\ 
    \rowcolor{Gray}
  TB3MS & 0.22 & 0.05 & 0.12 & 0.32 & -1.70 & 1.00 & 36704.00 & 0.22 & 0.08 & 0.07 & 0.38 & -0.04 & 1.00 & 32514.45 \\ 
  GS10 & -0.01 & 0.05 & -0.11 & 0.09 & -0.89 & 1.00 & 34050.67 & -0.07 & 0.07 & -0.20 & 0.06 & 0.56 & 1.00 & 33252.85 \\ 
  T10Y3MM & -0.05 & 0.25 & -0.55 & 0.45 & 0.69 & 1.00 & 32683.90 & -0.56 & 0.35 & -1.24 & 0.11 & -0.88 & 1.00 & 30362.71 \\ 
  T10YFFM & 0.16 & 0.26 & -0.34 & 0.67 & -0.02 & 1.00 & 32591.76 & 0.59 & 0.35 & -0.09 & 1.27 & -0.52 & 1.00 & 30614.83 \\ 
  M1SL & -0.02 & 0.10 & -0.21 & 0.17 & 0.92 & 1.00 & 38850.01 & -0.15 & 0.13 & -0.41 & 0.10 & 1.59 & 1.00 & 37152.35 \\ 
    \rowcolor{Gray}
  MICH & 0.28 & 0.07 & 0.14 & 0.42 & -1.92 & 1.00 & 34027.75 & 0.26 & 0.09 & 0.09 & 0.44 & -0.71 & 1.00 & 32478.63 \\ 
  PPIACO & \cellcolor{Gray}0.40 & \cellcolor{Gray}0.09 & \cellcolor{Gray}0.22 & \cellcolor{Gray}0.58 & \cellcolor{Gray}0.37 & \cellcolor{Gray}1.00 & \cellcolor{Gray}34506.33 & 0.05 & 0.14 & -0.23 & 0.32 & 0.01 & 1.00 & 30494.14 \\ 
  DJIA & -0.08 & 0.05 & -0.18 & 0.02 & -1.08 & 1.00 & 33014.20 & \cellcolor{Gray}-0.14 &\cellcolor{Gray} 0.06 & \cellcolor{Gray}-0.25 &\cellcolor{Gray} -0.03 & \cellcolor{Gray}-1.39 & \cellcolor{Gray}1.00 & \cellcolor{Gray}38072.04 \\ 
  NAPMPMI & 0.07 & 0.07 & -0.07 & 0.21 & -0.15 & 1.00 & 31349.18 & 0.15 & 0.09 & -0.02 & 0.31 & 0.15 & 1.00 & 31888.31 \\ 
  NAPMSDI & -0.14 & 0.08 & -0.29 & 0.01 & 0.44 & 1.00 & 29530.97 & -0.14 & 0.08 & -0.30 & 0.02 & 0.75 & 1.00 & 34495.51 \\ 
  OILPRICE & \cellcolor{Gray}0.17 &\cellcolor{Gray} 0.07 &\cellcolor{Gray} 0.02 &\cellcolor{Gray} 0.31 &\cellcolor{Gray} 1.61 &\cellcolor{Gray} 1.00 &\cellcolor{Gray} 35937.46 & 0.15 & 0.11 & -0.07 & 0.38 & 0.34 & 1.00 & 30461.90 \\ 
  GASPRICE & \cellcolor{Gray}0.22 & \cellcolor{Gray}0.08 & \cellcolor{Gray}0.06 & \cellcolor{Gray}0.39 & \cellcolor{Gray}-0.42 & \cellcolor{Gray}1.00 & \cellcolor{Gray}35307.84 & -0.17 & 0.14 & -0.44 & 0.09 & -2.18 & 1.00 & 28798.47 \\ 
\bottomrule
\end{tabular}}
\label{tab:US_inflation_quantile_regression_025}
\end{center}
\end{table}
%
\begin{table}[!ht]
\caption{Inflation dataset. Summary statistics of the posterior distribution from the quantile regression model at $\tau=0.50$, based on two inflation measures: CPIAUCSL on the left, and CPILFESL on the right.}
\begin{center}
\tabcolsep=2.0mm
\resizebox{0.95\columnwidth}{!}{\begin{tabular}{rrrrrrrr|rrrrrrr}
\toprule
&\multicolumn{7}{c}{CPIAUCSL} & \multicolumn{7}{c}{CPILFESL}\\
\cmidrule(lr){1-1}\cmidrule(lr){2-8}\cmidrule(lr){9-15}
 & mean & sd & \multicolumn{2}{c}{HPD$_{95\%}$} & CD & RHAT & ESS & mean & sd &\multicolumn{2}{c}{HPD$_{95\%}$} & CD & RHAT & ESS \\ 
\midrule
intercept & 0.02 & 0.04 & -0.05 & 0.09 & 0.08 & 1.00 & 38653.86 & -0.05 & 0.05 & -0.15 & 0.05 & -0.48 & 1.00 & 37258.09 \\ 
  UNRATE & -0.09 & 0.07 & -0.24 & 0.05 & -0.89 & 1.00 & 40000.00 & -0.12 & 0.10 & -0.31 & 0.07 & -1.47 & 1.00 & 40000.00 \\ 
  EC & 0.12 & 0.17 & -0.22 & 0.45 & -0.23 & 1.00 & 39872.73 & -0.29 & 0.26 & -0.80 & 0.22 & 0.59 & 1.00 & 38059.51 \\ 
  PRFI & 0.06 & 0.06 & -0.06 & 0.17 & 0.54 & 1.00 & 38831.28 & 0.02 & 0.08 & -0.14 & 0.18 & -1.05 & 1.00 & 38065.15 \\ 
  GDPC1 & -0.20 & 0.09 & -0.37 & -0.04 & 0.95 & 1.00 & 40000.00 & -0.05 & 0.13 & -0.29 & 0.20 & 0.80 & 1.00 & 40000.00 \\ 
  HOUST & 0.05 & 0.06 & -0.07 & 0.18 & 0.59 & 1.00 & 40000.00 & 0.00 & 0.09 & -0.18 & 0.19 & -0.81 & 1.00 & 40000.00 \\ 
  USPRIV & -0.12 & 0.19 & -0.49 & 0.26 & -1.02 & 1.00 & 40000.00 & 0.61 & 0.30 & 0.02 & 1.19 & 0.61 & 1.00 & 39141.42 \\ 
  \rowcolor{Gray}  
  TB3MS & 0.16 & 0.05 & 0.05 & 0.26 & -0.41 & 1.00 & 38450.20 & 0.17 & 0.07 & 0.02 & 0.31 & 1.43 & 1.00 & 39790.61 \\ 
  GS10 & -0.05 & 0.05 & -0.15 & 0.04 & 0.17 & 1.00 & 39462.73 & -0.02 & 0.07 & -0.16 & 0.11 & 0.55 & 1.00 & 37062.16 \\ 
  T10Y3MM & -0.02 & 0.22 & -0.46 & 0.41 & 0.36 & 1.00 & 40000.00 & -0.13 & 0.33 & -0.78 & 0.51 & -0.77 & 1.00 & 38312.90 \\ 
  T10YFFM & 0.11 & 0.22 & -0.33 & 0.54 & 1.79 & 1.00 & 40000.00 & 0.22 & 0.32 & -0.41 & 0.86 & -0.78 & 1.00 & 39778.29 \\ 
  M1SL & -0.01 & 0.09 & -0.18 & 0.16 & -0.16 & 1.00 & 40000.00 & -0.30 & 0.13 & -0.56 & -0.04 & -0.88 & 1.00 & 40000.00 \\ 
  MICH & \cellcolor{Gray}0.33 &\cellcolor{Gray} 0.08 & \cellcolor{Gray}0.18 & \cellcolor{Gray}0.48 & \cellcolor{Gray}-0.15 & \cellcolor{Gray}1.00 & \cellcolor{Gray}35640.61 & 0.18 & 0.10 & -0.02 & 0.37 & 2.00 & 1.00 & 36551.59 \\ 
  PPIACO & \cellcolor{Gray}0.49 & \cellcolor{Gray}0.09 & \cellcolor{Gray}0.32 & \cellcolor{Gray}0.66 & \cellcolor{Gray}-0.02 & \cellcolor{Gray}1.00 &\cellcolor{Gray}40000.00 & 0.01 & 0.15 & -0.29 & 0.29 & -0.90 & 1.00 & 34526.25 \\ 
  DJIA & -0.07 & 0.05 & -0.16 & 0.02 & 0.03 & 1.00 & 40000.00 & -0.02 & 0.07 & -0.15 & 0.11 & 0.87 & 1.00 & 37409.04 \\ 
  NAPMPMI & \cellcolor{Gray}0.13 & \cellcolor{Gray}0.06 & \cellcolor{Gray}0.01 & \cellcolor{Gray}0.25 & \cellcolor{Gray}-2.40 & \cellcolor{Gray}1.00 & \cellcolor{Gray}39937.07 & 0.04 & 0.09 & -0.14 & 0.22 & 0.58 & 1.00 & 36820.82 \\ 
  NAPMSDI & -0.09 & 0.07 & -0.21 & 0.04 & -1.45 & 1.00 & 35055.20 & -0.10 & 0.09 & -0.28 & 0.08 & 0.40 & 1.00 & 36562.28 \\ 
  OILPRICE & 0.09 & 0.07 & -0.06 & 0.23 & 0.74 & 1.00 & 38309.45 & 0.10 & 0.10 & -0.11 & 0.30 & 0.38 & 1.00 & 37895.93 \\ 
  GASPRICE & 0.18 & 0.10 & -0.02 & 0.38 & 0.82 & 1.00 & 33636.35 & -0.08 & 0.12 & -0.31 & 0.15 & -0.58 & 1.00 & 37705.93 \\ 
\bottomrule
\end{tabular}}
\label{tab:US_inflation_quantile_regression_050}
\end{center}
\end{table}
%
\begin{table}[!ht]
\caption{Inflation dataset. Summary statistics of the posterior distribution from the quantile regression model at $\tau=0.75$, based on two inflation measures: CPIAUCSL on the left, and CPILFESL on the right.}
\begin{center}
\tabcolsep=2.0mm
\resizebox{0.95\columnwidth}{!}{\begin{tabular}{rrrrrrrr|rrrrrrr}
\toprule
&\multicolumn{7}{c}{CPIAUCSL} & \multicolumn{7}{c}{CPILFESL}\\
\cmidrule(lr){1-1}\cmidrule(lr){2-8}\cmidrule(lr){9-15}
 & mean & sd & \multicolumn{2}{c}{HPD$_{95\%}$} & CD & RHAT & ESS & mean & sd &\multicolumn{2}{c}{HPD$_{95\%}$} & CD & RHAT & ESS \\ 
\midrule
  \rowcolor{Gray}
intercept & 0.29 & 0.04 & 0.21 & 0.36 & 1.74 & 1.00 & 32951.31 & 0.38 & 0.06 & 0.27 & 0.49 & -0.82 & 1.00 & 33013.25 \\ 
  UNRATE & -0.07 & 0.08 & -0.24 & 0.09 & -1.40 & 1.00 & 33340.50 & -0.03 & 0.13 & -0.28 & 0.22 & 0.56 & 1.00 & 32308.72 \\ 
  EC & 0.09 & 0.18 & -0.26 & 0.43 & 0.31 & 1.00 & 32594.47 & -0.24 & 0.24 & -0.72 & 0.24 & 1.49 & 1.00 & 36266.75 \\ 
  PRFI & 0.06 & 0.06 & -0.05 & 0.18 & -0.05 & 1.00 & 33969.31 & -0.09 & 0.08 & -0.25 & 0.06 & 0.40 & 1.00 & 37365.58 \\ 
  GDPC1 & \cellcolor{Gray}-0.20 & \cellcolor{Gray}0.09 & \cellcolor{Gray}-0.38 & \cellcolor{Gray}-0.02 & \cellcolor{Gray}0.14 & \cellcolor{Gray}1.00 & \cellcolor{Gray}34128.23 & -0.08 & 0.11 & -0.30 & 0.15 & -1.08 & 1.00 & 40000.00 \\ 
  HOUST & 0.01 & 0.07 & -0.12 & 0.14 & 0.58 & 1.00 & 35638.48 & 0.19 & 0.13 & -0.06 & 0.45 & 1.07 & 1.00 & 28518.78 \\ 
  USPRIV & -0.02 & 0.21 & -0.43 & 0.40 & 0.46 & 1.00 & 32972.46 & \cellcolor{Gray}0.67 & \cellcolor{Gray}0.30 & \cellcolor{Gray}0.07 & \cellcolor{Gray}1.25 & \cellcolor{Gray}2.05 & \cellcolor{Gray}1.00 & \cellcolor{Gray}35589.55 \\ 
  TB3MS & 0.07 & 0.05 & -0.03 & 0.16 & -0.56 & 1.00 & 36405.17 & -0.01 & 0.08 & -0.17 & 0.15 & -1.02 & 1.00 & 34322.40 \\ 
  GS10 & -0.07 & 0.05 & -0.17 & 0.03 & -0.58 & 1.00 & 31448.55 & -0.06 & 0.07 & -0.20 & 0.07 & -0.96 & 1.00 & 35307.54 \\ 
  T10Y3MM & -0.13 & 0.23 & -0.57 & 0.32 & -1.04 & 1.00 & 34202.66 & -0.11 & 0.33 & -0.75 & 0.56 & 0.32 & 1.00 & 35372.95 \\ 
  T10YFFM & 0.12 & 0.22 & -0.31 & 0.56 & 0.48 & 1.00 & 34789.85 & 0.24 & 0.32 & -0.40 & 0.87 & -2.06 & 1.00 & 36434.65 \\ 
  M1SL & -0.07 & 0.10 & -0.28 & 0.13 & -0.83 & 1.00 & 35216.03 &\cellcolor{Gray} -0.34 & \cellcolor{Gray}0.15 & \cellcolor{Gray}-0.65 &\cellcolor{Gray} -0.05 &\cellcolor{Gray} -1.02 &\cellcolor{Gray} 1.00 & \cellcolor{Gray}36117.69 \\ 
   \rowcolor{Gray} 
  MICH & 0.35 & 0.08 & 0.20 & 0.50 & 0.15 & 1.00 & 30922.58 & 0.26 & 0.11 & 0.05 & 0.47 & 0.13 & 1.00 & 32392.04 \\ 
  PPIACO & \cellcolor{Gray}0.36 & \cellcolor{Gray}0.11 &\cellcolor{Gray} 0.14 & \cellcolor{Gray}0.57 & \cellcolor{Gray}0.74 & \cellcolor{Gray}1.00 & \cellcolor{Gray}31005.83 & 0.05 & 0.15 & -0.25 & 0.34 & -0.26 & 1.00 & 31820.80 \\ 
  DJIA & -0.07 & 0.04 & -0.16 & 0.01 & 0.78 & 1.00 & 36925.69 & 0.06 & 0.06 & -0.06 & 0.18 & 1.22 & 1.00 & 39011.13 \\ 
  NAPMPMI & \cellcolor{Gray}0.14 & \cellcolor{Gray}0.06 & \cellcolor{Gray}0.01 & \cellcolor{Gray}0.26 & \cellcolor{Gray}-0.90 & \cellcolor{Gray}1.00 &\cellcolor{Gray}34007.31 & 0.03 & 0.09 & -0.15 & 0.20 & 2.63 & 1.00 & 36032.01 \\ 
  NAPMSDI & -0.06 & 0.06 & -0.17 & 0.05 & 0.07 & 1.00 & 35420.04 & -0.14 & 0.10 & -0.33 & 0.05 & -0.98 & 1.00 & 32717.83 \\ 
  OILPRICE & 0.09 & 0.07 & -0.04 & 0.23 & -0.26 & 1.00 & 35893.35 & -0.04 & 0.09 & -0.22 & 0.15 & -1.32 & 1.00 & 39985.40 \\ 
  GASPRICE & \cellcolor{Gray}0.29 & \cellcolor{Gray}0.11 &\cellcolor{Gray} 0.08 & \cellcolor{Gray}0.50 & \cellcolor{Gray}-0.09 & \cellcolor{Gray}1.00 & \cellcolor{Gray}31218.77 & 0.05 & 0.12 & -0.19 & 0.27 & 3.38 & 1.00 & 34787.51 \\ 
\bottomrule
\end{tabular}}
\label{tab:US_inflation_quantile_regression_075}
\end{center}
\end{table}
%
\begin{table}[!h]
\caption{Inflation dataset. Summary statistics of the posterior distribution from the quantile regression model at $\tau=0.90$, based on two inflation measures: CPIAUCSL on the left, and CPILFESL on the right.}
\begin{center}
\tabcolsep=2.0mm
\resizebox{0.95\columnwidth}{!}{\begin{tabular}{rrrrrrrr|rrrrrrr}
\toprule
&\multicolumn{7}{c}{CPIAUCSL} & \multicolumn{7}{c}{CPILFESL}\\
\cmidrule(lr){1-1}\cmidrule(lr){2-8}\cmidrule(lr){9-15}
 & mean & sd & \multicolumn{2}{c}{HPD$_{95\%}$} & CD & RHAT & ESS & mean & sd &\multicolumn{2}{c}{HPD$_{95\%}$} & CD & RHAT & ESS \\ 
\midrule
  \rowcolor{Gray}
intercept & 0.53 & 0.04 & 0.46 & 0.60 & 1.09 & 1.00 & 25970.01 & 0.85 & 0.06 & 0.74 & 0.97 & -0.10 & 1.00 & 26911.66 \\ 
  UNRATE & -0.07 & 0.09 & -0.25 & 0.12 & 0.14 & 1.00 & 21762.61 & 0.05 & 0.16 & -0.26 & 0.37 & 1.96 & 1.00 & 19969.68 \\ 
  EC & 0.19 & 0.20 & -0.19 & 0.58 & -0.21 & 1.00 & 21379.51 & -0.00 & 0.22 & -0.44 & 0.42 & -1.76 & 1.00 & 29158.31 \\ 
  PRFI & 0.10 & 0.06 & -0.01 & 0.21 & -1.48 & 1.00 & 23906.81 & -0.14 & 0.09 & -0.33 & 0.03 & 0.88 & 1.00 & 23996.31 \\ 
  GDPC1 & \cellcolor{Gray}-0.25 & \cellcolor{Gray}0.12 & \cellcolor{Gray}-0.49 & \cellcolor{Gray}-0.01 & \cellcolor{Gray}-0.35 & \cellcolor{Gray}1.00 & \cellcolor{Gray}19555.18 & -0.07 & 0.15 & -0.37 & 0.21 & 0.24 & 1.00 & 25092.79 \\ 
  HOUST & -0.03 & 0.08 & -0.18 & 0.12 & -0.07 & 1.00 & 22848.97 & 0.25 & 0.13 & -0.00 & 0.51 & -0.24 & 1.00 & 21035.48 \\ 
  USPRIV & -0.08 & 0.26 & -0.60 & 0.43 & -0.32 & 1.00 & 18817.75 & 0.36 & 0.36 & -0.33 & 1.07 & 0.66 & 1.00 & 22110.69 \\ 
  TB3MS & 0.00 & 0.06 & -0.12 & 0.13 & 0.20 & 1.00 & 20198.44 & \cellcolor{Gray}-0.24 & \cellcolor{Gray}0.08 & \cellcolor{Gray}-0.40 & \cellcolor{Gray}-0.09 & \cellcolor{Gray}-0.99 &\cellcolor{Gray} 1.00 & \cellcolor{Gray}28065.60 \\ 
  GS10 & -0.03 & 0.07 & -0.16 & 0.10 & 0.89 & 1.00 & 17799.89 & -0.07 & 0.08 & -0.22 & 0.09 & -1.05 & 1.00 & 25642.69 \\ 
  T10Y3MM & -0.13 & 0.24 & -0.60 & 0.33 & -1.07 & 1.00 & 24488.63 & -0.26 & 0.41 & -1.03 & 0.55 & 1.34 & 1.00 & 22170.82 \\ 
  T10YFFM & 0.10 & 0.23 & -0.34 & 0.57 & -0.03 & 1.00 & 24509.86 & 0.33 & 0.40 & -0.47 & 1.07 & 0.08 & 1.00 & 22799.52 \\ 
  M1SL & -0.05 & 0.13 & -0.30 & 0.20 & 0.43 & 1.00 & 19451.49 & -0.23 & 0.21 & -0.66 & 0.14 & 0.06 & 1.00 & 18862.21 \\ 
   \rowcolor{Gray} 
  MICH & 0.40 & 0.08 & 0.25 & 0.55 & -1.06 & 1.00 & 21317.26 & 0.30 & 0.11 & 0.09 & 0.53 & -0.87 & 1.00 & 23578.84 \\ 
  PPIACO & 0.29 & 0.15 & -0.00 & 0.58 & 1.16 & 1.00 & 16884.43 & \cellcolor{Gray}0.30 &\cellcolor{Gray} 0.12 & \cellcolor{Gray}0.07 & \cellcolor{Gray}0.53 & \cellcolor{Gray}0.01 & \cellcolor{Gray}1.00 & \cellcolor{Gray}29624.79 \\ 
  DJIA & -0.08 & 0.05 & -0.18 & 0.02 & -0.51 & 1.00 & 22684.78 & 0.09 & 0.07 & -0.05 & 0.23 & -0.27 & 1.00 & 25383.50 \\ 
  NAPMPMI & 0.10 & 0.09 & -0.07 & 0.27 & 0.03 & 1.00 & 16899.49 & 0.08 & 0.12 & -0.16 & 0.30 & -0.17 & 1.00 & 19874.40 \\ 
  NAPMSDI & -0.06 & 0.07 & -0.19 & 0.07 & -0.26 & 1.00 & 19733.80 & -0.22 & 0.12 & -0.45 & 0.02 & -0.15 & 1.00 & 19727.79 \\ 
  OILPRICE & 0.09 & 0.09 & -0.08 & 0.26 & 0.14 & 1.00 & 21259.95 & -0.18 & 0.13 & -0.43 & 0.07 & -0.85 & 1.00 & 21676.71 \\ 
  GASPRICE & \cellcolor{Gray}0.33 & \cellcolor{Gray}0.15 & \cellcolor{Gray}0.04 & \cellcolor{Gray}0.62 & \cellcolor{Gray}-0.97 & \cellcolor{Gray}1.00 & \cellcolor{Gray}15340.46 & 0.15 & 0.11 & -0.07 & 0.36 & 1.93 & 1.00 & 27521.01 \\ 
\bottomrule
\end{tabular}}
\label{tab:US_inflation_quantile_regression_090}
\end{center}
\end{table}
%
\newpage
\clearpage
\subsection{BDQMA: summary of the relevant parameters}
\label{sec:US_infl_additional_results_relevant_parameters}
%
\noindent In this section, for the US inflation dataset, we report two tables summarizing the most relevant regression parameters, e.g. those parameters whose inclusion probabilities are greater than $0.7$ for at least one period (quarter).
\begin{table}[!h]
\setlength{\tabcolsep}{5 pt}
\caption{Inflation dataset. Summary of the relevant covariates for the response variable \qmo\textsf{CPIAUCSL}\qmc. Here, by relevant covariates we mean those having inclusion probability (as estimated by the BDQMA model for quantile and Gaussian regression) larger than $0.7$ for at least one quarter.} 
\begin{center}\resizebox{0.5\columnwidth}{!}{\begin{tabular}{llcccccc}\\
\toprule
&&\multirow{2}{*}{Mean regression}&\multicolumn{5}{c}{Quantile regression}\\
Id& Name&& $\tau=0.10$ & $\tau=0.25$ & $\tau=0.50$ & $\tau=0.75$ & $\tau=0.90$\\
\cmidrule(lr){1-1}\cmidrule(lr){2-2}\cmidrule(lr){3-3}\cmidrule(lr){4-4}\cmidrule(lr){5-5}\cmidrule(lr){6-6}\cmidrule(lr){7-7}\cmidrule(lr){8-8}
%
4	&	\textsf{UNRATE} 	& &&&&&\\
\rowcolor{Gray}
5	&	\textsf{EC}		&	&&&&\checkmark&\checkmark\\
6	&	\textsf{PRFI}		&	\\
\rowcolor{Gray}
7	&	\textsf{GDPC1}		&	&&\checkmark&\checkmark&\cellcolor{Gray}&\cellcolor{Gray}\\
8	&	\textsf{HOUST}		&	\\	
9	&	\textsf{USPRIV}	&	\\
\rowcolor{Gray}
10	&	\textsf{TB3MS}		&	& \checkmark&\checkmark&\checkmark&\cellcolor{Gray}&\cellcolor{Gray}\\
11	&	\textsf{GS10}		&	\\
12	&	\textsf{T10Y3MM}	&	\\
13	&	\textsf{T10YFFM}	&	\\
\rowcolor{Gray}
14	&	\textsf{M1SL}		&	& \checkmark&\checkmark&\checkmark&\cellcolor{Gray}&\cellcolor{Gray}\\
\rowcolor{Gray}
15	&	\textsf{MICH}		&\checkmark	&  
\checkmark&\checkmark&\checkmark&\checkmark&\checkmark\\
\rowcolor{Gray}
16	&	\textsf{PPIACO}	& \checkmark& \checkmark	&\checkmark&\checkmark&\checkmark&\checkmark\\
\rowcolor{Gray}
17	&	\textsf{DJIA}		&	&&\checkmark&\checkmark&\checkmark&\cellcolor{Gray}\\
18	&	\textsf{PMI}		&	\\
19	&	\textsf{NAPMSDI}	&	\\
\rowcolor{Gray}
20	&	\textsf{OILPRICE}	&	&&\checkmark&\checkmark&\cellcolor{Gray}&\cellcolor{Gray}\\
\rowcolor{Gray}
21	&	\textsf{GASPRICE}	&	& \checkmark&\checkmark&\checkmark&\checkmark&\checkmark\\
%
\bottomrule 
\end{tabular}} 
\label{tab:table_US_inflation_data_CPIAUCSL_ss} 
\end{center} 
\end{table} 
%
\begin{table}[!h]
\setlength{\tabcolsep}{5 pt}
\caption{Inflation dataset. Summary of the relevant covariates for the response variable \qmo\textsf{CPILFESL}\qmc. Here, by relevant covariates we mean those having inclusion probability (as estimated by the BDQMA model for Gaussian and quantile regression) larger than $0.7$ for at least one quarter.} 
\begin{center}\resizebox{0.5\columnwidth}{!}{\begin{tabular}{llcccccc}\\
\toprule
&&\multirow{2}{*}{Mean regression}&\multicolumn{5}{c}{Quantile regression}\\
Id& Name&& $\tau=0.10$ & $\tau=0.25$ & $\tau=0.50$ & $\tau=0.75$ & $\tau=0.90$\\
\cmidrule(lr){1-1}\cmidrule(lr){2-2}\cmidrule(lr){3-3}\cmidrule(lr){4-4}\cmidrule(lr){5-5}\cmidrule(lr){6-6}\cmidrule(lr){7-7}\cmidrule(lr){8-8}
%
\rowcolor{Gray}
4	&	\textsf{UNRATE} 	&\checkmark &\checkmark&\checkmark&\checkmark&\checkmark&\checkmark\\
\rowcolor{Gray}
5	&	\textsf{EC}		&	&&\checkmark&\checkmark&\checkmark&\checkmark\\
6	&	\textsf{PRFI}		&	\\
7	&	\textsf{GDPC1}		&	&&&&&\\
\rowcolor{Gray}
8	&	\textsf{HOUST}		&\checkmark & \cellcolor{Gray}& \cellcolor{Gray}& \cellcolor{Gray}& \cellcolor{Gray}& \cellcolor{Gray}	\\	
\rowcolor{Gray}
9	&	\textsf{USPRIV}	&&\checkmark&&&&\checkmark	\\
\rowcolor{Gray}
10	&	\textsf{TB3MS}		&\checkmark	& \checkmark&\checkmark&\checkmark&\checkmark&\checkmark\\
11	&	\textsf{GS10}		&	\\
\rowcolor{Gray}
12	&	\textsf{T10Y3MM}	&&\checkmark &\checkmark&\checkmark&\checkmark&\checkmark	\\
\rowcolor{Gray}
13	&	\textsf{T10YFFM}	&&\checkmark	&\checkmark&\checkmark&\checkmark&\checkmark\\
\rowcolor{Gray}
14	&	\textsf{M1SL}		&\checkmark	& \checkmark&\checkmark&\checkmark&\checkmark&\checkmark\\
\rowcolor{Gray}
15	&	\textsf{MICH}		&\checkmark	&  
\checkmark&\checkmark&\checkmark&\checkmark&\checkmark\\
16	&	\textsf{PPIACO}	& & &&&&\\
\rowcolor{Gray}
17	&	\textsf{DJIA}		&	&&\checkmark&\checkmark&&\cellcolor{Gray}\\
18	&	\textsf{PMI}		&	\\
19	&	\textsf{NAPMSDI}	&	\\
\rowcolor{Gray}
20	&	\textsf{OILPRICE}	&	&\checkmark&\checkmark&\checkmark&\checkmark&\cellcolor{Gray}\\
21	&	\textsf{GASPRICE}	&	& &&&&\\
%
\bottomrule 
\end{tabular}} 
\label{tab:table_US_inflation_data_CPILFESL_ss} 
\end{center} 
\end{table} 
%
\newpage
\subsection{BDQMA: results for CPILFESL}
\label{sec:US_infl_additional_results_CPILFESL}
\begin{figure}[t]
\begin{center}
\resizebox*{0.9\textwidth}{!}{\subfigure[$\tau=0.10$]{\label{fig:CPILFESL_beta010_relevant}
\includegraphics[trim={0 1cm 0.5cm 1.5cm},clip,width=0.8\textwidth]{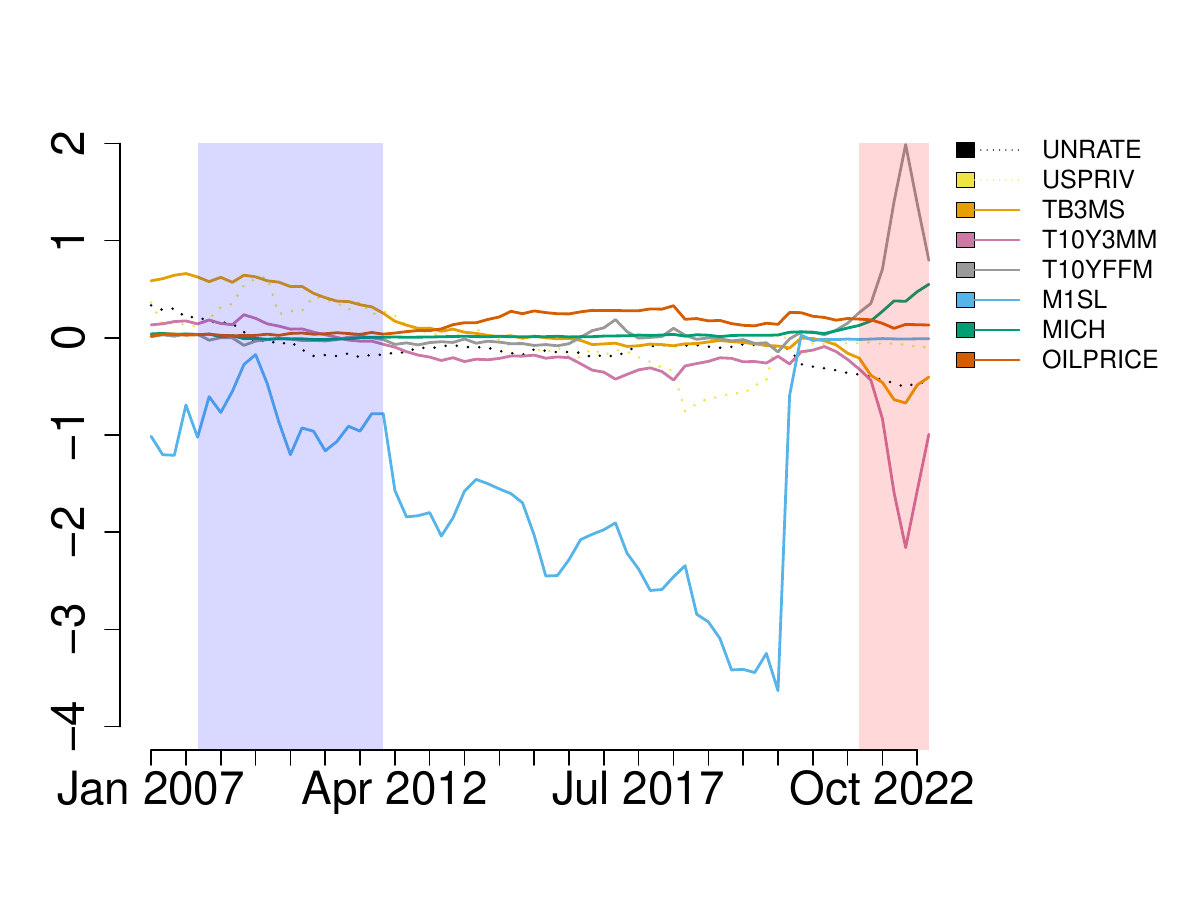}}
\subfigure[$\tau=0.25$]{\label{fig:CPILFESL_beta025_relevant}
\includegraphics[trim={0 1cm 0.5cm 1.5cm},clip,width=0.8\textwidth]{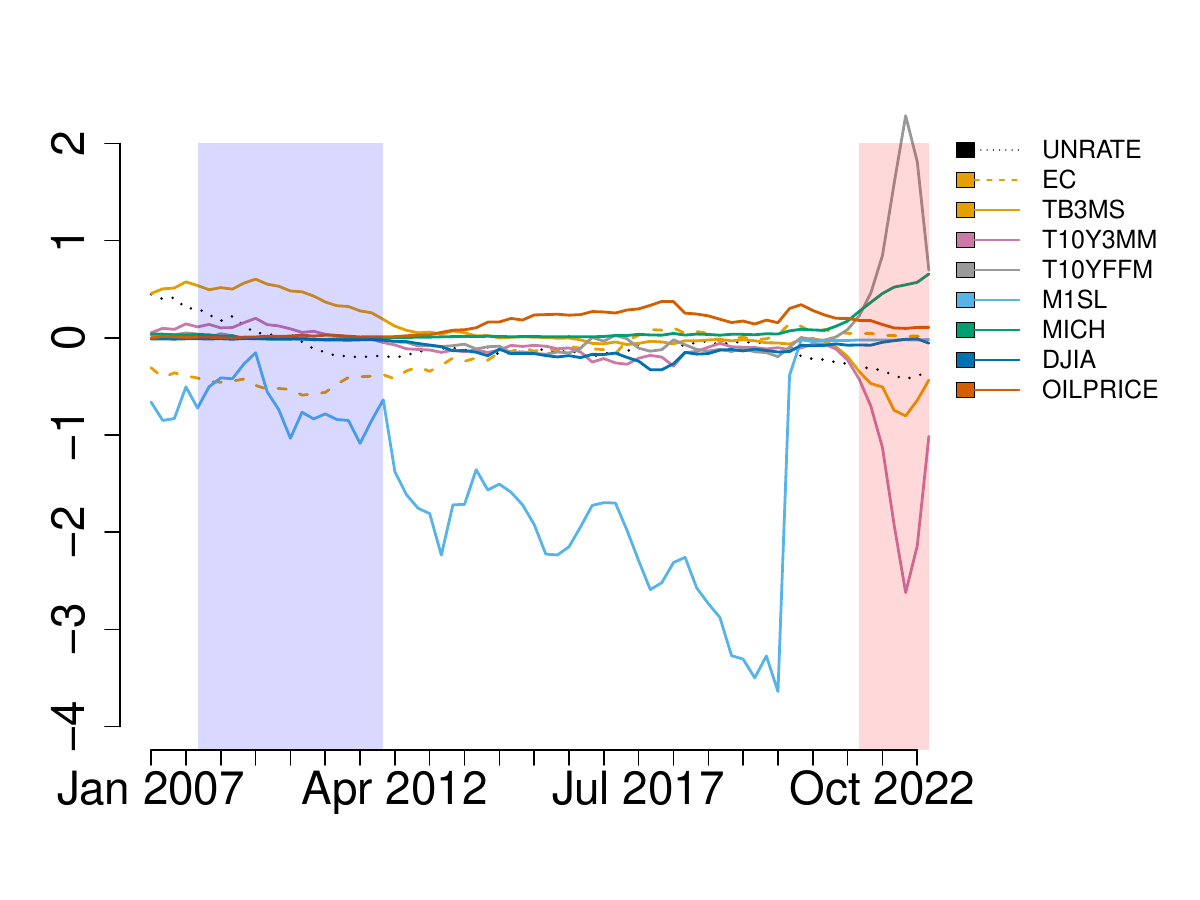}}}
\resizebox*{0.9\textwidth}{!}{\subfigure[$\tau=0.50$]{\label{fig:CPILFESL_beta050_relevant}
\includegraphics[trim={0 1cm 0.5cm 1.5cm},clip,width=0.8\textwidth]{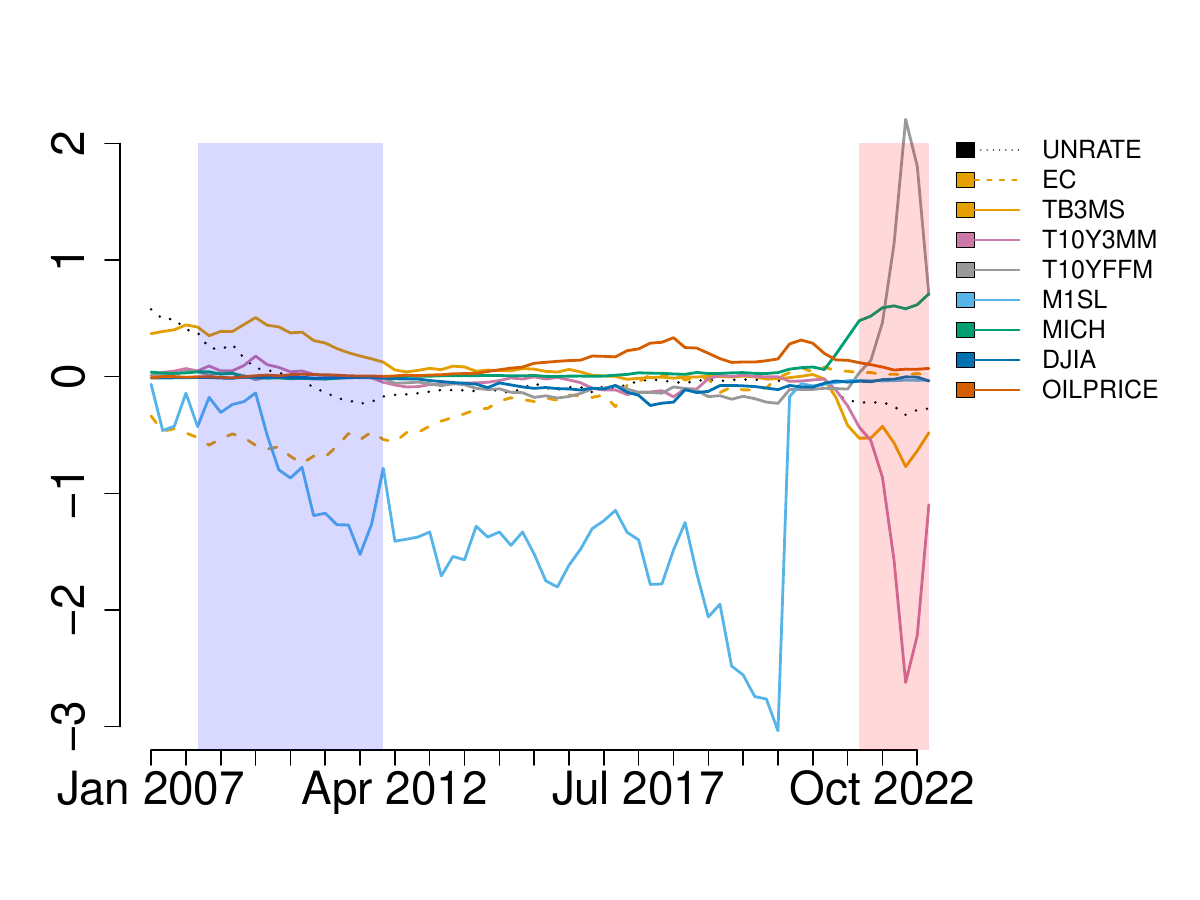}}
\subfigure[$\tau=0.75$]{\label{fig:CPILFESL_beta075_relevant}
\includegraphics[trim={0 1cm 0.5cm 1.5cm},clip,width=0.8\textwidth]{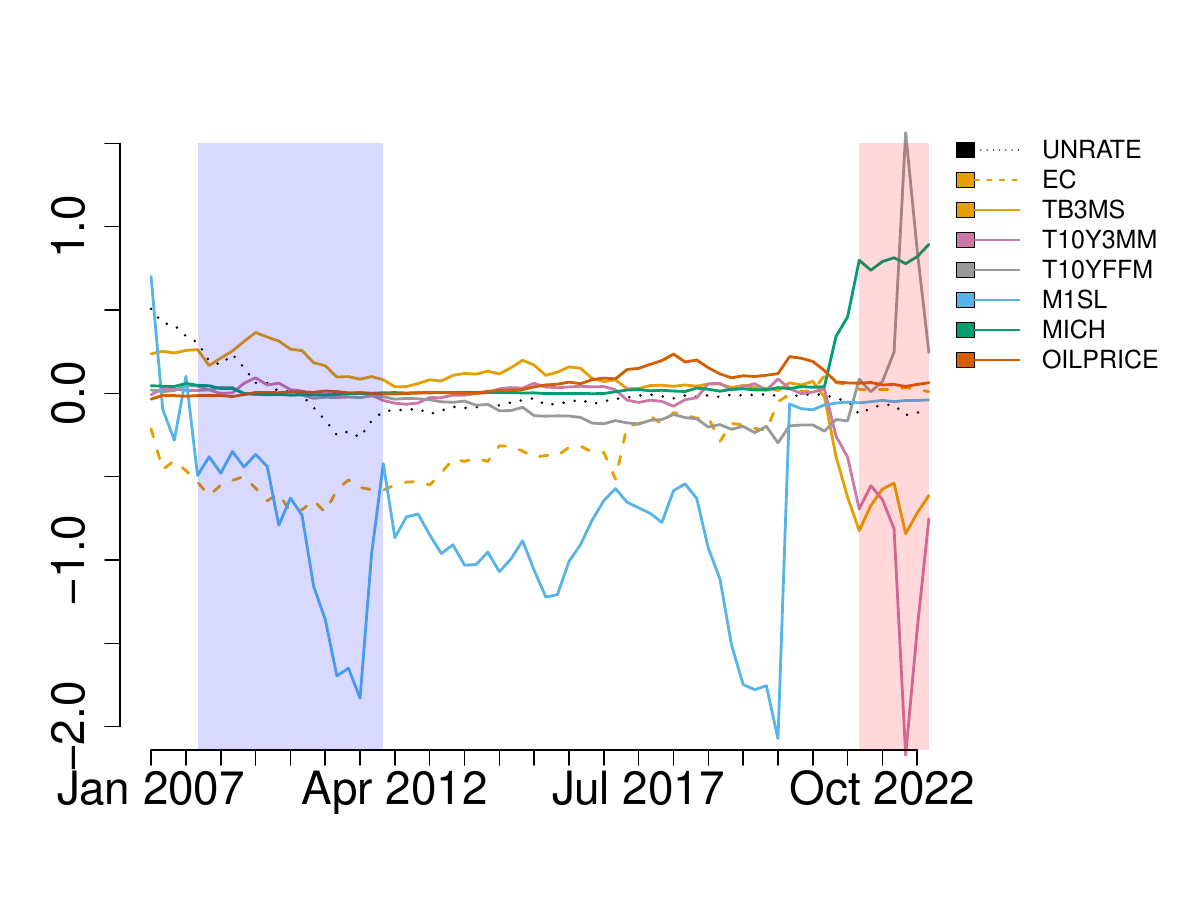}}}
\resizebox*{0.9\textwidth}{!}{\subfigure[$\tau=0.90$]{\label{fig:CPILFESL_beta090_relevant}
\includegraphics[trim={0 1cm 0.5cm 1.5cm},clip,width=0.8\textwidth]{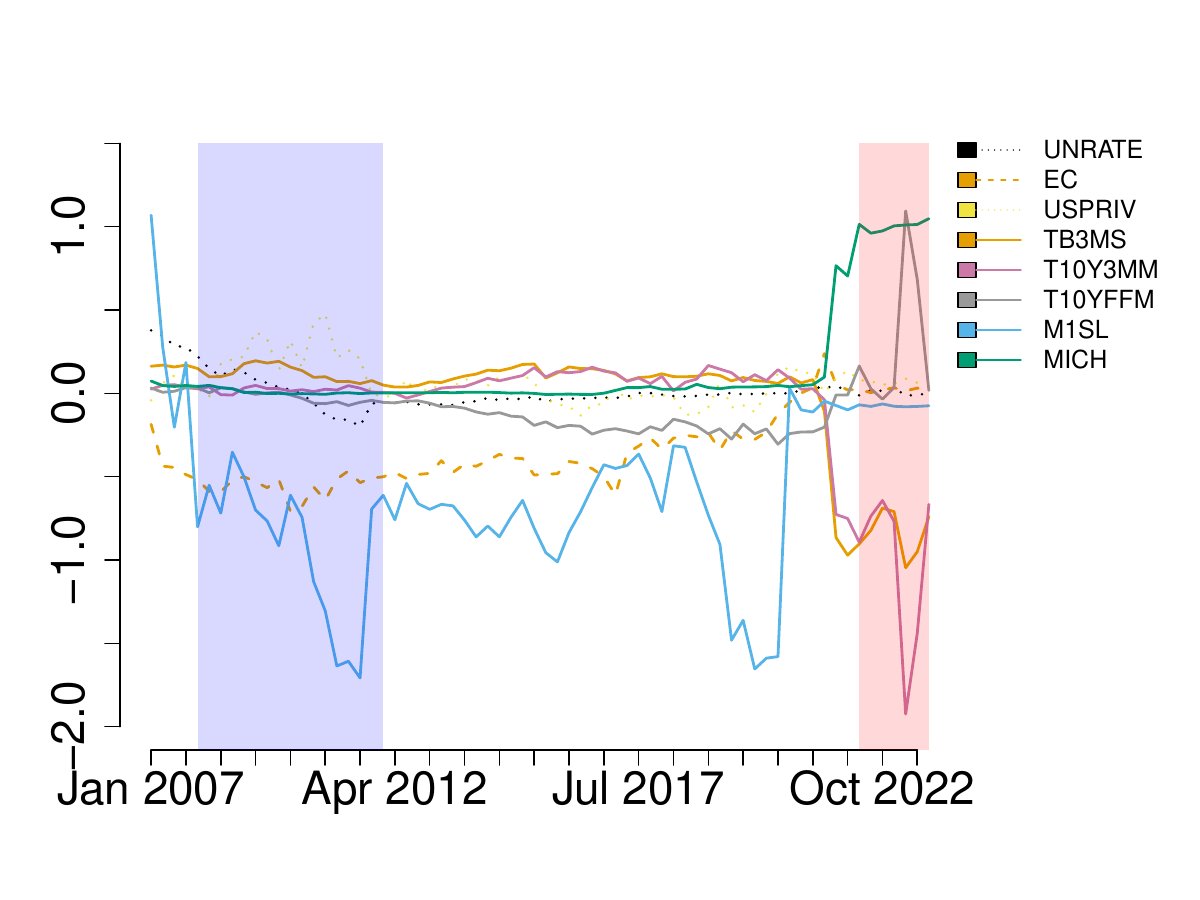}}
\subfigure[Mean regression]{\label{fig:CPILFESL_beta_relevant}
\includegraphics[trim={0 1cm 0.5cm 1.5cm},clip,width=0.8\textwidth]{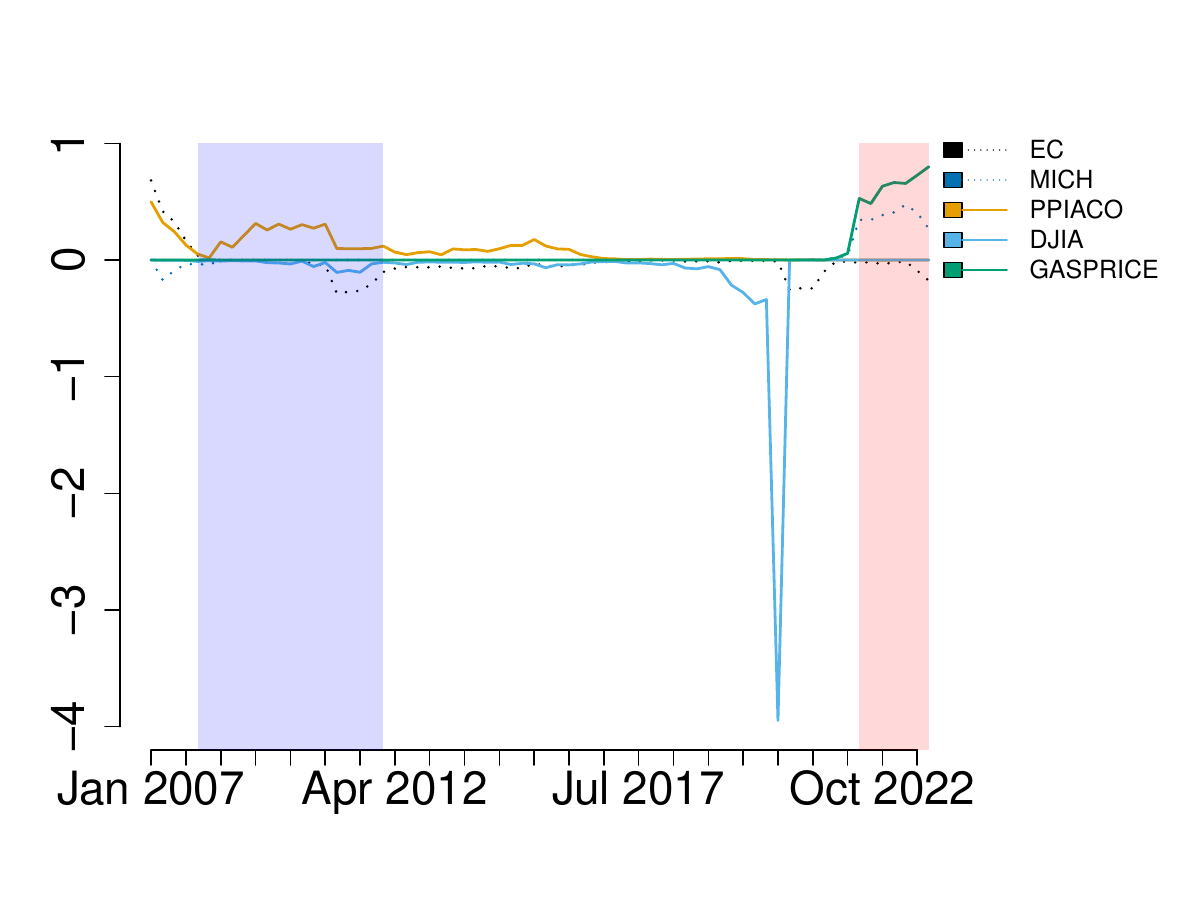}}}
\caption{US inflation data. Sequential estimates of the regression parameters by DMA (for the mean regression \ref{fig:CPILFESL_beta_relevant}) and BDQMA (for the quantile regression \ref{fig:CPILFESL_beta010_relevant}-\ref{fig:CPILFESL_beta090_relevant}) for the CPILFESL (Consumer Price Index For All Urban Consumers: all item excluding food and energy, US city average). For each quantile level $\tau$ the corresponding figure only reports those parameters having inclusion probability larger or equal to $0.7$ for at least one quarter. See Table \ref{tab:table_US_inflation_data_CPILFESL_ss} for a summary of the relevant covariate. The shaded areas identify two significant periods: the US Great Financial Crisis from 2007-Q4 to 2011-Q4 {\it (blue)}, and the Russian-Ukraine crisis {\it (red)}.}
\label{fig:inflation_RegDyn_CPILFESL}
\end{center}
\end{figure}

\begin{figure}[h]
\begin{center}
\resizebox*{0.9\textwidth}{!}{\subfigure[$\tau=0.10$]{\label{fig:CPILFESL_incl_prob010}
\includegraphics[trim={0 1cm 0.5cm 1.5cm},clip,width=0.8\textwidth]{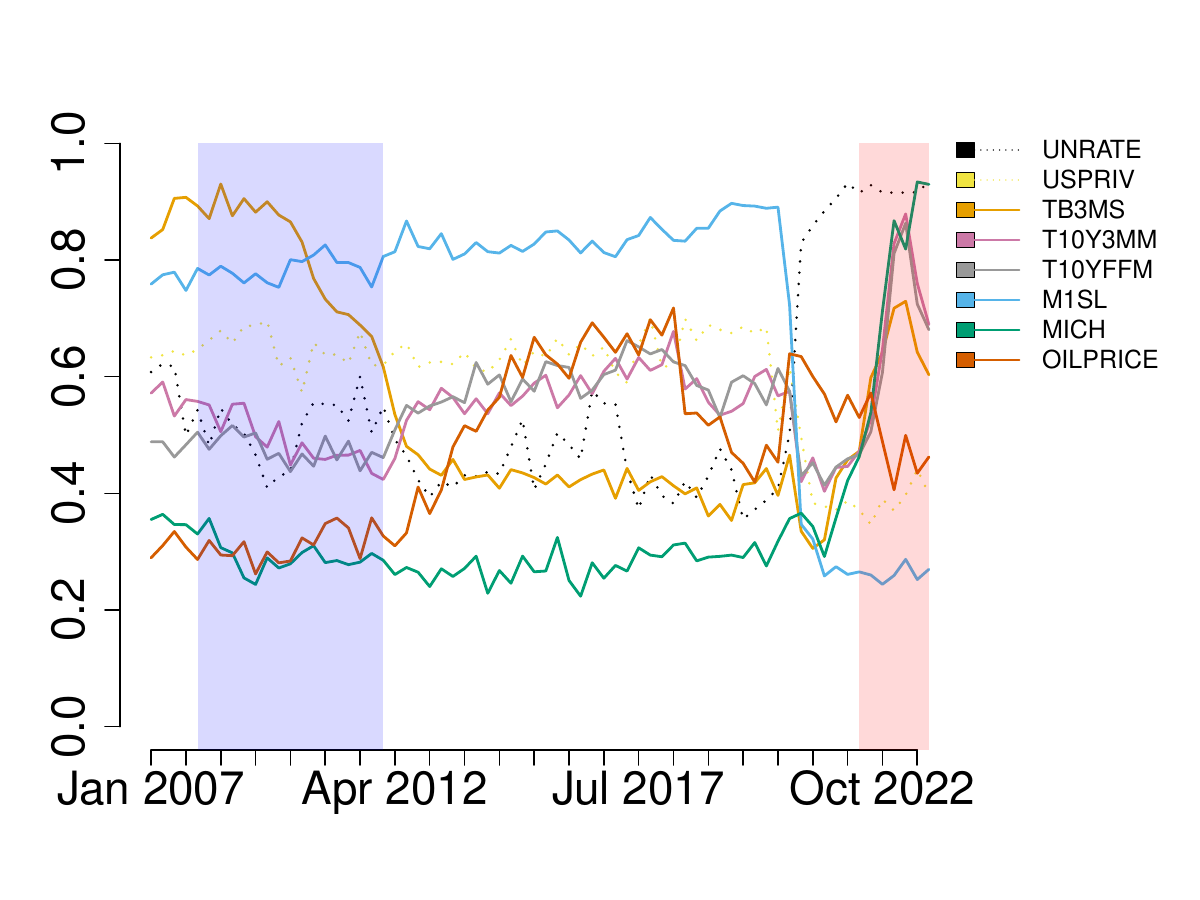}}
\subfigure[$\tau=0.25$]{\label{fig:CPILFESL_incl_prob025}
\includegraphics[trim={0 1cm 0.5cm 1.5cm},clip,width=0.8\textwidth]{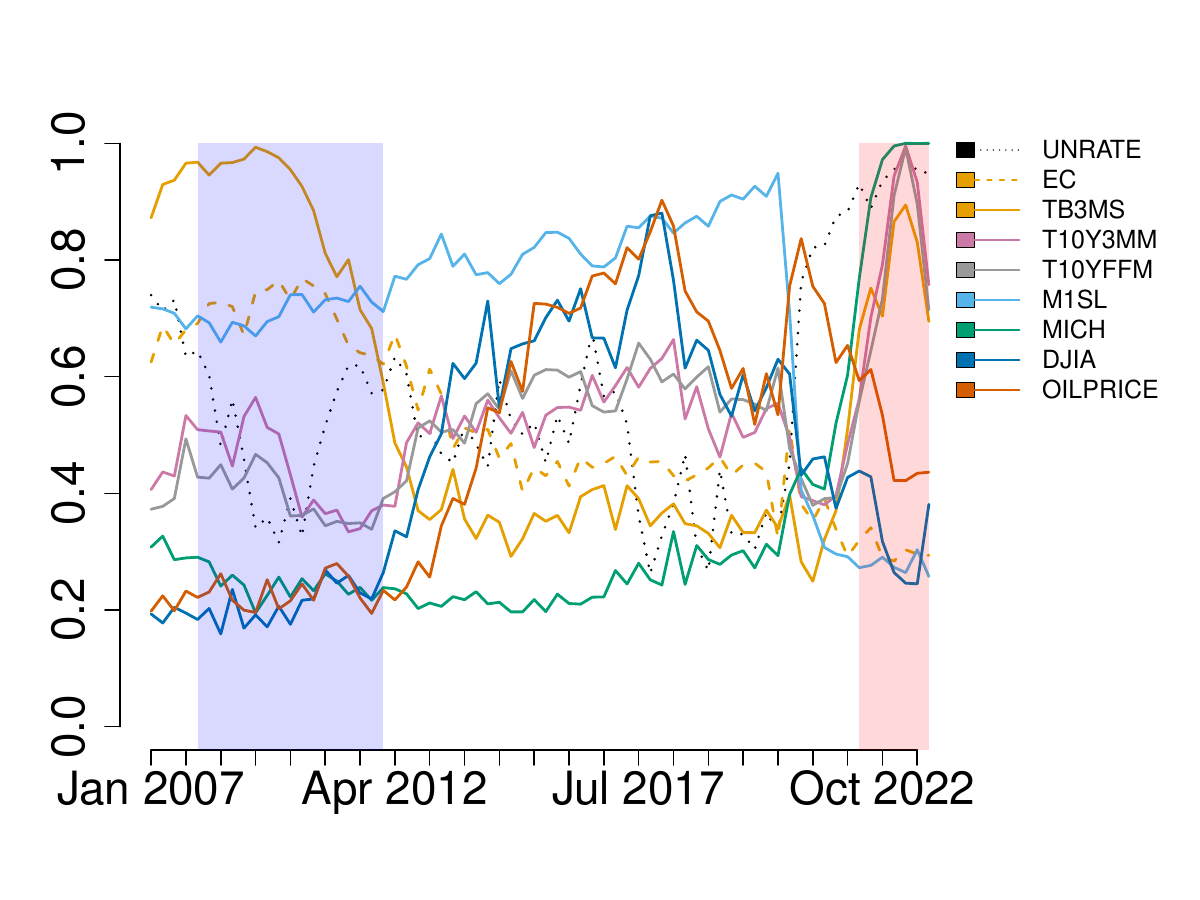}}}
\resizebox*{0.9\textwidth}{!}{\subfigure[$\tau=0.50$]{\label{fig:CPILFESL_incl_prob050}
\includegraphics[trim={0 1cm 0.5cm 1.5cm},clip,width=0.8\textwidth]{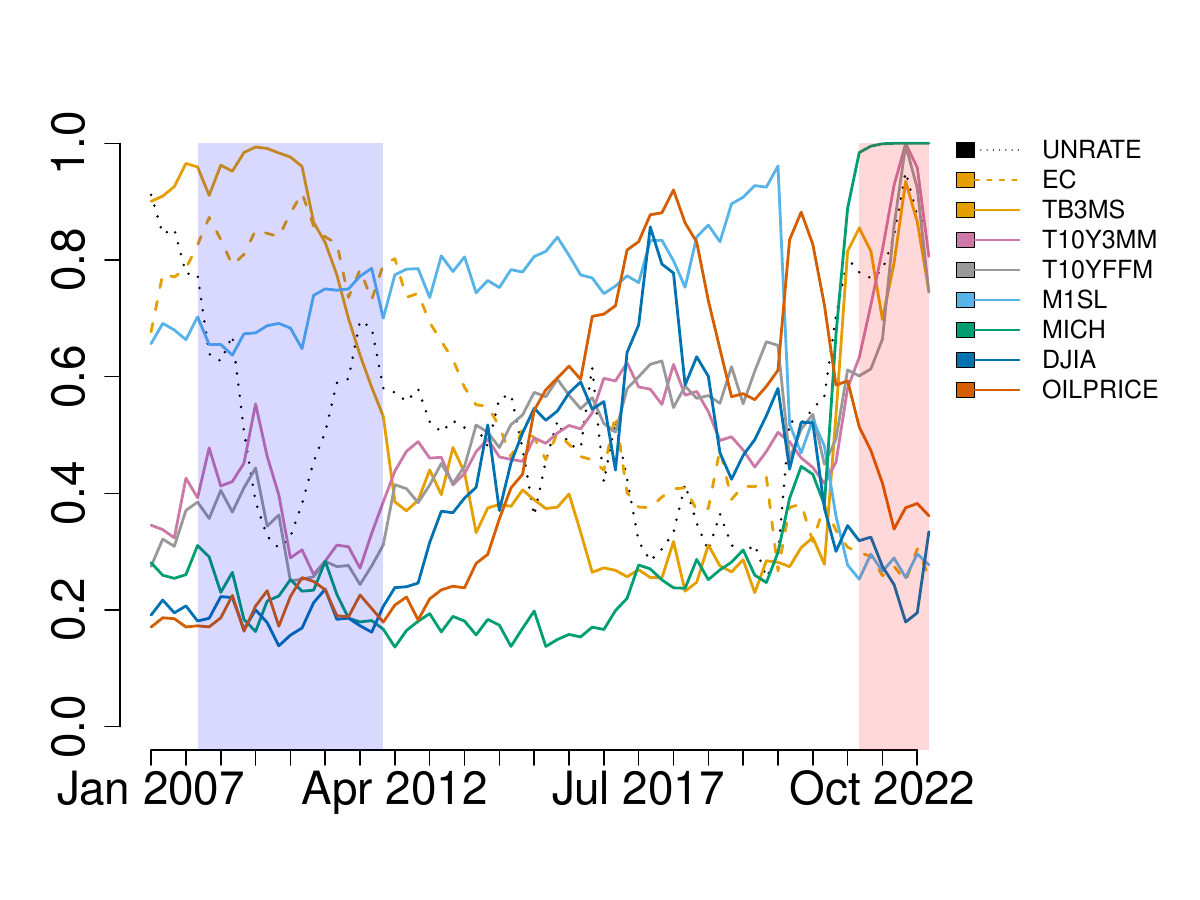}}
\subfigure[$\tau=0.75$]{\label{fig:CPILFESL_incl_prob075}
\includegraphics[trim={0 1cm 0.5cm 1.5cm},clip,width=0.8\textwidth]{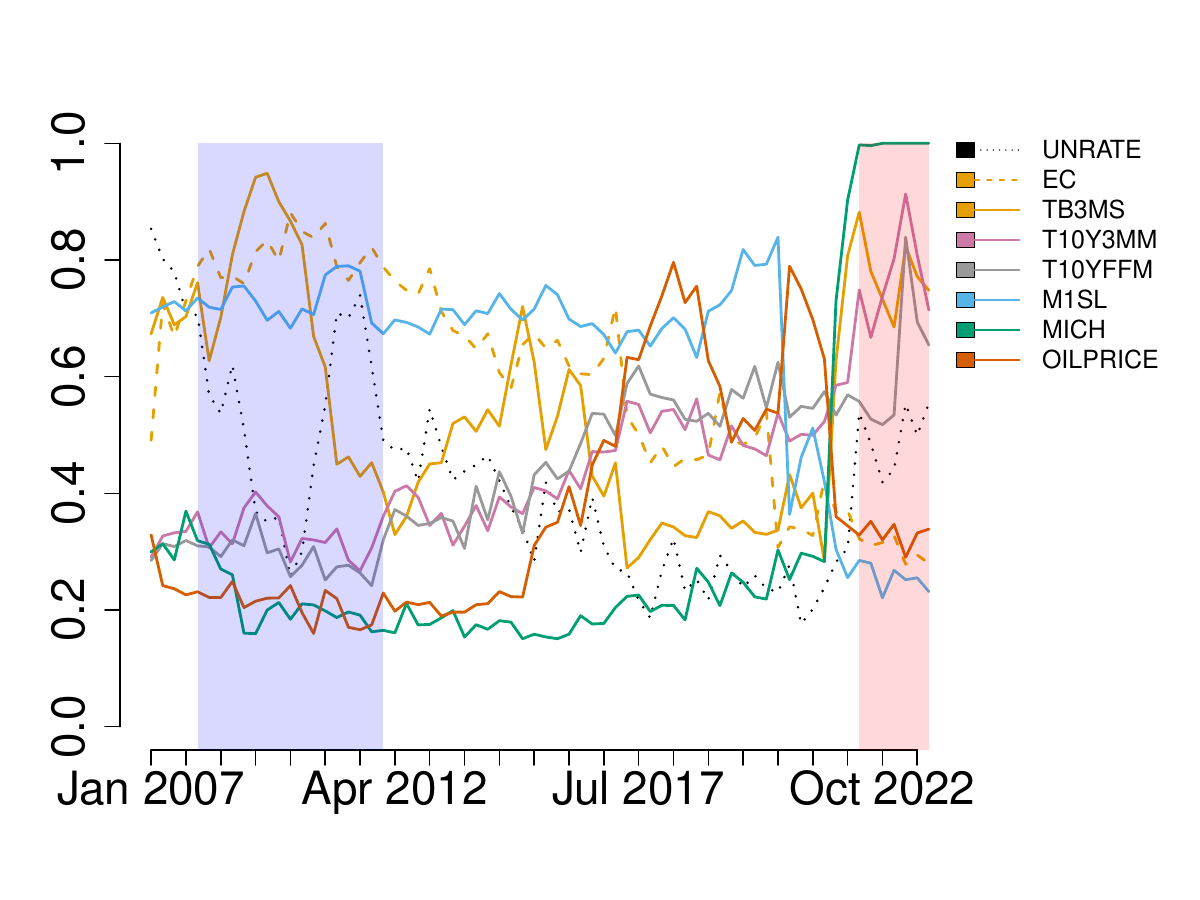}}}
\resizebox*{0.9\textwidth}{!}{\subfigure[$\tau=0.90$]{\label{fig:CPILFESL_incl_prob090}
\includegraphics[trim={0 1cm 0.5cm 1.5cm},clip,width=0.8\textwidth]{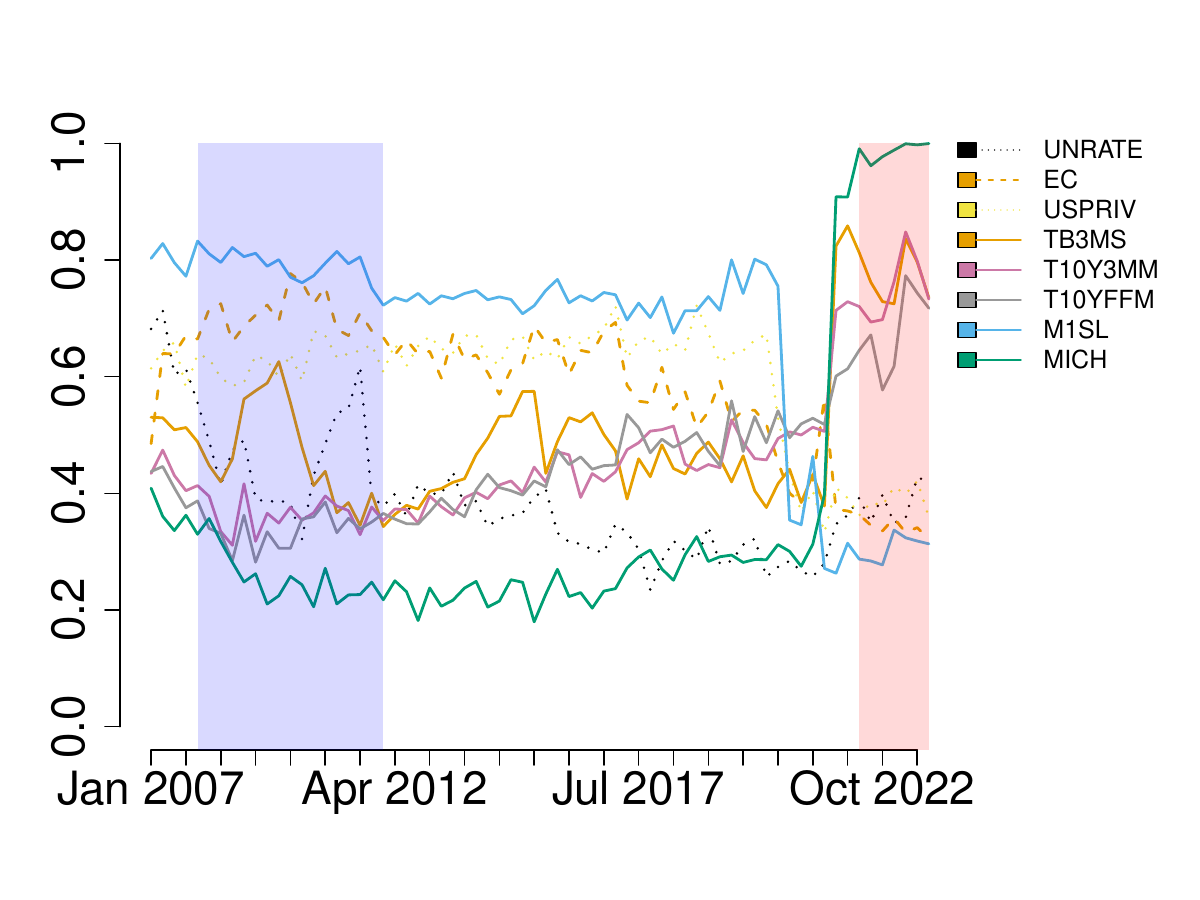}}
\subfigure[Mean regression]{\label{fig:CPILFESL_incl_prob}
\includegraphics[trim={0 1cm 0.5cm 1.5cm},clip,width=0.8\textwidth]{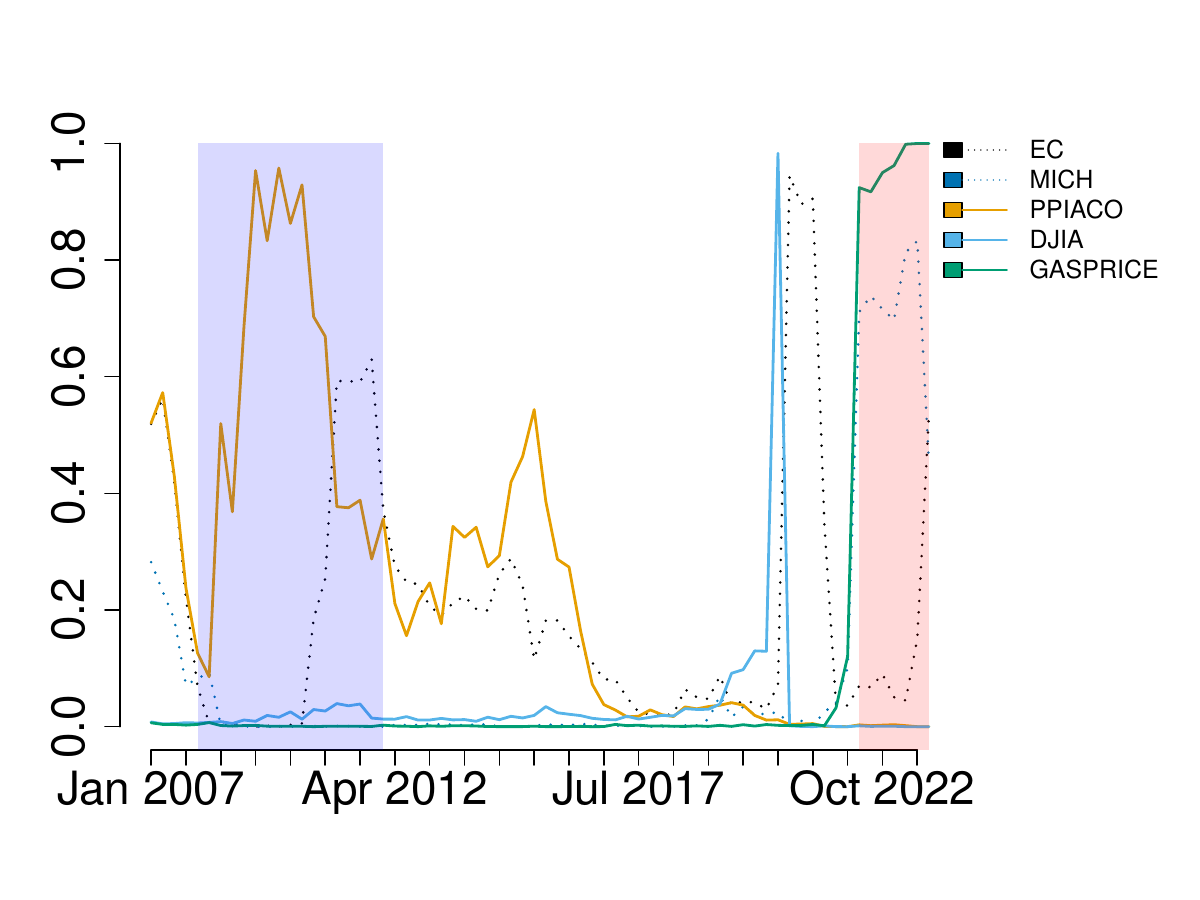}}}
\caption{US inflation dataset. Sequential update of the predicted inclusion
probabilities $\protect\pi _{t|t-1}$ by DMA (for the mean regression \ref{fig:CPILFESL_incl_prob}) and BDQMA (for the quantile regression \ref{fig:CPILFESL_incl_prob010}-\ref{fig:CPILFESL_incl_prob090}) for the CPILFESL (Consumer Price Index For All Urban Consumers: all item excluding food and energy, US city average). For each quantile level $\tau$ the corresponding figure only reports those parameters having inclusion probability larger or equal to $0.7$ for at least one quarter. See Table \ref{tab:table_US_inflation_data_CPILFESL_ss} for a summary of the relevant covariate. The shaded areas identify two significant periods: the US Great Financial Crisis from 2007-Q4 to 2011-Q4 {\it (blue)}, and the Russian-Ukraine crisis {\it (red)}.}
\label{fig:inflation_InclProb_CPILFESL}
\end{center}
\end{figure}

%
\noindent 
In this section, we present the time-varying parameters (Fig. \ref{fig:inflation_RegDyn_CPILFESL}) and inclusion probabilities (Fig. \ref{fig:inflation_InclProb_CPILFESL})  for the US inflation dataset, derived from our BDQMA model across all quantile levels $\tau=\left(0.10, 0.25, 0.5, 0.75, 0.90\right)$, using the CPILFESL inflation definition. These figures highlight the most significant regression parameters, specifically those with inclusion probabilities exceeding $0.7$ for at least one period (quarter).\newline

\noindent A visual inspection of Fig. \ref{fig:inflation_InclProb_CPILFESL} reveals two key insights: first, there are significantly fewer regressors deemed important when examining the mean compared to the quantiles. This discrepancy may be attributed to the mean's inherent tendency to smooth out variations and nuances present in the data, leading to the exclusion of significant predictors that contribute to specific quantile behavior. Second, the relevant predictors exhibit remarkable consistency across different quantiles, indicating that certain factors consistently influence inflation regardless of its distribution. Moreover, the patterns of both inclusion probabilities and parameter estimates follow a similar trajectory across quantiles, suggesting stable relationships throughout the analyzed period. Interestingly, a greater number of regressors are found to be significant when employing the CPILFESL definition of inflation compared to CPIAUCSL, highlighting CPILFESL's heightened sensitivity to a broader array of predictors. Additionally, it is noteworthy that financial variables such as the three-month Treasury Bill rate (TB3MS) and the two spreads (T10Y3MM and T10YFFM) consistently emerge as relevant regressors in the CPILFESL model. In contrast, these variables have not been considered significant in explaining the CPIAUCSL, as evidenced in Fig. \ref{fig:inflation_RegDyn} and Fig. \ref{fig:inflation_InclProb} in the main paper. This suggests that the CPILFESL definition may capture more complex dynamics related to financial conditions that influence inflation trends, further emphasizing the importance of selecting appropriate measures for inflation analysis.

%
%
%

\section{Additional results: US real estate analysis}
\label{sec:US_realestate_additional_results}
%
\subsection{Gaussian and quantile regression on the whole sample}
\label{sec:US_realestate_additional_results_static}
%
\noindent We report the results for the Gaussian and quantile regression models estimated on the whole sample: posterior mean (mean), standard deviation (sd), $95\%$ high posterior credible sets (HPD$_{95\%}$) and MCMC diagnostics (Geweke test of \cite{geweke.1992} (CD), Rhat statistic of \cite{gelman_rubin.1992} (RHAT) and effective sample size, ESS).  Statistics are based on $10,000$ draws, after $25,000$ burn-in draws, for $4$ independent chains. Shaded areas highlight those parameters for which the HPD$_{95\%}$ does not include zero. The lagged endogenous variable coefficients $(\phi_1,\phi_2)$ and the scale parameter $\sigma$ are not reported since they are nuisance parameters. The results for \qmo REITMKT\qmcsp are reported and the one for \qmo REITTBL\qmcsp and \qmo REITSM\qmcsp are similar and available upon request.
\begin{table}[h]
\setlength{\tabcolsep}{5 pt}
\caption{Real estate dataset (REITMKT). Posterior mean of the
significant parameters (selected by using the 95\% HPD intervals) of
the Gaussian and quantile regression model, estimated on the whole sample for
different quantile levels (columns). HPD$_{95\%}$ are reported in italics. The estimates for the Gaussian regression and for quantile regression, for each quantile level $\tau=(0.1,0.25,0.5,0.75,0.9)$, are reported in Tables \ref{tab:US_realestate_mean_regression}--\ref{tab:US_realestate_quantile_regression_090}.}
\begin{center}\resizebox{0.70\columnwidth}{!}{\begin{tabular}{lcccccc}\\
\toprule
&\multirow{2}{*}{Mean regression}&\multicolumn{5}{c}{Quantile regression}\\
Name&& $\tau=0.10$ & $\tau=0.25$ & $\tau=0.50$ & $\tau=0.75$ & $\tau=0.90$\\
\cmidrule(lr){1-1}\cmidrule(lr){2-2}\cmidrule(lr){3-3}\cmidrule(lr){4-4}\cmidrule(lr){5-5}\cmidrule(lr){6-6}\cmidrule(lr){7-7}
	\multirow{2}{*}{\textsf{intercept}} 	& &$-1.06$&$-0.56$&&$0.57$&$1.07$\\
	& &$\it(-1.14,  -0.98)$&  $\it(-0.65,  -0.47)$&& $\it(0.47,  0.66)$&$\it(0.99,  1.16)$\\
	\multirow{2}{*}{\textsf{TERM}} 	&$-0.18$ &$-0.30$&&&&\\
	&$\it(-0.36,-0.01)$ &$\it(-0.49 , -0.12)$ &&&&\\
	\multirow{2}{*}{\textsf{PREM}}		& &$-0.32308$&&&&\\
	& &$\it(-0.50,  -0.14)$&&&&\\
	\multirow{2}{*}{\textsf{TBILL}}		&$-0.45$ & $-0.65$ &$-0.40$&&&\\
	&$\it(-0.75,  -0.15)$ &$\it(-0.89,  -0.40)$&$\it(-0.70,  -0.09)$&&&\\
	\multirow{2}{*}{\textsf{MKTPE}}		& &&&&&\\
	& &&&&&\\
	\multirow{2}{*}{\textsf{MKTYLD}}	& $-0.12$ &$-0.13$&$-0.15$&$-0.17$&$-0.17$ &$-0.13$\\
	&$\it(-0.22,  -0.02)$ & $\it(-0.22,  -0.03)$& $\it(-0.24, -0.05)$& $\it(-0.26,  -0.07)$&$\it(-0.29,  -0.07)$&$\it(-0.26,  -0.015)$\\
	\multirow{2}{*}{\textsf{MKTMOM}}	&$-0.36$ &$-0.28$& $-0.27$ &$-0.34$&$-0.41$&$-0.46$ \\
	& $\it(-0.52,  -0.19)$ &$\it(-0.44, -0.11)$&$\it(-0.44,  -0.11)$&$\it(-0.48,  -0.19)$&$\it(-0.57 , -0.25)$&$\it(-0.62,  -0.30)$\\
	\multirow{2}{*}{\textsf{REITYLD}}		&$0.45$ &$0.56$&$0.32$&&&$0.30$\\
	&$\it(0.21,0.69)$ &$\it(0.37 , 0.75)$&$\it(0.07  ,0.57)$&&&$\it(0.08,  0.53)$\\
	\multirow{2}{*}{\textsf{REITMOM}}		&$0.78132$ &$0.73$&$0.66$&$0.74$&$0.67$&$0.59$\\
	& $\it(0.63,  0.93)$&  $\it(0.59,  0.86)$&$\it(0.52,  0.80)$&$\it(0.60,  0.87)$&$\it(0.51,  0.84)$&$\it(0.45,  0.73)$\\
	\multirow{2}{*}{\textsf{CONST}}	& &&&&&\\
	& &&&&&\\
	\multirow{2}{*}{\textsf{MBASE}}	& &$-0.21$&&&&\\
	& & $\it(-0.35,  -0.065)$&&&&\\
	\multirow{2}{*}{\textsf{INFLAT}}		&$-0.29$ &&&&&\\
	&$\it(0.14,  -0.57)$ &&&&&\\
	\multirow{2}{*}{\textsf{INFLAT\_LSE}}		& &&&&&\\
	& &&&&&\\
	\multirow{2}{*}{\textsf{INDPRD}}	& &&&&&\\
	& &&&&&\\
	\multirow{2}{*}{\textsf{CONSUM}}		& &&&&&\\
	& &&&&&\\
	\multirow{2}{*}{\textsf{PMI}}		& &&&&&\\
	& &&&&&\\
	\multirow{2}{*}{\textsf{DLEAD}}	& &&&&&\\
	& &&&&&\\
	\multirow{2}{*}{\textsf{MICH}}	&$0.20$  &$0.19$&&&&$0.19$\\
	&$\it(0.04,  0.36)$ &$\it(0.08,  0.30)$&&&&$\it(0.04,  0.35)$\\
	\multirow{2}{*}{\textsf{OILPRICE}}& &$0.16$&&&$0.24$&\\
	& &$\it(0.00,  0.31)$&&&$\it(0.01,  0.48)$&\\
	\multirow{2}{*}{\textsf{GASPRICE}}& &&&&&\\
	& &&&&&\\
\bottomrule 
\end{tabular}} 
\label{tab:RE_static_qreg}
\end{center} 
\end{table} 
%
%
\begin{table}[!t]
\caption{Real estate dataset (REITMKT). Summary statistics of the posterior distribution from the Gaussian regression model.}
\begin{center}
\tabcolsep=2.0mm
\resizebox{0.6\columnwidth}{!}{\begin{tabular}{rrrrrrrr}
\toprule
 & mean & sd & \multicolumn{2}{c}{HPD$_{95\%}$} & CD & RHAT & ESS \\ 
 \midrule
 intercept & 0.00 & 0.04 & -0.09 & 0.09 & -0.15 & 1.00 & 40000.00 \\ 
  \rowcolor{Gray}TERM & -0.18 & 0.09 & -0.36 & -0.01 & 0.32 & 1.00 & 40000.00 \\ 
  PREM & -0.20 & 0.12 & -0.43 & 0.04 & -0.77 & 1.00 & 40000.00 \\ 
  \rowcolor{Gray}TBILL & -0.45 & 0.15 & -0.75 & -0.15 & -0.32 & 1.00 & 40000.00 \\ 
  MKTPE & -0.04 & 0.06 & -0.16 & 0.08 & -0.22 & 1.00 & 40000.00 \\ 
  \rowcolor{Gray}MKTYLD & -0.12 & 0.05 & -0.22 & -0.02 & 0.15 & 1.00 & 40000.00 \\ 
  \rowcolor{Gray}MKTMOM & -0.36 & 0.08 & -0.52 & -0.20 & -0.26 & 1.00 & 40000.00 \\ 
 \rowcolor{Gray} REITYLD & 0.45 & 0.13 & 0.21 & 0.70 & 0.18 & 1.00 & 40000.00 \\ 
  \rowcolor{Gray}REITMOM & 0.78 & 0.07 & 0.64 & 0.93 & 0.07 & 1.00 & 40000.00 \\ 
  CONST & 0.00 & 0.05 & -0.09 & 0.10 & 0.14 & 1.00 & 40000.00 \\ 
  MBASE & -0.03 & 0.06 & -0.16 & 0.09 & 0.26 & 1.00 & 40000.00 \\ 
 \rowcolor{Gray} INFLAT & -0.29 & 0.14 & -0.57 & -0.01 & -0.29 & 1.00 & 40000.00 \\ 
  INFLAT\_LFE & 0.07 & 0.09 & -0.10 & 0.24 & -0.69 & 1.00 & 40000.00 \\ 
  INDPRO & -0.02 & 0.10 & -0.20 & 0.17 & 0.38 & 1.00 & 40000.00 \\ 
  CUMSUM & -0.01 & 0.11 & -0.22 & 0.20 & -0.53 & 1.00 & 40000.00 \\ 
  DLEAD & -0.11 & 0.11 & -0.34 & 0.11 & -0.75 & 1.00 & 40000.00 \\ 
 \rowcolor{Gray} MICH & 0.20 & 0.08 & 0.04 & 0.37 & -0.35 & 1.00 & 40000.00 \\ 
  OILPRICE & 0.17 & 0.10 & -0.02 & 0.37 & -0.43 & 1.00 & 40000.00 \\ 
  GASPRICE & -0.03 & 0.11 & -0.24 & 0.19 & 0.05 & 1.00 & 40000.00 \\ 
\bottomrule
\end{tabular}}
\label{tab:US_realestate_mean_regression}
\end{center}
\end{table}
%
\begin{table}[!h]
\caption{Real estate dataset (REITMKT). Summary statistics of the posterior distribution of the quantile regression model for $\tau=0.10$.}
\begin{center}
\tabcolsep=2.0mm
\resizebox{0.6\columnwidth}{!}{\begin{tabular}{rrrrrrrr}
\toprule
 & mean & sd & \multicolumn{2}{c}{HPD$_{95\%}$} & CD & RHAT & ESS  \\ 
\midrule
  \rowcolor{Gray}
intercept & -1.06 & 0.04 & -1.14 & -0.98 & -0.44 & 1.00 & 24025.88 \\ 
  \rowcolor{Gray}
  TERM & -0.30 & 0.09 & -0.48 & -0.11 & 2.33 & 1.00 & 19935.18 \\ 
    \rowcolor{Gray}
  PREM & -0.32 & 0.09 & -0.51 & -0.14 & -2.13 & 1.00 & 28816.86 \\ 
    \rowcolor{Gray}
  TBILL & -0.64 & 0.13 & -0.89 & -0.39 & -0.63 & 1.00 & 25238.84 \\ 
  MKTPE & -0.06 & 0.05 & -0.16 & 0.04 & 0.68 & 1.00 & 26133.83 \\ 
    \rowcolor{Gray}
  MKTYLD & -0.12 & 0.05 & -0.22 & -0.03 & 0.02 & 1.00 & 22793.72 \\ 
    \rowcolor{Gray}
  MKTMOM & -0.28 & 0.08 & -0.44 & -0.12 & 0.32 & 1.00 & 21587.79 \\ 
    \rowcolor{Gray}
  REITYLD & 0.56 & 0.10 & 0.37 & 0.75 & -0.45 & 1.00 & 27545.87 \\ 
    \rowcolor{Gray}
  REITMOM & 0.73 & 0.07 & 0.59 & 0.86 & 0.70 & 1.00 & 24671.42 \\ 
  CONST & 0.09 & 0.05 & -0.01 & 0.19 & -0.04 & 1.00 & 21463.95 \\ 
    \rowcolor{Gray}
  MBASE & -0.21 & 0.07 & -0.35 & -0.06 & -0.13 & 1.00 & 19402.35 \\ 
  INFLAT & -0.18 & 0.13 & -0.44 & 0.07 & -0.33 & 1.00 & 22650.68 \\ 
  INFLAT\_LFE & 0.01 & 0.07 & -0.13 & 0.15 & -0.30 & 1.00 & 25191.55 \\ 
  INDPRO & -0.09 & 0.09 & -0.27 & 0.09 & -0.30 & 1.00 & 21866.25 \\ 
  CUMSUM & 0.16 & 0.09 & -0.01 & 0.33 & 0.25 & 1.00 & 26039.98 \\ 
  DLEAD & -0.07 & 0.09 & -0.24 & 0.10 & -0.82 & 1.00 & 29057.38 \\ 
    \rowcolor{Gray}
  MICH & 0.19 & 0.05 & 0.09 & 0.30 & 1.00 & 1.00 & 31222.92 \\ 
    \rowcolor{Gray}
  OILPRICE & 0.16 & 0.08 & 0.00 & 0.31 & -0.32 & 1.00 & 28569.42 \\ 
  GASPRICE & -0.18 & 0.10 & -0.37 & 0.01 & 0.31 & 1.00 & 25543.91 \\ 
\bottomrule
\end{tabular}}
\label{tab:US_realestate_quantile_regression_010}
\end{center}
\end{table}
%
\begin{table}[!h]
\caption{Real estate dataset (REITMKT). Summary statistics of the posterior distribution of the quantile regression model for $\tau=0.25$.}
\begin{center}
\tabcolsep=2.0mm
\resizebox{0.6\columnwidth}{!}{\begin{tabular}{rrrrrrrr}
\toprule
 & mean & sd & \multicolumn{2}{c}{HPD$_{95\%}$} & CD & RHAT & ESS  \\ 
\midrule

  \rowcolor{Gray}
intercept & -0.56 & 0.04 & -0.65 & -0.47 & -0.09 & 1.00 & 27840.86 \\ 
  TERM & -0.11 & 0.09 & -0.28 & 0.06 & -1.56 & 1.00 & 28813.99 \\ 
  PREM & -0.23 & 0.12 & -0.46 & 0.01 & 0.27 & 1.00 & 29067.11 \\ 
    \rowcolor{Gray}
  TBILL & -0.40 & 0.16 & -0.70 & -0.09 & 1.42 & 1.00 & 28127.68 \\ 
  MKTPE & 0.01 & 0.06 & -0.11 & 0.13 & 1.04 & 1.00 & 27677.34 \\ 
    \rowcolor{Gray}
  MKTYLD & -0.15 & 0.05 & -0.24 & -0.05 & 0.90 & 1.00 & 30158.10 \\ 
    \rowcolor{Gray}
  MKTMOM & -0.27 & 0.08 & -0.44 & -0.11 & 0.94 & 1.00 & 28788.18 \\ 
    \rowcolor{Gray}
  REITYLD & 0.32 & 0.13 & 0.07 & 0.57 & 0.20 & 1.00 & 27062.58 \\ 
    \rowcolor{Gray}
  REITMOM & 0.67 & 0.07 & 0.53 & 0.80 & 2.78 & 1.00 & 30108.25 \\ 
  CONST & 0.07 & 0.05 & -0.02 & 0.17 & -1.00 & 1.00 & 27285.73 \\ 
  MBASE & -0.07 & 0.07 & -0.20 & 0.07 & -2.26 & 1.00 & 27966.47 \\ 
  INFLAT & -0.19 & 0.13 & -0.45 & 0.07 & -0.23 & 1.00 & 31197.94 \\ 
  INFLAT\_LFE & 0.05 & 0.08 & -0.11 & 0.22 & -0.31 & 1.00 & 30560.65 \\ 
  INDPRO & 0.05 & 0.08 & -0.12 & 0.21 & 1.35 & 1.00 & 31407.32 \\ 
  CUMSUM & -0.00 & 0.09 & -0.19 & 0.18 & -0.90 & 1.00 & 32995.79 \\ 
  DLEAD & -0.12 & 0.11 & -0.33 & 0.09 & 0.76 & 1.00 & 30493.51 \\ 
  MICH & 0.12 & 0.07 & -0.01 & 0.26 & 0.08 & 1.00 & 31703.86 \\ 
  OILPRICE & 0.09 & 0.10 & -0.10 & 0.29 & -0.50 & 1.00 & 29289.08 \\ 
  GASPRICE & -0.03 & 0.10 & -0.24 & 0.17 & -0.92 & 1.00 & 29412.36 \\ 
\bottomrule
\end{tabular}}
\label{tab:US_realestate_quantile_regression_025}
\end{center}
\end{table}

\begin{table}[!h]
\caption{Real estate dataset (REITMKT). Summary statistics of the posterior distribution of the quantile regression model for $\tau=0.50$.}
\begin{center}
\tabcolsep=2.0mm
\resizebox{0.6\columnwidth}{!}{\begin{tabular}{rrrrrrrr}
\toprule
 & mean & sd & \multicolumn{2}{c}{HPD$_{95\%}$} & CD & RHAT & ESS  \\ 
\midrule

intercept & -0.00 & 0.04 & -0.08 & 0.08 & 0.48 & 1.00 & 31815.12 \\ 
  TERM & -0.05 & 0.08 & -0.21 & 0.11 & 0.06 & 1.00 & 32303.28 \\ 
  PREM & -0.05 & 0.12 & -0.28 & 0.18 & -0.46 & 1.00 & 30865.66 \\ 
  TBILL & -0.22 & 0.15 & -0.52 & 0.08 & -0.11 & 1.00 & 30953.40 \\ 
  MKTPE & 0.02 & 0.05 & -0.08 & 0.13 & 0.06 & 1.00 & 34005.54 \\ 
   \rowcolor{Gray}
  MKTYLD & -0.17 & 0.05 & -0.26 & -0.08 & 0.02 & 1.00 & 34219.22 \\ 
   \rowcolor{Gray}
  MKTMOM & -0.34 & 0.07 & -0.48 & -0.19 & -1.13 & 1.00 & 33222.78 \\ 
  REITYLD & 0.17 & 0.12 & -0.07 & 0.41 & 0.07 & 1.00 & 30795.38 \\ 
   \rowcolor{Gray}
  REITMOM & 0.74 & 0.07 & 0.60 & 0.88 & -0.49 & 1.00 & 32111.00 \\ 
  CONST & 0.06 & 0.04 & -0.02 & 0.14 & -0.11 & 1.00 & 33282.37 \\ 
  MBASE & 0.01 & 0.06 & -0.11 & 0.13 & -0.79 & 1.00 & 31457.04 \\ 
  INFLAT & -0.14 & 0.14 & -0.42 & 0.14 & -0.50 & 1.00 & 31459.13 \\ 
  INFLAT\_LFE & 0.05 & 0.09 & -0.13 & 0.22 & 0.82 & 1.00 & 31100.54 \\ 
  INDPRO & 0.02 & 0.09 & -0.16 & 0.20 & 0.43 & 1.00 & 31430.16 \\ 
  CUMSUM & -0.10 & 0.11 & -0.31 & 0.11 & 1.03 & 1.00 & 31885.93 \\ 
  DLEAD & -0.10 & 0.11 & -0.32 & 0.13 & -0.33 & 1.00 & 31255.84 \\ 
  MICH & 0.12 & 0.10 & -0.07 & 0.31 & 1.37 & 1.00 & 28053.91 \\ 
  OILPRICE & 0.16 & 0.10 & -0.02 & 0.35 & 0.74 & 1.00 & 32253.69 \\ 
  GASPRICE & -0.10 & 0.10 & -0.30 & 0.10 & -1.60 & 1.00 & 31398.28 \\ 
\bottomrule
\end{tabular}}
\label{tab:US_realestate_quantile_regression_050}
\end{center}
\end{table}

\begin{table}[!h]
\caption{Real estate dataset (REITMKT). Summary statistics of the posterior distribution of the quantile regression model for $\tau=0.75$.}
\begin{center}
\tabcolsep=2.0mm
\resizebox{0.6\columnwidth}{!}{\begin{tabular}{rrrrrrrr}
\toprule
 & mean & sd & \multicolumn{2}{c}{HPD$_{95\%}$} & CD & RHAT & ESS  \\ 
\midrule
  \rowcolor{Gray}
intercept & 0.57 & 0.05 & 0.48 & 0.66 & -1.15 & 1.00 & 27619.58 \\ 
  TERM & -0.10 & 0.10 & -0.30 & 0.09 & -0.52 & 1.00 & 25859.08 \\ 
  PREM & 0.01 & 0.14 & -0.27 & 0.29 & 1.55 & 1.00 & 25553.32 \\ 
  TBILL & -0.13 & 0.17 & -0.46 & 0.19 & 1.41 & 1.00 & 27019.31 \\ 
  MKTPE & 0.05 & 0.07 & -0.09 & 0.18 & -0.81 & 1.00 & 25378.38 \\ 
    \rowcolor{Gray}
  MKTYLD & -0.18 & 0.06 & -0.29 & -0.07 & -0.23 & 1.00 & 28368.37 \\ 
    \rowcolor{Gray}
  MKTMOM & -0.41 & 0.08 & -0.57 & -0.25 & 0.44 & 1.00 & 28642.91 \\ 
  REITYLD & 0.20 & 0.15 & -0.09 & 0.48 & 0.57 & 1.00 & 24204.99 \\ 
    \rowcolor{Gray}
  REITMOM & 0.68 & 0.08 & 0.52 & 0.84 & -2.21 & 1.00 & 26386.37 \\ 
  CONST & -0.01 & 0.05 & -0.10 & 0.08 & -1.09 & 1.00 & 28859.86 \\ 
  MBASE & 0.02 & 0.06 & -0.10 & 0.14 & 1.87 & 1.00 & 29637.38 \\ 
  INFLAT & -0.12 & 0.15 & -0.42 & 0.17 & 0.85 & 1.00 & 28352.56 \\ 
  INFLAT\_LFE & -0.04 & 0.10 & -0.24 & 0.16 & 0.69 & 1.00 & 24322.02 \\ 
  INDPRO & -0.04 & 0.10 & -0.23 & 0.16 & -0.81 & 1.00 & 28364.99 \\ 
  CUMSUM & -0.13 & 0.11 & -0.36 & 0.09 & 0.80 & 1.00 & 28348.77 \\ 
  DLEAD & -0.02 & 0.13 & -0.27 & 0.24 & 0.61 & 1.00 & 25846.89 \\ 
  MICH & 0.19 & 0.10 & -0.00 & 0.39 & -0.24 & 1.00 & 24846.10 \\ 
    \rowcolor{Gray}
  OILPRICE & 0.24 & 0.12 & 0.00 & 0.47 & -1.25 & 1.00 & 24994.57 \\ 
  GASPRICE & -0.09 & 0.12 & -0.32 & 0.14 & -0.94 & 1.00 & 25541.09 \\ 
\bottomrule
\end{tabular}}
\label{tab:US_realestate_quantile_regression_075}
\end{center}
\end{table}

\begin{table}[!h]
\caption{Real estate dataset (REITMKT). Summary statistics of the posterior distribution of the quantile regression model for $\tau=0.90$.}
\begin{center}
\tabcolsep=2.0mm
\resizebox{0.6\columnwidth}{!}{\begin{tabular}{rrrrrrrr}
\toprule
 & mean & sd & \multicolumn{2}{c}{HPD$_{95\%}$} & CD & RHAT & ESS  \\ 
\midrule

  \rowcolor{Gray}
intercept & 1.07 & 0.04 & 1.00 & 1.15 & -0.38 & 1.00 & 24059.27 \\ 
  TERM & -0.09 & 0.08 & -0.26 & 0.07 & -0.20 & 1.00 & 23705.00 \\ 
  PREM & 0.01 & 0.10 & -0.19 & 0.20 & 0.51 & 1.00 & 28299.57 \\ 
  TBILL & -0.09 & 0.13 & -0.35 & 0.17 & 1.33 & 1.00 & 25580.08 \\ 
  MKTPE & 0.04 & 0.06 & -0.08 & 0.17 & -0.67 & 1.00 & 21905.86 \\ 
    \rowcolor{Gray}
  MKTYLD & -0.14 & 0.06 & -0.26 & -0.02 & 0.23 & 1.00 & 19413.51 \\ 
    \rowcolor{Gray}
  MKTMOM & -0.47 & 0.08 & -0.62 & -0.31 & -0.32 & 1.00 & 23689.98 \\ 
    \rowcolor{Gray}
  REITYLD & 0.30 & 0.11 & 0.08 & 0.53 & -0.09 & 1.00 & 24173.38 \\ 
    \rowcolor{Gray}
  REITMOM & 0.60 & 0.07 & 0.46 & 0.74 & -0.12 & 1.00 & 22028.82 \\ 
  CONST & -0.07 & 0.04 & -0.15 & 0.02 & 0.46 & 1.00 & 23569.10 \\ 
  MBASE & 0.04 & 0.06 & -0.08 & 0.17 & 0.24 & 1.00 & 23598.83 \\ 
  INFLAT & -0.20 & 0.12 & -0.44 & 0.03 & 0.52 & 1.00 & 26804.73 \\ 
  INFLAT\_LFE & -0.06 & 0.08 & -0.21 & 0.10 & 0.35 & 1.00 & 24886.73 \\ 
  INDPRO & -0.10 & 0.09 & -0.28 & 0.08 & -0.36 & 1.00 & 22803.85 \\ 
  CUMSUM & 0.00 & 0.09 & -0.18 & 0.19 & 2.75 & 1.00 & 26377.52 \\ 
  DLEAD & 0.03 & 0.10 & -0.17 & 0.22 & -0.05 & 1.00 & 25799.93 \\ 
    \rowcolor{Gray}
  MICH & 0.19 & 0.08 & 0.04 & 0.35 & 0.06 & 1.00 & 24275.79 \\ 
  OILPRICE & 0.17 & 0.10 & -0.02 & 0.37 & -0.65 & 1.00 & 22650.47 \\ 
  GASPRICE & -0.01 & 0.09 & -0.18 & 0.16 & -2.24 & 1.00 & 26734.99 \\ 
\bottomrule
\end{tabular}}
\label{tab:US_realestate_quantile_regression_090}
\end{center}
\end{table}

%
\begin{figure}[!h]
\begin{center}
\subfigure[$\tau=0.10$]{\label{fig:REITMKT_incl_prob010_LOESS}
\includegraphics[trim={0 1cm 0cm 1cm},clip,width=0.45\textwidth]{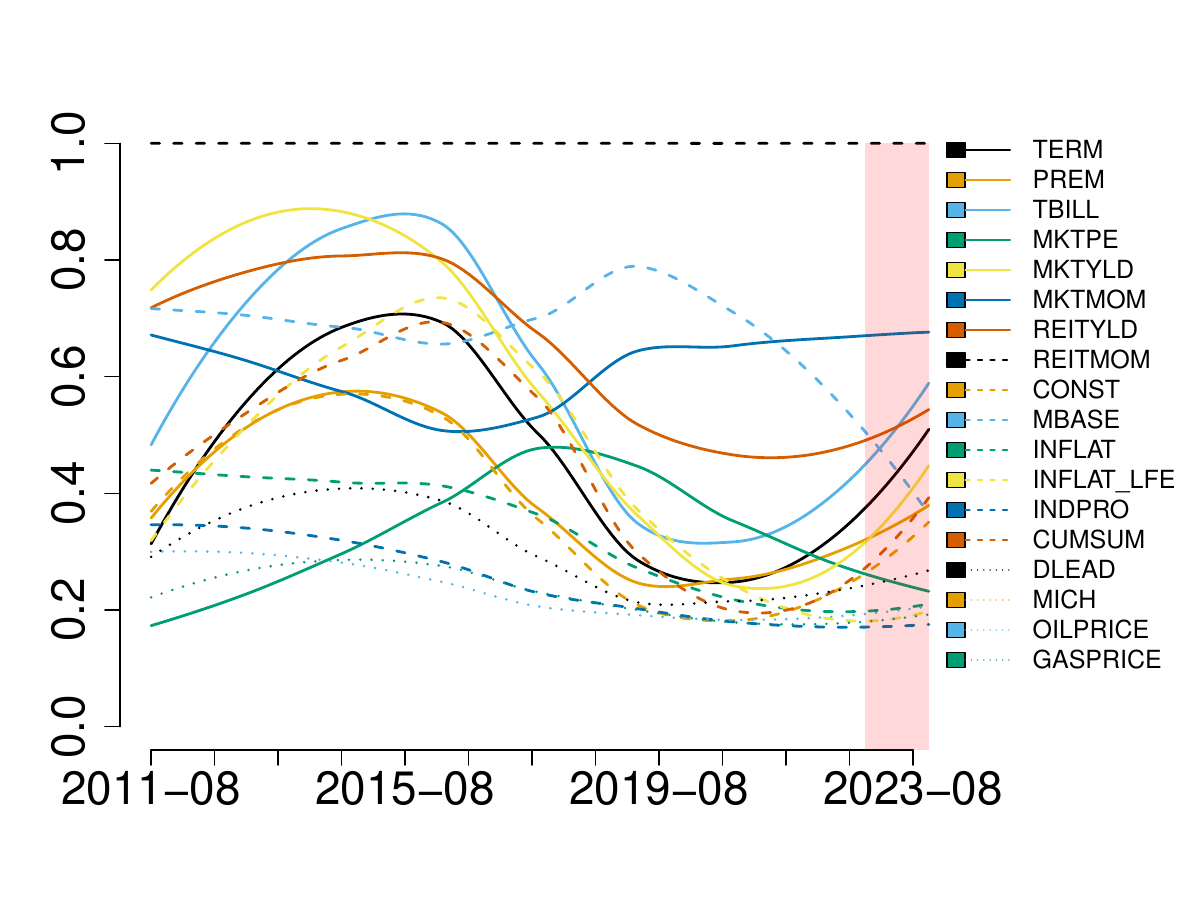}}\qquad
\subfigure[$\tau=0.25$]{\label{fig:REITMKT_incl_prob025_LOESS}
\includegraphics[trim={0 1cm 0cm 1cm},clip,width=0.45\textwidth]{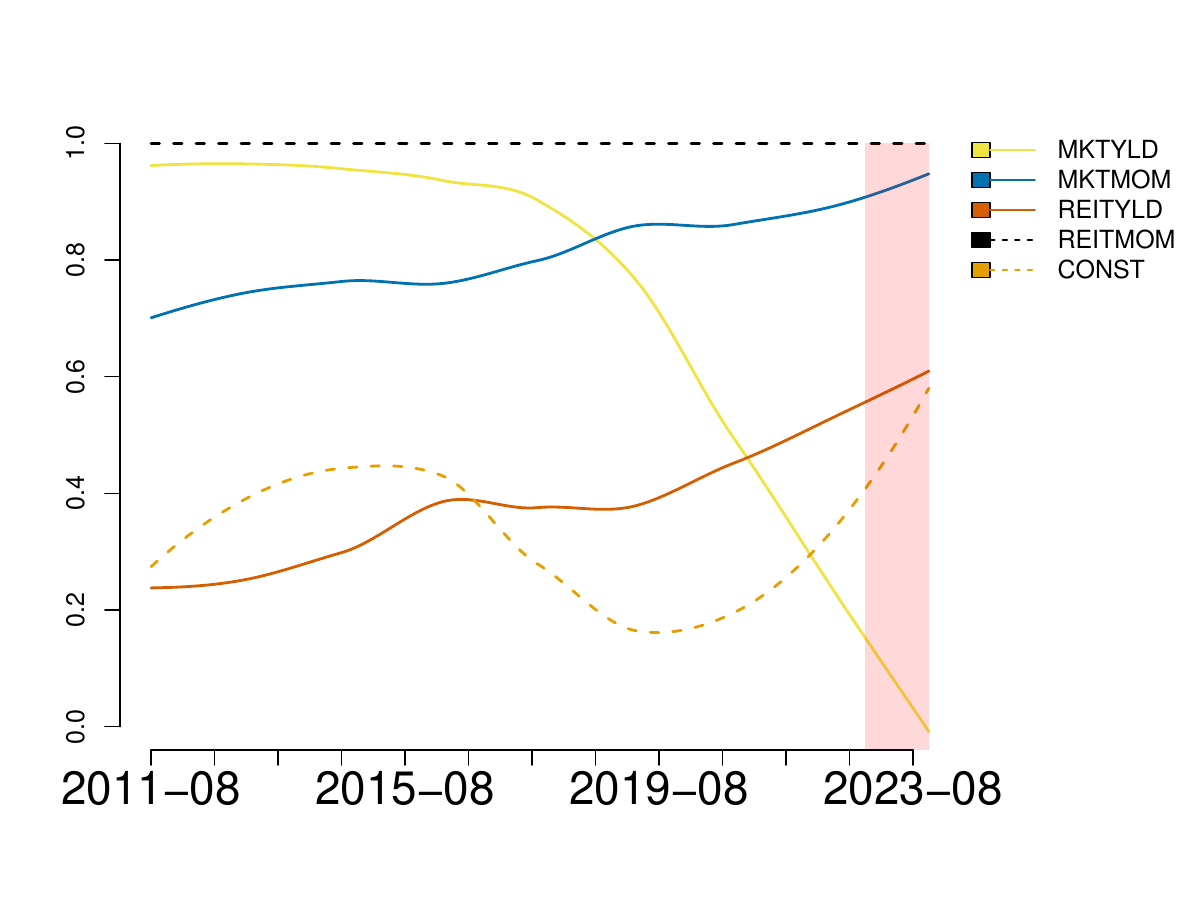}}
\subfigure[$\tau=0.50$]{\label{fig:REITMKT_incl_prob050_LOESS}
\includegraphics[trim={0 1cm 0cm 1cm},clip,width=0.45\textwidth]{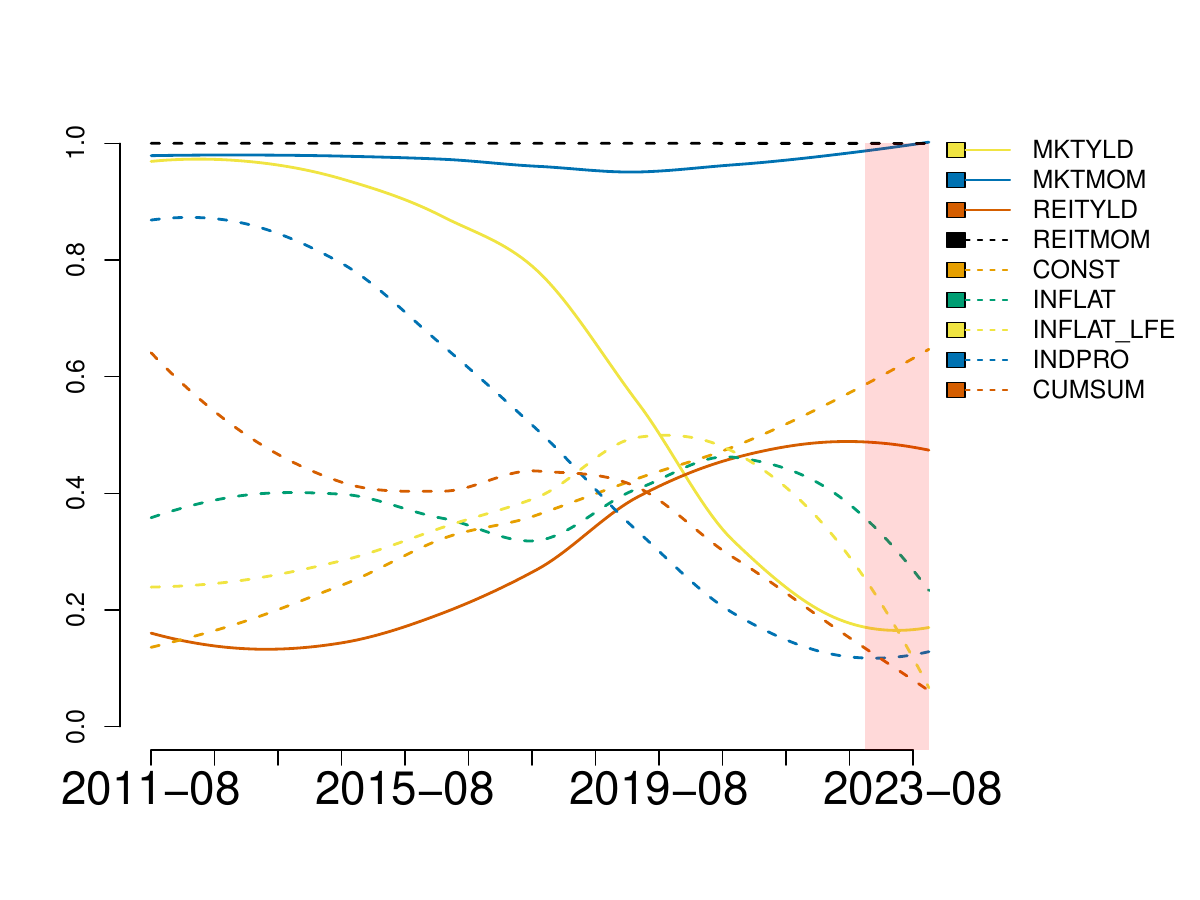}}\qquad
\subfigure[$\tau=0.75$]{\label{fig:REITMKT_incl_prob075_LOESS}
\includegraphics[trim={0 1cm 0cm 1cm},clip,width=0.45\textwidth]{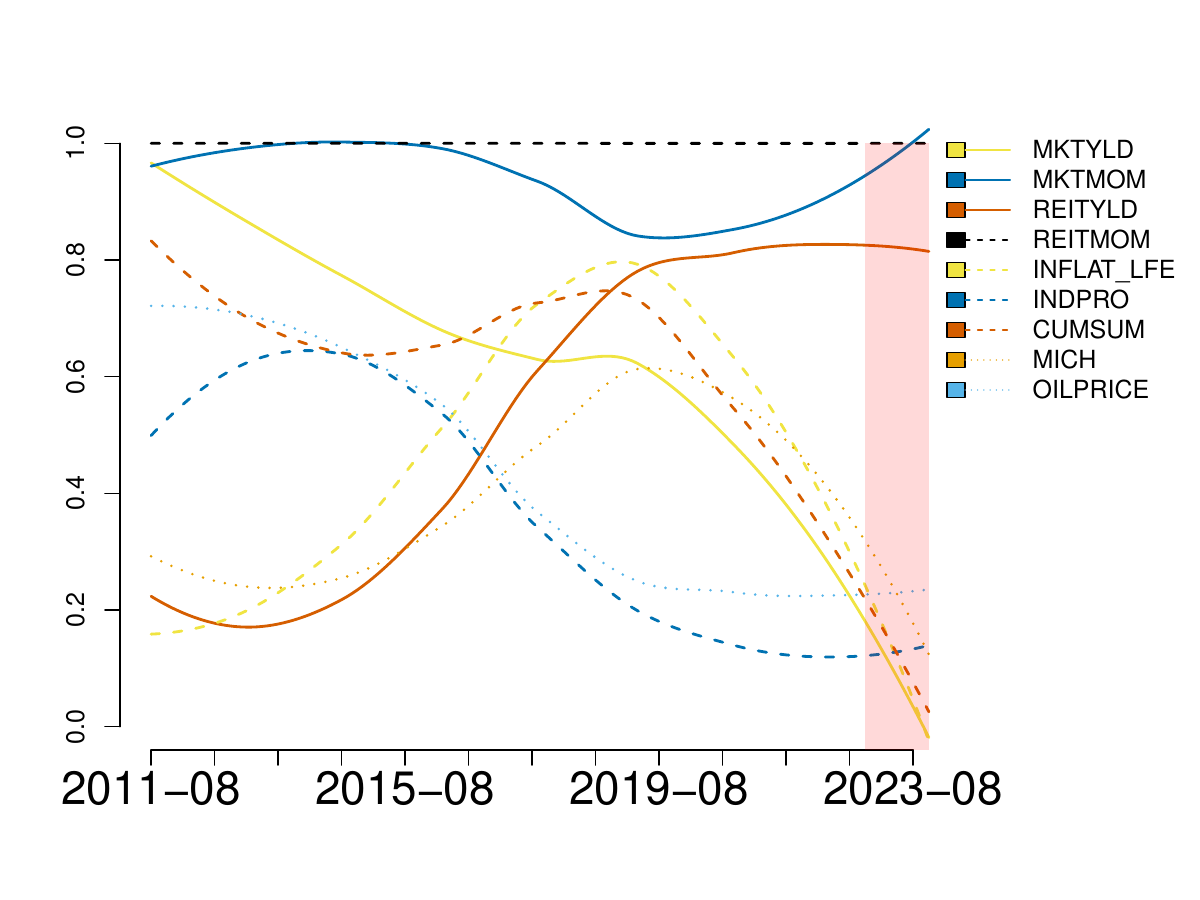}}
\subfigure[$\tau=0.90$]{\label{fig:REITMKT_incl_prob090_LOESS}
\includegraphics[trim={0 1cm 0cm 1cm},clip,width=0.45\textwidth]{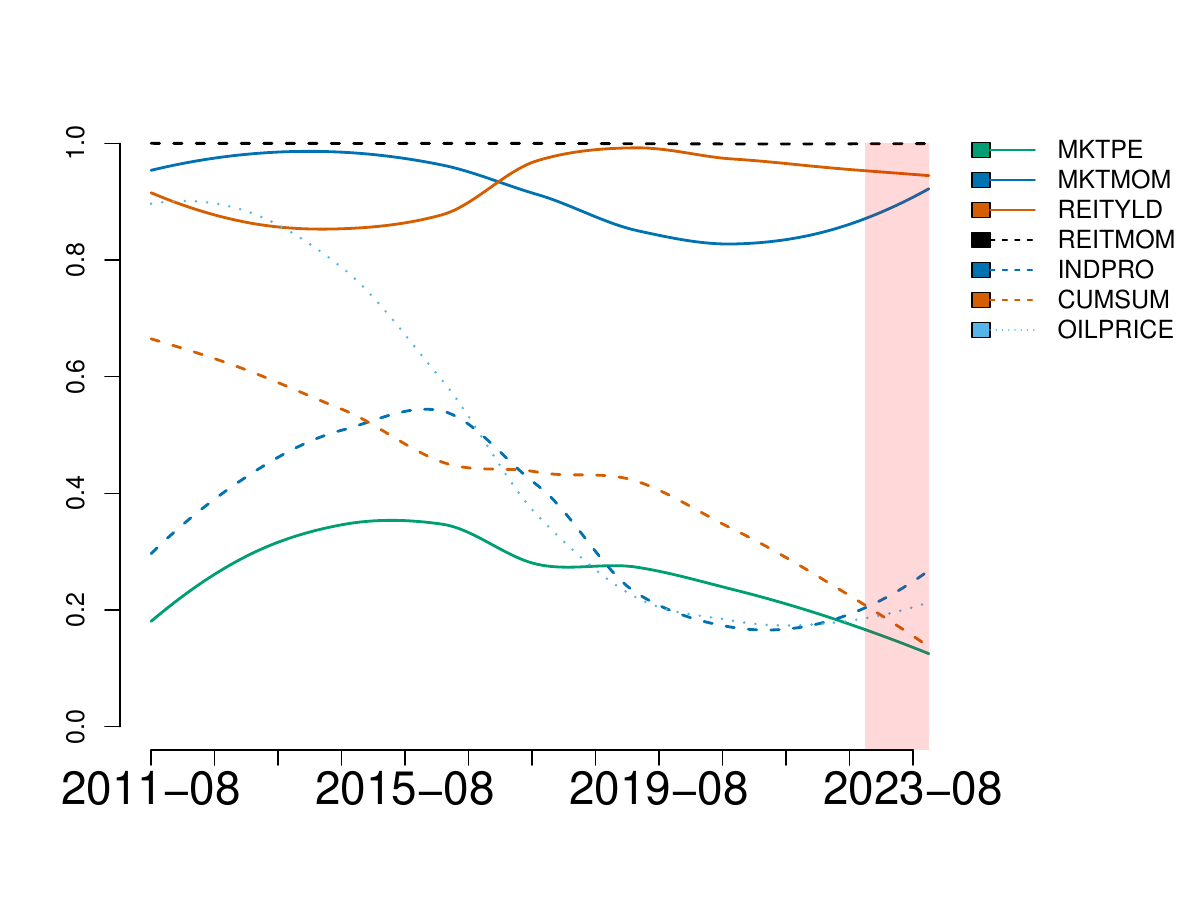}}\qquad
\subfigure[Mean regression]{\label{fig:REITMKT_incl_prob_LOESS}
\includegraphics[trim={0 1cm 0cm 1cm},clip,width=0.45\textwidth]{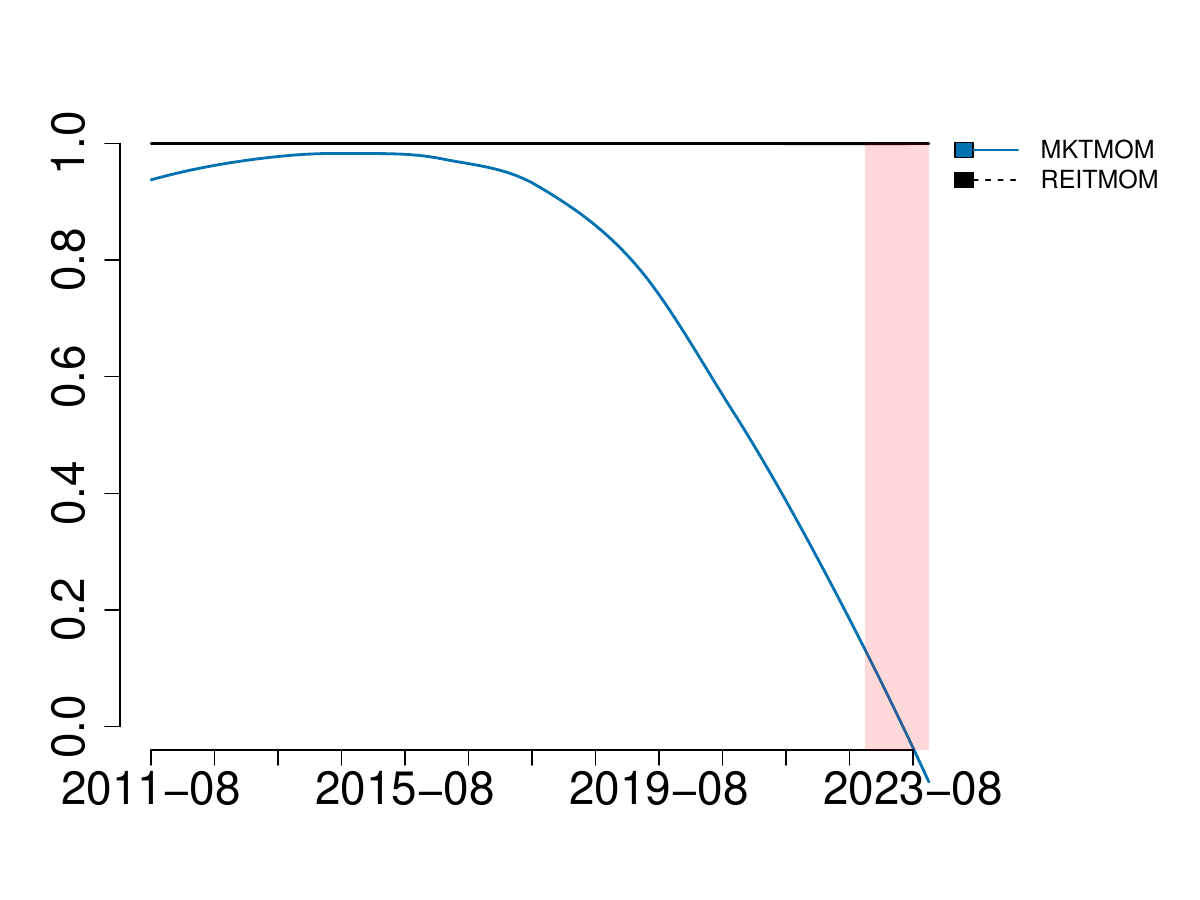}}
\caption{US real estate data. Sequential update of the predicted inclusion
probabilities $\protect\pi _{t|t-1}$ by DMA (for the mean regression \ref{fig:REITMKT_incl_prob}) and BDQMA (for the quantile regression \ref{fig:REITMKT_incl_prob010}-\ref{fig:REITMKT_incl_prob090}) for the REITMKT (monthly NAREIT equity return in excess of the monthly return on the S\&P 500 stock index), filtered by local polynomial regression to remove noise. For each quantile level $\tau$ the corresponding figure only reports those parameters having inclusion probability larger or equal to $0.7$ for at least one quarter. See Table \ref{tab:table_US_inflation_data_CPIAUCSL_ss} for a summary of the relevant covariates. The shaded area identifies the period from 2022-02 to 2023-10 of the Russian-Ukraine war.}
\label{fig:DQMA_RealEstate_data_results_LOESS}
\end{center}
\end{figure}
%
\newpage
\clearpage
\subsection{BDQMA: summary of the relevant parameters}
\label{sec:US_re_additional_results_relevant_parameters}
%
\noindent In this section, for the US real estate dataset (REITMKT), we report the table summarizing the most relevant regression parameters, e.g. those parameters whose inclusion probabilities are greater than $0.7$ for at least one period (month).
\begin{table}[!ht]
\setlength{\tabcolsep}{5 pt}
\caption{Real estate dataset. Summary of the relevant covariates for the response variable \qmo\textsf{REITMKT}\qmc. Here, by relevant covariates we mean those having inclusion probability (as estimated by the BDQMA model for Gaussian and quantile regression) larger than $0.7$ for at least one quarter.} 
\begin{center}\resizebox{0.7\columnwidth}{!}{\begin{tabular}{llcccccc}\\
\toprule
&&\multirow{2}{*}{Mean regression}&\multicolumn{5}{c}{Quantile regression}\\
Id& Name&& $\tau=0.10$ & $\tau=0.25$ & $\tau=0.50$ & $\tau=0.75$ & $\tau=0.90$\\
\cmidrule(lr){1-1}\cmidrule(lr){2-2}\cmidrule(lr){3-3}\cmidrule(lr){4-4}\cmidrule(lr){5-5}\cmidrule(lr){6-6}\cmidrule(lr){7-7}\cmidrule(lr){8-8}
%
    \rowcolor{Gray}
5	&	\textsf{TERM} 		&\checkmark &\checkmark&&&&	\\
    \rowcolor{Gray}
6	&	\textsf{PREM} 		&&\checkmark	&&&&\\
    \rowcolor{Gray}
7	&	\textsf{TBILL} 		&\checkmark &\checkmark&\checkmark&&&	\\
8	&	\textsf{MKTPE} 	&	\\
    \rowcolor{Gray}
9	&	\textsf{MKTYLD} 	&\checkmark &\checkmark&\checkmark&\checkmark&\checkmark&	\checkmark\\
    \rowcolor{Gray}
10	&	\textsf{MKTMOM}	&\checkmark 
&\checkmark&\checkmark&\checkmark&\checkmark&\checkmark	\\
    \rowcolor{Gray}
11	&	\textsf{REITYLD} 	&\checkmark &\checkmark&\checkmark&&&\checkmark	\\
    \rowcolor{Gray}
12	&	\textsf{REITMOM} 	&\checkmark &\checkmark&\checkmark&\checkmark&\checkmark&\checkmark	\\
13	& 	\textsf{CONST} 	&	\\
    \rowcolor{Gray}
14	&	\textsf{MBASE} 	&	&\checkmark&&&&\\
    \rowcolor{Gray}
15	&	\textsf{INFLAT} 		&\checkmark &&&&&	\\
16	&	\textsf{INFLAT\_LFE} 		&	\\
17	&	\textsf{INDPRD} 	&	\\
18	&	\textsf{CONSUM} 	&	\\
19	&	\textsf{DLEAD} 	&	 \\
    \rowcolor{Gray}
20	&	\textsf{MICH} 		&\checkmark &\checkmark&&&&\checkmark	\\
    \rowcolor{Gray}
21	&	\textsf{OILPRICE}&	&\checkmark	&&&\checkmark&\\		
22	&	\textsf{GASPRICE}&		\\	
%
\bottomrule 
\end{tabular}} 
\label{tab:table_US_re_data_ss} 
\end{center} 
\end{table} 
%

\newpage
\bibliographystyle{apalike}
\bibliography{biblio}

\end{document}